\documentclass[12pt,reqno]{amsart}

\usepackage{mathrsfs}
\usepackage{amsmath}
\usepackage{amsfonts}
\usepackage{amssymb}
\usepackage[all,cmtip]{xy}
\usepackage{hyperref}
\usepackage{enumerate}  
\usepackage{enumitem}
\usepackage{ulem}
\usepackage[all]{xy}
\usepackage[toc,page]{appendix}
\usepackage{tikz-cd}
\usepackage[all,cmtip]{xy}
\usepackage[toc,page]{appendix}
\usepackage{amsmath, amsthm, amsfonts, amssymb}
\usepackage{graphicx}
\usepackage{subcaption}
\usepackage{tikz-cd}
\usepackage[utf8]{inputenc}
\usepackage[english]{babel}
\usepackage{dsfont}
\usepackage{xcolor}
\usepackage{arydshln}
\usepackage{faktor}
\usepackage{mathtools}
\usepackage{enumitem}
\usepackage{scalerel}
\usepackage{booktabs}

\def\endpf{\hbox{\vrule height1.5ex width.5em}}

\usepackage{dsfont}
\newcommand\xdashmapsto[2][]{\mathrel{\mapstochar\xdashrightarrow[#1]{#2}}}

\makeatletter
\newcommand*{\da@rightarrow}{\mathchar"0\hexnumber@\symAMSa 4B }
\newcommand*{\da@leftarrow}{\mathchar"0\hexnumber@\symAMSa 4C }
\newcommand*{\xdashrightarrow}[2][]{%
	\mathrel{%
		\mathpalette{\da@xarrow{#1}{#2}{}\da@rightarrow{\,}{}}{}%
	}%
}
\newcommand{\xdashleftarrow}[2][]{%
	\mathrel{%
		\mathpalette{\da@xarrow{#1}{#2}\da@leftarrow{}{}{\,}}{}%
	}%
}
\newcommand{\xdashdownarrow}[2][]{%
	\mathrel{%
		\mathpalette{\da@xarrow{#1}{#2}\da@downarrow{}{}{\,}}{}%
	}%
}
\newcommand*{\da@xarrow}[7]{%
	\sbox0{$\ifx#7\scriptstyle\scriptscriptstyle\else\scriptstyle\fi#5#1#6\m@th$}%
	\sbox2{$\ifx#7\scriptstyle\scriptscriptstyle\else\scriptstyle\fi#5#2#6\m@th$}%
	\sbox4{$#7\dabar@\m@th$}%
	\dimen@=\wd0 %
	\ifdim\wd2 >\dimen@
	\dimen@=\wd2 %
	\fi
	\count@=2 %
	\def\da@bars{\dabar@\dabar@}%
	\@whiledim\count@\wd4<\dimen@\do{%
		\advance\count@\@ne
		\expandafter\def\expandafter\da@bars\expandafter{%
			\da@bars
			\dabar@ 
		}%
	}%
	\mathrel{#3}%
	\mathrel{%
		\mathop{\da@bars}\limits
		\ifx\\#1\\%
		\else
		_{\copy0}%
		\fi
		\ifx\\#2\\%
		\else
		^{\copy2}%
		\fi
	}%
	\mathrel{#4}%
}
\makeatother

\usepackage{graphicx}

\newcommand\tikzmark[1]{%
  \tikz[remember picture,overlay]\coordinate (#1);}

\newcommand{\underbracedmatrixl}[2]{%
  \left(\;
  \smash[b]{\underbrace{
    \begin{matrix}#1\end{matrix}
  }_{#2}}
  \;\right.
  \vphantom{\underbrace{\begin{matrix}#1\end{matrix}}_{#2}}
}

\newcommand{\underbracedmatrixll}[2]{%
  \left(\;\hspace{-.27in}
  \smash[b]{\underbrace{
    \begin{matrix}#1\end{matrix}
  }_{#2}}
  \;\right.
  \vphantom{\underbrace{\begin{matrix}#1\end{matrix}}_{#2}}
}

\newcommand{\underbracedmatrixr}[2]{%
  \left. \;
  \smash[b]{\underbrace{
    \begin{matrix}#1\end{matrix}
  }_{#2}}
  \;\right)
  \vphantom{\underbrace{\begin{matrix}#1\end{matrix}}_{#2}}
}

\newcommand{\underbracedmatrixrr}[2]{%
  \left. \;
  \smash[b]{\underbrace{
    \begin{matrix}#1\end{matrix}
  }_{#2}}
  \;\hspace{-.32in}\right)
  \vphantom{\underbrace{\begin{matrix}#1\end{matrix}}_{#2}}
}

\newcommand{\underbracedmatrix}[2]{%
  \left. \;
  \smash[b]{\underbrace{
    \begin{matrix}#1\end{matrix}
  }_{#2}}
  \;\right.
  \vphantom{\underbrace{\begin{matrix}#1\end{matrix}}_{#2}}
}
 
\newcommand{\bigslant}[2]{{\raisebox{.2em}{$#1$}\left/\raisebox{-.2em}{$#2$}\right.}}

\newtheorem{theorem}{Theorem}[section]
\newtheorem{lemma}[theorem]{Lemma}
\newtheorem{corollary}[theorem]{Corollary}
\newtheorem{proposition}[theorem]{Proposition}

\theoremstyle{definition}
\newtheorem{question}[theorem]{Question}
\newtheorem{definition}[theorem]{Definition}
\newtheorem{example}[theorem]{Example}
\newtheorem{remark}[theorem]{Remark}

\newtheorem{conjecture}[theorem]{Conjecture}
\newtheorem{definitionlemma}[theorem]{Definition-Lemma}
\date{}

\makeatletter
\let\@wraptoccontribs\wraptoccontribs
\makeatother

\begin{document}

\title[Canonical blow-ups of  Grassmann manifolds]{\bf  Canonical blow-ups of  Grassmann manifolds}

\author[Hanlong Fang]{Hanlong Fang\\ \\With an appendix by Kexing Chen}

\address{School of Mathematical Sciences, Peking University, Beijing 100871, China. }
\email{hlfang$ @$pku.edu.cn}

\vspace{3cm} \maketitle

\begin{abstract}
   We introduce  certain canonical blow-ups $\mathcal T_{s,p,n}$, as well as their distinct submanifolds $\mathcal M_{s,p,n}$, of Grassmann manifolds $G(p,n)$ by partitioning the Pl\"ucker coordinates with respect to a parameter $s$. Various geometric aspects of $\mathcal T_{s,p,n}$ and $\mathcal M_{s,p,n}$ are studied, for instance, the smoothness,  the holomorphic symmetries, the (semi-)positivity of the anti-canonical bundles, the existence of K\"ahler-Einstein metrics, the functoriality, etc. 
   In particular, we introduce the notion of homeward compactification, of which  $\mathcal T_{s,p,n}$  are examples, as a generalization of the wonderful compactification. Lastly, a generalization of $\mathcal T_{s,p,n}$ according to  vector-valued parameters $\overline s$ is given, and open questions are raised.
\end{abstract}

\tableofcontents

\section{Introduction}

Let $G(p,n)$, $0<p<n$, be the Grassmann manifold  consisting of complex $p$-planes in the complex $n$-space. 
Each point $x\in G(p,n)$ one to one corresponds to an equivalence class of $p\times n$ non-degenerate matrices, where the equivalence relation is induced by the matrix multiplication from the left by the elements of the general linear group $GL(p,\mathbb C)$. A matrix representative  $\widetilde x$ of $x$ is a matrix in the corresponding equivalence class.

Define an index set $\mathbb I_{p,n}$ by
\begin{equation}\label{ipq}
\mathbb I_{p,n}:=\big\{(i_1,\cdots,i_p)\in\mathbb Z^p\big|1\leq i_p<i_{p-1}<\cdots <\cdots<i_1\leq n\}.
\end{equation}
Denote by $[\cdots ,z_I,\cdots]_{I\in\mathbb I_{p,n}}$ the homogeneous coordinates for the complex projective space $\mathbb {CP}^{N_{p,n}}$ where $N_{p,n}=\frac{n!}{p! (n-p)!}-1$. For each index $I=(i_1,\cdots,i_p)\in\mathbb I_{p,n}$  and  a matrix representative  $\widetilde x$ of  $x\in G(p,n)$, denote by  $P_I(\widetilde x)$  the determinant of the submatrix of $\widetilde x$ consisting of the $i_1^{th}, \cdots,i_p^{th}$ 
columns. The Pl\"ucker embedding of  $G(p,n)$ into $\mathbb {CP}^{N_{p,n}}$  can be given  by
\begin{equation}\label{plucker}
\begin{split}
e:G(p,n)&\longrightarrow\mathbb {CP}^{N_{p,n}}\\
x&\mapsto [\cdots ,P_I(\widetilde x), \cdots]_{I\in\mathbb I_{p,n}}
\end{split}.
\end{equation}

This paper relies on the following interesting partition of the index set $\mathbb I_{p,n}$. Fix an integer $s$, and define index sets  $\mathbb I_{s,p,n}^{k}$, $0\leq k\leq p$, by
\begin{equation}\label{ikk}
\mathbb I_{s,p,n}^{k}:=\big\{(i_1,\cdots,i_p)\in\mathbb Z^p\big|{\footnotesize 1\leq i_p<\cdots<i_{k+1}\leq s\,;s+1\leq i_{k}<i_{k-1}<\cdots<i_1\leq n}\}.
\end{equation}
Then there is a partition
\begin{equation}
\mathbb I_{p,n}=\bigsqcup_{k=0}^p\mathbb I_{s,p,n}^{k}\,\,.
\end{equation}

Consider the following linear subspace of $\mathbb {CP}^{N_{p,n}}$ for $0\leq k\leq p$,
\begin{equation}\label{subl}
    \big\{[\cdots ,z_I,\cdots]_{I\in\mathbb I_{p,n}}\in\mathbb {CP}^{N_{p,n}}\big|z_I=0\,,\,\,\forall I\notin\mathbb I_{s,p,n}^k \big\}\,.
\end{equation}
It is clear that the subspace (\ref{subl}) is isomorphic to the complex projective space $\mathbb {CP}^{N^k_{s,p,n}}$
where $N^k_{s,p,n}$ is the cardinal number of the set $\mathbb I_{s,p,n}^{k}$ minus $1$; by a slight abuse of notation, we denote (\ref{subl}) by $\mathbb {CP}^{N^k_{s,p,n}}$ and its homogeneous coordinates by  $[\cdots ,z_I,\cdots]_{I\in\mathbb I^k_{s,p,n}}$. Recall the following projection (rational) map  $F_s^k$ by dropping the coordinates whose indices are not in $\mathbb I_{s,p,n}^k$.
\begin{equation}\label{FK}
\begin{split}
F_s^k:\mathbb {CP}^{N_{p,n}}&\dashrightarrow\mathbb {CP}^{N^k_{s,p,n}} \\
[\cdots ,z_I,\cdots]_{I\in\mathbb I_{p,n}}&\xdashmapsto{}[\cdots, z_I,\cdots]_{I\in\mathbb I^k_{s,p,n}}.
\end{split}
\end{equation}
We make the convention that  $\mathbb{CP}^{N^k_{s,p,n}}$ is a point and $F_s^k$ is the trivial map when $N^k_{s,p,n}=0,-1$.

Now we define a rational map  \begin{equation}
   \mathcal K_{s,p,n}:=(e, f_s^0,\cdots,f_s^p) \,,
\end{equation}
where  $e$ is the Pl\"ucker embedding in (\ref{plucker}) and $f_s^k:=F_s^k\circ e$ for $0\leq k\leq p$\,. More precisely, 
\begin{equation}\label{bpp}
\begin{split}
\mathcal K_{s,p,n}&:\,\,G(p,n)\dashrightarrow \mathbb {CP}^{N_{p,n}}\times\mathbb {CP}^{N^0_{s,p,n}}\times\cdots\times\mathbb {CP}^{N^p_{s,p,n}}\\
x&\xdashmapsto{} \big( [\cdots ,P_I(\widetilde x),\cdots]_{I\in\mathbb I_{p,n}},[\cdots,P_I(\widetilde x),\cdots]_{I\in\mathbb I^0_{s,p,n}},\cdots,[\cdots,P_I(\widetilde x),\cdots]_{I\in\mathbb I^p_{s,p,n}}\big)
\end{split}.
\end{equation}

\begin{definition}\label{ypt}
Assume that $0<p<n$ and $0<s<n$. Let $\mathcal T_{s,p,n}$ be the scheme-theoretic closure  of the birational image of $G(p,n)$ under  $\mathcal K_{s,p,n}$ in $\mathbb {CP}^{N_{p,n}}\times\mathbb {CP}^{N^0_{s,p,n}}\times\cdots\times\mathbb {CP}^{N^p_{s,p,n}}$. We call $\mathcal T_{s,p,n}$ the {\it canonical blow-up of $G(p,n)$} with respect to the parameter $s$. 
\end{definition}

\begin{example}\label{de}
Denote by $[x_1,\cdots, x_s,y_1,\cdots,y_{n-s}]$ the homogeneous coordinates for the projective space $\mathbb{CP}^{n-1}$. Then $\mathcal T_{s,1,n}$ is the blow-up of $\mathbb{CP}^{n-1}$ along the union of the disjoint linear subspaces $\mathbb{CP}^{s-1}$ and $\mathbb{CP}^{n-s-1}$; $\mathcal K_{s,1,n}$ is given by 
\begin{equation}
\mathcal K_{s,1,n}([x_1,\cdots, x_s,y_1,\cdots,y_{n-s}])=\left([x_1,\cdots, x_s,y_1,\cdots,y_{n-s}], [x_1,\cdots, x_s],[y_1,\cdots,y_{n-s}]\right).
\end{equation}
\end{example}

We  have

\begin{proposition}\label{tspnsmooth}
	$\mathcal T_{s,p,n}$ is  smooth. 
\end{proposition}

Denote by $Z(n,\mathbb C)$ the center of the general linear group $GL(n,\mathbb C)$. Define a subgroup $GL(s,\mathbb C)\times GL(n-s,\mathbb C)$  of $GL(n,\mathbb C)$ by
\begin{equation}\label{sns}
GL(s,\mathbb C)\times GL(n-s,\mathbb C):=\left\{\left.\left(
\begin{matrix}
g_1&0\\
0&g_2\\
\end{matrix}\right)\right\vert_{} g_1\in GL(s,\mathbb C), g_2\in GL(n-s,\mathbb C)\right\}.
\end{equation}
It is clear that the quotient group $\bigslant{GL(s,\mathbb C)\times GL(n-s,\mathbb C)} {Z(n,\mathbb C)}$ is a subgroup of the projective linear group $PGL(n,\mathbb C)$.

By a result of Chow (\cite{C}), the holomorphic automorphism group of the Grassmannian $G(p,n)$ is  $PGL(n,\mathbb C)$  when $n\neq 2p$, and  $PGL(n,\mathbb C)\rtimes\mathbb Z/2\mathbb Z$ when $n=2p$ (the discrete symmetry is given by the so-called dual map). 

We determine the symmetry group of $\mathcal T_{s,p,n}$ as follows.

\begin{proposition}\label{auto1}
The holomorphic automorphism group ${\rm Aut}(\mathcal T_{s,p,n})$ of $\mathcal T_{s,p,n}$ is
\begin{equation}
   \bigslant{GL(s,\mathbb C)\times GL(n-s,\mathbb C)} {Z(n,\mathbb C)}\,\,,
\end{equation}
except for the following cases.
\begin{enumerate}[label={\rm(\Alph*)}]
	\item {\rm (USD\footnote{USD stands for a discrete symmetry which turns $\mathcal T_{s,p,n}$ upside down.} case)}. If $n=2s$ and $s\neq p$,
	\begin{equation}
     {\rm Aut}(\mathcal T_{s,p,2s})=\left(\bigslant{GL(s,\mathbb C)\times GL(n-s,\mathbb C)} {Z(n,\mathbb C)}\right)\rtimes\mathbb Z/2\mathbb Z\,\,.
    \end{equation} 
	\item {\rm (DUAL\footnote{DUAL stands for the discrete symmetry  induced by the dual map of $G(p,2p)$.} case)}. If $n=2p$ and $s\neq p$,
	\begin{equation}
      {\rm Aut}(\mathcal T_{s,p,2p})=\left(\bigslant{GL(s,\mathbb C)\times GL(n-s,\mathbb C)} {Z(n,\mathbb C)}\right)\rtimes\mathbb Z/2\mathbb Z\,\,.
    \end{equation} 
	\item {\rm (USD+DUAL\footnote{The two types of discrete symmetries in (A) and (B) coexist in (C).} case)}.  If $n=2s=2p$ and $p\geq 2$,
	\begin{equation}
      {\rm Aut}(\mathcal T_{p,p,2p})=\left(\bigslant{GL(s,\mathbb C)\times GL(n-s,\mathbb C)} {Z(n,\mathbb C)}\right)\rtimes\mathbb Z/2\mathbb Z\rtimes\mathbb Z/2\mathbb Z\,\,.
    \end{equation} 
    \item {\rm (Degenerate cases)}. 
    \begin{enumerate}[label={\rm(\alph*).}]
    \item  $\mathcal T_{1,1,2}\cong\mathbb {CP}^1$ and  ${\rm Aut}(\mathcal T_{1,1,2})=PGL(2,\mathbb C)$. 
    \item For $m\geq 2$, $\mathcal T_{m,1,m+1}\cong\mathcal T_{m,m,m+1}\cong\mathcal T_{1,1,m+1}\cong\mathcal T_{1,m,m+1}$; their holomorphic automorphism groups are isomorphic to the following subgroup of $PGL(n,\mathbb C)$. \begin{equation}\label{parabolic} \left\{\left. \left( \begin{matrix} V&\eta\\ 0&1\\ \end{matrix}\right)\right\vert_{} \begin{matrix} V\in GL(n-1,\mathbb C) \,,\,\,\eta \,\,{\rm is \,\,a\,\,} (n-1)\times 1\,\,{\rm matrix}\\ \end{matrix} \right\}\,. \end{equation}
    \end{enumerate}
\end{enumerate}   
\end{proposition}

We note that  the geometry of $\mathcal T_{s,p,n}$ is largely encoded in a subvariety $\mathcal M_{s,p,n}$.
We next give a  definition of $\mathcal M_{s,p,n}$ when $p\leq s$ (see Definition \ref{mspn} for the general case).
First notice that the inverse rational map $\mathcal K_{s,p,n}^{-1}:\mathcal T_{s,p,n}\dashrightarrow G(p,n)$ has a regular extension $R_{s,p,n}:\mathcal T_{s,p,n}\rightarrow G(p,n)$ (see Definition-Lemma \ref{rspn}). Let $ \mathcal M_{s,p,n} :=R_{s,p,n}^{-1}(G(p,s))$ be the scheme-theoretic inverse image of $G(p,s)$ under the holomorphic map $R_{s,p,n}$, where $G(p,s)\subset G(p,n)$ is a sub-grassmannian consisting of the points with the following matrix representatives.
\begin{equation}
\widetilde x=\underbracedmatrixl{
* &* &\cdots& * \\
* &* &\cdots& *\\
\vdots&\vdots&\ddots&\vdots\\
* &* &\cdots&* \\ }{s\,\rm columns}\hspace{-.22in}
\begin{matrix}
  &\hfill\tikzmark{e}\\
  \\
  \\
  \\
  &\hfill\tikzmark{f}\end{matrix}\hspace{-.1in}
  \begin{matrix}
  &\hfill\tikzmark{a}\\
  \\
  \\
  \\
  &\hfill\tikzmark{b}\end{matrix}\,
  \underbracedmatrixr{ 0  &0& \cdots & 0\\
 0  &0 & \cdots&0\\
  \vdots&\vdots&\ddots&\vdots \\
 0  &0 & \cdots& 0\\} {(n-s) \,\rm columns}\,.
  \tikz[remember picture,overlay]   \draw[dashed,dash pattern={on 4pt off 2pt}] ([xshift=0.5\tabcolsep,yshift=7pt]a.north) -- ([xshift=0.5\tabcolsep,yshift=-2pt]b.south);\tikz[remember picture,overlay]   \draw[dashed,dash pattern={on 4pt off 2pt}] ([xshift=0.5\tabcolsep,yshift=7pt]e.north) -- ([xshift=0.5\tabcolsep,yshift=-2pt]f.south);
\end{equation}

We have that
\begin{proposition}\label{msmooth}
$\mathcal M_{s,p,n}$ is smooth. 
\end{proposition}

\begin{proposition}\label{mauto} The holomorphic automorphism group ${\rm Aut}(\mathcal M_{s,p,n})$ of $\mathcal M_{s,p,n}$ is 
\begin{equation}
    PGL(s,\mathbb C)\times PGL(n-s,\mathbb C)\,,
\end{equation} 
except for the following cases. 
\begin{enumerate}[label={\rm(\Alph*)}]
\item {\rm (Usd\footnote{Usd stands for a discrete symmetry induced by a symmetry turns $\mathcal T_{s,p,n}$ upside down.} case)}. If $n=2s$ and $p\neq s$, 
\begin{equation}
  {\rm Aut}(\mathcal M_{s,p,2s})=\left(PGL(s,\mathbb C)\times PGL(n-s,\mathbb C)\right) \rtimes\mathbb Z/2\mathbb Z\,. 
\end{equation}
\item {\rm (Dual\footnote{Dual stands for the discrete symmetry induced by the dual map of $G(p,2p)$.}  case)}. If $n=2p$ and $p\neq s$, 
\begin{equation}
  {\rm Aut}(\mathcal M_{s,p,2p})=\left(PGL(s,\mathbb C)\times PGL(n-s,\mathbb C)\right) \rtimes\mathbb Z/2\mathbb Z\,. 
\end{equation}
\item {\rm (Usd+Dual\footnote{The two types of discrete symmetries in (A) and (B) coexist in (C).} case)}. If $n=2s=2p$ and $p\geq 3$,
\begin{equation}
  {\rm Aut}(\mathcal M_{p,p,2p})=\left(PGL(s,\mathbb C)\times PGL(n-s,\mathbb C)\right) \rtimes\mathbb Z/2\mathbb Z\rtimes\mathbb Z/2\mathbb Z\,. 
\end{equation}
\item {\rm (Degenerate cases)}. \begin{enumerate}[label={\rm(\alph*).}]
      \item  $\mathcal M_{1,1,2}$ is a point.
      \item  $\mathcal M_{2,2,4}\cong\mathbb {CP}^3$ and ${\rm Aut}(\mathcal M_{2,2,4})=PGL(4,\mathbb C)$\,.
\end{enumerate}
\end{enumerate}
 \end{proposition}

\begin{remark}\label{dualusd}
The above discrete symmetries exist in a more general setting as the following isomorphisms   (see Definitions \ref{tiso}, \ref{USD1}, \ref{musd}, and \ref{mdual}).
\begin{equation}
\begin{split}
     &{\rm DUAL}:\,\mathcal T_{s,p,n}\rightarrow\mathcal T_{s,n-p,n}\,\,\,\,{\rm and}\,\,\,\,
     {\rm USD}:\mathcal T_{s,p,n}\rightarrow\mathcal T_{n-s,p,n}\,;\\
     &  {\rm Dual}:\mathcal M_{s,p,n}\rightarrow \mathcal M_{s,n-p,n} \,\,\,\,{\rm and}\,\,\,\,  {\rm Usd}:\mathcal M_{s,p,n}\rightarrow\mathcal M_{n-s,p,n}\,.
\end{split}
\end{equation}
We can thus assume that $2p\leq n\leq 2s$ without loss of  generality.
\end{remark}

When $s=p$, the manifold
$\mathcal M_{s,p,n}$ is isomorphic to a classical object studied in algebraic geometry called the variety of complete collineations.  Study (\cite {Stu}), Severi (\cite{Sev1,Sev2}), and Van der Waerden (\cite{Van}) studied the complete conics in $\mathbb P^2$ from the perspective of enumerative geometry.  Semple (\cite{Se1, Se2}), Alguneid (\cite {Al}), and  Tyrrell (\cite{Ty}) studied the complete collineations in $\mathbb P^n$ of higher dimensions.  Vainsencher (\cite{Va}) established the smoothness and proved that  $\mathcal M_{p,p,n}$ is a wonderful variety in the sense of   De Concini and Procesi (\cite{DP}). Brion (\cite{Br6}) and Massarenti (\cite{Mas}) determined the holomorphic automorphism groups of  $\mathcal M_{p,p,2p}$ and $\mathcal M_{p,p,n}$ respectively; the connected component of the identity element of the holomorphic automorphism group of  $\mathcal M_{s,p,n}$  follows from the general result of Brion (\cite{Br6}) on wonderful varieties.

We show that $\mathcal M_{s,p,n}$  is a relative version of the variety of complete collineations.
\begin{proposition}\label{fiber}
$\mathcal M_{s,p,n}$ has locally holomorphically trivial fibration structures as follows. 
\begin{enumerate}[label=\rm(\Alph*)]
    \item If $p<n-s$ and $p<s$, $\mathcal M_{s,p,n}$ is a fibration with the base $G(p,s)$ and the fiber $\mathcal M _{p,p,n-s+p}$, and a fibration with the base $G(p,n-s)$  and the fiber  $\mathcal M_{p,p,s+p}$.
    \item If $n-s<p<s$, $\mathcal M_{s,p,n}$ is a fibration with the base  $G(p,s)$ and the fiber $\mathcal M _{n-s,n-s,n-s+p}$, and a fibration with the base   $G(s+p-n,s)$ and the fiber  $\mathcal M _{n-s,n-s,2n-s-p}$.
    \item If $p=n-s<s$, $\mathcal M_{s,p,n}$ is a  fibration with the base   $G(p,s)$ and the fiber  $\mathcal M _{p,p,2p}$.
\end{enumerate}
\end{proposition}

To have a better understanding of  $\mathcal T_{s,p,n}$ and $\mathcal M_{s,p,n}$, we  briefly recall the notion of spherical varieties, and in particular the wonderful compactification.

\begin{definition} Let $G$ be a connected reductive group. An irreducible normal $G$-variety $X$ is called {\it spherical} if a Borel subgroup of $G$ has an open orbit on $X$. \end{definition}

\begin{definition}[wonderful compactification]\label{wond}
Let $G$ be a connected reductive group  and $X$ an irreducible algebraic $G$-variety. Then $X$ is a {\it wonderful variety} if the following holds.
\begin{enumerate}[label={\rm(\Alph*)}]
\item $X$ is smooth and complete.
\item $G$ has an open orbit in $X$ whose complement is the union of $m$ smooth prime divisors $D_i$, $1\leq i\leq m$, with normal crossing. We call each $D_i$, $1\leq i\leq m$, a $G$-stable divisor of $X$.
\item The intersection of the divisors $D_i$ is nonempty. Each $G$-orbit one to one corresponds to the variety $X_I$ defined by \begin{equation}
X_I:=\left(\,\bigcap_{i\in I}D_i\right)\mathbin{\scaleobj{2.1}{\backslash}} \left(\bigcup_{\substack{1\leq j\leq m\\\,\,j\notin I}}D_j\right),
\end{equation}
where $I$ is a subset of $\{1,2,\cdots,m\}$. Moreover, the closure of each $G$-orbit in $X$ is smooth.
\end{enumerate}
If $X^0$ is a $G$-invariant open subvariety of the wonderful variety $X$, we call $X$ a {\it wonderful compactification} of $X^0$.
\end{definition}

One can verify that $\mathcal T_{s,p,n}$ and $\mathcal M_{s,p,n}$ are spherical $GL(s,\mathbb C)\times GL(n-s,\mathbb C)$-varieties. 
Moreover, we prove

\begin{proposition}\label{mwond}
  $\mathcal M_{s,p,n}$ is a wonderful variety.
\end{proposition}
\begin{proposition}\label{gwond}
The complement of the open
$GL(s,\mathbb C)\times GL(n-s,\mathbb C)$-orbit in $\mathcal T_{s,p,n}$ is a simple normal crossing divisor consisting of $2r$  
smooth, irreducible divisors as follows.
\begin{equation}    
D^-_1, 
D^-_2,\cdots,D^-_r,D^+_1,D^+_2,\cdots,D^+_r\,.
\end{equation}
The following holds under a certain rearrangement of the indices of the above divisors. 
\begin{enumerate}[label={\rm(\Alph*)}]
\item $D^-_1=\mathcal M_{s,p,n}$.
\item Each $GL(s,\mathbb C)\times GL(n-s,\mathbb C)$-orbit of $\mathcal T_{s,p,n}$ one to one corresponds to the  quasi-projective variety $X_{(I^-,I^+)}$ defined  by
\begin{equation}\label{inrule}
X_{(I^-,I^+)}:=\left(\bigcap_{\,\,i\in I^-}D^-_i\mathbin{\scaleobj{1.2}{\bigcap}} \bigcap_{\,\,i\in I^+}D^+_i\right)\mathbin{\scaleobj{2.1}{\backslash}} \left(\bigcup_{\substack{1\leq j\leq r\\\,\,j\notin I^-}}D^-_j\mathbin{\scaleobj{1.2}{\bigcup}}\,\bigcup_{\substack{1\leq j\leq r\\\,\,\,j\notin I^+}}D^+_j\right),
\end{equation} 
where $I^-,I^+$ are subsets of $\{1,2,\cdots,r\}$
such that 
\begin{equation}
\min(I^-)+\min(I^+)\geq r+2.
\end{equation}
Here we make the convention that $\min(\emptyset)=+\infty$.
Moreover, the closure of each $GL(s,\mathbb C)\times GL(n-s,\mathbb C)$-orbit in $\mathcal T_{s,p,n}$ is smooth.
\item There is a $GL(s,\mathbb C)\times GL(n-s,\mathbb C)$-equivariant flat map $\mathcal P_{s,p,n}:\mathcal T_{s,p,n}\rightarrow D_1^-$ such that the following  holds.
\begin{enumerate}[label={\rm(\alph*)}]
\item $\mathcal P_{s,p,n}$ is a  retraction, that is, the restriction of $\mathcal P_{s,p,n}$ on $D^-_1$ is the identity map.
\item The restriction of $\mathcal P_{s,p,n}$ on $D^+_1$ is an isomorphism.
\item  $\mathcal P_{s,p,n}(D^-_i)=\mathcal P_{s,p,n}(D^+_{r+2-i})$ for $2\leq i\leq r$.
\end{enumerate}
\item  $D_1^-$ is a wonderful variety with the $GL(s,\mathbb C)\times GL(n-s,\mathbb C)$-stable divisors $\check D_2,\cdots,\check D_r$, where $\check D_i:=\mathcal P_{s,p,n}(D^-_i)$,  $2\leq i\leq r$.
\item When $n=2s$ or $n=2p$ there is a holormorphic automoprhism  $\sigma$ of $\mathcal T_{s,p,n}$ such that $\sigma\circ \sigma$ is the identity map and $\sigma(D^{\pm}_i)=D^{\mp}_i$ for $1\leq i\leq r$.
\end{enumerate}
\end{proposition}
\begin{remark}
The generic fiber of $\mathcal P_{s,p,n}$ is a smooth rational curve of degree $r$ with respect to the Pl\"ucker embedding of $G(p,n)$. There is a moduli-interpretation of  $\mathcal P_{s,p,n}:\mathcal T_{s,p,n}\rightarrow\mathcal M_{s,p,n}$ in the sense of Fujiki (\cite{Fu}). 
\end{remark} 
\begin{remark} Kausz (\cite{Ka}) constructed $\mathcal T_{p,p,2p}$ as the modular compactification of the reductive group $GL(p,\mathbb C)$, and observed  Property (B) in Proposition \ref{gwond}. Rittatore (\cite{Ri1,Ri2}) formulated and classified the toroidal embeddings of reductive groups.
\end{remark}

\begin{definition}\label{rank}
 We call the integer $r$ in Proposition \ref{gwond} the {\it rank} of $\mathcal T_{s,p,n}$.
\end{definition}
\begin{remark}
We can show that \begin{equation}
    r=\min\{s,n-s,p,n-p\}.
\end{equation} We make the convention that  $r$ is always referred to the above quantity in this paper.
\end{remark}

We thus propose a generalization of the wonderful compactification as follows. 
\begin{definition}[homeward compactification] \label{wondm}
Let $G$ be a connected reductive group and $\mathcal T$ an irreducible algebraic $G$-variety.  Then $\mathcal T$ is called a {\it homeward variety}   if the following holds.
	\begin{enumerate}[label={\rm(\Alph*)}]
		\item $\mathcal T$ is smooth and complete.
		\item $G$ has an open orbit in $\mathcal T$ whose complement is a simple normal crossing divisor. Denote the irreducible components of the  complement by $D_1, D_2,\cdots,D_l$.
		\item Each $G$-orbit one to one corresponds to a quasi-projective variety $X_{I}$, provided that $X_{I}$ is non-empty, as follows.
        \begin{equation}
		X_{I}=\left(\bigcap_{\,\,i\in I}D_i\right)\mathbin{\scaleobj{2.1}{\backslash}} \left(\bigcup_{\substack{1\leq j\leq l\\\,\,j\notin I}}D_j\right),
		\end{equation} 
		where $I$ is a subset of $\{1,2,\cdots,l\}$.
	    The closure of each $G$-orbit in $\mathcal T$ is smooth.
		\item There is a $G$-stable subvariety $\mathcal M$ of $\mathcal T$ and a $G$-equivariant flat map $\mathcal P:\mathcal T\rightarrow\mathcal M$  such that the following holds.
		\begin{enumerate}[label={\rm(\alph*)}]
		\item $\mathcal P$ is a retraction, that is, the restriction of $\mathcal P$ on $\mathcal M$ is the identity map.
		\item $\mathcal M$ is a wonderful $G$-variety, and the  set of the $G$-stable divisors is given by
		\begin{equation}
		    \left\{\,\mathcal P(D_i)\,\big|\,1\leq i\leq l\,\right\}\mathbin{\scaleobj{1.2}{\backslash}} \left\{\mathcal M\right\}\,.
		\end{equation} 
		We call such $\mathcal M$ a {\it home} of $\mathcal T$.
		\end{enumerate}
	\end{enumerate}
	If $\mathcal T^0$ is a $G$-invariant open subvariety of $\mathcal T$, we call $\mathcal T$ a  {\it homeward compactification} of $\mathcal T^0$.
\end{definition}

\begin{definition}[homeward compactification on a roll] \label{gnwond}
$\mathcal T$ is called a {\it homeward variety on a roll} if $\mathcal T$ is a homeward variety with the following additional property.
	\begin{enumerate}[label={\rm(\Alph*)}]
	\setcounter{enumi}{4}
		\item There is an automoprhism  $\sigma$ of $\mathcal T$ such that $\sigma(\mathcal M)\neq\mathcal M$ where $\mathcal M$ is a home of $\mathcal T$.
		\end{enumerate}
\end{definition}
\begin{remark}
It is clear that wonderful varieties and Cartesion products of homeward varieties are homeward, and that a Cartesion product of homeward varieties is a homeward variety on a roll if one of its factor is a homeward variety on a roll.
\end{remark}

Combining Propositions \ref{auto1} and \ref{gwond}, we have that
\begin{theorem}\label{toctt}
$\mathcal T_{s,p,n}$ is a homeward variety.  $\mathcal T_{s,p,n}$ is a homeward variety on a roll if and only if $n=2s$ or $n=2p$.
\end{theorem}

Next, we consider the differential-geometric nature of $\mathcal T_{s,p,n}$ and $\mathcal M_{s,p,n}$.  Recall that a line bundle $L$
on a projective variety $X$ of dimension $n$ is called big if its highest self-intersection number $(L^n)$ is positive, and called numerical effective (or nef) if the intersection number $(L\cdot C)$ is non-negative for any curve $C$ on $X$.

We establish the (semi-)positivity of  the anti-canonical bundles as follows.
\begin{theorem}\label{mfano}
The anti-canonical bundle $-K_{\mathcal M_{s,p,n}}$  of $\mathcal M_{s,p,n}$ is ample.
\end{theorem}
\begin{theorem}\label{fano}
The anti-canonical bundle $-K_{\mathcal T_{s,p,n}}$ of $\mathcal T_{s,p,n}$ is big and numerical effective. $-K_{\mathcal T_{s,p,n}}$ is ample if and only if the rank $r\leq 2$. 
\end{theorem} 
\begin{remark} The Fanoness of $\mathcal M_{p,p,2p}$ was first derived by De Concini and Procesi (\cite{DP}).
\end{remark}

Mabuchi (\cite{Mab})  established  a combinatorial criterion for the existence of K\"ahler-Einstein metrics on toric Fano manifolds in terms of the  barycenter of the moment polytope.  Delcroix  (\cite{De1,De2}) largely generalized it  to spherical Fano  manifolds. 
Applying Delcroix's criterion, we have

\begin{proposition}\label{tke}
There are K\"ahler-Einstein metrics on a Fano manifold $\mathcal T_{s,p,n}$ ($r=1, 2$), if and only if $n=2s$ or $n=2p$ (or equivalently, $\mathcal T_{s,p,n}$ is a homeward compactification on a roll).
\end{proposition}

\begin{proposition}\label{mr13} There are K\"ahler-Einstein metrics on $\mathcal M_{s,p,n}$ if and only if there are K\"ahler-Einstein metrics on $\mathcal M_{p,s,n}$.
\end{proposition} 

\begin{proposition}\label{mr14}  There are K\"ahler-Einstein metrics on $\mathcal M_{s,p,n}$ if $|s+p-n|\gg r$, $|s-p|\gg r$, and there are K\"ahler-Einstein metrics on $\mathcal M_{r,r,2r}$.
\end{proposition} 

We reduce the criterion for $\mathcal M_{s,p,n}$ to  check the positivity of certain integrals over a polytope (see (\ref{criterion})). Based on the reduction, we  verify for $\mathcal M_{s,p,n}$  of small ranks that
\begin{corollary}\label{mr15}
There are K\"ahler-Einstein metrics on $\mathcal M_{p,p,2p}$ when $p\leq 5$,
and $\mathcal M_{s,p,n}$ when $r\leq 2$. 
\end{corollary}

We make the following conjecture.
\begin{conjecture}
There are K\"ahler-Einstein metrics on the manifolds  $\mathcal M_{s,p,n}$.
\end{conjecture}
\begin{remark}
The results for $\mathcal M_{2,2,4}$, $\mathcal M_{3,3,6}$, and $\mathcal T_{s,1,n}$ are derived in \cite{De2,De3}
\end{remark}
\smallskip

In what follows, we generalize  the rational map $\mathcal K_{s,p,n}$ in a natural way. Let  $\overline s=(s_1,\cdots,s_t)$ be an integer-valued vector such that $t\geq 2$ and $\sum_{l=1}^ts_l=n$. Define an index set $\mathbb K^{\overline s}$ by
\begin{equation}\label{kkk}
\mathbb K^{\overline s}:=\left\{(k_1,\cdots,k_t)\in\mathbb Z^t\,\,\rule[-.15in]{0.01in}{.4in}\,0\leq k_1\leq s_1,\,0\leq k_2\leq s_2,\,\cdots,0\leq k_t\leq s_t\,\,{\rm and}\,\,\sum_{i=1}^tk_i=p\right\}.
\end{equation}
Denote the elements of $\mathbb K^{\overline s}$ by $K_0,K_1,\cdots,K_{L_{\overline s}}$ where $L_{\overline s}$ is the cardinal number of $\mathbb K^{\overline s}$ minus $1$. Define index sets  $\mathbb I_{\overline s,p,n}^{K}$ for $K=(k_1,\cdots,k_t)\in\mathbb K^{\overline s}$ by
\begin{equation}\label{ikkk}
\mathbb I_{\overline{s},p,n}^{K}:= \left\{
(i_1,\cdots,i_p)\in\mathbb Z^p\,\,\rule[-.28in]{0.01in}{.65in}\,\,\footnotesize\begin{matrix}
	s_1+\cdots+s_{t-1}+1\leq i_{k_t}<\cdots<i_{1}\leq s_1+s_2+\cdots+s_{t}\,;\\
	s_1+\cdots+s_{t-2}+1\leq i_{k_{t-1}+k_{t}}<\cdots<i_{k_{t}+1}\leq s_1+\cdots+s_{t-1}\,;\\ 
	\vdots\,\\
	1\leq i_{k_1+\cdots+k_{t-1}+k_t}<\cdots<i_{k_2+\cdots+k_{t}+1}\leq s_1
	\end{matrix}\right\}.
\end{equation}
Let  $N_{\overline{s},p,n}^{K}$ be the cardinal number of $\mathbb I_{\overline{s},p,n}^{K}$ minus $1$. 
Similarly to (\ref{FK}), we can define  a projection (rational) map $F_{\overline s}^K$ from $\mathbb {CP}^{N_{p,n}}$ to $\mathbb {CP}^{N_{\overline{s},p,n}^{K}}$ by  dropping the homogeneous coordinates whose indices are not in  $\mathbb I_{\overline{s},p,n}^{K}$. 

Define $\mathcal K_{{\overline s},p,n}:=\left(e, F_{\overline s}^{K_1}\circ e,\cdots,F_{\overline s}^{K_{L^{\overline s}}}\circ e\right)$. More precisely, 
\begin{equation}\label{r}
\begin{split}
\mathcal K_{\overline s,p,n}:&\,\, G(p,n)\dashrightarrow \mathbb {CP}^{N_{p,n}}\times\mathbb {CP}^{N_{\overline{s},p,n}^{K_1}}\times\cdots\times\mathbb {CP}^{N_{\overline{s},p,n}^{K_{L_{\overline s}}}}\\
x&\xdashmapsto{} \big( [\cdots ,P_I(\widetilde x),\cdots]_{I\in\mathbb I_{p,n}},[\cdots,P_I(\widetilde x),\cdots]_{I\in\mathbb I^{K_1}_{\overline s,p,n}},\cdots,[\cdots,P_I(\widetilde x),\cdots]_{I\in\mathbb I_{\overline{s},p,n}^{K_{L_{\overline s}}}}\big)
\end{split}\,.
\end{equation}

\begin{definition}\label{multspn}
    Denote by $\mathcal T_{{\overline s},p,n}$ the closure of the birational image of $G(p,n)$ under $\mathcal K_{\overline s,p,n}$. 
\end{definition}

Notice that in general $\mathcal T_{\overline s,p,n}$  is neither spherical nor smooth; its singularities are defined by equations of a combinatorial nature (see Remark \ref{rmon}). We can define the subvariety $\mathcal M_{\overline s,p,n}$ and the holomorphic map $\mathcal P_{\overline s,p,n}:\mathcal T_{\overline s,p,n}\rightarrow \mathcal M_{\overline s,p,n}$ in a similar way, and would like to ask 
\begin{question}
Is $\mathcal P_{\overline s,p,n}$ a flat map?
\end{question}

Motivated by Sato (\cite{S1,S2}), we end up the paper by constructing compatible embeddings among $\mathcal T_{s,p,n}$ and $\mathcal M_{s,p,n}$.
\medskip

We now briefly describe the organization of the paper and the basic ideas for the proof. To establish the smoothness, we generalize the Van der Waerden representation, which is essentially Gaussian elimination, to give explicit coordinate charts for  $\mathcal T_{s,p,n}$. Through the general theory of Bialynichi-Birula for algebraic $\mathbb C^*$-actions on projective manifolds,  we define submanifolds $\mathcal M_{s,p,n}$ as the source of a natural $\mathbb C^*$-action on $\mathcal T_{s,p,n}$; various fibration structures and isomorphisms of $\mathcal M_{s,p,n}$ follow from an observation that the source and the sink of $\mathcal T_{s,p,n}$ are isomorphic. Viewing $\mathcal M_{s,p,n}$ as a moduli space in the sense of Fujiki, we derive a natural flat map  $\mathcal P_{s,p,n}:\mathcal T_{s,p,n}\rightarrow\mathcal M_{s,p,n}$.  The fundamental work of Brion  on spherical varieties is another main technical tool used extensively in this paper. We determine the cone of effective divisors/curves of $\mathcal T_{s,p,n}$ and $\mathcal M_{s,p,n}$. 
Combined with the fibration structures, the holomorphic automorphism groups  follow. The invariants in Brion's theory can be computed conveniently in the Van der Waerden representation. By Kleiman's ampleness criterion, the (semi-)positivity of the anti-canonical bundles is a consequence of the calculation of the intersection numbers.   

The organization of the paper is as follows. In Section \ref{iterated}, we realize $\mathcal T_{s,p,n}$ as iterated blow-ups of $G(p,n)$. In  Section \ref{group}, we  introduce a group actions on $\mathcal T_{s,p,n}$. In Section \ref{isomt}, we define USD and DUAL isomorphisms of $\mathcal T_{s,p,n}$. We generalize the Van der Waerden representation and prove Propositions \ref{tspnsmooth}, \ref{msmooth} in Section \ref{vander}.
Section \ref{foliation} is devoted to the theory of Bia{\l}ynicki-Birula and its application to  $\mathcal T_{s,p,n}$. We define in  Section \ref{basicp} the subvariety $\mathcal M_{s,p,n}$ and the map $\mathcal P_{s,p,n}:\mathcal T_{s,p,n}\rightarrow\mathcal M_{s,p,n}$;   in Section \ref{basicif}, we introduce the Usd and Dual isomorphisms of $\mathcal M_{s,p,n}$ and the fibration structures on  $\mathcal M_{s,p,n}$;  in Section \ref{basicid}, we prove  Propositions \ref{fiber}, \ref{mwond}, \ref{gwond} based on the geometry of the foliation on $\mathcal T_{s,p,n}$ induced by the $\mathbb C^*$-action. In Section \ref{curanddivg}, we investigate the invariant divisors of $\mathcal T_{s,p,n}$ and $\mathcal M_{s,p,n}$ and give a geometric interpretation; in Section \ref{curanddivi},  we study $T$-invariant curves of $\mathcal T_{s,p,n}$ and compute the intersection numbers. Theorems \ref{mfano},  \ref{fano}, and Propositions \ref{auto1}, \ref{mauto} are proved in Section \ref{curanddivp} and Section \ref{sym} respectively. In Section \ref{ke}, we address   the existence of K\"ahler-Einstein metrics on $\mathcal T_{s,p,n}$ and $\mathcal M_{s,p,n}$ by proving Propositions \ref{tke}, \ref{mr13}, \ref{mr14}. Sections \ref{vps}, \ref{sat} are devoted to the study of the rational map $\mathcal K_{\overline s,p,n}$ and the functoriality of $\mathcal K_{s,p,n}$ respectively. 

We attach five appendices to the paper. In Appendix \ref{section:cover}, we prove that the Van der Waerden representation is a holomorphic atlas, and  define the intermediate Van der Waerden representation.
In Appendix \ref{section:projbc},
we assign certain local coordinate charts for the projective bundles associated with the source and the sink in $G(p,n)$ which are helpful to understand the fibration structures of $\mathcal M_{s,p,n}$ and used in the proof of Lemma \ref{fibpre}. In Appendix \ref{section:rigidmc},
we classify the induced automorphisms of the Picard group of $\mathcal T_{s,p,n}$. In Appendix 
\ref{section:rigidtc}, we determine the automorphism of $\mathcal M_{s,p,n}$ under the assumption that the induced automorphisms of the Picard group is the identiy; we give an alternative proof in Appendix
\ref{section:rigidtc2}.
We include certain numerical calculations for Delcroix's criterion in Appendix \ref{section:nc}.
\medskip

{\noindent\bf Acknowledgement.}  The author appreciates greatly Prof.~Huang for suggesting the topic. He is very grateful to Kexing Chen for writing
Appendix~\ref{section:nc}, to Prof.~Brion and Prof.~Delcroix for answering  questions with patience, and to  Shi Chen,  Xiaocheng Li, Zhan Li, Xu Yang, and Chuyu Zhou for helpful discussions. The author finally thanks Hyeyeong Kim and Yue Ni for leading  wonderful journeys in his life over the past few years.

\section{Basic properties of \texorpdfstring{$\mathcal T_{s,p,n}$}{rr}}\label{basict}

\subsection{Iterated blow-ups}\label{iterated}

In this subsection, we will realize $\mathcal T_{s,p,n}$ as iterated blow-ups of $G(p,n)$, and show that $\mathcal K_{s,p,n}^{-1}$  has a regular extension. 

\begin{definition}\label{sk}
For $0\leq k\leq p$, define a subscheme $S_k\subset G(p,n)$ by
\begin{equation}
S_k:=\big\{ x\in G(p,n)\,\big| P_{I}(\widetilde x)=0\,\,\forall{I\in \mathbb I_{s,p,n}^k}\,\, {\rm where}\,\,\widetilde x\,\,{\rm \,\,is\,\,a\,\,matrix\,\,representative\,\,of}\,\,x \big\}\,. 
\end{equation}
Denote by  $\mathcal S_k\subset \mathcal O_{G(p,n)}$, $0\leq k\leq p$, the ideal sheaf of $S_k$.
\end{definition}

The iterated blow-ups is defined by
\begin{definition}[iterated blow-ups]\label{calt}
Let $\sigma$ be a permutation of $\{0,1,\cdots,p\}$. Let $Y_0^{\sigma}$ be the blow-up of $G(p,n)$ along $S_{\sigma(0)}$; denote the blow-up map by $g_0^{\sigma}:Y^{\sigma}_0\rightarrow G(p,n)$. Inductively, for  $0\leq i\leq p-1$ define $Y^{\sigma}_{i+1}$ to be the blow-up of $Y^{\sigma}_i$ along the  scheme-theoretic inverse image  $(g^{\sigma}_0\circ g^{\sigma}_{1}\circ\cdots\circ g^{\sigma}_i)^{-1}(S_{\sigma(i+1)})$; denote the corresponding blow-up map by $g^{\sigma}_{i+1}:Y^{\sigma}_{i+1}\rightarrow Y^{\sigma}_{i}$. 	Graphically, we have
\vspace{-.08in}
\begin{equation}\label{sblow}
\begin{tikzcd}
&Y^{\sigma}_p\ar{r}{g^{\sigma}_p}&Y^{\sigma}_{p-1}\ar{r}{g^{\sigma}_{p-1}}&\cdots\ar{r}{g^{\sigma}_1}&Y^{\sigma}_0\ar{r}{g^{\sigma}_0}&G(p,n) \\
&&(g^{\sigma}_0\circ\cdots\circ g^{\sigma}_{p-1})^{-1}(S_{\sigma(p)})\ar[hook]{u}&\cdots&(g^{\sigma}_0)^{-1}(S_{\sigma(1)})\ar[hook]{u}&S_{\sigma(0)}\ar[hook]{u}\\
\end{tikzcd}.\vspace{-20pt}
\end{equation}
\end{definition}

We shall show that the construction of $Y_p^{\sigma}$ is independent of the choice of the permutation. 
\begin{lemma}\label{inb}
Let $\sigma_1$ and $\sigma_2$ be  permutations of $\{0,1,\cdots,p\}$. There exists a biholomorphic map $r_{\sigma_1\sigma_2}:Y_p^{\sigma_1}\rightarrow Y_p^{\sigma_2}$ such that the following diagram commutes.
\vspace{-0.05in}
\begin{equation}
\begin{tikzcd}
&Y^{\sigma_1}_{p}\arrow{dr}[swap]{\hspace{-.03in}(g^{\sigma_1}_0\circ\cdots\circ g^{\sigma_1}_{p})}\arrow{rr}{r_{\sigma_1\sigma_2}}&&Y^{\sigma_2} _{p} \arrow{ld}{\hspace{-0.03in}(g^{\sigma_2}_0\circ\cdots\circ g^{\sigma_2}_{p})} \\
&&\,\,G(p,n)\,\,&\\
\end{tikzcd}\vspace{-20pt}\,.
\end{equation}
\end{lemma}

The idea for the proof of Lemma \ref{inb} is to embed the iterated blow-ups into the  blow-ups of the ambient projective space $\mathbb {CP}^{N_{p,n}}$ defined as follows.

Let $L_k$, $0\leq k\leq p$,  be the linear subspace  of $\mathbb {CP}^{N_{p,n}}$   defined by \begin{equation}
L_k:= \big\{[\cdots,z_I,\cdots]\in\mathbb {CP}^{N_{p,n}}\big|z_I=0,\,\,\forall I\in\mathbb I^k_{s,p,n}\big\}; 
\end{equation} or equivalently,  $L_k$ is the locus of indeterminacy of the projection $F^k_s$ defined by (\ref{FK}). 

For each permutation $\sigma$  of $\{0,1,\cdots,p\}$, we can define  iterated blow-ups of  $\mathbb {CP}^{N_{p,n}}$ in the same way as  Definition \ref{calt}.  Let $\widetilde Y_0^{\sigma}$ be the blow-up of  $\mathbb {CP}^{N_{p,n}}$  along $L_{\sigma(0)}$, and denote the blow-up map by $G_0^{\sigma}$. Inductively, for  $0\leq i\leq p-1$ define $\widetilde Y^{\sigma}_{i+1}$ to be the blow-up of $\widetilde Y^{\sigma}_i$ along $(G^{\sigma}_0\circ G^{\sigma}_{1}\circ\cdots\circ G^{\sigma}_i)^{-1}(L_{\sigma(i+1)})$, and denote the corresponding blow-up map by $G^{\sigma}_{i+1}$. Then we have 
\vspace{-.08in}
\begin{equation}\label{lblow}
\small
\begin{tikzcd}
&\widetilde  Y^{\sigma}_p\ar{r}{G^{\sigma}_p}&\widetilde Y^{\sigma}_{p-1}\ar{r}{G^{\sigma}_{p-1}}&\cdots\ar{r}{G^{\sigma}_1}&\widetilde Y^{\sigma}_0\ar{r}{G^{\sigma}_0}&\mathbb {CP}^{N_{p,n}} \\
&&(G^{\sigma}_0\circ\cdots\circ G^{\sigma}_{p-1})^{-1}(L_{\sigma(p)})\ar[hook]{u}&\cdots&(G^{\sigma}_0)^{-1}(L_{\sigma(1)})\ar[hook]{u}&L_{\sigma(0)}\ar[hook]{u}\\
\end{tikzcd}.\vspace{-20pt}
\end{equation}

\begin{lemma}\label{piv}
Let $\sigma_1$ and $\sigma_2$ be  permutations of $\{0,1,\cdots,p\}$. There exists a biholomorphic map  $\widetilde r_{\sigma_1\sigma_2}:\widetilde Y_p^{\sigma_1}\rightarrow\widetilde Y_p^{\sigma_2}$ such that the following diagram commutes.	
\vspace{-0.05in}
\begin{equation}\label{cpin}
\begin{tikzcd}
&\widetilde Y^{\sigma_1}_{p}\arrow{dr}[swap]{\hspace{-.05in}(G^{\sigma_1}_0\circ\cdots\circ G^{\sigma_1}_{p})}\arrow{rr}{\widetilde r_{\sigma_1\sigma_2}}&&\widetilde Y^{\sigma_2} _{p} \arrow{ld}{\hspace{-0.05in}(G^{\sigma_2}_0\circ\cdots\circ G^{\sigma_2}_{p})} \\
&&\,\,\mathbb {CP}^{N_{p,n}}\,\,&\\
\end{tikzcd}\vspace{-20pt}\,.
\end{equation}
\end{lemma}

{\noindent\bf Proof of Lemma \ref{piv}.} It is clear that there exists a unique, natural birational map $\widetilde r_{\sigma_1\sigma_2}:\widetilde Y_p^{\sigma_1}\dashrightarrow\widetilde Y_p^{\sigma_2}$ such that (\ref{cpin}) commutes.	Hence, it suffices to prove Lemma \ref{piv} locally. 

We will   give the Proj construction of the blow-ups in the follows. For each index $I^*\in\mathbb I_{p,n}$, let $\mathcal A_{I^*}\subset\mathbb {CP}^{N_{p,n}}$ be the affine open set defined by $z_{I^*}\neq 0$. Denote the structure ring of $\mathcal A_{I^*}$ by $\mathbb C[\cdots,x_I,\cdots]_{I\in\mathbb I_{p,n}}$, where $x_I:=\frac{z_I}{z_{I^*}}$ is the inhomogeneous coordinate (we make the convention that $x_{I^*}=1$). Then, 
\begin{equation}\label{b0}
\begin{split}
  (G^{\sigma_1}_{0})^{-1}&(\mathcal A_{I^*})
  \cong {\rm Proj}\,\mathbb
 C[\cdots,x_I,\cdots][\cdots,t_I,\cdots]_{I\in\mathbb I_{s,p,n}^{\sigma_1(0)}}\big/
 \big(x_{I^{\prime}}t_{I^{\prime\prime}}-x_{I^{\prime\prime}}t_{I^{\prime}}\big)_{I^{\prime},\,I^{\prime\prime}\in\mathbb I_{s,p,n}^{\sigma_1(0)}}\\
 &={\rm Spec}\,\mathbb C[\cdots,x_I,\cdots]_{I\notin\mathbb I_{s,p,n}^{\sigma_1(0)}}\times{\rm Proj}\,\mathbb
 C\frac{[\cdots,x_I,\cdots]_{I\in\mathbb I_{s,p,n}^{\sigma_1(0)}}[\cdots,t_I,\cdots]_{I\in\mathbb I_{s,p,n}^{\sigma_1(0)}}} { \big(x_{I^{\prime}}t_{I^{\prime\prime}}-x_{I^{\prime\prime}}t_{I^{\prime}}\big)_{I^{\prime},\,I^{\prime\prime}\in\mathbb I_{s,p,n}^{\sigma_1(0)}}}.\\
 \end{split}
\end{equation}
Notice that here we view $\mathbb C\frac{[\cdots,x_I,\cdots]_{I\in\mathbb I_{s,p,n}^{\sigma_1(0)}}[\cdots,t_I,\cdots]_{I\in\mathbb I_{s,p,n}^{\sigma_1(0)}}} { \big(x_{I^{\prime}}t_{I^{\prime\prime}}-x_{I^{\prime\prime}}t_{I^{\prime}}\big)_{I^{\prime},\,I^{\prime\prime}\in\mathbb I_{s,p,n}^{\sigma_1(0)}}}$ as a graded ring such that $x_I$ is of degree zero and $t_I$ is of degree one.

Inductively, we can conclude that $( G^{\sigma_1}_0\circ\cdots\circ G^{\sigma_1}_{p})^{-1}(\mathcal A_{I^*})$ is isomorphic to
\begin{equation}\label{proj1}
\begin{split}
&\prod_{i=0}^p{\rm Proj}\,\mathbb C\frac{[\cdots,x_I,\cdots]_{I\in\mathbb I_{s,p,n}^{\sigma_1(i)}}[\cdots,t_I,\cdots]_{I\in\mathbb I_{s,p,n}^{\sigma_1(i)}}}{ \big(x_{I^{\prime}}t_{I^{\prime\prime}}-x_{I^{\prime\prime}}t_{I^{\prime}}\big)_{I^{\prime},\,I^{\prime\prime}\in\mathbb I_{s,p,n}^{\sigma_1(i)}}}\,\,\,.\\
\end{split}
\end{equation}
By a similar argument, 
$(G^{\sigma_2}_0\circ\cdots\circ G^{\sigma_2}_{p})^{-1}(\mathcal A_{I^*})$ is isomorphic to
\begin{equation}\label{proj2}
\begin{split}
&\prod_{i=0}^p{\rm Proj}\,\mathbb C\frac{[\cdots,x_I,\cdots]_{I\in\mathbb I_{s,p,n}^{\sigma_2(i)}}[\cdots,t_I,\cdots]_{I\in\mathbb I_{s,p,n}^{\sigma_2(i)}}}{ \big(x_{I^{\prime}}t_{I^{\prime\prime}}-x_{I^{\prime\prime}}t_{I^{\prime}}\big)_{I^{\prime},\,I^{\prime\prime}\in\mathbb I_{s,p,n}^{\sigma_2(i)}}}\,\,\,.\\
\end{split}
\end{equation}
It is clear that (\ref{proj1}) and (\ref{proj2}) are isomorphic naturally over  $\mathbb {CP}^{N_{p,n}}$.

We complete the proof of Lemma \ref{piv}.\,\,\,\,\,$\endpf$
\medskip

\begin{lemma}\label{res}
For $0\leq k\leq p$, $Y^{\sigma}_k$ is the strict transformation of $G(p,n)$ under the blow-up $G^{\sigma}_0\circ \cdots\circ G^{\sigma}_{k}:\widetilde Y^{\sigma}_k\rightarrow \mathbb {CP}^{N_{p,n}}$; or equivalently, there is a commutative diagram as follows.
\vspace{-.09in}
\begin{equation}\label{fucgG}
\begin{tikzcd}
&Y^{\sigma}_p\ar[hook]{d}\ar{r}{g^{\sigma}_p}&Y^{\sigma}_{p-1}\ar[hook]{d}\ar{r}{g^{\sigma}_{p-1}}&\cdots\ar{r}{g^{\sigma}_1}&Y^{\sigma}_0\ar[hook]{d}\ar{r}{g^{\sigma}_0}&G(p,n)\ar[hook]{d}{e} \\
&\widetilde Y^{\sigma}_p\ar{r}{G^{\sigma}_p}&\widetilde Y^{\sigma}_{p-1}\ar{r}{G^{\sigma}_{p-1}}&\cdots\ar{r}{G^{\sigma}_{1}}&\widetilde Y^{\sigma}_0\ar{r}{G^{\sigma}_{0}}&\mathbb{CP}^{N_{p,n}}\\
\end{tikzcd}.\vspace{-20pt}
\end{equation}
\end{lemma}
{\noindent\bf Proof of  Lemma \ref{res}.}  Recall the blow-up closure lemma as follows  (see (22.2.5) in \cite{V}). Consider a fibered diagram
\begin{equation}\label{bcl}
\begin{tikzcd}
&W\arrow[hook]{d}\arrow[hook]{r}&Z\arrow[hook]{d} \\
&\widetilde W\arrow[hook]{r}&\widetilde Z\\
\end{tikzcd}\,,\vspace{-20pt}
\end{equation}
where $W,Z,\widetilde W,\widetilde Z$ are projective schemes and the maps are closed embeddings. Then the strict transformation of $Z$ under the blow-up of $\widetilde Z$ along $\widetilde W$ is isomorphic to the blow-up of $Z$ along $W$. 

We  prove Lemma \ref{res} by an induction on $k$. By definition $S_k$, $0\leq k\leq p$, is  the scheme-theoretic intersection of $L_k$ and $G(p,n)$. Then  we have the following  fibered diagram.
\vspace{-0.05in}
\begin{equation}
\begin{tikzcd}
&S_{\sigma(0)}\arrow[hook]{d}\arrow[hook]{r}&G(p,n)\ar[hook]{d} \\
&L_{\sigma(0)}\arrow[hook]{r}&\mathbb{CP}^{N_{p,n}}\\
\end{tikzcd}.\vspace{-20pt}
\end{equation}
Applying the blow-up closure lemma, we can derive the following  commutative diagram and thus complete the proof when $k=0$.
\vspace{-0.06in}
\begin{equation}
\begin{tikzcd}
&Y^{\sigma}_0\ar[hook]{d}\ar{r}&G(p,n)\ar[hook]{d} \\
&\widetilde Y^{\sigma}_0\ar{r}&\mathbb{CP}^{N_{p,n}}\\
\end{tikzcd}.\vspace{-20pt}
\end{equation}

Suppose  that Lemma \ref{res}  is true for $k=m<p$. Consider the case $k=m+1$.
By assumption we have the following 
commutative diagram.
\vspace{-0.06in}
\begin{equation}\label{cf}
\begin{tikzcd}
&(g^{\sigma}_0\circ \cdots\circ g^{\sigma}_{m})^{-1}(S_{\sigma(m+1)})\ar[hook]{r}&Y^{\sigma}_m\ar[hook]{d}\ar{r}{g^{\sigma}_0\circ \cdots\circ g^{\sigma}_{m}}&[2em]G(p,n)\ar[hook]{d} \\
&(G^{\sigma}_0\circ \cdots\circ G^{\sigma}_{m})^{-1}(L_{\sigma(m+1)})\ar[hook]{r}&\widetilde Y^{\sigma}_m\ar{r}{G^{\sigma}_0\circ \cdots\circ G^{\sigma}_{m}}&\mathbb{CP}^{N_{p,n}}\\
\end{tikzcd}\,.\vspace{-20pt}
\end{equation}
We will complete ($\ref{cf}$)
and derive a fibered diagram as follows.
\vspace{-0.06in}
\begin{equation}\label{cf2}
\begin{tikzcd}
&(g^{\sigma}_0\circ \cdots\circ g^{\sigma}_{m})^{-1}(S_{\sigma(m+1)})\arrow[hook]{d}\arrow[hook]{r}&Y^{\sigma}_m \arrow[hook]{d} \\
&(G^{\sigma}_0\circ \cdots\circ G^{\sigma}_{m})^{-1}(L_{\sigma(m+1)})\arrow[hook]{r}&\widetilde Y^{\sigma}_m\\
\end{tikzcd}\,.\vspace{-20pt}
\end{equation}
It suffices to compute locally. Let $A,A/I,B,C$ are the structure rings of $\mathbb{CP}^{N_{p,n}}, G(p,n), \widetilde Y^{\sigma}_m, Y^{\sigma}_m$,  in certain affine subsets, respectively. Let $a_1,\cdots, a_{{N^{\sigma(m+1)}_{s,p,n}+1}}$ be the elements of $A$ corresponding to $z_I$, $I\in\mathbb I_{s,p,n}^{\sigma(m+1)}$, and  $\pi({a_1}),\cdots, \pi(a_{{N^{\sigma(m+1)}_{s,p,n}}+1})$ the images of  $a_1,\cdots, a_{{N^{\sigma(m+1)}_{s,p,n}}+1}$ under the quotient map $\,\pi:A\rightarrow A/I$. Dual to (\ref{cf}) we have the following commutative diagram of structure rings.
\begin{equation}\label{scf}
\begin{tikzcd}
&\bigslant{C}{\left(h(\pi({a_1})),\cdots,h\big(\pi({a_{{N^{\sigma(m+1)}_{s,p,n}}+1}})\big)\right)}&C\ar{l}&A/I\ar{l}[swap]{h} \\
&\bigslant{B}{\left(\widetilde h(a_1),\cdots,\widetilde h(a_{{N^{\sigma(m+1)}_{s,p,n}}+1})\right)}\ar[dashed]{u}{\mathbin{\scaleobj{1.4}{\iota}}\,\,}&B\ar{l}\ar{u}&A\ar{l}[swap]{\widetilde h}\ar{u}{\mathbin{\scaleobj{1.4}{\pi}}}\\
\end{tikzcd}.\vspace{-20pt}
\end{equation}
It is clear that there exists a natural map $\iota$ completing diagram (\ref{scf}), and, moreover,
\begin{equation}
C\otimes_B\bigslant{B}{\left(\widetilde h(a_1),\cdots,\widetilde h(a_{{N^{\sigma(m+1)}_{s,p,n}}+1})\right)}=\bigslant{C}{\left(h(\pi({a_1})),\cdots,h\big(\pi({a_{{N^{\sigma(m+1)}_{s,p,n}}+1}})\big)\right)}.
\end{equation}
By applying the blow-up closure lemma to (\ref{cf2}), we conclude  Lemma \ref{res} for $k=m+1$. 

This completes the proof of Lemma \ref{res}.
$\endpf$
\medskip

{\noindent\bf Proof of  Lemma  \ref{inb}.} This is a direct consequence of Lemmas \ref{piv} and \ref{res}.\,\,\,\,\,$\endpf$

\medskip

Next we will relate the iterated blow-up to the rational map $\mathcal K_{s,p,n}$ and the canonical blow-up $\mathcal T_{s,p,n}$. For convenience, we  set $\sigma$ to  the identity permutation and omit the supscrit $\sigma$ in the following.
\begin{lemma}\label{cpc}  
There is an isomorphism $\nu:Y_p\rightarrow \mathcal T_{s,p,n}$ such that the following diagram commutes.
\vspace{-0.05in}
\begin{equation}
\begin{tikzcd}
&Y_{p}\arrow{dr}[swap]{\hspace{-.03in}(g_0\circ\cdots\circ g_{p})}\arrow{rr}{\nu}&&\mathcal T_{s,p,n} \arrow[dashed]{ld}{\hspace{-0.03in}\mathcal K^{-1}_{s,p,n}} \\
&&\,\,G(p,n)\,\,&\\
\end{tikzcd}\vspace{-20pt}\,.
\end{equation}
\end{lemma}
{\noindent\bf Proof of  Lemma \ref{cpc}.} Denote by $\mathcal Y_{k}$  the blow-up of $\mathbb {CP}^{N_{p,n}}\times\mathbb {CP}^{N^0_{s,p,n}}\times\cdots\times\mathbb {CP}^{N^{k-1}_{s,p,n}}$ along  $L_{k}\times\mathbb {CP}^{N^0_{s,p,n}}\times\cdots\times\mathbb {CP}^{N^{k-1}_{s,p,n}}$ and $\mathcal G_k$ the blow-up map for $0\leq k\leq p$. Define a rational map  $K_{k}:\mathbb {CP}^{N_{p,n}}\times\mathbb {CP}^{N^0_{s,p,n}}\times\cdots\times\mathbb {CP}^{N^{k-1}_{s,p,n}}\dashrightarrow\mathbb {CP}^{N_{p,n}}\times\mathbb {CP} ^{N^0_{s,p,n}}\times\cdots\times\mathbb {CP}^{N^{k}_{s,p, n}}$
by
\begin{equation}\label{Kk2}
\begin{split}
&K_{k}\big([\cdots,z_I,\cdots]_{I\in\mathbb I_{p,n}}\times[\cdots,t_I,\cdots]_{I\in\mathbb I_{s,p,n}^0}\times\cdots\times[\cdots,t_I,\cdots]_{I\in\mathbb I_{s,p,n}^{k-1}}\big)\\
=[\cdots&,z_I,\cdots]_{I\in\mathbb I_{p,n}}\times[\cdots,t_I,\cdots]_{I\in\mathbb I_{s,p,n}^0}\times\cdots\times[\cdots,t_I,\cdots]_{I\in\mathbb I_{s,p,n}^{k-1}}\times[\cdots,z_I,\cdots]_{I\in\mathbb I_{s,p,n}^{k}}
\end{split};
\end{equation} 
that is $K_{k}:=\left({\rm id}\big|_{\mathbb {CP}^{N_{p,n}}}\,,\,\, {\rm id}\big|_{\mathbb {CP}^{N^0_{s,p,n}}}\,,\,\, \cdots,\,\, {\rm id}\big|_{\mathbb {CP}^{N^{k-1}_{s,p,n}}}\,,\,\, F_s^k\right)$.
Denote by $pr_k$ the projection map 
\begin{equation}
pr_k:\mathbb {CP}^{N_{p,n}}\times\mathbb {CP}^{N^0_{s,p,n}}\times\cdots\times\mathbb {CP}^{N^k_{s,p,n}}\rightarrow \mathbb {CP}^{N_{p,n}}\times\mathbb {CP}^{N^0_{s,p,n}}\times\cdots\times\mathbb {CP}^{N^{k-1}_{s,p,n}}\,.    
\end{equation}

Similarly to Lemma \ref{res}, to prove Lemma \ref{cpc} it suffices to embed $\mathcal Y_k$ into $\mathbb {CP}^{N_{p,n}}\times\mathbb {CP}^{N^0_{s,p,n}}\times\cdots\times\mathbb {CP}^{N^k_{s,p,n}}$, $0\leq k\leq p$,  such that $\mathcal G_k=pr_k\big|_{\mathcal Y_k}$ and  $(K_{k}\circ {\mathcal G}_k)(x)=x$ for a generic point $x\in\mathcal Y_k$.

We will prove by induction. By (\ref{b0}) $\mathcal Y_0$ is isomorphic to the subvariety of $\mathbb{CP}^{N_{p,n}}\times \mathbb{CP}^{N_{s,p,n}^0}$ defined by the ideal generated by $\big\{z_{I^{\prime}}t_{I^{\prime\prime}}-z_{I^{\prime\prime}}t_{I^{\prime}}\big\}_{I^{\prime},\,\,I^{\prime\prime}\in\mathbb I_{s,p,n}^{0}}$, where $[\cdots,t_I,\cdots]_{I\in\mathbb I_{s,p,n}^{0}}$ and  $[\cdots, z_I, \cdots]_{I\in\mathbb I_{p,n}}$  are the homogeneous coordinates for $\mathbb{CP}^{N_{s,p,n}^0}$ and $\mathbb{CP}^{N_{p,n}}$ respectively. Therefore, we can embed $\mathcal Y_0$ into $\mathbb{CP}^{N_{p,n}}\times \mathbb{CP}^{N_{s,p,n}^0}$ by identifying it with the above subvariety. Moreover, by computing in an open set $\mathcal A_{I_0^{\prime}}\subset\mathbb{CP}^{N_{p,n}}$ defined by  $z_{I_0^{\prime}}\neq 0$ where $I_0^{\prime}=(p,p-1,\cdots, 1)$, it is clear that  $(K_{0}\circ \mathcal G_0)(x)=x$ for a generic point $x\in\mathcal Y_0$.

Assume the induction hypothesis for $k=m\leq p-1$. Similar to the case $k=0$, we can embed $\mathcal Y_{m+1}$  into $\mathbb {CP}^{N_{p,n}}\times\mathbb {CP}^{N^0_{s,p,n}}\times\cdots\times\mathbb {CP}^{N^{m+1}_{s,p,n}}$.
Let $\mathcal A_{I^{\prime}_{m+1}}\subset\mathbb{CP}^{N_{p,n}}$ be an open set defined by  $z_{I^{\prime}_{m+1}}\neq 0$, where $I^{\prime}_{m+1}=(s+m+1,s+m,\cdots,s+m+2-p)\in\mathbb I_{s,p,n}^{m+1}$.  Computing in the scheme-theoretic inverse image of $\mathcal A_{I^{\prime}_k}$ under the map $(pr_{m+1}\circ pr_{m}\circ\cdots\circ pr_{0})$,  we can verify that $(K_{m+1}\circ \mathcal G_{m+1})(x)=x$ for a generic point $x\in\mathcal  Y_{m+1}$.

We complete the proof of Lemma \ref{cpc}.
\,\,\,\,$\endpf$

\begin{definitionlemma}\label{rspn}
Denote by $R_{s,p,n}:\mathcal T_{s,p,n}\rightarrow G(p,n)$ the regular extension of $\mathcal K^{-1}_{s,p,n}$\,. We call $R_{s,p,n}:\mathcal T_{s,p,n}\rightarrow G(p,n)$ {\it  the canonical blow-up map} or {\it  the canonical blow-up} for short.
\end{definitionlemma}
{\noindent\bf Proof of Definition-Lemma \ref{rspn}.} By Lemma \ref{cpc}, we can embed the iterated blow-up $Y_p$ into $\mathbb {CP}^{N_{p,n}}\times\mathbb {CP}^{N^0_{s,p,n}}\times\cdots\times\mathbb {CP}^{N^{p}_{s,p,n}}$. Then $\mathcal K^{-1}_{s,p,n}$ has a regular extension which is isomorphic to the blow-up map $(g_0\circ\cdots\circ g_{p})$.\,\,\,\,$\endpf$
\begin{remark}
It is clear that $R_{s,p,n}$ is induced by the projection 
of $\mathcal T_{s,p,n}$ to the first factor $\mathbb {CP}^{N_{p,n}}$ of the ambient space $\mathbb {CP}^{N_{p,n}}\times\mathbb {CP}^{N^0_{s,p,n}}\times\cdots\times\mathbb {CP}^{N^p_{s,p,n}}$.
\end{remark}

\subsection{Extended group actions  on \texorpdfstring{$\mathcal T_{s,p,n}$}{dd}} \label{group}

In this subsection, we will lift certain group actions from $G(p,n)$ to $\mathcal T_{s,p,n}$.
 
Recall that the general linear group $GL(n,\mathbb C)$ has a natural action on the Grassmannian $G(p,n)$  by  the matrix multiplication from the right. The subgroup $GL(s,\mathbb C)\times GL(n-s,\mathbb C)$ in (\ref{sns}) has an induced action on $G(p,n)$; we denote it by $\delta_{s,p,n}$. 

We shall lift the above $GL(s,\mathbb C)\times GL(n-s,\mathbb C)$-action to $\mathcal T_{s,p,n}$ in the following.

Recall that the $PGL(n,\mathbb C)$-action of the Grassmannian $G(p,n)$ is induced by the projective linear transformation of the ambient complex projective space $\mathbb {CP}^{N_{p,n}}$.  Moreover, we can define a linear representation $\rho:GL(n,\mathbb C)\rightarrow GL(N_{p,n}+1,\mathbb C)$ through the Pl\"ucker coordinates $P_I$ as follows. For each index $I=(i_1,\cdots,i_p)\in\mathbb I_{p,n}$  denote by $E_I$ the $p\times n$ matrix such that its submatrix consisting of the $i_1^{th}, \cdots,i_p^{th}$ columns is the identity matrix, and that all other columns are zero. Let $(\cdots,z_I,\cdots)_{I\in\mathbb I_{p,n}}$ be the coordinates for the complex Euclidean space $\mathbb C^{N_{p,n}+1}$; define  vectors $e_I\in\mathbb C^{N_{p,n}+1}$ by $z_I=1$ and $z_{I^{\prime}}=0$ for all $I^{\prime}\in \mathbb I_{p,n}$ and $I^{\prime}\neq I$. Define
\begin{equation}
    \rho(g)\left((\cdots,z_I,\cdots)_{\mathbb I_{p,n}}\right):=\sum_{\widetilde I\in\mathbb I_{p,n}}\left(\sum_{I\in\mathbb I_{p,n}}z_I\cdot P_{\widetilde I}\left(E_I\cdot g\right)\right)e_{\widetilde I}\,.
\end{equation}
Here $E_I\cdot g$ is the matrix multiplication.

Denote by $\widetilde \rho$ the restriction of $\rho$ to the subgroup $GL(s,\mathbb C)\times GL(n-s,\mathbb C)$. Define linear subspaces $M_k$, $0\leq k\leq p$, of $\mathbb C^{N_{p,n}+1}$ as follows  (see (\ref{subl}) as well). 
\begin{equation}
 M_k:=   \left\{(\cdots ,z_I,\cdots)_{I\in\mathbb I_{p,n}}\in\mathbb {CP}^{N_{p,n}}\big|z_I=0\,,\,\,\forall I\notin\mathbb I_{s,p,n}^k \right\}\,.
\end{equation}
It is clear that $\oplus_{k=0}^pM_k=\mathbb C^{N_{p,n}+1}$ and that $M_k$, $0\leq k\leq p$, are invariant subspaces of  the representation $\widetilde \rho$. Then $\widetilde \rho$ has a decomposition
\begin{equation}
    \widetilde \rho=\rho_0+\rho_1+\cdots+\rho_p\,,
\end{equation}
where $\rho_k$,$0\leq k\leq p$,  is a representation as follows.
\begin{equation} \rho_k:GL(s,\mathbb C)\times GL(n-s,\mathbb C)\rightarrow GL(M_k)\cong GL(N^k_{s,p,n}+1,\mathbb C)\,.
\end{equation} 
\begin{lemma}\label{gg}
		There exists a unique $GL(s,\mathbb C)\times GL(n-s,\mathbb C)$-action $\Delta_{s,p,n}$ on $\mathcal T_{s,p,n}$ such that the following diagram commutes for each element $g\in GL(s,\mathbb C)\times GL(n-s,\mathbb C)$.
	\begin{equation}\label{t1}
	\begin{tikzcd} &\mathcal T_{s,p,n}  \arrow{r}{\Delta_{s,p,n}(g)}&[2em]\mathcal T_{s,p,n}\arrow{d}{{R_{s,p,n}}} \\ &G(p,n)\arrow[leftarrow]{u}{R_{s,p,n}} \arrow{r}{\delta_{s,p,n}(g)}&G(p,n)\\ \end{tikzcd}.\vspace{-20pt} \end{equation}
\end{lemma}
{\noindent\bf Proof of Lemma \ref{gg}.}
The uniqueness follows from the fact that $R_{s,p,n}$ is a birational map. 

To prove the existence, we define a $GL(s,\mathbb C)\times GL(n-s,\mathbb C)$-action $\widehat\Delta_{s,p,n}$ on the ambient space $G(p,n)\times\mathbb {CP}^{N^0_{s,p,n}}\times\cdots\times\mathbb {CP}^{N^p_{s,p,n}}$ as follows. For each $g\in GL(s,\mathbb C)\times GL(n-s,\mathbb C)$,
\begin{equation}\label{p3} \begin{split} \widehat\Delta_{s,p,n}(g)&\left(x,[\cdots,t_I,\cdots]_{I\in\mathbb I_{s,p,n}^0}, \cdots, [\cdots,t_I,\cdots]_{I\in\mathbb I_{s,p,n}^p}\right)\\ &:=\left(\delta_{s,p,n}(g)(x), \rho_0(g)\big([\cdots,t_I,\cdots]_{I\in\mathbb I_{s,p,n}^0}\big), \cdots, \rho_p(g)\left([\cdots,t_I,\cdots]_{I\in\mathbb I_{s,p,n}^p}\right)\right). \end{split} \end{equation} 

Next, we will show that $\mathcal T_{s,p,n}$ is invarinat under $\widehat\Delta_{s,p,n}$.  By the definition of $\rho_k$, we can verify that for a generic $x\in G(p,n)$
\begin{equation}\label{p4} 
\begin{split} \mathcal K_{s,p,n} &\left(\delta_{s,p,n}(g)(x)\right)=\left(\delta_{s,p,n}(g)(x),[\cdots, P_I(\delta_{s,p,n}(g)(x)),\cdots]_{I\in\mathbb I_{s,p,n}^0},\cdots,\right.\\
&\,\,\,\,\,\,\,\,\,\,\,\,\,\,\,\,\,\,\,\,\,\,\,\,\,\,\,\,\,\,\,\,\,\,\,\,\,\,\,\,\,\,\,\,\,\,\,\,\,\,\,\,\,\,\,\,\,\,\,\,\,\,\,\,\left.[\cdots,P_I(\delta_{s,p,n}(g)(x)),\cdots]_{I\in\mathbb I_{s,p,n}^p}\right)\\ 
&=\left(\delta_{s,p,n}(g)(x), \rho_0(g)([\cdots,P_I(x),\cdots]_{I\in\mathbb I_{s,p,n}^0}), \cdots, \rho_p(g)(\left[\cdots,P_I(x),\cdots\right]_{I\in\mathbb I_{s,p,n}^p})\right)\\ &=\widehat\Delta_{s,p,n}(g)\left(\mathcal K_{s,p,n}(x)\right)\,.\end{split} \end{equation} Therefore, $\widehat\Delta_{s,p,n}(g)\left(\mathcal T_{s,p,n}\right)=\mathcal T_{s,p,n}$ for each $g\in GL(s,\mathbb C)\times GL(n-s,\mathbb C)$. 

We complete the proof of Lemma \ref{gg} by defining $\Delta_{s,p,n}:=\widehat\Delta_{s,p,n}\big|_{\mathcal T_{s,p,n}}$.\,\,\,$\endpf$
\medskip

Define an algebraic $\mathbb C^*$-action $\psi_{s,p,n}$ on $G(p,n)$  by
\begin{equation}\label{grac} \psi_{s,p,n}(\lambda):=\left(\begin{matrix} I_{s\times s}&0\\ 0&\lambda\cdot I_{(n-s)\times (n-s)}\\ \end{matrix} \right),\,\,\lambda\in\mathbb C^*.
\end{equation}
By Lemma \ref{gg}, we have a unique lifting $\Psi_{s,p,n}$ of $\psi_{s,p,n}$ from $G(p,n)$ to $\mathcal T_{s,p,n}$. 

\begin{lemma}\label{exc}
Let $\psi_{s,p,n}$ be  the $\mathbb C^*$-action in (\ref{grac}). There is a unique equivariant $\mathbb C^*$-action $\Psi_{s,p,n}$ on $\mathcal T_{s,p,n}$ such that the following diagram commutes for all $\lambda\in\mathbb C^*$.
\vspace{-.03in}
\begin{equation}\label{t2}
\begin{tikzcd} &\mathcal T_{s,p,n}  \arrow{r}{\Psi_{s,p,n}(\lambda)}&[2em]\mathcal T_{s,p,n}\arrow{d}{{R_{s,p,n}}} \\ &G(p,n)\arrow[leftarrow]{u}{R_{s,p,n}} \arrow{r}{\psi_{s,p,n}(\lambda)}&G(p,n)\\ \end{tikzcd}.\vspace{-20pt} \end{equation}
\end{lemma}
{\noindent\bf Proof of Lemma \ref{exc}.} Define an  $\mathbb C^*$-action $\widehat{\Psi}_{s,p,n}$ on $G(p,n)\times\mathbb {CP}^{N^0_{s,p,n}}\times\cdots\times\mathbb {CP}^{N^m_{s,p,n}}$ by \begin{equation}
 \begin{split}
 \widehat{\Psi}_{s,p,n}(\lambda)&\big(x,[\cdots,t_I,\cdots]_{I\in\mathbb I_{s,p,n}^0}, \cdots, [\cdots,t_I,\cdots]_{I\in\mathbb I_{s,p,n}^p}\big)\\
 &=\big(\psi_{s,p,n}(\lambda)(x),[\cdots,t_I,\cdots]_{I\in\mathbb I_{s,p,n}^0}, \cdots, \left[\cdots,t_I,\cdots\right]_{I\in\mathbb I_{s,p,n}^p}\big),\,\,\,\, \lambda\in\mathbb C^*.
 \end{split}
 \end{equation} 
 
We can complete the proof by setting ${\Psi}_{s,p,n}:=\widehat{\Psi}_{s,p,n}\big|_{\mathcal T_{s,p,n}}$.\,\,\,$\endpf$

\subsection{Isomorphisms among \texorpdfstring{$\mathcal T_{s,p,n}$}{rr}} \label{isomt}
In this subsection, we will introduce certain isomorphisms among $\mathcal T_{s,p,n}$ of various parameters.

Recall that every $p$-dimensional subspace $W$ of an $n$-dimensional space  $V$ determines an $(n-p)$-dimensional quotient space $V/W$ of $V$. 
Taking the dual yields an inclusion of $(V/W)^*$ in $V^*$;  in terms of the Grassmannian, this gives a canonical isomorphism
\begin{equation}\label{dual}
  {\empty}^*:\,\,\,\,Gr(p, V)\cong Gr(n-p, V^*).  
\end{equation}

In the following, we will extend (\ref{dual}) to an isomorphism between $\mathcal T_{s,p,n}$ and $\mathcal T_{s,n-p,n}$. For each point $x\in Gr(p, V)\cong G(p,n)$, let $\widetilde x$ be a matrix representative of $x$, and $\widetilde{x^*}$ a matrix representative of $x^*$, where $x^*\in Gr(n-p, V^*)\cong G(n-p,n)$ is the image of $x$ under ${\empty}^*$. Then by definition,
\begin{equation}
    \widetilde x\cdot \left(\widetilde{x^*}\right)^T=0\,,
\end{equation}
where  $\left(\widetilde{x^*}\right)^T$ is the transpose of the matrix $\widetilde {x^*}$,  $0$ is the $p\times (n-p)$ zero matrix, and the product is the matrix multiplication. 
Hence, in terms of matrix representatives, we  can write   (\ref{dual})  explicitly as follows.
\begin{equation}\label{gdual}
\begin{split}
{\empty}^*:\,\,\,\,G(p, n)&\cong G(n-p, n)\\
\widetilde x&\mapsto \widetilde {x^*}
\end{split}\,\,,\,\,\,\,\,\,\,\,\,\,\,\,\widetilde x\cdot \left(\widetilde{x^*}\right)^T=0\,.
\end{equation}

For $0\leq k\leq p$ and  $I=(i_1,\cdots, i_p)\in\mathbb I^k_{s,p,n}$, define $I^*=(i^*_1,\cdots,i^*_{n-p})\in\mathbb I^{s-k}_{s,n-p,n}$ such that
\begin{equation}
    \left\{i_1,\cdots, i_p, i^*_1,\cdots,i^*_{n-p}\right\}=\left\{1,2,\cdots,n \right\}.
\end{equation}
Denote by $\left[\cdots,z^k_I,\cdots\right]_{I\in\mathbb I^k_{s,p,n}}$ and $\left[\cdots,z^{s-k}_{I^*},\cdots\right]_{I^*\in\mathbb I^{s-k}_{s,n-p,n}}$ the homogeneous coordinates for $\mathbb {CP}^{N^{k}_{s,p,n}}$ and  $\mathbb {CP}^{N^{s-k}_{s,n-p,n}}$ respectively.
Define an isomorphism $g_k:\mathbb {CP}^{N^{k}_{s,p,n}}\rightarrow\mathbb {CP}^{N^{s-k}_{s,n-p,n}}$ by 
\begin{equation}
\begin{split}
\left[\cdots,z^{s-k}_{I^*},\cdots\right]_{I^*\in\mathbb I^{s-k}_{s,n-p,n}}&=g_k\left(\left[\cdots,z^k_I,\cdots\right]_{I\in\mathbb I^k_{s,p,n}}\right)=\big[\cdots,\underset{\substack{\uparrow\\z_{I^*}}}{(-1)^{\sigma(I)}\cdot z^k_{I},}\cdots\big]_{I^*\in\mathbb I^{s-k}_{s,n-p,n}}\,,
\end{split}
\end{equation}
where $\sigma(I)$ is the signature of the following permutation.
\begin{equation}
     \left(\begin{matrix}
		i_1&i_2&\cdots&i_p&i^*_1&i^*_2&\cdots&i^*_{n-p}\\
		1&2&\cdots&p&p+1&p+2&\cdots&n\\
	\end{matrix}\right).
\end{equation}

Based on the above construction,  we will define an isomorphism 
\begin{equation}
    \widetilde{{\rm DUAL}}:\,\,G(p,n)\times\mathbb {CP}^{N^{0}_{s,p,n}}\times\cdots\times\mathbb {CP}^{N^p_{s,p,n}}\rightarrow G(n-p,n)\times\mathbb {CP}^{N^{0}_{s,n-p,n}}\times\cdots\times\mathbb {CP}^{N^{n-p}_{s,n-p,n}}\,.
\end{equation}
Let $a=\left(x,\left[\cdots,z^0_I,\cdots\right]_{I\in\mathbb I^0_{s,p,n}}, \cdots, \left[\cdots,z^p_I,\cdots\right]_{I\in\mathbb I^p_{s,p,n}}\right)\in G(p,n)\times\mathbb {CP}^{N^{0}_{s,p,n}}\times\cdots\times\mathbb {CP}^{N^p_{s,p,n}}$ be an arbitrary point. We define the $G(n-p,n)$-component of $\widetilde{{\rm DUAL}}(a)$ to be $x^*$ the image under the dual map of the Grassmannian. Recall that when $n-s<k\leq n-p$ or $0\leq k<n-s-p$, $\mathbb I^{k}_{s,n-p,n}=\emptyset$ and hence the projective space $\mathbb {CP}^{N^k_{s,n-p,n}}$ consists of a single point $c_k$ by convention; define the $\mathbb {CP}^{N^k_{s,n-p,n}}$-component of $\widetilde{{\rm DUAL}}(a)$ to be $c_k$.
Define the $\mathbb {CP}^{N^k_{s,n-p,n}}$-component of $\widetilde{{\rm DUAL}}(a)$ by
\begin{equation}
    \left[\cdots,z^{k}_{I^*}\left(\widetilde{{\rm DUAL}}(a)\right),\cdots\right]_{I^*\in\mathbb I^{k}_{s,n-p,n}}=g_{s-k}\left(\left[\cdots,z^{s-k}_I,\cdots\right]_{I\in\mathbb I^{s-k}_{s,p,n}}\right)\,,
\end{equation}
when $\max\{0,n-s-p\}\leq k\leq \min\{n-s,n-p\}$.

\begin{definitionlemma}\label{tiso}
The restriction of $\widetilde{{\rm DUAL}}$  induces an isomorphism DUAL from  $\mathcal T_{s,p,n}$ to $\mathcal T_{s,n-p,n}$ such that the following diagram commutes.
\vspace{-.03in}
\begin{equation}\label{DUAL}
\begin{tikzcd} &\mathcal T_{s,p,n}\arrow{r}{\rm DUAL}&\mathcal T_{s,n-p,n}\arrow{d}{R_{s,n-p,n}} \\ &G(p,n) \arrow{r}{{\empty}^*}\arrow[leftarrow]{u}{R_{s,p,n}}&G(n-p,n)\\ \end{tikzcd}\,\,.\vspace{-20pt} \end{equation}
\end{definitionlemma}
{\bf\noindent Proof of Lemma \ref{tiso}.}
It suffices to prove that the image of $\mathcal T_{s,p,n}$ under $\widetilde{{\rm DUAL}}$ is contained in $\mathcal T_{s,n-p,n}$. 

Consider the following open set of $G(p,n)$,
\begin{equation}
   U:=\left\{\left.  \left(I_{p\times p} \hspace{-0.13in}\begin{matrix}
  &\hfill\tikzmark{c1}\\
  &\hfill\tikzmark{d1}
  \end{matrix}\,\,\, A \right)\right\vert_{}
   A\,\,{\rm  is\,\, a\,\,} p\times (n-p)\,\,{\rm matrix}\,\right\}.
  \tikz[remember picture,overlay]   \draw[dashed,dash pattern={on 4pt off 2pt}] ([xshift=0.5\tabcolsep,yshift=7pt]c1.north) -- ([xshift=0.5\tabcolsep,yshift=-2pt]d1.south);
\end{equation}
Denote by  $U^*$ the image of $U$ under the dual map ${\empty}^*$. Then,
\begin{equation}
   U^*=\left\{\left.  \left(-A^T \hspace{-0.13in}\begin{matrix}
  &\hfill\tikzmark{c2}\\
  &\hfill\tikzmark{d2}
  \end{matrix}\,\,\, I_{(n-p)\times (n-p)} \right)\right\vert_{}
   A\,\,{\rm  is\,\, a\,\,} p\times (n-p)\,\,{\rm matrix}\,\right\};
  \tikz[remember picture,overlay]   \draw[dashed,dash pattern={on 4pt off 2pt}] ([xshift=0.5\tabcolsep,yshift=7pt]c2.north) -- ([xshift=0.5\tabcolsep,yshift=-2pt]d2.south);
\end{equation}
moreover, the dual map takes the following form in terms of the matrix representatives,
\begin{equation}\label{dtrans}
  { \left(I_{p\times p} \hspace{-0.13in}\begin{matrix}
  &\hfill\tikzmark{c1}\\
  &\hfill\tikzmark{d1}
  \end{matrix}\,\,\, A \right)}^*=\left(-A^T \hspace{-0.13in}\begin{matrix}
  &\hfill\tikzmark{c2}\\
  &\hfill\tikzmark{d2}
  \end{matrix}\,\,\, I_{(n-p)\times (n-p)} \right).
  \tikz[remember picture,overlay]   \draw[dashed,dash pattern={on 4pt off 2pt}] ([xshift=0.5\tabcolsep,yshift=7pt]c1.north) -- ([xshift=0.5\tabcolsep,yshift=-2pt]d1.south);
  \tikz[remember picture,overlay]   \draw[dashed,dash pattern={on 4pt off 2pt}] ([xshift=0.5\tabcolsep,yshift=7pt]c2.north) -- ([xshift=0.5\tabcolsep,yshift=-2pt]d2.south);
\end{equation}

Computing the images of $U$ and $U^*$ under the birational maps $\mathcal K_{s,p,n}$ and  $\mathcal K_{s,n-p,n}$ respectively, we can conclude Lemma \ref{tiso}. \,\,\,$\endpf$
\medskip

In a similar way, we will introduce another  isomorphism in the following. 
Let $\mathfrak E=(e_{ij})$ be an $n\times n$ anti-diagonal matrix defined by
\begin{equation}
e_{ij}=\left\{
    \begin{array}{cc}
        1 & {\rm for\,\,}j=n+1-i \\
        0 & {\rm for\,\,}j\neq n+1-i
    \end{array}\right..
\end{equation}
For $0\leq k\leq p$ and $I=(i_1,\cdots, i_p)\in\mathbb I^k_{s,p,n}$, denote by  index $I^*$ the following index
\begin{equation}
 I^*:=(n+1-i_p,n+1-i_{p-1},\cdots,n+1-i_2,n+1-i_1)\in\mathbb I^{p-k}_{n-s,p,n}\,.
\end{equation}
Let  $[\cdots,z^k_I,\cdots]_{I\in\mathbb I^k_{s,p,n}}$ and $[\cdots,z^{p-k}_{I^*},\cdots]_{I^*\in\mathbb I^{p-k}_{n-s,p,n}}$ be the homogeneous coordinates for $\mathbb {CP}^{N^{k}_{s,p,n}}$ and  $\mathbb {CP}^{N^{p-k}_{n-s,p,n}}$ respectively.
For $0\leq k\leq p$, define an isomorphism $g_k:\mathbb {CP}^{N^{k}_{s,p,n}}\rightarrow\mathbb {CP}^{N^{p-k}_{n-s,p,n}}$ by assigning
$z^{p-k}_{I^*}=z^k_{I}$.
Then, we can define an isomorphism 
\begin{equation}
    \widetilde{\rm USD}:\,\,G(p,n)\times\mathbb {CP}^{N^{0}_{s,p,n}}\times\cdots\times\mathbb {CP}^{N^p_{s,p,n}}\rightarrow G(p,n)\times\mathbb {CP}^{N^{0}_{n-s,p,n}}\times\cdots\times\mathbb {CP}^{N^p_{n-s,p,n}}\,
\end{equation}
as follows.  Let $a=\left(x,\left[\cdots,z^0_I,\cdots\right]_{I\in\mathbb I^0_{s,p,n}}, \cdots, \left[\cdots,z^p_I,\cdots\right]_{I\in\mathbb I^p_{s,p,n}}\right)$ be an arbitrary point of $G(p,n)\times\mathbb {CP}^{N^{0}_{s,p,n}}\times\cdots\times\mathbb {CP}^{N^p_{s,p,n}}$.
Define the $G(p,n)$-component of $\widetilde{{\rm USD}}(a)$ to be $\mathfrak E\cdot x$. When $s<k\leq p$ or $0\leq k<s+p-n$, we define the $\mathbb {CP}^{N^k_{n-s,p,n}}$-component of $\widetilde{{\rm USD}}(a)$ to be the single point of $\mathbb {CP}^{N^k_{n-s,p,n}}$.
When $\max\{0,s+p-n\}\leq k\leq \min\{p,s\}$,  we define the $\mathbb {CP}^{N^k_{n-s,p,n}}$-component of $\widetilde{{\rm USD}}(a)$ by
\begin{equation}
\left[\cdots,z^{k}_{I^*}\left(\widetilde{{\rm USD}}(a)\right),\cdots\right]_{I^*\in\mathbb I^{k}_{n-s,p,n}}=g_{p-k}\left(\left[\cdots,z^{p-k}_I,\cdots\right]_{I\in\mathbb I^{p-k}_{s,p,n}}\right)\,.
\end{equation}

Similarly, we have
\begin{definitionlemma}\label{USD1}
$\widetilde {\rm USD}$ induces an isomorphism ${\rm USD}:\mathcal T_{s,p,n}\rightarrow\mathcal T_{n-s,p,n}$ such that the following diagram commutes.
\vspace{-.05in}
\begin{equation}\label{USD}
\begin{tikzcd} &\mathcal T_{s,p,n}\arrow{r}{\rm USD}&\mathcal T_{n-s,p,n}\arrow{d}{R_{n-s,p,n}} \\ &G(p,n) \arrow{r}{\mathfrak E}\arrow[leftarrow]{u}{R_{s,p,n}}&G(p,n)\\ \end{tikzcd}\,\,.\vspace{-20pt} \end{equation}  
\end{definitionlemma}

\begin{remark}
${\rm DUAL}$ is an automorphism of $\mathcal T_{s,p,2p}$ and USD is an automorphism of $\mathcal T_{s,p,2s}$. \end{remark}

\section{Van der Waerden Representation} \label{vander}

In this section, we will introduce the Van der Waerden representation for $\mathcal T_{s,p,n}$ and establish the smoothness of $\mathcal T_{s,p,n}$ as an immediate corollary. It also lays the foundation for the later sections as an efficient computational tool.

For convenience, we make the following conventions in the remainder of this paper. Firstly,  we denote by $P_I(x)$ the functions $P_I(\widetilde x)$ by a slight abuse of notation (for they are only used in such a way that the choice of matrix representatives is irrelevant); call $P_I$ the  Pl\"ucker coordinate functions. 
Secondly, depending on the context we enumerate and compare the components of the homogeneous coordinates $[\cdots,z_I,\cdots]_{I\in\mathbb I^k_{s,p,n}}$ in convenient orders  without mentioning the order explicitly each time. 

By Remark \ref{dualusd}, we may assume that $2p\leq n\leq 2s$ in this section.
\subsection{Coordinate charts for \texorpdfstring{$R_{s,p,n}^{-1}(U_p)$}{ee} when \texorpdfstring{$p=n-s\leq s$}{ee}} \label{vanderp}
To illustrate the idea of the original Van der Waerden representation,  we will address the case $s=n-p$ separately in this subsection.  

Define an affine  open subset $U_p\subset G(p,n)$ by
\begin{equation}\label{up1}
   U_p:=\left\{\left.  \left(Z \hspace{-0.13in}\begin{matrix}
  &\hfill\tikzmark{c2}\\
  &\hfill\tikzmark{d2}
  \end{matrix}\,\,\, I_{p\times p} \right)\right\vert_{}
   Z\,\,{\rm  is\,\, a\,\,} p\times s\,\,{\rm matrix}\,\right\}\,.
  \tikz[remember picture,overlay]   \draw[dashed,dash pattern={on 4pt off 2pt}] ([xshift=0.5\tabcolsep,yshift=7pt]c2.north) -- ([xshift=0.5\tabcolsep,yshift=-2pt]d2.south);
\end{equation}
Equip it with  holomorphic coordinates
\begin{equation}\label{u0}
    Z:=\left(\begin{matrix}
z_{11}&z_{12}&\cdots &z_{1s}\\
z_{21}&z_{22}&\cdots &z_{2s}\\
\vdots&\vdots&\ddots&\vdots\\
z_{p1}&z_{p2}&\cdots &z_{ps}\\
\end{matrix}\right).
\end{equation}
Define an index set $\mathbb J_p$ by
\begin{equation}
\mathbb J_p:=\left\{\left.\left(
\begin{matrix}
i_1&i_2&\cdots&i_p\\
j_1&j_2&\cdots&j_p\\
\end{matrix}\right) \right\vert_{}\footnotesize\begin{matrix}
(i_1,i_2,\cdots,i_p)\,\, {\rm is\,\,a\,\,permutation\,\,of\,\,}(1,2,\cdots,\,p);\\1\leq j_k\leq s,\,1\leq k\leq p,\,{\rm and\,\,}j_{k_1}\neq j_{k_2}\,\,{\rm for\,\,} k_1\neq k_2
\end{matrix}\right\}.
\end{equation}

In the following, we will associate each index $\tau\in \mathbb J_p$ with  a quasi-projective subvariety $A^{\tau}$ $\left(\cong\mathbb C^{p(n-p)}\right)$ of $\mathbb {CP}^{N_{p,n}}\times\mathbb {CP}^{N^0_{s,p,n}}\times\cdots\times\mathbb {CP}^{N^p_{s,p,n}}$, and a holomorphic map  $J^{\tau}:\mathbb C^{p(n-p)}\rightarrow A^{\tau}$ which is compatible with $R_{s,p,n}$.

Associate each $\tau=\left(\begin{matrix} i_1&i_2&\cdots&i_p\\
j_1&j_2&\cdots&j_p\\
\end{matrix}\right)\in\mathbb J_p$  with a complex Euclidean space $\mathbb {C}^{p(n-p)}$ equipped with  holomorphic coordinates $\left(\overrightarrow B^1,\cdots,\overrightarrow B^p\right)$ where
\begin{equation}\label{bp}
\small
\begin{split}
&\overrightarrow  B^{1}:=\left(a_{i_1j_1},\xi^{(1)}_{i_11},\xi^{(1)}_{i_12},\cdots,\xi^{(1)}_{i_1(j_1-1)},\xi^{(1)}_{i_1(j_1+1)},\cdots,\xi^{(1)}_{i_1s},\xi^{(1)}_{1j_1},\xi^{(1)}_{2j_1},\cdots,\xi^{(1)}_{(i_1-1)j_1},\xi^{(1)}_{(i_1+1)j_1},\cdots,\xi^{(1)}_{pj_1}\right)\\
&\overrightarrow  B^{2}:=\left(a_{i_2j_2},\xi^{(2)}_{i_21},\xi^{(2)}_{i_22},\cdots,\xi^{(2)}_{i_2(j_1-1)},\widehat{\xi^{(2)}_{i_2j_1}},\xi^{(2)}_{i_2(j_1+1)},\cdots,\xi^{(2)}_{i_2(j_2-1)},\widehat{\xi^{(2)}_{i_2j_2}},\xi^{(2)}_{i_2(j_2+1)},\cdots,\xi^{(2)}_{i_2s},\right.\\
&\,\,\,\,\,\,\,\,\,\,\,\,\,\left.\cdots,\xi^{(2)}_{1j_2},\cdots,\xi^{(2)}_{(i_1-1)j_2},\widehat{\xi^{(2)}_{i_1j_2}},\xi^{(2)}_{(i_1+1)j_2},\cdots,\xi^{(2)}_{(i_2-1)j_2},\widehat{\xi^{(2)}_{i_2j_2}},\xi^{(2)}_{(i_2+1)j_2},\cdots,\xi^{(2)}_{pj_2}\right)\\
&\,\,\,\,\,\,\,\,\,\,\,\,\,\,\,\,\,\,\,\,\,\,\,\,\,\,\,\vdots\\
&\overrightarrow B^{k}:=\left(a_{i_kj_k},\xi^{(k)}_{i_k1},\xi^{(k)}_{i_k2},\cdots,\widehat{\xi^{(k)}_{i_kj_1}},\cdots,\widehat{\xi^{(k)}_{i_kj_2}},\cdots,\widehat{\xi^{(k)}_{i_kj_k}},\cdots,\xi^{(k)}_{i_ks},\right.\\
&\,\,\,\,\,\,\,\,\,\,\,\,\,\,\,\,\,\,\,\,\,\,\,\,\,\,\,\,\,\,\left.\xi^{(k)}_{1j_k},\xi^{(k)}_{2j_k},\cdots,\widehat{\xi^{(k)}_{i_1j_k}},\cdots,\widehat{\xi^{(k)}_{i_2j_k}},\cdots,\widehat{\xi^{(k)}_{i_kj_k}},\cdots,\xi^{(k)}_{pj_k}\right)\in\mathbb C^{n+1-2k}\,\,\,\,\,\,\,\,\,\,\,\,\,\,\,\,\,\,\,\,\,\,\,\,\,\,\,\,\,\,\,\,\,\\
\end{split}
\end{equation}
\begin{equation*}
\begin{split}    
&\,\,\,\,\,\,\,\,\,\,\,\,\,\,\,\,\,\,\,\,\,\,\,\,\,\,\,\vdots\\
&\overrightarrow  B^{p}:=\left(a_{i_pj_p},\xi^{(p)}_{i_p1},\cdots,\widehat{\xi^{(p)}_{i_pj_1}},\cdots,\widehat{\xi^{(p)}_{i_pj_2}},\cdots,\widehat{\xi^{(p)}_{i_pj_p}},\cdots,\xi^{(p)}_{i_ps}\right)\,\,.\,\,\,\,\,\,\,\,\,\,\,\,\,\,\,\,\,\,\,\,\,\,\,\,\,\,\,\,\,\,\,\,\,\,\,\,\,\,\,\,\,\,\,\,\,\,\,\,\,\,\,\,\,\,\,\,\,\,\,\,\,\,\,\,\,\,\,\,\,\,\,\,\,\,\,\,\\
\end{split}
\end{equation*}

Define a holomorphic map $\Gamma^{\tau}:\mathbb C^{p(n-p)}\rightarrow U_p$  by
\begin{equation}\label{ngamma}
\Gamma^{\tau}\left(\overrightarrow B^1,\cdots,\overrightarrow B^p\right):=\left( \sum_{k=1}^p\,\Xi_k^T\cdot\Omega_k\cdot\prod_{t=1}^{k}a_{i_tj_t}\hspace{-0.13in}\begin{matrix}
  &\hfill\tikzmark{c2}\\
  \\
  &\hfill\tikzmark{d2}
  \end{matrix}\,\,\,\,I_{p\times p} \right)\,,
\tikz[remember picture,overlay]   \draw[dashed,dash pattern={on 4pt off 2pt}] ([xshift=0.5\tabcolsep,yshift=7pt]c2.north) -- ([xshift=0.5\tabcolsep,yshift=-2pt]d2.south);
\end{equation} 
where $\Xi_k^T$ is the transpose of the function-valued vector $\Xi_k$,  and  $\Xi_k$, $\Omega_k$ are defined for $1\leq k\leq p$ as follows.  $\Xi_k:=
\left(v_1^k,\cdots,v_p^k\right)$ where
\begin{equation}\label{nxi}
     v_t^k=\left\{\begin{array}{ll}
    \xi^{(k)}_{tj_k} \,\,\,\,\,\,\,\,\,\,\,\,\,\,\,\,\,\,\,\,\,\,\,\,\,\,\,\,\,\, & t\in\{1,2,\cdots,p\}\backslash\{i_1,i_2,\cdots,i_k\}\,\,\,\,\,\,\,\,\,\,\,\,\,\,\,\,\, \,\,\,\,\,\,\,\,\,\,\,\,\,\,\,\,\,\\
    0& t\in\{i_1,i_2,\cdots,i_{k-1}\}\,\,\,\,\,\,\, \\
  
    1 \,\,\,\,\,\,\,\,\,\,\,\,\,\,\,\,\,\,\,\,\,\,\,\,\,\,\,\,\,\, &t=i_k\,\,\,\,\,\,\,\,\,\,\,\,\,\,\,\,\, \,\,\,\,\,\,\,\,\,\,\,\,\,\,\,\,\,\\
    \end{array}\right.;
\end{equation}
 $\Omega_k:=
\left(w_1^k,\cdots,w_s^k\right)$ where
\begin{equation}\label{nome}
     w_t^k=\left\{\begin{array}{ll}
    \xi^{(k)}_{i_kt} \,\,\,\,\,\,\,\,\,\,\,\,\,\,\,\,\,\,\,\,\,\,\,\,\,\,\,\,\,\, & t\in\{1,2,\cdots,s\}\backslash\{j_1,j_2,\cdots,j_k\}\,\,\,\,\,\,\,\,\,\,\,\,\,\,\,\,\, \,\,\,\,\,\,\,\,\,\,\,\,\,\,\,\,\,\\
    0& t\in\{j_1,j_2,\cdots,j_{k-1}\}\,\,\,\,\,\,\, \\
  
    1 \,\,\,\,\,\,\,\,\,\,\,\,\,\,\,\,\,\,\,\,\,\,\,\,\,\,\,\,\,\, &t=j_k\,\,\,\,\,\,\,\,\,\,\,\,\,\,\,\,\, \,\,\,\,\,\,\,\,\,\,\,\,\,\,\,\,\,\\
    \end{array}\right..
\end{equation}

Define a rational map  $J^{\tau}:\mathbb C^{p(n-p)}\dashrightarrow\mathbb {CP}^{N_{p,n}}\times\mathbb {CP}^{N^0_{s,p,n}}\times\cdots\times\mathbb {CP}^{N^p_{s,p,n}}$ by 
\begin{equation}\label{jjtp}
 J^{\tau}:=\mathcal K_{s,p,n}\circ \Gamma^{\tau}  \,, 
\end{equation} 
where $\mathcal K_{s,p,n}=\left(e, f_s^0,\cdots,f_s^p\right)$ is the rational map defined by (\ref{bpp}).

\begin{example}
The above definition of $\Gamma^{\tau}$ is essentially the process of Gaussian elimination. Let $X=G(3,6)$ and $\tau=\left(
\begin{matrix} 1&2&3\\ 1&2&3\\
\end{matrix}\right)$.  The coordinates $\left(\overrightarrow B^1,\overrightarrow B^2,\overrightarrow B^3\right)$ of  $\mathbb {C}^{9}$ are
\begin{equation}
\begin{split}
&\overrightarrow  B^1=\left(a_{11},\xi^{(1)}_{12},\xi^{(1)}_{13}\,,\,\,\xi^{(1)}_{21},\xi^{(1)}_{31}\right),\overrightarrow  B^2=\left(a_{22},\xi^{(2)}_{23},\xi^{(2)}_{32}\right)\,,\,\,\overrightarrow  B^3=\left(a_{33}\right).\\
\end{split}
\end{equation}
The holomorphic map $\Gamma^{\tau}:\mathbb C^{9}\rightarrow U_p$ is given by $\Gamma^{\tau}\left(\overrightarrow B^1,\overrightarrow B^2,\overrightarrow B^3\right)=$
\begin{equation}
\left( \begin{matrix}
a_{11}&a_{11}\xi^{(1)}_{12}&a_{11}\xi^{(1)}_{13}\\
a_{11}\xi^{(1)}_{21}&a_{11}(\xi^{(1)}_{21}\xi^{(1)}_{12}+a_{22})&a_{11}(\xi^{(1)}_{21}\xi^{(1)}_{13}+a_{22}\xi^{(2)}_{23})\\
a_{11}\xi^{(1)}_{31}&a_{11}(\xi^{(1)}_{31}\xi^{(1)}_{12}+a_{22}\xi^{(2)}_{32})&a_{11}(\xi^{(1)}_{31}\xi^{(1)}_{13}+a_{22}(\xi^{(2)}_{32}\xi^{(2)}_{23}+a_{33})))\\
\end{matrix}\hspace{-0.13in}\begin{matrix}
 &\hfill\tikzmark{c2}\\
 \\
 &\hfill\tikzmark{d2}
 \end{matrix}\,\,\,\,I_{3\times 3} \right)\,.
\tikz[remember picture,overlay]   \draw[dashed,dash pattern={on 4pt off 2pt}] ([xshift=0.5\tabcolsep,yshift=7pt]c2.north) -- ([xshift=0.5\tabcolsep,yshift=-2pt]d2.south);
\end{equation} 
\end{example}

We have the following key lemma for the Van der Waerden representation.
\begin{lemma}\label{em}
The rational map $J^{\tau}$ defined by (\ref{jjtp}) is a holomorphic embedding of $\mathbb C^{p(n-p)}$.
\end{lemma} 
{\noindent\bf Proof of  Lemma \ref{em}.} One can show that different choices of $\tau$ can be related through the $GL(s,\mathbb C)\times GL(n-s,\mathbb C)$-action of $\mathcal T_{s,p,n}$, which only results in different signs. Without loss of generality, we can assume in the following that 
\begin{equation}
\tau= \left(
\begin{matrix}
p&p-1&\cdots&1\\
s&s-1&\cdots&s-p+1\\
\end{matrix}\right).    
\end{equation}

We will split the proof into two parts. We first show that $J^{\tau}$ has a holomorphic extension over $\mathbb C^{p(n-p)}$, and then the extension is injective.
\smallskip

{\bf\noindent Step 1.} Since  $e\circ\Gamma^{\tau}$ is holomorphic on $\mathbb C^{p(n-p)}$, to prove that $J^{\tau}$ has a holomorphic extension over $\mathbb C^{p(n-p)}$, it suffices to establish this for the rational maps $f_s^k\circ \Gamma^{\tau}:\mathbb C^{p(n-p)}\dashrightarrow\mathbb {CP}^{N^k_{s,p,n}}$, $0\leq k\leq p$.

When $k=p$ it is trivial for $\mathbb {CP}^{N_{s,p,n}^p}$ consists of a single point by convention.  

When $k=p-1$,  $\left(f_s^{p-1}\circ \Gamma^{\tau}\right)\left(\overrightarrow B^1,\cdots,\overrightarrow B^p\right)$ takes the following form, in terms of the homogeneous coordinates for $\mathbb {CP}^{N^{p-1}_{s,p,n}}$.
\begin{equation*}
\begin{split}
\bigg[&\cdots,P_{I}\left(\Gamma^{\tau}\left(\overrightarrow B^1,\cdots,\overrightarrow B^p \right)\right),\cdots\bigg]_{I\in\mathbb I_{s,p,n}^{p-1}}\\
=&\bigg[\,\,a_{ps},\,a_{ps}\xi^{(1)}_{p(s-1)},\,\cdots,\, a_{ps}\xi^{(1)}_{p1},\, (-1)\cdot a_{ps}\xi^{(1)}_{(p-1)s}, (-1)\cdot a_{ps}\left(\xi^{(1)}_{(p-1)s}\xi^{(1)}_{p(s-1)}+a_{(p-1)(s-1)}\right),\,\\
&(-1)\cdot a_{ps}\left(\xi^{(1)}_{(p-1)s}\xi^{(1)}_{p(s-2)}+a_{(p-1)(s-1)}\xi^{(2)}_{(p-1)(s-2)}\right),\cdots,\\
&\,\,\,\,\,\,\,\,\,\,\,\,\,\,\,\,\,\,\,\,\,\,\,\,\,\,\,\vdots\\
&(-1)^{k-1}\cdot a_{ps}\xi^{(1)}_{(p+1-k)s},\,(-1)^{k-1}\cdot a_{ps}\left(\xi^{(1)}_{(p+1-k)s}\xi^{(1)}_{p(s-1)}+a_{(p-1)(s-1)}\xi^{(2)}_{(p+1-k)(s-1)}\right),\,\\
&(-1)^{k-1}\cdot a_{ps}\left(\xi^{(1)}_{(p+1-k)s}\xi^{(1)}_{p(s-2)}+a_{2(s-1)}\left(\xi^{(2)}_{(p+1-k)(s-1)}\xi^{(2)}_{(p-1)(s-2)}+a_{(p-2)(s-2)}\xi^{(3)}_{(p+1-k)(s-2)}\right)\right),\\
&\,\,\,\,\,\,\,\,\,\,\,\,\,\,\,\,\,\,\,\,\,\,\,\,\,\,\,\vdots\\
&(-1)^{p-1}\cdot a_{ps}\xi^{(1)}_{1s},\,(-1)^{p-1}\cdot a_{ps}\left(\xi^{(1)}_{1s}\xi^{(1)}_{p(s-1)}+a_{(p-1)(s-1)}\xi^{(2)}_{1(s-1)}\right),\,\cdots\bigg]\\
=&\bigg[\,\,1,\,\xi^{(1)}_{p(s-1)},\,\cdots,\, \xi^{(1)}_{p1},\, (-1)\cdot \xi^{(1)}_{(p-1)s}, (-1)\cdot\left(\xi^{(1)}_{(p-1)s}\xi^{(1)}_{p(s-1)}+a_{(p-1)(s-1)} \right),\\
&(-1)\cdot \left(\xi^{(1)}_{(p-1)s}\xi^{(1)}_{p(s-2)}+a_{(p-1)(s-1)}\xi^{(2)}_{(p-1)(s-2)}\right),\cdots\\
&\,\,\,\,\,\,\,\,\,\,\,\,\,\,\,\,\,\,\,\,\,\,\,\,\,\,\,\vdots\\
\,\,\,\,&(-1)^{k-1}\cdot \xi^{(1)}_{(p+1-k)s},\,(-1)^{k-1}\cdot \left(\xi^{(1)}_{(p+1-k)s}\xi^{(1)}_{p(s-1)}+a_{(p-1)(s-1)}\xi^{(2)}_{(p+1-k)(s-1)}\right),\,\\
&(-1)^{k-1}\cdot\left(\xi^{(1)}_{(p+1-k)s}\xi^{(1)}_{p(s-2)}+a_{(p-1)(s-1)}(\xi^{(2)}_{(p+1-k)(s-1)}\xi^{(2)}_{(p-1)(s-2)}+a_{(p-2)(s-2)}\xi^{(3)}_{(p+1-k)(s-2)})\right),\\
\end{split}
\end{equation*}
\begin{equation}\label{f1}
\begin{split}
&\,\,\,\,\,\,\,\,\,\,\,\,\,\,\,\,\,\,\,\,\,\,\,\,\,\,\,\vdots\\
&(-1)^{p-1}\cdot \xi^{(1)}_{1s},\,(-1)^{p-1}\cdot \left(\xi^{(1)}_{1s}\xi^{(1)}_{p(s-1)}+a_{(p-1)(s-1)}\xi^{(2)}_{1(s-1)}\right),\,\cdots\bigg]\,.\,\,\,\,\,\,\,\,\,\,\,\,\,\,\,\,\,\,\,\,\,\,\,\,\,\,\,\,\,\,\,\,\,\,\,\,\,\,\,\,\,\,\,\,\,\,\,\,\,\,\,\\
\end{split}
\end{equation}
It is clear that (\ref{f1}) gives a holomorphic map from $\mathbb C^{p(n-p)}$ to $\mathbb {CP}^{N^{p-1}_{s,p,n}}$ which extends $f_s^{p-1}\circ \Gamma^{\tau}$.

In the following, we will establish the holomorphic extension for  $f_s^k\circ \Gamma^{\tau}$, $0\leq k\leq p-1$, by a similar argument.

We first introduce certain special indices $I_k,I_k^*,  I_{\mu\nu}^k, I_{\mu\nu}^{k*}\in\mathbb I^k_{s,p,n}$ as follows. For $0\leq k\leq p$, define
\begin{equation}\label{I_k}
    I_k:=(s+k,s+k-1,\cdots,s-p+k+1)\,;
\end{equation}
for $1\leq k\leq p-1$, define
\begin{equation}\label{I*k}
    I^*_{k}:=(s+k+1,s+k-1,s+k-2,\cdots,s-p+k+3,s-p+k+2,s-p+k)\,;
\end{equation}
for $0\leq k\leq p$, $s-p+k+1\leq \mu\leq s$, and $1\leq \nu\leq s-p+k$ define
\begin{equation}\label{Ikmn}
    I_{\mu\nu}^k:=(s+k,s+k-1,\cdots,,\widehat \mu,\cdots, \nu)\,;
\end{equation}
for $0\leq k\leq p$, $s+1\leq \mu\leq s+k$, and $s+k+1\leq \nu\leq n$ define
\begin{equation}\label{Ikmn*}
    I^{k*}_{\mu\nu}:=(\nu,s+k,s+k-1,\cdots,\widehat \mu,\cdots, s-p+k+1)\,.
\end{equation}

\medskip

{\noindent\bf Claim I.} For $0\leq k\leq p$,
\begin{equation}\label{hol}
P_{I_k}\left(\Gamma^{\tau}\left(\overrightarrow B^1,\cdots,\overrightarrow B^p\right)\right)=(-1)^{k(p-k)}\cdot \prod_{t=1}^{p-k} a^{p-k+1-t}_{(p+1-t)(s+1-t)}\,.
\end{equation}

{\noindent\bf Claim II.} For each $I\in\mathbb I^k_{s,p,n}$, $0\leq k\leq p$, there is a polynomial $Q_{I}$ in variables $(\overrightarrow B^1,\cdots,\overrightarrow B^p)$ such that \begin{equation}\label{hol2}
P_{I}\left(\Gamma^{\tau}\left(\overrightarrow B^1,\cdots,\overrightarrow B^p\right)\right)=\left(\,\prod_{t=1}^{p-k} a^{p-k+1-t}_{(p+1-t)(s+1-t)}\right)\cdot Q_{I}\left(\overrightarrow B^1,\cdots,\overrightarrow B^p\right).
\end{equation}

{\noindent\bf Proof of  Claim I.} For $1\leq i\leq p-k$, denote by  $e_i$ the submatrix  formed by the $(s-p+k+1)^{th}$, $(s-p+k+2)^{th}$, $\cdots$, $(s-i+1)^{th}$  columns and  the $(k+1)^{th}$, $(k+2)^{th}$, $\cdots$,  $(p+1-i)^{th}$ rows of the following $p\times s$ matrix.
\begin{equation}
  \sum_{j=i}^{p-k}\left(\Xi_j^T\cdot\Omega_j\cdot\prod_{t=1}^{j}a_{(p+1-t)(s+1-t)}\right)\,.
\end{equation} 
It is clear that 
\begin{equation}
    P_{I_k}\left(\Gamma^{\tau}\left(\overrightarrow B^1,\cdots,\overrightarrow B^p\right)\right)=(-1)^{k(p-k)}\cdot\det e_1\,.
\end{equation}

For $1\leq i\leq p-k-1$, define a $(p-k+1-i)\times(p-k+1-i)$ matrix $A_i$ by
\begin{equation}
A_i:=\left(
\begin{matrix}
&1&0&\cdots&0&-\xi^{(i)}_{(k+1)(s+1-i)}\\
&0&1&\cdots&0&-\xi^{(i)}_{(k+2)(s+1-i)}\\
&\vdots&\vdots&\ddots&\vdots&\vdots\\
&0&0&\cdots&1&-\xi^{(i)}_{(p-i)(s+1-i)}\\
&0&0&\cdots&0 &1\\
\end{matrix}
\right)\,.
\end{equation}
Then,
\begin{equation}
  {\rm det}\, e_i={\rm det}\, (A_i\cdot e_i)=\left(\prod_{t=1}^{i} a_{(p+1-t)(s+1-t)}\right)\cdot{\rm det}\, e_{i+1}\,,\,\,\,\,1\leq i\leq p-k-1.  
\end{equation} 
Noticing that 
\begin{equation}\label{last}
    e_{p-k}=\prod_{t=1}^{p-k} a_{(p+1-t)(s+1-t)}\,,
\end{equation}  
we complete the proof of Claim I. \,\,\,\,$\endpf$
\medskip

{\noindent\bf Proof of  Claim II.} Define a $p\times s$ matrix by
\begin{equation}\label{ee}
    E:=C_1+C_2+\cdots+C_p
\end{equation} where
\begin{equation}
 C_m:= \left(\prod_{t=1}^{m}a_{(p+1-t)(s+1-t)}\right)\cdot\Xi_m^T\cdot\Omega_m\,,\,\,\,1\leq m\leq p\,.
\end{equation} 
One can show that  $P_{I}\left(\Gamma^{\tau}\left(\overrightarrow B^1,\cdots,\overrightarrow B^p \right)\right)$ is the determinant of a certain $(p-k)\times (p-k)$ submatrix $\mathfrak e$ of $E$; assume that $\mathfrak e$ is formed by the  $\alpha_1^{th},\cdots,\alpha_{p-k}^{th}$ columns and the $\beta_1^{th},\cdots,\beta_{p-k}^{th}$ rows of $E$ such that $1\leq \alpha_1<\cdots< \alpha_{p-k}\leq s$ and $1\leq \beta_1<\cdots< \beta_{p-k}\leq p$.

For $1\leq m\leq p$, denote by $C_m^{\mathfrak e}$ the submatrix  of $C_m$ formed by the  $\beta_1^{th},\cdots,\beta_{p-k}^{th}$ rows. Notice that $C_m^{\mathfrak e}$ is a $(p-k)\times s$ matrix of rank $1$ and that its entries are monomials divided by $\prod_{t=1}^{m} a_{(p+1-t)(s+1-t)}$. 

Write $C_m^{\mathfrak e}=\left(\eta^{(1)}_m,\cdots,\eta_m^{(s)}\right)$ as column vectors. Then,
\begin{equation}\label{expa}
\begin{split}
{\rm det}\,\mathfrak e&=\sum_{m_1=1}^p\sum_{m_2=1}^p\cdots\sum_{m_{p-k}=1}^p{\rm det}(\eta_{m_1}^{(\alpha_1)},\eta_{m_2}^{(\alpha_2)},\cdots,\eta_{m_{p-k}}^{(\alpha_{p-k})})\\
&=\sum_{\substack{1\leq m_1,\,m_2,\,\cdots,\,m_{p-k}\leq p\\ m_1,\,m_2,\,\cdots,\,m_{p-k}\,{\rm are\, distinct}}}{\rm det}(\eta_{m_1}^{(\alpha_1)},\eta_{m_2}^{(\alpha_2)},\cdots,\eta_{m_{p-k}}^{(\alpha_{p-k})})\,.\\
\end{split}
\end{equation}
It is easy to verify that each term in the last line of (\ref{expa}) is a polynomial in $\left(\overrightarrow B^1,\cdots,\overrightarrow B^p \right)$ which is divided by $\prod_{t=1}^{p-k} a^{p-k+1-t}_{(p+1-t)(s+1-t)}$.

We thus complete the proof of Claim II. \,\,\,\,$\endpf$
\smallskip

By Claims I and II we can conclude that the rational map $J^{\tau}$ defined by $(\ref{jjtp})$ has a holomorphic extension over $\mathbb {C}^{p(n-p)}$. We denote the extension also by $J^{\tau}$ by a slight abuse of notation.

\smallskip

{\bf\noindent Step 2.}
To show that  $J^{\tau}:\mathbb C^{p(n-p)}\rightarrow \mathbb {CP}^{N_{p,n}} 
\times\mathbb {CP}^{N^0_{s,p,n}}\times\cdots\times\mathbb {CP}^{N^p_{s,p,n}}$
is an embedding, we make the following claims.
\medskip

{\noindent\bf Claim III.} For $0\leq k\leq p-1$, $\mu=s-p+k+1$, and $1\leq \nu\leq s-p+k$,
\begin{equation}\label{hol3}
P_{I^k_{\mu\nu}}\left(\Gamma^{\tau}\left(\overrightarrow B^1,\cdots,\overrightarrow B^p\right)\right)=(-1)^{k(p-k)}\cdot \left(\prod_{t=1}^{p-k} a^{p-k+1-t}_{(p+1-t)(s+1-t)}\right)\cdot\xi^{(p-k)}_{(k+1)\nu}\,.
\end{equation}

{\noindent\bf Claim III$^{\prime}$.} For $1\leq k\leq p-1$, $s+1\leq \mu\leq s+k$, and $\nu=s+k+1$,
\begin{equation}\label{hol4}
P_{I_{\mu\nu}^{k*}}\left(\Gamma^{\tau}\left(\overrightarrow B^1,\cdots,\overrightarrow B^p\right)\right)=(-1)^{k(p-k)+s+k+1-\mu}\cdot\left(\prod_{t=1}^{p-k} a^{p-k+1-t}_{(p+1-t)(s+1-t)}\right)\cdot\xi^{(p-k)}_{(\mu-s)(s-p+k+1)}\,,\,\,\,
\end{equation}
\begin{equation}\label{hol5}
\small
\begin{split}
 P_{I^*_{k}}\left(\Gamma^{\tau}\left(\overrightarrow B^1,\cdots,\overrightarrow B^p\right)\right)=&(-1)^{k(p-k)+1}\left(\prod_{t=1}^{p-k} a^{p-k+1-t}_{(p+1-t)(s+1-t)}\right)\left(a_{k(s-p+k)}+\xi^{(p-k)}_{k(s-p+k+1)}\xi^{(p-k)}_{(k+1)(s-p+k)}\right)\,.\,\,\,  \,\,\, 
\end{split}
\end{equation}

{\noindent\bf Proof of  Claims III, III$^{\prime}$.}   We will sketch the proof for it is the same as in Claim I.

For the case ${I^k_{\mu\nu}}$, we can proceed in the same way as in Claim I, and the only difference is that in the last step (\ref{last}) we have $\left(\prod_{t=1}^{p-k} a_{(p+1-t)(s+1-t)}\right)\cdot\xi^{(p-k)}_{(k+1)\mu}$ instead of $\prod_{t=1}^{p-k} a_{(p+1-t)(s+1-t)}$.
This concludes  (\ref{hol3}). 

For the case $I_{\mu\nu}^{k*}$, let $e_1$ be the submatrix of $E$ (defined by (\ref{ee})) formed by the $(s-p+k+1)^{th}$, $\cdots$, $(s-1)^{th}$, $s^{th}$ columns and the $(\nu-s)^{th},(k+2)^{th},(k+3)^{th},\cdots,p^{th}$  rows. Then 
$P_{I_{\mu\nu}^{k*}}\left(\Gamma^{\tau}\left(\overrightarrow B^1,\cdots,\overrightarrow B^p\right)\right)=(-1)^{k(p-k)+s+k+1-\nu}\cdot{\rm det}\,  e_1$.
Compute the determinant inductively;  notice that in the last step (\ref{last}) we have $\left(\prod_{t=1}^{p-k} a_{(p+1-t)(s+1-t)}\right)\cdot\xi^{(p-k)}_{(\nu-s)(s-p+k+1)}$. 
This completes the proof of (\ref{hol4}).

For the case $I^*_k$, let $e_1$ be the submatrix of $E$ formed by the $(s-p+k)^{th}$, $(s-p+k+2)^{th}$, $(s-p+k+3)^{th}$, $\cdots$, $s^{th}$ columns and the $k^{th},(k+2)^{th},(k+3)^{th},\cdots,p^{th}$ rows. Then 
$P_{I^*_{k}}=(-1)^{k(p-k)+1}\cdot{\rm det}\, e_1$. We can conclude (\ref{hol5}) similarly by derive in the last step that 
$\left(\prod_{t=1}^{p-k} a_{(p+1-t)(s+1-t)}\right)\cdot\left(a_{k(s-p+k)}+\xi^{(p-k)}_{k(s-p+k+1)}\cdot\xi^{(p-k)}_{(k+1)(s-p+k)}\right)$.
\,\,\,\,$\endpf$
\smallskip

Applying Claims I,  III, III$^{\prime}$, we have the following expression for $(f_s^k\circ\Gamma^{\tau})(\overrightarrow B^1,\cdots,\overrightarrow B^p)$ in terms of the homogeneous coordinates for  $\mathbb {CP}^{N^k_{s,p,n}}$. When $1\leq k\leq p-1$, 
\begin{equation*}
\begin{split}
\bigg[&\cdots,P_{I}\big(\Gamma^{\tau}(\overrightarrow B^1,\cdots,\overrightarrow B^p)\big),\cdots\bigg]_{I\in\mathbb I_{s,p,n}^k}\,\,\,\,\,\,\,\,\,\,\,\,\,\,\,\,\,\,\,\,\,\,\,\,\,\,\,\,\,\,\,\,\,\,\,\,\,\,\,\,\,\,\,\,\,\,\,\,\,\,\,\,\,\,\,\,\,\,\,\,\,\,\,\,\,\,\,\,\,\,\,\,\,\,\,\,\,\,\,\,\,\,\,\,\,\,\,\,\,\,\,\,\,\,\,\,\,\,\,\,\,\,\,\,\,\,\,\,\,\,\,\,\,\,\,\,\,\,\,\,\,\,\,\,\,\,\,\,\,\,\,\,\,\,\,\,\,\,\,\,\\
\end{split}
\end{equation*}
\begin{equation}
\small
\begin{split}
=&\bigg[\,(-1)^{k(p-k)} \prod_{t}^{p-k} a^{p-k+1-t}_{(p+1-t)(s+1-t)},\,\cdots,\,(-1)^{k(p-k)} \xi^{(p-k)}_{(k+1)(s-p+k)}\prod_{t=1}^{p-k} a^{p-k+1-t}_{(p+1-t)(s+1-t)},\\
&\,(-1)^{k(p-k)} \xi^{(p-k)}_{(k+1)(s-p+k-1)}\prod_{t=1}^{p-k} a^{p-k+1-t}_{(p+1-t)(s+1-t)},\cdots,(-1)^{k(p-k)} \xi^{(p-k)}_{(k+1)1}\prod_{t=1}^{p-k} a^{p-k+1-t}_{(p+1-t)(s+1-t)},\\
&\,\,\,\,\,\,\,\,\,\,\,\,\,\,\,\,\,\,\,\,\,\,\,\,\,\,\,\,\,\,\,\,\,\,\,\,\,\,\,\,\cdots\\
&(-1)^{k(p-k)+1} \xi^{(p-k)}_{k(s-p+k+1)}\prod_{t=1}^{p-k} a^{p-k+1-t}_{(p+1-t)(s+1-t)},\cdots,(-1)^{k(p-k)+2} \xi^{(p-k)}_{(k-1)(s-p+k+1)}\prod_{t=1}^{p-k} a^{p-k+1-t}_{(p+1-t)(s+1-t)},\cdots,\\
&\,\,\,\,\,\,\,\,\,\,\,\,\,\,\,\,\,\,\,\,\,\,\,\,\,\,\,\,\,\,\,\,\,\,\,\,\,\,\,\,\cdots\\
&(-1)^{k(p-k)+k} \xi^{(p-k)}_{1(s-p+k+1)}\prod_{t=1}^{p-k} a^{p-k+1-t}_{(p+1-t)(s+1-t)},\cdots,\\
&(-1)^{k(p-k)+1} \left(a_{k(s-p+k)}+\xi^{(p-k)}_{k(s-p+k+1)}\cdot\xi^{(p-k)}_{(k+1)(s-p+k)}\right)\prod_{t=1}^{p-k} a^{p-k+1-t}_{(p+1-t)(s+1-t)},\cdots\bigg]\,\,\,\,\,\,\,\,\,\,\,\,\,\,\,\,\,\,\,\,\,\,\\
\end{split}
\end{equation}
\begin{equation}
\begin{split}
=&\bigg[\,1,\,\cdots,\, \xi^{(p-k)}_{(k+1)(s-p+k)},\,\xi^{(p-k)}_{(k+1)(s-p+k-1)},\,\cdots, \xi^{(p-k)}_{(k+1)1},\cdots,(-1) \xi^{(p-k)}_{k(s-p+k+1)},\cdots,\,\\
&(-1)^{2} \xi^{(p-k)}_{(k-1)(s-p+k+1)},\cdots,(-1)^{k} \xi^{(p-k)}_{1(s-p+k+1)},\cdots,-\left(a_{k(s-p+k)}+\xi^{(p-k)}_{k(s-p+k+1)}\xi^{(p-k)}_{(k+1)(s-p+k)}\right),\cdots\bigg]\,.\\
\end{split}
\end{equation}
When $k=0$ and $s\geq p+1$,
\begin{equation}
\bigg[\cdots,P_{I}\left(\Gamma^{\tau}\left(\overrightarrow B^1,\cdots,\overrightarrow B^p\right)\right),\cdots\bigg]_{I\in\mathbb I_{s,p,n}^0}\\
=\big[\,1,\,\cdots,\,\xi^{(p)}_{1(s-p)},\,\xi^{(p)}_{1(s-p-1)},\,\cdots,\,\xi^{(p)}_{11},\,\cdots]\,.
\end{equation}

Moreover, recalling the  coordinate $z_{1s}$ for $U_p$ defined by (\ref{u0}), we have that
\begin{equation}
    z_{ps}\left(\left(e\circ\Gamma^{\tau}\right)\left(\overrightarrow B^1,\cdots,\overrightarrow B^p\right)\right)=a_{ps}\,.
\end{equation}
Then it is easy to verify that $J^{\tau}=\big(e\circ\Gamma^{\tau}, f_s^0\circ\Gamma^{\tau}, \cdots, f_s^k\circ\Gamma^{\tau},\cdots, f_s^p\circ\Gamma^{\tau}\big)$  is an embedding.
\smallskip

We complete the proof of Lemma \ref{em}.
\,\,\,\,$\endpf$
\medskip

Let $A^{\tau}$, $\tau\in\mathbb J_p$, be the image of $\mathbb C^{p(n-p)}$ under the holomorphic map $J^{\tau}$.  By a slight abuse of notation, we denote the inverse holomorphic map by $(J^{\tau})^{-1}:A^{\tau}\rightarrow \mathbb C^{p(n-p)}$.   It is clear that  $\left\{\left(A^{\tau},(J^{\tau})^{-1}\right)\right\}_{\tau\in\mathbb J_p}$  is a system of coordinate charts of $R^{-1}_{s,p,n}(U_p)$.

\begin{definition}\label{vande0}
We call the above defined system of coordinate charts $\left\{\left(A^{\tau},(J^{\tau})^{-1}\right)\right\}_{\tau\in\mathbb J_p}$ the {\it Van der Waerden representation} of $R^{-1}_{s,p,n}(U_p)$.
\end{definition}

Next, we will show that that  the union of $A^{\tau}$ is $R^{-1}_{s,p,n}(U_p)$; hence the Van der Waerden representation $\left\{\left(A^{\tau},(J^{\tau})^{-1}\right)\right\}_{\tau\in\mathbb J_p}$  is a holomorphic atlas for $R^{-1}_{s,p,n}(U_p)$.

\begin{lemma}\label{coor2} $\bigcup_{\tau\in\mathbb J_p}A^{\tau}= R_{s,p,n}^{-1}(U_p)$.
\end{lemma}
{\bf\noindent Proof of Lemma \ref{coor2}.} See Appendix \ref{section:coverp}.\,\,\,\,$\endpf$
\medskip

We thus have the following corollary. 
\begin{corollary}\label{0smoo}
	$R_{s,p,n}^{-1}(U_p)$ is smooth.
\end{corollary}

\subsection{Coordinate charts for \texorpdfstring{$R_{s,p,n}^{-1}(U_l)$}{ee} } \label{vanderl}
In this subsection, we will introduce  the Van der Waerden representation for  $R_{s,p,n}^{-1}(U_l)$ provided that $p\leq s$. Here $U_l$, $0\leq l\leq p$, is an affine open subset of $G(p,n)$ defined by
\begin{equation}\label{ul}
U_l:=\left\{\,\,\,\,\underbracedmatrixll{Z\\Y}{s-p+l\,\,\rm columns}
  \hspace{-.45in}\begin{matrix}
  &\hfill\tikzmark{a}\\
  &\hfill\tikzmark{b}  
  \end{matrix} \,\,\,\,\,
  \begin{matrix}
  0\\
I_{(p-l)\times(p-l)}\\
\end{matrix}\hspace{-.11in}
\begin{matrix}
  &\hfill\tikzmark{c}\\
  &\hfill\tikzmark{d}
  \end{matrix}\hspace{-.11in}\begin{matrix}
  &\hfill\tikzmark{g}\\
  &\hfill\tikzmark{h}
  \end{matrix}\,\,\,\,
\begin{matrix}
I_{l\times l}\\
0\\
\end{matrix}\hspace{-.11in}
\begin{matrix}
  &\hfill\tikzmark{e}\\
  &\hfill\tikzmark{f}\end{matrix}\hspace{-.3in}\underbracedmatrixrr{X\\W}{(n-s-l)\,\,\rm columns}\,\,\,\,\right\}
  \tikz[remember picture,overlay]   \draw[dashed,dash pattern={on 4pt off 2pt}] ([xshift=0.5\tabcolsep,yshift=7pt]a.north) -- ([xshift=0.5\tabcolsep,yshift=-2pt]b.south);\tikz[remember picture,overlay]   \draw[dashed,dash pattern={on 4pt off 2pt}] ([xshift=0.5\tabcolsep,yshift=7pt]c.north) -- ([xshift=0.5\tabcolsep,yshift=-2pt]d.south);\tikz[remember picture,overlay]   \draw[dashed,dash pattern={on 4pt off 2pt}] ([xshift=0.5\tabcolsep,yshift=7pt]e.north) -- ([xshift=0.5\tabcolsep,yshift=-2pt]f.south);\tikz[remember picture,overlay]   \draw[dashed,dash pattern={on 4pt off 2pt}] ([xshift=0.5\tabcolsep,yshift=7pt]g.north) -- ([xshift=0.5\tabcolsep,yshift=-2pt]h.south);
\end{equation}
and equipped  with the following holomorphic coordinates.
\begin{equation}\label{ulx}
    Z:=\left(\begin{matrix}
    z_{11}&\cdots&z_{1(s-p+l)}\\
    \vdots&\ddots&\vdots\\
    z_{l1}&\cdots&z_{l(s-p+l)}\\
    \end{matrix}\right)\,\,,\,\,\,\,X:=\left(\begin{matrix}x_{1(s+l+1)}&\cdots&x_{1n}\\ \vdots&\ddots&\vdots\\ x_{l(s+l+1)}&\cdots &x_{ln}\\
\end{matrix}\right)\,,\,\,\,\,\,\,\,\,\,\,\,\,\,\,\,\,\,\,\,\,\,\,\,\,\,\,\,\,\,\,\,\,\,\,\,\,\,\,\,\,\,\,\,\,\,\,\,
\end{equation}
\begin{equation}
\,\,Y:=\left(\begin{matrix}
y_{(l+1)1}&\cdots& y_{(l+1)(s-p+l)}\\
\vdots&\ddots&\vdots\\
y_{p1}&\cdots& y_{p(s-p+l)}\\
\end{matrix}  \right)\,\,,\,\,\,\,W:=\left(\begin{matrix}
w_{(l+1)(s+l+1)}&\cdots& w_{(l+1)n}\\
\vdots&\ddots&\vdots\\
w_{p(s+l+1)}&\cdots& w_{pn}\\
\end{matrix}  \right) \,.
\end{equation}
\begin{remark}
We make the convention that, in the remainder of this paper, the open set $U_l$, $0\leq l\leq p$, is always referred to (\ref{ul}).
\end{remark}
\begin{remark}
  We make the convention that the matrices of dimension zero are omitted in (\ref{ul}). For instance,  when $l=p$ and $p=n-s$, the definition of $U_p$ in (\ref{ul}) coincides with that given in Section \ref{vanderp} by omitting the matrices $W,X,Y$ and $I_{(p-l)\times(p-l)}$. When $l=0$ by omitting matrices $X,Z$ and $I_{l\times l}$ in (\ref{ul}), we have
\begin{equation}
  U_0:=  \left\{\left. \left(Y\,\hspace{-0.15in}\begin{matrix}
  &\hfill\tikzmark{c1}\\
  &\hfill\tikzmark{d1}
  \end{matrix}\,\,\,\, I_{p\times p}\hspace{-0.13in}\begin{matrix}
  &\hfill\tikzmark{c2}\\
  &\hfill\tikzmark{d2}
  \end{matrix}\hspace{-.11in}\begin{matrix}
  &\hfill\tikzmark{g}\\
  &\hfill\tikzmark{h}
  \end{matrix}\,\,\,\, W\right)\right\vert_{}\footnotesize\begin{matrix}
  \,Y\,\,{\rm  is\,\, a\,\,} p\times (s-p)\,\,{\rm matrix}\,;\\
  W\,\,{\rm  is\,\, a\,\,} p\times (n-s)\,\,{\rm matrix}\,\\
  \end{matrix}\right\}\,.
  \tikz[remember picture,overlay]   \draw[dashed,dash pattern={on 4pt off 2pt}] ([xshift=0.5\tabcolsep,yshift=7pt]c1.north) -- ([xshift=0.5\tabcolsep,yshift=-2pt]d1.south);\tikz[remember picture,overlay]   \draw[dashed,dash pattern={on 4pt off 2pt}] ([xshift=0.5\tabcolsep,yshift=7pt]c2.north) -- ([xshift=0.5\tabcolsep,yshift=-2pt]d2.south);\tikz[remember picture,overlay]   \draw[dashed,dash pattern={on 4pt off 2pt}] ([xshift=0.5\tabcolsep,yshift=7pt]g.north) -- ([xshift=0.5\tabcolsep,yshift=-2pt]h.south);
\end{equation}
\end{remark}
\begin{remark}
We can adapt the Van der Waerden representation for $R_{s,p,n}^{-1}(U_l)$, $l=0,p$, either from the general construction to be introduced in the following by the above convention, or from that given in Section \ref{vanderp}. For simplicity, we will not  address it separately.
\end{remark}
\begin{remark}\label{dilation}
From the dynamic point of view, the Van der Waerden representation for $R_{s,p,n}^{-1}(U_l)$ when $l=0,p$ is different from that when $1\leq l\leq p-1$. When $l=0,p$, there is one dilation parameter, while there are two  when $1\leq l\leq p-1$. In Section \ref{foliation} we will discuss this  phenomenon in detail for the purpose of understanding the foliation structure on $\mathcal T_{s,p,n}$.
\end{remark}

In the following, we will construct the Van der Waerden representation for $R_{s,p,n}^{-1}(U_l)$. 

For $0\leq l\leq p$, define an index set $\mathbb J_l$ by \begin{equation}
\mathbb J_l:=\left\{\left(
\begin{matrix}
i_1&i_2&\cdots&i_{p-l}&\cdots&i_p\\
j_1&j_2&\cdots&j_{p-l}&\cdots&j_p\\
\end{matrix}\right)\rule[-.38in]{0.01in}{.82in}\,\,\footnotesize\begin{matrix}
l+1\leq \, i_k\,\leq p\,\,{\rm for}\,\,1\leq\,k\,\leq p-l\,;\,\,\\
1\leq \, i_k\,\leq l\,\,{\rm for}\,\,p-l+1\leq\,k\,\leq p\,;\,\,\\
s+l+1\leq \, j_k\,\leq n\,\,{\rm for}\,\,1\leq\,k\,\leq p-l\,;\,\,\\
1\leq \, j_k\,\leq s-p+l\,\,{\rm for}\,\,p-l+1\leq\,k\,\leq p\,;\,\,\\
i_{k_1}\neq i_{k_2}\,\,{\rm and\,\,}j_{k_1}\neq j_{k_2}\,\,{\rm for\,\,} k_1\neq k_2.\end{matrix}\right\}.
\end{equation}
Associate each $\tau=\left(\begin{matrix}
i_1&i_2&\cdots&i_p\\
j_1&j_2&\cdots&j_p\\
\end{matrix}\right)\in\mathbb J_l$ with a complex Euclidean space $\mathbb {C}^{p(n-p)}$ equipped with holomorphic coordinates  $\left(\widetilde X,\widetilde Y,\overrightarrow B^1,\cdots,\overrightarrow B^p\right)$ defined as follows.
\begin{equation}\label{ulu}
\widetilde X:=\left(\begin{matrix}
x_{1(s+l+1)}&\cdots &x_{1n}\\
\vdots&\ddots&\vdots\\
x_{l(s+l+1)}&\cdots &x_{ln}\\
\end{matrix}\right)\,\,\,{\rm and}\,\,\,\widetilde Y:=\left(\begin{matrix}
y_{(l+1)1}&\cdots& y_{(l+1)(s-p+l)}\\
\vdots&\ddots&\vdots\\
y_{p1}&\cdots& y_{p(s-p+l)}\\
\end{matrix}  \right)  \,;\,
\end{equation}
for $1\leq k\leq p-l$,
\begin{equation}\label{qbp}
\begin{split}
&\overrightarrow B^{k}:=\left(b_{i_{k}j_{k}},\xi^{(k)}_{i_{k}(s+l+1)},\xi^{(k)}_{i_{k}(s+l+2)},\cdots,\widehat{\xi^{(k)}_{i_{k}j_1}},\cdots,\widehat{\xi^{(k)}_{i_{k}j_2}},\cdots,\widehat{\xi^{(k)}_{i_{k}j_{k}}},\cdots,\xi^{(k)}_{i_{k}n},\right.\\
&\,\,\,\,\,\,\,\,\,\,\,\,\,\,\,\,\,\,\,\,\,\,\,\,\,\,\,\,\,\,\left.\xi^{(k)}_{(l+1)j_{k}},\xi^{(k)}_{(l+2)j_{k}},\cdots,\widehat{\xi^{(k)}_{i_1j_{k}}},\cdots,\widehat{\xi^{(k)}_{i_2j_{k}}},\cdots,\widehat{\xi^{(k)}_{i_{k}j_{k}}},\cdots,\xi^{(k)}_{pj_{k}}\right)\,;\\
\end{split}
\end{equation}
for $p-l+1\leq k\leq p$,
\begin{equation}\label{qbpb}
\begin{split}
&\overrightarrow B^{k}:=\left(a_{i_{k}j_{k}},\xi^{(k)}_{i_{k}1},\xi^{(k)}_{i_{k}2},\cdots,\widehat{\xi^{(k)}_{i_{k}j_{p-l+1}}},\cdots,\widehat{\xi^{(k)}_{i_{k}j_{p-l+2}}},\cdots,\widehat{\xi^{(k)}_{i_{k}j_{k}}},\cdots,\xi^{(k)}_{i_{k}(s-p+l)},\right.\\
&\,\,\,\,\,\,\,\,\,\,\,\,\,\,\,\,\,\,\,\,\,\,\,\,\,\,\,\,\,\,\left.\xi^{(k)}_{1j_{k}},\xi^{(k)}_{2j_{k}},\cdots,\widehat{\xi^{(k)}_{i_{p-l+1}j_{k}}},\cdots,\widehat{\xi^{(k)}_{i_{p-l+2}j_{k}}},\cdots,\cdots,\widehat{\xi^{(k)}_{i_{k}j_{k}}},\cdots,\xi^{(k)}_{lj_{k}}\right)\,.\,\,\,\,\,\,\,\,\,\,\,\,\,\\
\end{split}
\end{equation}
Define a holomorphic map $\Gamma_l^{\tau}:\mathbb C^{p(n-p)}\rightarrow U_l$ by
\begin{equation}\label{ws}
\begin{split}
&\,\,\,\,\,\,\,\,\,\,\,\,\,\,\,\,\,\,\,\,\,\,\,\,\,\,\,\,\,\,\,\,\,\,\,\,\,\,\,\,\,\,\,\,\,\,\,\,\,\,\,\,\,\,\,\,\,\,\,\,\,\,\,\,\,\Gamma_l^{\tau}\left(\widetilde X,\widetilde Y,\overrightarrow B^1,\cdots,\overrightarrow B^p\right):=\\
&\left(
\begin{matrix}
\sum\limits_{k=p-l+1}^p\left(\prod\limits_{t=p-l+1}^{k}a_{i_{t}j_t}\right)\cdot\Xi_k^T\cdot\Omega_k &0_{l\times(p-l)}&I_{l\times l}&\widetilde X\\ \widetilde Y&I_{(p-l)\times(p-l)}&0_{(p-l)\times l}&\sum\limits_{k=1}^{p-l}\left(\prod\limits_{t=1}^{k}b_{i_tj_t}\right)\cdot\Xi_k^T\cdot\Omega_k\\
\end{matrix}\right)\,,      \end{split}
\end{equation}
where $\Xi_k$ and $\Omega_k$ are defined as follows. For $1\leq k\leq p-l$, $\Xi_k:=
\left(v_{l+1}^k,\cdots,v_p^k\right)$ where
\begin{equation}\label{w3}
     v_t^k=\left\{\begin{array}{ll}
    \xi^{(k)}_{tj_k} \,\,\,\,\,\,\,\,\,\,\,\,\,\,\,\,\,\,\,\,\,\,\,\,\,\,\,\,\,\, & t\in\{l+1,l+2,\cdots,p\}\backslash\{i_1,i_2,\cdots,i_k\}\,\,\,\,\,\,\,\,\,\,\,\,\,\,\,\,\, \,\,\,\,\,\,\,\,\,\,\,\,\,\,\,\,\,\\
    0& t\in\{i_1,i_2,\cdots,i_{k-1}\}\,\,\,\,\,\,\, \\
  
    1 \,\,\,\,\,\,\,\,\,\,\,\,\,\,\,\,\,\,\,\,\,\,\,\,\,\,\,\,\,\, &t=i_k\,\,\,\,\,\,\,\,\,\,\,\,\,\,\,\,\, \,\,\,\,\,\,\,\,\,\,\,\,\,\,\,\,\,\\
    \end{array}\right.,
\end{equation}
and $\Omega_k:=
\left(w_{s+l+1}^k,\cdots,w_n^k\right)$ where
\begin{equation}\label{w4}
     w_t^k=\left\{\begin{array}{ll}
    \xi^{(k)}_{i_kt} \,\,\,\,\,\,\,\,\,\,\,\,\,\,\,\,\,\,\,\,\,\,\,\,\,\,\,\,\,\, & t\in\{s+l+1,s+l+2,\cdots,n\}\backslash\{j_1,j_2,\cdots,j_k\}\,\,\\
    0& t\in\{j_1,j_2,\cdots,j_{k-1}\}\,\,\,\,\,\,\, \\
  
    1 \,\,\,\,\,\,\,\,\,\,\,\,\,\,\,\,\,\,\,\,\,\,\,\,\,\,\,\,\,\, &t=j_k\,\,\,\,\,\,\,\,\,\\
    \end{array}\right..
\end{equation}
For $p-l+1\leq k\leq p$,  $\Xi_k:=
\left(v_{1}^k,\cdots,v_l^k\right)$ where
\begin{equation}\label{w5}
     v_t^k=\left\{\begin{array}{ll}
    \xi^{(k)}_{tj_k} \,\,\,\,\,\,\,\,\,\,\,\,\,\,\,\,\,\,\,\,\,\,\,\,\,\,\,\,\,\, & t\in\{1,2,\cdots,l\}\backslash\{i_{p-l+1},i_{p-l+2},\cdots,i_k\}\,\,\,\,\,\,\,\,\,\,\,\,\,\,\,\,\, \,\,\,\,\,\,\,\,\,\,\,\,\,\,\,\,\,\\
    0& t\in\{i_{p-l+1},i_{p-l+2},\cdots,i_{k-1}\}\,\,\,\,\,\,\, \\
  
    1 \,\,\,\,\,\,\,\,\,\,\,\,\,\,\,\,\,\,\,\,\,\,\,\,\,\,\,\,\,\, &t=i_k\,\,\,\,\,\,\,\,\,\,\,\,\,\,\,\,\, \,\,\,\,\,\,\,\,\,\,\,\,\,\,\,\,\,\\
    \end{array}\right.,
\end{equation}
and $\Omega_k:=
\left(w_{1}^k,\cdots,w_{s-p+l}^k\right)$ where
\begin{equation}\label{w2}
     w_t^k=\left\{\begin{array}{ll}
    \xi^{(k)}_{i_kt} \,\,\,\,\,\,\,\,\,\,\,\,\,\,\,\,\,\,\,\,\,\,\,\,\,\,\,\,\,\, & t\in\{1,2,\cdots,s-p+l\}\backslash\{j_{p-l+1},j_{p-l+2},\cdots,j_k\} \,\,\,\,\,\,\,\,\,\,\,\,\,\,\,\,\,\\
    0& t\in\{j_{p-l+1},j_{p-l+2},\cdots,j_{k-1}\}\,\,\,\,\,\,\, \\
  
    1 \,\,\,\,\,\,\,\,\,\,\,\,\,\,\,\,\,\,\,\,\,\,\,\,\,\,\,\,\,\, &t=j_k\,\,\,\,\,\,\,\,\,\,\,\,\,\\
    \end{array}\right..
\end{equation}

We can define a rational map  $J_l^{\tau}:\mathbb C^{p(n-p)}\dashrightarrow\mathbb {CP}^{N_{p,n}}\times\mathbb {CP}^{N^0_{s,p,n}}\times\cdots\times\mathbb {CP}^{N^p_{s,p,n}}$ by 
\begin{equation}\label{jl}
 J_l^{\tau}:=\mathcal K_{s,p,n}\circ \Gamma_l^{\tau} \,.  
\end{equation}

\begin{example}
Consider $G(4,8)$, with $s=4$, $l=2$, and $\tau=\left(
\begin{matrix}
3&4&1&2\\
7&8&1&2\\
\end{matrix}\right)\in\mathbb J_2$. The holomorphic coordinates $\left(\widetilde X,\widetilde Y,\overrightarrow B^1,\overrightarrow B^2, \overrightarrow B^3, \overrightarrow B^4\right)$ of $\mathbb {C}^{16}$ are	\begin{equation}
\begin{split}
&\widetilde X=\left(x_{17},x_{18},x_{27},x_{28}\right),\,\,\,\,\,\,\widetilde Y=\left(y_{31},y_{32},y_{41},y_{42}\right),\\
&\overrightarrow  B^1=\left(b_{37},\xi^{(1)}_{38},\xi^{(1)}_{47}\right),\,\,\,\,\,\,\overrightarrow  B^2=\left(b_{48}\right),\,\,\,\,\,\, B^3=\left(a_{11},\xi^{(3)}_{12},\xi^{(3)}_{21}\right),\,\,\,\,\,\,\overrightarrow  B^4=\left(a_{22}\right).\\
\end{split}
\end{equation}
The holomorphic map $\Gamma_2^{\tau}:\mathbb C^{16}\rightarrow U_2$ is given by $\Gamma^{\tau}\left(\widetilde X,\widetilde Y,\overrightarrow B^1,\cdots,\overrightarrow B^4\right)=$
\begin{equation}
\left(
\begin{matrix}
a_{11}&a_{1}\cdot\xi^{(3)}_{12}&0&0&1&0&x_{17}&x_{18}\\
a_{11}\cdot\xi^{(3)}_{21}&a_{11}\cdot(\xi^{(3)}_{12}\cdot\xi^{(3)}_{21}+a_{22})&0&0&0&1&x_{27}&x_{28}\\
y_{31}&y_{32}&1&0&0&0&b_{37}&b_{37}\cdot\xi^{(1)}_{38}\\
y_{41}&y_{42}&0&1&0&0&b_{37}\cdot\xi^{(1)}_{47}&b_{37}\cdot(\xi^{(1)}_{47}\cdot\xi^{(1)}_{38}+b_{48})\\
	\end{matrix}
	\right).
	\end{equation}
\end{example}

Similarly to Lemma \ref{em}, we have the following key lemma.
\begin{lemma}\label{qem}
The rational map $J_l^{\tau}$ defined by (\ref{jl}) is a holomorphic embedding of $\mathbb C^{p(n-p)}$.
\end{lemma} 
{\noindent\bf Proof of  Lemma \ref{qem}.}
 Without loss of generality, we can assume  that 
\begin{equation}
\tau= \left(
\begin{matrix}
l+1&l+2&\cdots&p&l&l-1&\cdots&1\\
s+l+1&s+l+2&\cdots&s+p&s-p+l&s-p+l-1&\cdots&s-p+1\\
\end{matrix}\right).    
\end{equation}

{\noindent\bf Claim I.} Assume that $0\leq k\leq p$. Let $I_k\in\mathbb I^k_{s,p,n}$ be defined by (\ref{I_k}). Then, 
\begin{equation}\label{a1}
P_{I_k}\left(\Gamma_l^{\tau}\left(\widetilde X,\widetilde Y,\overrightarrow B^1,\cdots,\overrightarrow B^p\right)\right)=\left\{
\begin{array}{ll}
(-1)^{k(p-k)} & {\rm when}\,\,k=l\\
(-1)^{k(p-k)}\cdot \prod\limits_{t=1}^{l-k} a^{l-k+1-t}_{(l+1-t)(s-p+l+1-t)}& {\rm when}\,\,k<l \\
(-1)^{k(p-k)}\cdot \prod\limits_{t=1}^{k-l}b^{k-l+1-t}_{(l+t)(s+l+t)}&{\rm when}\,\,k>l\\
\end{array}\right..
\end{equation}

{\noindent\bf Claim II.} For each $I\in\mathbb I_{s,p,n}^k$, $0\leq k\leq p$, there is a polynomial $Q_{I}$  in variables $\left(\widetilde X,\widetilde Y,\overrightarrow B^1,\cdots,\right.$ $\left.\overrightarrow B^p\right)$ such that
\begin{equation}\label{ghol2}
\small
P_{I}\left(\Gamma_l^{\tau}\left(\widetilde X,\widetilde Y,\overrightarrow B^1,\cdots,\overrightarrow B^p\right)\right)=P_{I_k}\left(\Gamma_l^{\tau}\left(\widetilde X,\widetilde Y,\overrightarrow B^1,\cdots,\overrightarrow B^p\right)\right)\cdot Q_{I}\left(\widetilde X,\widetilde Y,\overrightarrow B^1,\cdots,\overrightarrow B^p \right).
\end{equation} 

{\noindent\bf Proof of  Claims I, II.} 
We can conclude Claim I in a similar way as Lemma \ref{em} by applying Gaussian elimination. 

When $k=l$, Claim II holds trivially. We hence assume that $k\neq l$ in the following. Define
\begin{equation} E_1:=\sum\limits_{k=1}^{p-l}\left(\prod_{t=1}^{k}b_{i_tj_t}\right)\cdot\Xi_k^T\cdot\Omega_k\,\,\,\,{\rm and}\,\, E_2:=\sum\limits_{k=p-l+1}^p\left(\prod_{t=p-l+1}^{k}a_{i_{t}j_t}\right)\cdot\Xi_k^T\cdot\Omega_k\,.
\end{equation} 
It is easy to verify that  $P_{I}$,  $I\in\mathbb I_{s,p,n}^k$, can be expanded as \begin{equation}
P_{I}\left(\Gamma_l^{\tau}\left(\widetilde X,\widetilde Y,\overrightarrow B^1,\cdots,\overrightarrow B^p\right)\right)=\sum_{\,\,{\rm finitely\,\,many\,\,indices\,\,}\mathfrak a}a_{\mathfrak a}\left(\widetilde X,\widetilde Y\right)\cdot H_{\mathfrak a}\left(\overrightarrow B^1,\cdots,\overrightarrow B^p\right)
\end{equation}
such that the following holds.
\begin{enumerate}[label=(\alph*).]
    \item $a_{\mathfrak a}$ is a polynomial in variables $\widetilde X$ and $\widetilde Y$\,.
    \item $H_{\mathfrak a}\left(\overrightarrow B^1,\cdots,\overrightarrow B^p\right)$ is a product of the determinants of a certain submatrix $e^\mathfrak a_{1}$ of  $E_1$ and  a certain submatrix $e^\mathfrak a_{2}$ of $E_2$\,.
    \item  When $k>l$, $e^\mathfrak a_{1}$ is of rank at least $k-l$; when $k<l$,  $e^\mathfrak a_{2}$ is of rank at least $l-k$.
\end{enumerate} 
Then by Gaussian elimination, we can prove (\ref{ghol2}) in the same way as Claim II in Lemma $\ref{em}$.

We complete the proof.  \,\,\,\,$\endpf$
\smallskip

Applying Claims I and II, we can conclude that the rational map  $J_l^{\tau}$ extends to a holomorphic map from $\mathbb C^{p(n-p)}$ to $\mathbb {CP}^{N_{p,n}}\times\mathbb {CP}^{N^0_{s,p,n}}\times\cdots\times\mathbb {CP}^{N^p_{s,p,n}}$ in the same way as Lemma \ref{em}.
\smallskip

Next we will show that  $J_l^{\tau}$ 
is an embedding. 
\smallskip

{\noindent\bf Claim III.}  For $0\leq k\leq l-1$, $\mu=s-p+k+1$,  and $1\leq \nu\leq s-p+k$, let $I^k_{\mu\nu}\in\mathbb I^k_{s,p,n}$ be defined by (\ref{Ikmn}). Then,
\begin{equation}\label{thol2}
P_{I^k_{\mu\nu}}\left(\Gamma_l^{\tau}\left(\widetilde X, \widetilde Y,\cdots,\overrightarrow B^p\right)\right)=(-1)^{k(p-k)} \left(\prod_{t=1}^{l-k} a^{l-k+1-t}_{(l+1-t)(s-p+l+1-t)}\right)\xi^{(p-k)}_{(k+1)\nu}\,.
\end{equation}
\smallskip

{\noindent\bf Claim III$^{\prime}$.}  For $1\leq k\leq l-1$, $s+1\leq \mu\leq s+k$, and $\nu=s+k+1$, let   $I_{\mu\nu}^{k*}, I^*_k\in\mathbb I^k_{s,p,n}$ be defined by (\ref{Ikmn*}) and (\ref{I*k}) respectively. Then,
\begin{equation}\label{thol3}
\small
P_{I_{\mu\nu}^{k*}}\left(\Gamma_l^{\tau}\left(\widetilde X, \widetilde Y, \cdots,\overrightarrow B^p\right)\right)=(-1)^{k(p-k)+s+k+1-\mu}\left(\prod_{t=1}^{l-k} a^{l-k+1-t}_{(l+1-t)(s-p+l+1-t)}\right)\xi^{(p-k)}_{(\mu-s)(s-p+k+1)}\,,
\end{equation}
\begin{equation}
\begin{split}
&P_{I^*_{k}}\left(\Gamma_l^{\tau}\left(\widetilde X, \widetilde Y,\cdots,\overrightarrow B^p\right)\right)=(-1)^{k(p-k)+1} \left(\prod_{t=1}^{l-k} a^{l-k+1-t}_{(l+1-t)(s-p+l+1-t)}\right)\\
&\,\,\,\,\,\,\,\,\,\,\,\,\,\,\,\,\,\,\,\,\,\,\,\,\,\,\,\,\,\,\,\,\,\,\,\,\,\,\,\,\,\,\,\,\,\,\,\,\,\,\,\,\,\,\,\,\,\,\,\,\,\,\,\,\,\,\,\,\,\,\,\,\,\,\,\,\,\,\,\,\,\,\,\,\,\,\,\,\,\,\,\,\,\,\,\,\,\,\,\,\,\,\cdot\left(a_{k(s-p+k)}+\xi^{(p-k)}_{(k+1)(s-p+k)}\cdot\xi^{(p-k)}_{k(s-p+k+1)}\right)\,.
\end{split}
\end{equation}
\smallskip

{\noindent\bf Claim III$^{\prime\prime}$.} 
For $l+1\leq k\leq p$, $\mu=s+k$, and $s+k+1\leq \nu\leq n$, 
\begin{equation}\label{thol4}
\begin{split}
&P_{I_{\mu\nu}^{k*}}\left(\Gamma_l^{\tau}\left(\widetilde X, \widetilde Y,\cdots,\overrightarrow B^p\right)\right)=(-1)^{k(p-k)}\left(\prod_{t=1}^{k-l}b^{k-l+1-t}_{(l+t)(s+l+t)}\right)\xi^{(k-l)}_{k\nu}\,.\\
\end{split}
\end{equation}

{\noindent\bf Claim III$^{\prime\prime\prime}$.} 
For $l+1\leq k\leq p-1$,  $s-p+k+1\leq \mu\leq s$, and $\nu=s-p+k$, 
\begin{equation}\label{thol5}
P_{I_{\mu\nu}^k}\left(\Gamma_l^{\tau}\left(\widetilde X, \widetilde Y,\cdots,\overrightarrow B^p\right)\right)=(-1)^{k(p-k)+\mu-s+p-k}\left(\prod_{t=1}^{k-l}b^{k-l+1-t}_{(l+t)(s+l+t)}\right)\xi^{(k-l)}_{(\mu-s+p)(s+k)}\,,
\end{equation}
\begin{equation}\label{thol6}
\begin{split}
&P_{I^*_{k}}\left(\Gamma_l^{\tau}\left(\widetilde X, \widetilde Y,\cdots,\overrightarrow B^p\right)\right)=(-1)^{k(p-k)+1}\left(\prod_{t=1}^{k-l}b^{k-l+1-t}_{(l+t)(s+l+t)}\right)\\
&\,\,\,\,\,\,\,\,\,\,\,\,\,\,\,\,\,\,\,\,\,\,\,\,\,\,\,\,\,\,\,\,\,\,\,\,\,\,\,\,\,\,\,\,\,\,\,\,\,\,\,\,\,\,\,\,\,\,\,\,\,\,\,\,\,\,\,\,\,\,\,\,\,\,\,\,\,\,\,\,\,\,\,\,\,\,\,\,\,\,\,\,\,\,\,\,\,\,\,\,\,\,\cdot\left(b_{(k+1)(s+k+1)}+\xi^{(k-l)}_{k(s+k+1)}\cdot\xi^{(k-l)}_{(k+1)(s+k)}\right)\,.
\end{split}
\end{equation}

{\noindent\bf Claim III$^{\prime\prime\prime\prime}$.} 
For $k=l$,  $s-p+k+1\leq \mu\leq s$, and $1\leq \nu\leq s-p+k$, 
\begin{equation}
P_{I_{\mu\nu}^k}\left(\Gamma_l^{\tau}\left(\widetilde X, \widetilde Y,\cdots,\overrightarrow B^p\right)\right)=(-1)^{k(p-k)+\mu-s+p-l}y_{(\mu-s+p)\nu}\,.
\end{equation}
For $k=l$,  $s+1\leq \mu\leq s+k$, and $s+k+1\leq \nu\leq n$, 
\begin{equation}
P_{I_{\mu\nu}^{k*}}\left(\Gamma_l^{\tau}\left(\widetilde X, \widetilde Y,\cdots,\overrightarrow B^p\right)\right)=(-1)^{k(p-k)+\mu-s-l}x_{(\mu-s)\nu}\,.
\end{equation}

{\noindent\bf Proof of  Claims III, III$^{\prime}$, III$^{\prime\prime}$, III$^{\prime\prime\prime}$,  and III$^{\prime\prime\prime\prime}$.} The proof is the same as that of Claim III in Lemma \ref{em}. We omit it here for simplicity. 
\,\,\,\,$\endpf$
\medskip

Thanks to (\ref{a1}), (\ref{thol2}), (\ref{thol3}),  (\ref{thol4}), (\ref{thol5}), (\ref{thol6}), we can write the rational map $\left(f_s^k\circ\Gamma_l^{\tau}\right)$ in terms of the homogeneous coordinates for $\mathbb {CP}^{N^k_{s,p,n}}$ as follows. 
When $1\leq k\leq l-1$,
\begin{equation}
\begin{split}
\bigg[&\cdots,P_{I}\left(\Gamma_l^{\tau}\left(\widetilde X, \widetilde Y,\overrightarrow B^1,\cdots,\overrightarrow B^p\right)\right),\cdots\bigg]_{I\in\mathbb I_{s,p,n}^k}\\
=&\bigg[1,\cdots,\xi^{(p-k)}_{(k+1)1},\xi^{(p-k)}_{(k+1)2},\cdots,\xi^{(p-k)}_{(k+1)(s-p+k)},\cdots,\\
&\cdots, (-1)^{k}\cdot\xi^{(p-k)}_{1(s-p+k+1)},(-1)^{k-1}\cdot\xi^{(p-k)}_{2(s-p+k+1)},\cdots,(-1)\cdot\xi^{(p-k)}_{k(s-p+k+1)},\cdots,\\
&\cdots,(-1)\cdot\left(a_{k(s-p+k)}+\xi^{(p-k)}_{(k+1)(s-p+k)}\cdot\xi^{(p-k)}_{k(s-p+k+1)}\right),\cdots\bigg]\,;\\
\end{split}
\end{equation}
when $k=l$,
\begin{equation}
\begin{split}
\bigg[&\cdots,P_{I}\left(\Gamma_l^{\tau}\left(\widetilde X, \widetilde Y,\overrightarrow B^1,\cdots,\overrightarrow B^p\right)\right),\cdots\big]_{I\in\mathbb I_{s,p,n}^k}\\
=&\bigg[1,\cdots,-y_{(l+1)1},\cdots, -y_{(l+1)(s-p+l)},(-1)^2y_{(l+2)1},\cdots,(-1)^2y_{(l+2)(s-p+l)},\cdots,(-1)^{p-l}y_{p1}, \cdots,  \\
&\,\,\,\,\,\, (-1)^{p-l}y_{p(s-p+l)},\cdots, (-1)^{l-1}x_{1(s+l+1)},\cdots, (-1)^{l-1}x_{1n}, \cdots,(-1)^{l-2}x_{2(s+l+1)},\cdots,\\
&\,\,\,\,\,\, (-1)^{l-2}x_{2n},\cdots,(-1)^0x_{l(s+l+1)},\cdots,(-1)^0x_{ln},\cdots\bigg]\,;\\
\end{split}
\end{equation}
when $l+1\leq k\leq p-1$,
\begin{equation}
\begin{split}
\bigg[&\cdots,P_{I}\left(\Gamma_l^{\tau}\left(\widetilde X, \widetilde Y,\overrightarrow B^1,\cdots,\overrightarrow B^p\right)\right),\cdots\big]_{I\in\mathbb I_{s,p,n}^k}\\
=&\bigg[1,\cdots,\xi^{(k-l)}_{k(s+k+1)},\xi^{(k-l)}_{k(s+k+2)},\cdots,\xi^{(k-l)}_{kn},\cdots,\\
&\cdots, (-1)\cdot\xi^{(k-l)}_{(k+1)(s+k)},(-1)^{2}\cdot\xi^{(k-l)}_{(k+2)(s+k)},\cdots,(-1)^{p-k}\cdot\xi^{(k-l)}_{p(s+k)},\cdots,\\
&\cdots,(-1)\cdot\left(b_{(k+1)(s+k+1)}+\xi^{(k-l)}_{k(s+k+1)}\cdot\xi^{(k-l)}_{(k+1)(s+k)}\right),\cdots\bigg]\,.\\
\end{split}
\end{equation}
When $k=0$ and $s\geq p+1$,
\begin{equation}
\begin{split}
\bigg[&\cdots,P_{I}\left(\Gamma_l^{\tau}\left(\widetilde X, \widetilde Y,\overrightarrow B^1,\cdots,\overrightarrow B^p\right)\right),\cdots\bigg]_{I\in\mathbb I_{s,p,n}^k}=\big[1,\cdots,\xi^{(p)}_{11},\xi^{(p)}_{12},\cdots,\xi^{(p)}_{1(s-p)},\cdots\big]\,;\\
\end{split}
\end{equation}
when $k=p$ and $n-s\geq p+1$,
\begin{equation}
\begin{split}
\bigg[&\cdots,P_{I}\left(\Gamma_l^{\tau}\left(\widetilde X, \widetilde Y,\overrightarrow B^1,\cdots,\overrightarrow B^p\right)\right),\cdots\bigg]_{I\in\mathbb I_{s,p,n}^p}=\big[1,\cdots,\xi^{(p-l)}_{p(s+p+1)},\xi^{(p-l)}_{p(s+p+2)},\cdots,\xi^{(p-l)}_{pn},\cdots\big]\,.\\
\end{split}
\end{equation}
Recalling the affine coordinates $z_{l(s-p+l)}$ and $w_{(l+1)(s+l+1)}$  for $U_l$ in (\ref{ulx}), we have that
\begin{equation}
\begin{split}
&z_{l(s-p+l)}\left(\Gamma_l^{\tau}\left(\widetilde X, \widetilde Y,\overrightarrow B^1,\cdots,\overrightarrow B^p\right)\right)=a_{l(s-p+l)}\,,\\ 
&w_{(l+1)(s+l+1)}\left(\Gamma_l^{\tau}\left(\widetilde X, \widetilde Y,\overrightarrow B^1,\cdots,\overrightarrow B^p\right)\right)=b_{(l+1)(s+l+1)}\,.\\ 
\end{split}
\end{equation}

Therefore, the holomorphic map
$J_l^{\tau}=\mathcal K_{s,p,n}\circ \Gamma_l^{\tau}$ 
is an embedding. This completes the proof of Lemma \ref{qem}.
\,\,\,\,$\endpf$
\smallskip

Let $A^{\tau}$, $\tau\in\mathbb J_l$, be the image of $\mathbb C^{p(n-p)}$ under the holomorphic map $J_l^{\tau}$.  By a slight abuse of notation, we denote the inverse holomorphic map by $(J_l^{\tau})^{-1}:A^{\tau}\rightarrow \mathbb C^{p(n-p)}$.   It is clear that  $\left\{\left(A^{\tau},(J_l^{\tau})^{-1}\right)\right\}_{\tau\in\mathbb J_p}$  is a system of coordinate charts of $R^{-1}_{s,p,n}(U_l)$.

\begin{definition}\label{vandel}
We call the above defined system of coordinate charts $\left\{\left(A^{\tau},(J_l^{\tau})^{-1}\right)\right\}_{\tau\in\mathbb J_l}$ the {\it Van der Waerden representation} of $R^{-1}_{s,p,n}(U_l)$.
\end{definition}

Similarly, we will show that  the union of $A^{\tau}$ is $R^{-1}_{s,p,n}(U_l)$; hence the Van der Waerden representation $\left\{\left(A^{\tau},(J_l^{\tau})^{-1}\right)\right\}_{\tau\in\mathbb J_l}$  is a holomorphic atlas for $R^{-1}_{s,p,n}(U_l)$.

\begin{lemma}\label{qcoor2} $\bigcup_{\tau\in\mathbb J_l}A^{\tau}= R_{s,p,n}^{-1}(U_l)$.
\end{lemma}
{\bf\noindent Proof of Lemma \ref{qcoor2}.} See Appendix \ref{section:coverl}.\,\,\,\,$\endpf$

\subsection{The case \texorpdfstring{$n-s<p$}{ff}} \label{vandern}
In this subsection,  we will introduce the Van der Waerden representation for $R_{s,p,n}^{-1}(U_l)$ when $n-s<p<s$, where $U_l$, $0\leq l\leq n-s$, is an affine open subset of $G(p,n)$ defined by (\ref{ul}).  Since the proof is the same as in Section \ref{vanderl}, we omit it here.

For $0\leq l\leq n-s$, define an index set $\mathbb J_l$ by
	\begin{equation}
 \left\{\left(
	\begin{matrix}
	i_1&\cdots&i_{n-s}\\
	j_1&\cdots&j_{n-s}\\
	\end{matrix}\right)\rule[-.35in]{0.01in}{.72in}\footnotesize\begin{matrix}
	(i_{n-s-l+1},\,i_{n-s-l+2},\,\cdots,\,i_{n-s})\,\, {\rm is\,\,a\,\,permutation\,\,of\,\,}(1,\,2,\,\cdots,\,\,l)\,;\\
	(j_1,\,\cdots,\,j_{n-s-l})\,\,{\rm is\,\,a\,\,permutations\,\,of\,\,}(s+l+1,\,s+l+2,\,\cdots,\,n)\,;\\
	1\leq\, j_t\,\leq s-p+l\,\,{\rm for}\,\,n-s-l+1\leq\,t\,\leq n-s\,,\,\,j_{t_1}\neq j_{t_2}\,\,{\rm for\,\,} t_1\neq t_2\,;\\
	l+1\leq\, i_t\,\leq p\,\,{\rm for}\,\,1\leq\,t\,\leq n-s-l\,,i_{t_1}\neq i_{t_2}\,\,{\rm for\,\,} t_1\neq t_2\\
	\end{matrix}\right\}.
	\end{equation}
Associate each $\tau=\left(\begin{matrix}
i_1&i_2&\cdots&i_{n-s}\\
j_1&j_2&\cdots&j_{n-s}\\
\end{matrix}\right)\in\mathbb J_l$ with a complex Euclidean space $\mathbb {C}^{p(n-p)}$ equipped with the holomorphic coordinates  $\left(\widetilde X,\widetilde Y,\overrightarrow B^1,\cdots,\overrightarrow B^{n-s}\right)$ defined as follows.
\begin{equation}\label{rulu}
\widetilde X:=\left(\begin{matrix}
x_{1(s+l+1)}&\cdots &x_{1n}\\
\vdots&\ddots&\vdots\\
x_{l(s+l+1)}&\cdots &x_{ln}\\
\end{matrix}\right)\,\,\,\,{\rm and}\,\,\,\,\,\widetilde Y:=\left(\begin{matrix}
y_{(l+1)1}&\cdots& y_{(l+1)(s-p+l)}\\
\vdots&\ddots&\vdots\\
y_{p1}&\cdots& y_{p(s-p+l)}\\
\end{matrix}  \right)\,;\,\,
\end{equation}
for $1\leq k\leq n-s-l$,
\begin{equation}\label{rqbp}
\begin{split}
&\overrightarrow B^{k}:=\left(a_{i_{k}j_{k}},\xi^{(k)}_{i_{k}(s+l+1)},\xi^{(k)}_{i_{k}(s+l+2)},\cdots,\widehat{\xi^{(k)}_{i_{k}j_1}},\cdots,\widehat{\xi^{(k)}_{i_{k}j_2}},\cdots,\widehat{\xi^{(k)}_{i_{k}j_{k}}},\cdots,\xi^{(k)}_{i_{k}n},\right.\\
&\,\,\,\,\,\,\,\,\,\,\,\,\,\,\,\,\,\,\,\,\,\,\,\,\,\,\,\,\left.\xi^{(k)}_{(l+1)j_{k}},\xi^{(k)}_{(l+2)j_{k}},\cdots,\widehat{\xi^{(k)}_{i_1j_{k}}},\cdots,\widehat{\xi^{(k)}_{i_2j_{k}}},\cdots,\cdots,\widehat{\xi^{(k)}_{i_{k}j_{k}}},\cdots,\xi^{(k)}_{pj_{k}}\right)\,;\\
\end{split}
\end{equation}
for $n-s-l+1\leq k\leq n-s$,
\begin{equation}\label{rqbpb}
\begin{split}
&\overrightarrow B^{k}:=\left(a_{i_{k}j_{k}},\xi^{(k)}_{i_{k}1},\xi^{(k)}_{i_{k}2},\cdots,\widehat{\xi^{(k)}_{i_{k}j_{n-s-l+1}}},\cdots,\widehat{\xi^{(k)}_{i_{k}j_{n-s-l+2}}},\cdots,\widehat{\xi^{(k)}_{i_{k}j_{k}}},\cdots,\xi^{(k)}_{i_{k}(s-p+l)},\right.\\
&\,\,\,\,\,\,\,\,\,\,\,\,\,\,\,\,\,\,\,\,\,\,\,\,\,\,\,\,\left.\xi^{(k)}_{1j_{k}},\xi^{(k)}_{2j_{k}},\cdots,\widehat{\xi^{(k)}_{i_{n-s-l+1}j_{k}}},\cdots,\widehat{\xi^{(k)}_{i_{n-s-l+2}j_{k}}},\cdots,\cdots,\widehat{\xi^{(k)}_{i_{k}j_{k}}},\cdots,\xi^{(k)}_{lj_{k}}\right)\,.\,\,\,\,\,\,\,\,\,\,\,\,\,\,\\
\end{split}
\end{equation}

Define a holomorphic map $\Gamma_l^{\tau}:\mathbb C^{p(n-p)}\rightarrow U_l$ by
\begin{equation}\label{gammalr}
\small
\begin{split}
&\,\,\,\,\,\,\,\,\,\,\,\,\,\,\,\,\,\,\,\,\,\,\,\,\,\,\,\,\,\,\,\,\,\,\,\,\,\,\,\,\,\,\,\,\,\,\,\,\,\,\,\,\,\,\,\,\,\,\,\,\,\,\,\,\,\Gamma_l^{\tau}\left(\widetilde X,\widetilde Y,\overrightarrow B^1,\cdots,\overrightarrow B^{n-s}\right):=\\
&\left(
\begin{matrix}
\sum\limits_{k=n-s-l+1}^{n-s}\left(\prod\limits_{t=n-s-l+1}^{k}a_{i_{t}j_t}\right)\cdot\Xi_k^T\cdot\Omega_k &0_{l\times(p-l)}&I_{l\times l}&\widetilde X\\ \widetilde Y&I_{(p-l)\times(p-l)}&0_{(p-l)\times l}&\sum\limits_{k=1}^{n-s-l}\left(\prod\limits_{t=1}^{k}b_{i_tj_t}\right)\cdot\Xi_k^T\cdot\Omega_k\\
\end{matrix}\right)\,.       
\end{split}
\end{equation}
where $\Xi_k$ and $\Omega_k$ are defined as follows. For $1\leq k\leq n-s-l$, $\Xi_k:=
\left(v_{l+1}^k,\cdots,v_p^k\right)$ where
\begin{equation}\label{lw1}
     v_t^k=\left\{\begin{array}{ll}
    \xi^{(k)}_{tj_k} \,\,\,\,\,\,\,\,\,\,\,\,\,\,\,\,\,\,\,\,\,\,\,\,\,\,\,\,\,\, & t\in\{l+1,l+2,\cdots,p\}\backslash\{i_1,i_2,\cdots,i_k\}\,\,\,\,\,\,\,\,\,\,\,\,\,\,\,\,\, \,\,\,\,\,\,\,\,\,\,\,\,\,\,\,\,\,\\
    0& t\in\{i_1,i_2,\cdots,i_{k-1}\}\,\,\,\,\,\,\, \\
  
    1 \,\,\,\,\,\,\,\,\,\,\,\,\,\,\,\,\,\,\,\,\,\,\,\,\,\,\,\,\,\, &t=i_k\,\,\,\,\,\,\,\,\,\,\,\,\,\,\,\,\, \,\,\,\,\,\,\,\,\,\,\,\,\,\,\,\,\,\\
    \end{array}\right.,
\end{equation}
and $\Omega_k:=
\left(w_{s+l+1}^k,\cdots,w_n^k\right)$ where
\begin{equation}\label{lw2}
     w_t^k=\left\{\begin{array}{ll}
    \xi^{(k)}_{i_kt} \,\,\,\,\,\,\,\,\,\,\,\,\,\,\,\,\,\,\,\,\,\,\,\,\,\,\,\,\,\, & t\in\{s+l+1,s+l+2,\cdots,n\}\backslash\{j_1,j_2,\cdots,j_k\}\,\,\\
    0& t\in\{j_1,j_2,\cdots,j_{k-1}\}\,\,\,\,\,\,\, \\
  
    1 \,\,\,\,\,\,\,\,\,\,\,\,\,\,\,\,\,\,\,\,\,\,\,\,\,\,\,\,\,\, &t=j_k\,\,\,\,\,\,\,\,\,\\
    \end{array}\right..
\end{equation}
For $n-s-l+1\leq k\leq n-s$, $\Xi_k:=
\left(v_{1}^k,\cdots,v_l^k\right)$ where
\begin{equation}\label{lw3}
     v_t^k=\left\{\begin{array}{ll}
    \xi^{(k)}_{tj_k} \,\,\,\,\,\,\,\,\,\,\,\,\,\,\,\,\,\,\,\,\,\,\,\,\,\,\,\,\,\, & t\in\{1,2,\cdots,l\}\backslash\{i_{n-s-l+1},i_{n-s-l+2},\cdots,i_k\}\,\,\,\,\,\,\,\,\,\,\,\,\,\,\,\,\, \,\,\,\,\,\,\,\,\,\,\,\,\,\,\,\,\,\\
    0& t\in\{i_{n-s-l+1},i_{n-s-l+2},\cdots,i_{k-1}\}\,\,\,\,\,\,\, \\
  
    1 \,\,\,\,\,\,\,\,\,\,\,\,\,\,\,\,\,\,\,\,\,\,\,\,\,\,\,\,\,\, &t=i_k\,\,\,\,\,\,\,\,\,\,\,\,\,\,\,\,\, \,\,\,\,\,\,\,\,\,\,\,\,\,\,\,\,\,\\
    \end{array}\right.,
\end{equation}
and $\Omega_k:=
\left(w_{1}^k,\cdots,w_{s-p+l}^k\right)$ where
\begin{equation}\label{lw4}
     w_t^k=\left\{\begin{array}{ll}
    \xi^{(k)}_{i_kt} \,\,\,\,\,\,\,\,\,\,\,\,\,\,\,\,\,\,\,\,\,\,\,\,\,\,\,\,\,\, & t\in\{1,2,\cdots,s-p+l\}\backslash\{j_{n-s-l+1},j_{n-s-l+2},\cdots,j_k\} \,\,\,\,\,\,\,\,\,\,\,\,\,\,\,\,\,\\
    0& t\in\{j_{n-s-l+1},j_{n-s-l+2},\cdots,j_{k-1}\}\,\,\,\,\,\,\, \\
  
    1 \,\,\,\,\,\,\,\,\,\,\,\,\,\,\,\,\,\,\,\,\,\,\,\,\,\,\,\,\,\, &t=j_k\,\,\,\,\,\,\,\,\,\,\,\,\,\\
    \end{array}\right..
\end{equation}

We can define a holomorphic embedding $J_l^{\tau}:\mathbb C^{p(n-p)}\rightarrow\mathbb {CP}^{N_{p,n}}\times\mathbb {CP}^{N^0_{s,p,n}}\times\cdots\times\mathbb {CP}^{N^p_{s,p,n}}$ by \begin{equation}
J_l^{\tau}:=\mathcal K_{s,p,n}\circ \Gamma_l^{\tau} \,.  
\end{equation}
Let $A^{\tau}$ be the image of $\mathbb C^{p(n-p)}$ under $J_l^{\tau}$.  Moreover, we can show that $\left\{\left(A^{\tau},(J^{\tau})^{-1}\right)\right\}_{\tau\in\mathbb J_l}$  is a holomorphic atlas for $R^{-1}_{s,p,n}(U_l)$, $0\leq l\leq n-s$, similarly to Section \ref{vanderl}.
\medskip

We end up this section by showing that $\mathcal T_{s,p,n}$ is smooth.
\smallskip

{\noindent\bf Proof of  Proposition \ref{tspnsmooth}.}  
By Remark \ref{dualusd} we can assume that $2p\leq n\neq 2s$. Let $S$ be the set  of the singular points of $\mathcal T_{s,p,n}$. Suppose that $S$ is non-empty. Since $S$ is invariant under the $GL(s,\mathbb C)\times GL(n-s,\mathbb C)$-action, we can conclude that
 \begin{equation}
S\bigcap\left(\bigcup_{l=0}^rR_{s,p,n}^{-1}(U_l)\right)\neq\emptyset\,
\end{equation}
where $r=\min\{p,n-s\}$.
However, $R_{s,p,n}^{-1}(U_l)$, $0\leq l\leq r$ is smooth by  the Van der Waerden representation constructed in Sections \ref{vanderl}, \ref{vandern}. This is a a contradiction.

We thus complete the proof of Proposition \ref{tspnsmooth}.
\,\,\,\,$\endpf$

\section{Foliation on \texorpdfstring{$\mathcal T_{s,p,n}$}{hh}} \label{foliation}
In this section, we recall the Bia{\l}ynicki-Birula decomposition (\cite{Bi}), and then study dynamic behavior of the $\mathbb C^*$-actions $\psi_{s,p,n}$ and $\Psi_{s,p,n}$ on $G(p,n)$ and $\mathcal T_{s,p,n}$ respectively (see Section \ref{group} for the definition of $\psi_{s,p,n}$ and $\Psi_{s,p,n}$). 

Without loss of generality, we may assume that  $2p\leq n\leq 2s$.
\subsection{The Bia{\l}ynicki-Birula decomposition} \label{foliationbb}
Recall the following notation in \cite{Bi}.  Let $V$ be a $\mathbb C^*$-module. Denote by $V_0$ the $\mathbb C^*$-submodule composed of all vectors $v\in V$ such that  $\lambda v = v$ for $\lambda\in\mathbb C^*$; denote by $V^+$ (resp. $V^-$) the $\mathbb C^*$-submodule spanned by all vectors $v\in V$ such that for $\lambda\in\mathbb C^*$ the image of $v$ under $\lambda$ equals $\lambda^mV$, for some positive integer $m$ (resp. negative integer $m$). It follows from the above that then $V = V^0\oplus V^+\oplus V^-$.

We state the  Bia{\l}ynicki-Birula  decomposition theorem over $\mathbb C$ as follows.

\begin{theorem}[\cite{Bi}, Theorems 4.1, 4.2,  4.3]\label{bb} Let $\mathcal X$ be a complete smooth complex manifold with an algebraic action of $\mathbb C^*$.  Denote by $\mathcal X^{\mathbb C^*}$ the set of the fixed points of $\mathcal X$ under the $\mathbb C^*$-action; let  $\mathcal X^{\mathbb C^*}=\bigcup_{l=0}^r  (\mathcal X^{\mathbb C^*})_l$ be the decomposition of $\mathcal X^{\mathbb C^*}$ into connected components. Then, there exists a unique locally closed
$\mathbb C^*$-invariant decomposition of $\mathcal X$, 
\begin{equation}
  \mathcal X=\bigsqcup_{l=0}^r  (\mathcal X^+)_l\,\,\,\,\left(resp.\,\,\mathcal X=\bigsqcup_{l=0}^r  (\mathcal X^-)_l\right)  
\end{equation}
and morphisms $\gamma_l^+: \mathcal X_l^+\rightarrow  (\mathcal  X^{\mathbb C^*})_l$ (resp. $\gamma_l^-: \mathcal X_l^-\rightarrow  (\mathcal X^{\mathbb C^*})_l$),  $0\leq l\leq r$, such that the following holds.
\begin{enumerate}[label=(\alph*).]
    \item  $\mathcal X_l^+$ and $\mathcal X_l^-$ are smooth complex manifolds. $(\mathcal X^{\mathbb C^*})_l$ is a closed complex submanifold of $\mathcal X_l^+$ and $\mathcal X_l^-$. $\mathcal X_l^+\cap \mathcal X_l^- =(\mathcal X^{\mathbb C^*})_l$.
	\item  $\gamma_l^{\pm}$ is a $\mathbb C^*$-fibration over $(\mathcal X^{\mathbb C^*})_l$ such that  $\gamma_l^{\pm}\big|(\mathcal X^{\mathbb C^*})_l$  is the identity and $(\mathcal X_l^{\pm})^{\mathbb C^*}=(\mathcal X^{\mathbb C^*})_l$.
	\item For each $\mathfrak a\in(\mathcal X^{\mathbb C^*})_l$, the tangent space $T_{\mathfrak a}(\mathcal X_l^{\pm})=T_{\mathfrak a}(\mathcal X)^0\oplus T_{\mathfrak a}(\mathcal X)^{\pm}$, and the dimension of the fibration $\gamma_l^{\pm}$  equals $\dim T_{\mathfrak a}(\mathcal X)^{\pm}$.
\end{enumerate}	
We call $\mathcal X^+_l$ (resp. $\mathcal X^-_l$), $0\leq l\leq r$, a stable (resp. unstable) manifold of $\mathcal X$ with respect to the $\mathbb C^*$-action.
\end{theorem} 

\begin{definition}
We say that $(\mathcal X^{\mathbb C^*})_l$ has $p$ {\it positive directions} (resp. $n$ {\it negative directions}) if there is a point $\mathfrak a\in(\mathcal X^{\mathbb C^*})_l$ such that $\dim T_{\mathfrak a}(\mathcal X)^+=p$ (resp. $\dim T_{\mathfrak a}(\mathcal X)^-=n$). The definition is independent of the choice of the point $\mathfrak a$ by Theorem \ref{bb}.
\end{definition}
\begin{definition}[\cite{BS}]\label{ss}
	We call $\mathcal X_l^{\mathbb C^*}$ a {\it source} (resp. {\it sink}) if  it has no negative direction (resp. positive direction).  By the Bia{\l}ynicki-Birula decomposition $\mathcal X$ has exact one  source and one sink. 
\end{definition}

\begin{definition}\label{order}
Define a function $\Lambda^+:\mathcal X\rightarrow \mathcal X^{\mathbb C^*}$ (resp. $\Lambda^-:\mathcal X\rightarrow \mathcal X^{\mathbb C^*}$)  by
\begin{equation}\label{pl}
\Lambda^+(x)=\lim_{t\rightarrow 0}t\cdot x\,\,\,\,\,\,({\rm resp}.\,\,\Lambda^-(x)=\lim_{t\rightarrow\infty}t\cdot x).
\end{equation}
Call a component $(\mathcal X^{\mathbb C^*})_i$ {\it lower than} a component $(\mathcal X^{\mathbb C^*})_j$ if $i\neq j$ and there is a point $x\in\mathcal X$ such that
\begin{equation}
    \Lambda^+(x)\in (\mathcal X^{\mathbb C^*})_i\,\,\,\,\,\,{\rm and}\,\,\,\,\,\, \Lambda^-(x)\in (\mathcal X^{\mathbb C^*})_j\,;
\end{equation}
write $(\mathcal X^{\mathbb C^*})_i<(\mathcal X^{\mathbb C^*})_j$. 
\end{definition}

In what follows, we will derive an explicit Bia{\l}ynicki-Birula decomposition for Grassmann manifolds $G(p,n)$. Define an integer $r$ by
\begin{equation}\label{rank2}
    r=\min\left\{s,n-s,p,n-p\right\}.
\end{equation}
For $0\leq l\leq r$, define subsets $\mathcal V_{(p-l,l)}$, $\mathcal V_{(p-l,l)}^+$ and $\mathcal V_{(p-l,l)}^-$ of $G(p,n)$ in matrix representatives by
\begin{equation}\label{vani}
\begin{split}
&\mathcal V_{(p-l,l)} :=\left\{\left. \left(
\begin{matrix}
0&X\\
Y&0\\
\end{matrix}\right)\right\vert_{}\footnotesize\begin{matrix}
X\,\,{\rm is\,\,an\,\,}l\times (n-s)\,\,{\rm matrix\,\,of\,\,rank}\,\,l\,;\\
Y\,\,{\rm is\,\,a\,\,}(p-l)\times s\,\,{\rm matrix\,\,of\,\,rank}\,\,(p-l)\\
\end{matrix}
\right\}\,,
\end{split}
\end{equation}
and
\begin{equation}\label{stab}
\begin{split}
&\mathcal V_{(p-l,l)}^+:= \left\{\left.\left(
\begin{matrix}
0&X\\
Y&W\\
\end{matrix}\right)\right\vert_{}\footnotesize\begin{matrix}
X\,\,{\rm is\,\,an\,\,}l\times (n-s)\,\,{\rm matrix \,\,of\,\,rank\,\,}l\,;\\
Y\,\,{\rm is\,\,a\,\,}(p-l)\times s\,\,{\rm matrix \,\,of\,\,rank\,\,}(p-l)\,\\
\end{matrix}\right\},\\
&\mathcal V_{(p-l,l)}^-:= \left\{\left.\left(
\begin{matrix}
Z&X\\
Y&0\\
\end{matrix}\right)\right\vert_{}{\footnotesize\begin{matrix}
X\,\,{\rm is\,\,an\,\,}l\times (n-s)\,\,{\rm matrix \,\,of\,\,rank\,\,}l\,;\\
Y\,\,{\rm is\,\,a\,\,}(p-l)\times s\,\,{\rm matrix \,\,of\,\,rank\,\,}(p-l)\,\\
\end{matrix}}\right\}\,.    
\end{split}
\end{equation}

\begin{lemma}\label{EVS}
Let $\psi_{s,p,n}$ be the $\mathbb C^*$-action on $G(p,n)$ defined in Section \ref{group}. 
Then the connected components of the  set of the fixed points of $G(p,n)$ under  $\psi_{s,p,n}$ are
\begin{equation}
   \mathcal V_{(p,0)},\mathcal V_{(p-1,1)},\cdots,\mathcal V_{(p-r,r)}\,.
\end{equation}
The corresponding stable (resp. unstable) manifolds are 
\begin{equation}
   \mathcal V^+_{(p,0)},\mathcal V^+_{(p-1,1)},\cdots,\mathcal V^+_{(p-r,r)}\,\,(resp.\,\,\mathcal V^-_{(p,0)},\mathcal V^-_{(p-1,1)},\cdots,\mathcal V^-_{(p-r,r)}\,\,)\,.
\end{equation}
Moreover, for $0\leq l\leq r$,
\begin{equation}\label{schubert}
\overline{\mathcal V_{(p-l,l)}^+}\mathbin{\scaleobj{1.7}{\backslash}} \mathcal V_{(p-l,l)}^+=\bigsqcup_{k=l+1}^r\mathcal V_{(p-k,k)}^+\,\,\,\,\,\,\,\,{\rm and}\,\,\,\,\,\,\,\overline{\mathcal V_{(p-l,l)}^-}\mathbin{\scaleobj{1.7}{\backslash}} \mathcal V_{(p-l,l)}^-=\bigsqcup_{k=0}^{l-1}\mathcal V_{(p-k,k)}^-\,,
\end{equation}
where $\overline{\mathcal V_{(p-l,l)}^{\pm}}$ is the Zariski closure of $\mathcal V_{(p-l,l)}^{\pm}$.
\end{lemma}

{\noindent\bf Proof of Lemma \ref{EVS}.}
For simplicity, we only prove for the case $p\leq n-s$ since the proof when $n-s<p$ is the same.

It is clear that the point of $\mathcal V_{(p-l,l)}$, $0\leq l\leq r$, is fixed by $\psi_{s,p,n}$. In the following, we will show that each fixed point $x\in G(p,n)$ of $\psi_{s,p,n}$ belongs a certain $\mathcal V_{(p-l,l)}$.

Let $\widetilde x$ be a matrix representative of $x$ such that
\begin{equation}    \widetilde x=\Big( \widetilde Y \hspace{-.1in}\begin{matrix} &\hfill\tikzmark{c}\\ &\hfill\tikzmark{d} \end{matrix} \hspace{.12in} \widetilde X\Big)
\tikz[remember picture,overlay]   \draw[dashed,dash pattern={on 4pt off 2pt}] ([xshift=0.5\tabcolsep,yshift=6pt]c.north) -- ([xshift=0.5\tabcolsep,yshift=3pt]d.south);
\end{equation}
where $\widetilde Y$ is a $p\times s$ matrix of rank  $p-l$ for a certain $0\leq l\leq p$.

If  $l=p$,  $\widetilde Y$ is a zero matrix and hence $x\in\mathcal V_{(0,p)}$. If $l\leq p-1$,  by a certain $GL(p,\mathbb C)$-action from left we can assume that
\begin{equation}
    \widetilde x=\left(\begin{matrix}
    0&X\\
    Y&W\\
    \end{matrix}\right)
    \end{equation}
where  $X$ is a $l\times (n-s)$ matrix of rank $l$ and $Y$ is a $(p-l)\times s$ matrix of rank $p-l$. Since $\psi_{s,p,n}(\lambda)\cdot x=x$ for each $\lambda\in\mathbb C^*$, we can conclude that $W$ is a null matrix.
Therefore  $x\in\mathcal V_{(p-l,l)}$. 

We can verify that the submanifolds in $(\ref{stab})$ are the corresponding stable and unstable manifolds of $\mathcal V_{(p-l,l)}$ in the sense of Theorem \ref{bb}, where the fibration is given by projecting $W$ and $Z$ to null matrices.

Decomposing $\mathcal V_{(p-k,k)}$ into Schubert cells, we can show that (\ref{schubert}) holds.

We conclude Lemma \ref{EVS}.\,\,\,\,
$\endpf$

\begin{remark}
The integer $r$ defined in (\ref{rank2}) coincides with the rank $r$ of $\mathcal T_{s,p,n}$ defined in Definition \ref{rank}. We make the convention that in the remainder of the paper $r$ is always taken to be the rank of $\mathcal T_{s,p,n}$.
\end{remark}

\begin{remark}\label{less}
According to Definition \ref{ss}, $\mathcal V_{(p,0)}$ is the source and $\mathcal V_{(p-r,r)}$ is the sink. It is clear that
\begin{equation}
     \mathcal V_{(p,0)}<\mathcal V_{(p-1,1)}<\cdots<\mathcal V_{(p-l,l)}<\cdots<\mathcal V_{(p-r+1,r-1)}<\mathcal V_{(p-r,r)}.
 \end{equation}
\end{remark} 
\smallskip

Recall that the $\mathbb C^*$-action $\psi_{s,p,n}(\lambda)$, $\lambda\in\mathbb C^*$,  induces a family of biholomorphisms of $G(p,n)$. We can derive a holomorphic vector field $e_{s,p,n}$  on $G(p,n)$ by taking the differential of $\psi_{s,p,n}(\lambda)$  with respect to $\lambda\in\mathbb C^*$ at $\lambda=1$.

We have the following normal form of the vector field $e_{s,p,n}$ near fixed points.

\begin{lemma}\label{nf}
For each $\mathfrak a\in\mathcal V_{(p-l,l)}$, $0\leq l\leq r$, there is a  holomorphic coordinate chart for a certain neighborhood of $\mathfrak a$ such that  $e_{s,p,n}$ takes the following form in the corresponding local coordinates.
\begin{equation}
e_{s,p,n}=-\sum_{i=1}^l\sum_{j=1}^{s-p+l}z_{ij}\frac{\partial}{\partial z_{ij}}+\sum_{i=l+1}^{p}\sum_{j=s+l+1}^{n}w_{ij}\frac{\partial}{\partial w_{ij}}\,.
\end{equation}
\end{lemma}
{\noindent\bf Proof of Lemma \ref{nf}.}
Without loss of generality, we can assume that $\mathfrak a\in U_l$ where $U_l$ is defined by (\ref{ul}).  It is clear that $\psi_{s,p,n}(\lambda)$ takes the following form in terms of the local coordinates defined in (\ref{ulx}). 
\begin{equation}
\psi_{s,p,n}(\lambda)\left(\left(\,\,\,\begin{matrix}
Z&0&I_{l\times l}&X\\
Y&I_{(p-l)\times(p-l)}&0&W\\
\end{matrix}\right)\right)=\left(\left(\,\,\,\begin{matrix}
\frac{1}{\lambda}\cdot Z&0&I_{l\times l}&X\\
Y&I_{(p-l)\times(p-l)}&0&\lambda\cdot W\\
\end{matrix}\right)\right)\,.
\end{equation}
Taking the differential with respect to $\lambda$ at $\lambda=1$, we conclude Lemma \ref{nf}.\,\,\,\,$\endpf$
\medskip

Since the $\mathbb C^*$-action $\Psi_{s,p,n}(\lambda)$, $\lambda\in\mathbb C^*$, induces a family of biholomorphisms of $\mathcal T_{s,p,n}$, we can define a vector field $E_{s,p,n}$  on $\mathcal T_{s,p,n}$ by taking the differential of $\Psi_{s,p,n}(\lambda)$ at $\lambda=1$.

We have the following explicit  Bia{\l}ynicki-Birula decomposition for $\mathcal T_{s,p,n}$.

\begin{lemma}\label{11} 
There are $r+1$ connected components $\mathcal D_{(p-l,l)}$, $0\leq l\leq r$, of the  set of the fixed points of $\mathcal T_{s,p,n}$ under  the $\mathbb C^*$-action $\Psi_{s,p,n}$, such that  the following holds.
\begin{enumerate}[label=(\alph*).]
\item   $R_{s,p,n}(\mathcal D_{(p-l,l)})=\mathcal V_{(p-l,l)}$, $0\leq l\leq r$. 
		
\item $\mathcal D_{(p,0)}$ and $\mathcal D_{(p-r,r)}$ are smooth divisors of $\mathcal T_{s,p,n}$.  $\mathcal D_{(p-l,l)}$, $1\leq l\leq r-1$,  is a closed complex manifold of codimension two in $\mathcal T_{s,p,n}$.  
		
\item Denote the corresponding stable (resp. unstable) manifolds by  $\mathcal D^+_{(p-l,l)}$ (resp. $\mathcal D^-_{(p-l,l)}$). Then for $0\leq l\leq r$,
\begin{equation}\label{tschubert}
\overline{\mathcal D_{(p-l,l)}^+}\mathbin{\scaleobj{1.7}{\backslash}} \mathcal D_{(p-l,l)}^+=\bigsqcup_{k=l+1}^r\mathcal D_{(p-k,k)}^+\,\,\,\,\,\,\,\,{\rm and}\,\,\,\,\,\,\,\overline{\mathcal D_{(p-l,l)}^-}\mathbin{\scaleobj{1.7}{\backslash}} \mathcal D_{(p-l,l)}^-=\bigsqcup_{k=0}^{l-1}\mathcal D_{(p-k,k)}^-\,,
\end{equation}
	
\item The vector field $E_{s,p,n}$ takes the following form near the fixed points.
\begin{enumerate}[label=(\roman*)]
\item For each point $\mathfrak a\in\mathcal D_{(p,0)}$, there is a holomorphic coordinate chart with coordinates $(b,\cdots)$, such that $\mathfrak a=(0,0,\cdots,0)$ and  \begin{equation}
E_{s,p,n}=b\frac{\partial}{\partial b}\,.\,\,\,\,\,\,\,\,\,\,\,\,\,\,\,\,\,\,\,\,\,\,\,\end{equation}
		
\item For each point $\mathfrak a\in\mathcal D_{(p-l,l)}$, $1\leq l\leq r-1$, there is a holomorphic coordinate chart with coordinates $(b,a,\cdots)$, such that $\mathfrak a=(0,0,\cdots,0)$ and  \begin{equation}
E_{s,p,n}= -a\frac{\partial}{\partial a} +b\frac{\partial}{\partial b}\,.
\end{equation}
	
\item For each point $\mathfrak a\in\mathcal D_{(p-r,r)}$, there is a holomorphic coordinate chart with coordinates $(a,\cdots)$, such that $\mathfrak a=(0,0,\cdots,0)$ and  \begin{equation}
E_{s,p,n}=-a\frac{\partial}{\partial a}\,.\,\,\,\,\,\,\,\,\,\,\,\,\,\,\,\,\,\,\end{equation}
 \end{enumerate}
\end{enumerate}
\end{lemma}

{\noindent\bf Proof of Lemma \ref{11}.}
For $0\leq l\leq r$, let $\left\{\left(A^{\tau},(J_l^{\tau})^{-1}\right)\right\}_{\tau\in\mathbb J_l}$ be the Van der Waerden representation of $R_{s,p,n}^{-1}(U_l)$  where  $U_l\subset G(p,n)$ is  defined by (\ref{ul}). 

For each local coordinate chart $\left(A^{\tau},(J_l^{\tau})^{-1}\right)$, $\tau\in\mathbb J_l$, the $\mathbb C^*$-action $\Psi_{s,p,n}(\lambda)$ takes the following form, in terms of the holomorphic coordinates $\left(\widetilde X,\widetilde Y,\overrightarrow B^1,\cdots,\overrightarrow B^r\right)$ defined by (\ref{ulu}), (\ref{qbp}), (\ref{qbpb}) (or (\ref{rulu}), (\ref{rqbp}), (\ref{rqbpb}) when $n-s<p$). When $1\leq l\leq r-1$,
\begin{equation}\label{ddil}
\begin{split}
 \Psi_{s,p,n}(\lambda)&\left( \widetilde X,\widetilde Y,\overrightarrow B^1,\cdots,\overrightarrow B^r\right)=\left(\frac{1}{\lambda}\cdot a_{i_{r-l+1}j_{r-l+1}},\lambda\cdot b_{i_1j_1},\xi^{(1)}_{i_1(s+l+1)},\cdots,\widehat{\xi^{(1)}_{i_1j_1}},\cdots,\xi^{(1)}_{i_1n},\right.\\
&\xi^{(1)}_{(l+1)j_1},\cdots,\widehat{\xi^{(1)}_{i_1j_1}},\cdots,\xi^{(1)}_{pj_1},\xi^{(r-l+1)}_{i_{r-l+1}1},\cdots,\widehat{\xi^{(p-l+1)}_{i_{r-l+1}j_{r-l+1}}},\cdots,\xi^{(r-l+1)}_{i_{r-l+1}(s-p+l)},\\
&\left.\xi^{(r-l+1)}_{1j_{r-l+1}},\cdots,\widehat{\xi^{(r-l+1)}_{i_{r-l+1}j_{r-l+1}}},\cdots,\xi^{(r-l+1)}_{lj_{r-l+1}},\widetilde X,\widetilde Y,\overrightarrow B^2,\cdots,\overrightarrow B^{r-l},\overrightarrow B^{r-l+2},\cdots,\overrightarrow B^r\right).
\end{split}
\end{equation}
That is, $\Psi_{s,p,n}$ dilates $a_{i_{r-l+1}j_{r-l+1}}$ and $b_{i_ij_1}$ by $\frac{1}{\lambda}$ and $\lambda$ respectively. Similarly, when $l=0$, $\Psi_{s,p,n}$ dilates $b_{i_ij_1}$ by $\lambda$. When $l=r$, $\Psi_{s,p,n}$ dilates $a_{i_{r-l+1}j_{r-l+1}}$ by $\frac{1}{\lambda}$.

Applying the $GL(s,\mathbb C)\times GL(n-s,\mathbb C)$ of $\mathcal T_{s,p,n}$, Lemma \ref{11} follows similarly to Lemmas \ref{EVS} and \ref{nf}.\,\,\,\,\,$\endpf$
\medskip

We introduce certain important divisors of $\mathcal T_{s,p,n}$ in the following.
\begin{definition}\label{div}
For $1\leq k\leq r$, denote by $D_{k}^-$ (resp. $D_{k}^+$) the Zariski closure of the unstable (resp. stable) manifold $\mathcal D_{(p-k+1,k-1)}^-$  $\left({\rm resp}.\,\,\mathcal D_{(p-r+k-1,r-k+1)}^+\right)$  of $\mathcal T_{s,p,n}$. 
\end{definition}
\begin{remark}
It is clear that $D_1^-=\mathcal D_{(p,0)}$ and $D_1^+=\mathcal D_{(p-r,r)}$. We call $D_1^-$ and $D_1^+$ the  source and the sink respectively according to Definition \ref{ss}.
\end{remark}

\begin{lemma}\label{snc}
$D_{1}^+$, $D_{2}^+$, $\cdots$, $D_{r}^+$, $D_{1}^-$, $D_{2}^-$, $\cdots$, $D_{r}^-$ are distinct smooth divisors of $\mathcal T_{s,p,n}$. The following divisor is simple normal crossing. 
\begin{equation}
   D_{1}^++D_{2}^++\cdots+D_{r}^++D_{1}^-+D_{2}^-+\cdots+D_{r}^-\,.
\end{equation}
\end{lemma}

{\noindent\bf Proof of  Lemma \ref{snc}.} 
Let $\left(A^{\tau},(J_l^{\tau})^{-1}\right)$ be an arbitrary holomorphic coordinate chart in the Van der Waerden representation where $0\leq l\leq r$ and $\tau=\left(
	\begin{matrix}
	i_1&i_2&\cdots&i_{p-l}&\cdots&i_p\\
	j_1&j_2&\cdots&j_{p-l}&\cdots&j_p\\
	\end{matrix}\right)\in\mathbb J_l$.
Combining Lemma \ref{EVS} and Property (a) in Lemma \ref{11},  we can determine the stable and unstable manifolds explicitly as follows.
\begin{equation}
    \begin{split}
        &D_k^+\bigcap A^{\tau} =\left\{a_{i_{k}j_{k}}=0\right\}\,\,\,{\rm when}\,\, r-l+1\leq k\leq r\,;D_k^+\bigcap A^{\tau} =\emptyset\,\,\,{\rm when}\,\, 1\leq k\leq r-l\,.\\
        &D_k^-\bigcap A^{\tau} =\left\{b_{i_{k-l}j_{k-l}}=0\right\}\,\,\,{\rm when}\,\, l+1\leq k\leq r\,;D_k^-\bigcap A^{\tau} =\emptyset\,,\,\,{\rm when}\,\, 1\leq k\leq l\,.\\
    \end{split}
\end{equation}

Applying the $GL(s,\mathbb C)\times GL(n-s,\mathbb C)$-action of $\mathcal T_{s,p,n}$, Lemma \ref{snc} follows. \,\,\,$\endpf$

\subsection{Properties of the exceptional divisor} \label{foliationex}
In this subsection  we will relate the divisors $D_{k}^{\pm}$  to the exceptional divisor of the canonical blow-up $R_{s,p,n}:\mathcal T_{s,p,n}\rightarrow G(p,n)$.

We characterize the center of $R_{s,p,n}$ as follows.    Let $S_k\subset G(p,n)$ be the subscheme in Definition \ref{sk}, and $\overline{\mathcal V_{(p-k-1,k+1)}^{\pm}}$  the Zariski closure of $\mathcal V_{(p-k-1,k+1)}^{\pm}$ in $G(p,n)$. Under the convention that $\mathcal V_{(p-r-1,r+1)}^+=\mathcal V_{(p+1,-1)}^-=\emptyset$, we have that
\begin{lemma}\label{kl}For $0\leq k\leq r$, $S_k$ is the  scheme-theoretic union of  $\overline{\mathcal V_{(p-k-1,k+1)}^+}$ and  $\overline{\mathcal V_{(p-k+1,k-1)}^-}$. 
\end{lemma}
{\noindent\bf Proof of  Lemma \ref{kl}.}
It is clear that $\overline{\mathcal V_{(p-k-1,k+1)}^+}$ and  $\overline{\mathcal V_{(p-k+1,k-1)}^-}$ are Schubert cycles. Therefore,  we can conclude by Theorem 1.4 in \cite{BV} that  the ideal sheaf of  $\overline{\mathcal V_{(p-k-1,k+1)}^+}$  in $G(p,n)$ $\left({\rm resp}.\,\,\overline{\mathcal V_{(p-k+1,k-1)}^-}\,\,{\rm in}\,\,G(p,n)\right)$ is generated by  
\begin{equation}\label{pine}
  PL_k^+:=\bigcup_{0\leq j\leq k}\left\{P_I\big|I\in \mathbb I^j_{s,p,n}\right\}\,\,\,\,\left({\rm resp.}\,\,\, PL_k^-:=\bigcup_{k\leq j\leq r}\left\{P_I\big|I\in \mathbb I^j_{s,p,n}\right\} \right).  
\end{equation}

It is clear that $\overline{\mathcal V_{(p-k-1,k+1)}^+}$ and  $\overline{\mathcal V_{(p-k+1,k-1)}^-}$ are disjoint. Then, to prove Lemma \ref{kl}, it suffices to show that the ideal sheaf $\mathcal S_k$ of $S_k$ coincides with the ideal sheaf generated by $PL_k^+$ (resp. $PL_k^-$) locally near $\overline{\mathcal V_{(p-k-1,k+1)}^+}$ $\left({\rm resp}.\overline{\mathcal V_{(p-k+1,k-1)}^-}\right)$.

Let $x\in\overline{\mathcal V_{(p-k-1,k+1)}^+}$\,.  By  Lemma \ref{EVS}, there is an integer $k+1\leq l\leq r$ such that $x\in \mathcal V_{(p-l,l)}^+$. Since the variety $\overline{\mathcal V_{(p-k-1,k+1)}^+}$ and the ideal sheaf $\mathcal S_j$, $0\leq j\leq r$,  are invariant under the $GL(s,\mathbb C)\times GL(n-s,\mathbb C)$-action of $G(p,n)$, we may further assume that $x\in U_l$ where
\begin{equation}\label{xul}
U_l=\left\{\left(\,\,\,\begin{matrix}
Z&0&I_{l\times l}&X\\
Y&I_{(p-l)\times(p-l)}&0&W\\
\end{matrix}\right)\right\}\,.
\end{equation}

The restriction $\mathcal S_j\big|_{U_l}$, $0\leq j\leq k$,  is an ideal of the coordinate ring of $U_l$.
Applying the Laplace expansion repeatedly,
we can show that $\mathcal S_j\big|_{U_l}$, $0\leq j\leq k$,  is generated by the $(l-j)\times(l-j)$ minors of the matrix $Z$ in (\ref{xul}). Then,
\begin{equation}
\mathcal S_j\big|_{U_l}\subset \mathcal S_k\big|_{U_l}\,\,\,{\rm for}\,\,\,\,0\leq j\leq k\,.
\end{equation}

The same result holds for $\overline{\mathcal V_{(p-k+1,k-1)}^+}$. We thus conclude Lemma \ref{kl}. $\endpf$

\begin{lemma}\label{excep}
Let $E$ be the exceptional divisor of the canonical blow-up $R_{s,p,n}:\mathcal T_{s,p,n}\rightarrow G(p,n)$. When $p<s$ and $n-s\neq p$,
\begin{equation}
    E=D_{1}^++D_{2}^++\cdots+D_{r}^++D_{1}^-+D_{2}^-+\cdots+D_{r}^-\,;
\end{equation}
when $n-s=p<s$,
\begin{equation}
    E=D_{1}^++D_{2}^++\cdots+D_{r}^++D_{1}^-+D_{2}^-+\cdots+D_{r-1}^-\,;
\end{equation}
when $n-s=p=s$,
\begin{equation}
    E=D_{1}^++D_{2}^++\cdots+D_{r-1}^++D_{1}^-+D_{2}^-+\cdots+D_{r-1}^-\,.
\end{equation}
\end{lemma}

{\noindent\bf Proof of Lemma \ref{excep}.} By Lemma \ref{kl}, the support of the center of the canonical blow-up $R_{s,p,n}$ is 
\begin{equation}
     \left(\bigcup_{k=1}^r\overline{\mathcal V_{(p-k,k)}^+}\right)\bigcup\left(\bigcup_{k=0}^{r-1}\overline{\mathcal V_{(p-k,k)}^-}\right).
\end{equation}
Therefore, by Property (a) in Lemma \ref{11} we can show that the support of $E$ is contained in
\begin{equation}
    \left(\bigcup_{k=1}^{r}D_{k}^+\right)\bigcup\left(\bigcup_{k=1}^{r}D_{k}^-\right).
\end{equation}

In what follows, we will compute the multiplicity of $E$ along $D_{l}^{\pm}$ in certain  holomorphic coordinate charts. The proof is based on a case by case argument.
\smallskip

{\noindent\bf Case I ($p<s$ and $n-s\neq p$).} Recall the iterated blow-up in Definition \ref{calt} with the permutation $\sigma:=\sigma_0$ defined by 
\begin{equation}
  \sigma_0(k)=\left\{  \begin{matrix}
        &k+1\,\,\,\,\,&{\rm when}\,\, &0\leq k\leq p-1\,\\
        & 0\,\,\,\,\,&{\rm when}\,\, &k=p
    \end{matrix}\right..
\end{equation}Restricting $(\ref{sblow})$ to $U_0$, we have that
\vspace{-.05in}
\begin{equation}\label{cblow1}
\footnotesize
\begin{tikzcd}
&Y^{0}_p\ar{r}{g^{\sigma_0}_p}&Y^{0}_{p-1}\ar{r}{g^{\sigma_0}_{p-1}}&\cdots\ar{r}{g^{\sigma_0}_1}&Y^{0}_0\ar{r}{g^{\sigma_0}_0}&U_0 \\
&&(g^{\sigma_0}_0\circ\cdots\circ g^{\sigma_0}_{p-1})^{-1}(S_{\sigma_0(p)}\cap U_0)\ar[hook]{u}&\cdots&(g^{\sigma_0}_0)^{-1}(S_{\sigma_0(1)}\cap U_0)\ar[hook]{u}&S_{\sigma_0(0)}\cap U_0\ar[hook]{u}\\
\end{tikzcd}\vspace{-20pt}
\end{equation}
where for $0\leq i\leq p$,
\begin{equation}
   Y^{0}_{i}:=Y_i^{\sigma_0}\cap (g^{\sigma_0}_0\circ\cdots\circ g^{\sigma_0}_{i})^{-1}(U_0)\subset  G(p,n)\times\mathbb {CP}^{N^0_{s,p,n}}\times\cdots\times\mathbb {CP}^{N^{i}_{s,p,n}}\,\,\,\,. 
\end{equation}
It is clear that $g_j^{\sigma_0}$ is an isomorphism for $r\leq j\leq p$.

For $1\leq k\leq r$, define an affine quasi-projective variety $d^-_{k}$ by 
\begin{equation}
 d^-_{k}:= \left\{\left. \left(Y\,\hspace{-0.15in}\begin{matrix}
  &\hfill\tikzmark{c3}\\
  &\hfill\tikzmark{d3}
  \end{matrix}\,\,\,\, I_{p\times p}\hspace{-0.13in}\begin{matrix}
  &\hfill\tikzmark{c4}\\
  &\hfill\tikzmark{d4}
  \end{matrix}\hspace{-0.11in}\begin{matrix}
  &\hfill\tikzmark{c5}\\
  &\hfill\tikzmark{d5}
  \end{matrix}\,\,\,\,\, W\right)\right\vert_{}\begin{matrix}
  \,Y\,\,{\rm  is\,\, a\,\,} p\times (s-p)\,\,{\rm matrix}\,;\\
  W\,\,{\rm  is\,\, a\,\,} p\times (n-s)\,\,{\rm matrix\,\,of\,\,rank\,\,}k-1\,\\
  \end{matrix}\right\}\,.
  \tikz[remember picture,overlay]   \draw[dashed,dash pattern={on 4pt off 2pt}] ([xshift=0.5\tabcolsep,yshift=7pt]c3.north) -- ([xshift=0.5\tabcolsep,yshift=-2pt]d3.south);\tikz[remember picture,overlay]   \draw[dashed,dash pattern={on 4pt off 2pt}] ([xshift=0.5\tabcolsep,yshift=7pt]c4.north) -- ([xshift=0.5\tabcolsep,yshift=-2pt]d4.south);\tikz[remember picture,overlay]   \draw[dashed,dash pattern={on 4pt off 2pt}] ([xshift=0.5\tabcolsep,yshift=7pt]c5.north) -- ([xshift=0.5\tabcolsep,yshift=-2pt]d5.south);
\end{equation}
Then $d_{k}^-$ is a dense open subset of $\mathcal V_{(p-k+1,k-1)}^-$. Moreover, by the Van der Waerden representation, we can show that $\overline{R^{-1}_{s,p,n}(d_{k}^-)}=D_{k}^-$ where   $\overline{R^{-1}_{s,p,n}(d_{k}^-)}$ is the Zariski closure of the inverse image of $d_{k}^-$.

By Lemma \ref{kl}, $\overline{\mathcal V_{(p,0)}^-}$ is the center of the blow-up $g^{\sigma_0}_{0}$; $(g^{\sigma_0}_0\circ\cdots\circ g^{\sigma_0}_{i})^{-1}\left(\overline{\mathcal V_{(p-i,i)}^-}\right)$ is the center of the blow-up $g^{\sigma_0}_{i}$ for $1\leq i\leq r-1$.
Notice that $(g^{\sigma_0}_0\circ\cdots\circ g^{\sigma_0}_{k-2})^{-1}(d_{k}^-)\cong d_{k}^-$ for $2\leq k\leq r$. Therefore, $(g^{\sigma_0}_0\circ\cdots\circ g^{\sigma_0}_{k-2})^{-1}(d_{k}^-)$  is the center of the blow-up $g^{\sigma_0}_{k-1}$ in a certain smaller open subset of $U_0$.  Since $d_k^-$ is smooth, we can conclude that $E$ has multiplicity $1$ along $D_{k}^-$ for $1\leq k\leq r$.

Next we will compute the multiplicity of $E$ along $D_l^+$. Let  $\sigma_p$ be a permutation such that $\sigma_p(j)=p-j$ for $0\leq j\leq p$; consider the following iterated blow-ups by restricting (\ref{sblow}) to $U_p$. 
\vspace{-.05in}
\begin{equation}\label{csblow}
\footnotesize
\begin{tikzcd}
&Y^p_p\ar{r}{g^{\sigma_p}_p}&Y^{p}_{p-1}\ar{r}{g^{\sigma_p}_{p-1}}&\cdots\ar{r}{g^{\sigma_p}_1}&Y^{p}_0\ar{r}{g^{\sigma_p}_0}&U_p \\
&&(g^{\sigma_p}_0\circ\cdots\circ g^{\sigma_p}_{p-1})^{-1}(S_0\cap U_p)\ar[hook]{u}&\cdots&(g^{\sigma_p}_0)^{-1}(S_{p-1}\cap U_p)\ar[hook]{u}&S_{p}\cap U_p\ar[hook]{u}\\
\end{tikzcd}\vspace{-20pt}
\end{equation}
where for $0\leq i\leq p$,
\begin{equation}
   Y^p_{i}:=Y_i^{\sigma_p}\cap (g^{\sigma_p}_0\circ\cdots\circ g^{\sigma_p}_{i})^{-1}(U_p)\subset  G(p,n)\times\mathbb {CP}^{N^p_{s,p,n}}\times\cdots\times\mathbb {CP}^{N^{p-i}_{s,p,n}}\,\,\,\,. 
\end{equation}

For $1\leq k\leq r$, define an affine quasi-projective variety $d^+_{k}$  as follows. When $r=p$, 
\begin{equation}
 d^+_{k}:=  \left\{\left. \left(Z\,\hspace{-0.15in}\begin{matrix}
  &\hfill\tikzmark{c3}\\
  &\hfill\tikzmark{d3}
  \end{matrix}\hspace{-0.11in}\begin{matrix}
  &\hfill\tikzmark{c5}\\
  &\hfill\tikzmark{d5}
  \end{matrix}\,\,\,\, I_{p\times p}\hspace{-0.13in}\begin{matrix}
  &\hfill\tikzmark{c4}\\
  &\hfill\tikzmark{d4}
  \end{matrix}\,\,\,\,\, X\right)\right\vert_{}\begin{matrix}
  X\,\,{\rm  is\,\, a\,\,} p\times (n-s-p)\,\,{\rm matrix}\,;\\
  \,Z\,\,{\rm  is\,\, a\,\,} p\times s\,\,{\rm matrix\,\,of\,\,rank\,\,}k-1\\
  \end{matrix}\right\}\,;
  \tikz[remember picture,overlay]   \draw[dashed,dash pattern={on 4pt off 2pt}] ([xshift=0.5\tabcolsep,yshift=7pt]c3.north) -- ([xshift=0.5\tabcolsep,yshift=-2pt]d3.south);\tikz[remember picture,overlay]   \draw[dashed,dash pattern={on 4pt off 2pt}] ([xshift=0.5\tabcolsep,yshift=7pt]c4.north) -- ([xshift=0.5\tabcolsep,yshift=-2pt]d4.south);\tikz[remember picture,overlay]   \draw[dashed,dash pattern={on 4pt off 2pt}] ([xshift=0.5\tabcolsep,yshift=7pt]c5.north) -- ([xshift=0.5\tabcolsep,yshift=-2pt]d5.south);
\end{equation}
when $r=n-s<p$,
\begin{equation}
 d^+_{k}:=  \left\{\left. \left(\begin{matrix}
 Z\\Y
 \end{matrix}\,\hspace{-0.15in}\begin{matrix}
  &\hfill\tikzmark{c3}\\
  &\hfill\tikzmark{d3}
  \end{matrix}\,\,\,\, \begin{matrix}
      0\\I_{(p-r)\times (p-r)}
    \end{matrix}\hspace{-0.13in}\begin{matrix}
  &\hfill\tikzmark{c4}\\
  &\hfill\tikzmark{d4}
  \end{matrix}\hspace{-0.11in}\begin{matrix}
  &\hfill\tikzmark{c5}\\
  &\hfill\tikzmark{d5}
  \end{matrix}\,\,\,\,\, \begin{matrix}
  I_{r\times r}\\0
  \end{matrix}\right)\right\vert_{}\begin{matrix}
  Y\,\,{\rm  is\,\, a\,\,} (p-r)\times (n-p)\,\,{\rm matrix}\,;\\
  \,Z\,\,{\rm  is\,\, a\,\,} r\times (n-p)\,\,{\rm matrix\,\,of\,\,rank\,\,}k-1\\
  \end{matrix}\right\}\,.
  \tikz[remember picture,overlay]   \draw[dashed,dash pattern={on 4pt off 2pt}] ([xshift=0.5\tabcolsep,yshift=7pt]c3.north) -- ([xshift=0.5\tabcolsep,yshift=-2pt]d3.south);\tikz[remember picture,overlay]   \draw[dashed,dash pattern={on 4pt off 2pt}] ([xshift=0.5\tabcolsep,yshift=7pt]c4.north) -- ([xshift=0.5\tabcolsep,yshift=-2pt]d4.south);\tikz[remember picture,overlay]   \draw[dashed,dash pattern={on 4pt off 2pt}] ([xshift=0.5\tabcolsep,yshift=7pt]c5.north) -- ([xshift=0.5\tabcolsep,yshift=-2pt]d5.south);
\end{equation}

Similarly, for $1\leq k\leq r$ we can show that  $d_{k}^+$ is a dense open subset of $\mathcal V_{(p-r+k-1,r-k+1)}^+$\,;  $\overline{R^{-1}_{s,p,n}(d_{k}^+)}=D_{k}^+$; $(g^{\sigma_p}_0\circ\cdots\circ g^{\sigma_p}_{k-1})^{-1}(d_{k}^+)$  is the center of the blow-up $g_{k}^{\sigma_p}$ in a certain smaller open subset of $U_p$; $(g^{\sigma_p}_0\circ\cdots\circ g^{\sigma_p}_{k-1})^{-1}(d_{k}^+)\cong d_{k}^+$ is smooth. We can conclude that $E$ has multiplicity $1$ along $D_{k}^+$ for $1\leq k\leq r$.

\smallskip

{\noindent\bf Case II ($n-s=p<s$).}  Notice that $\overline{\mathcal V_{(p-r+1,r-1)}^-}$ is a divisor of $G(p,n)$. Therefore,  the blow-up $g^{\sigma_0}_{r}:Y^{0}_{r}\rightarrow Y^{0}_{r-1}$ in (\ref{cblow1}) is trivial, and $E$ has multiplicity $0$ along $D_{r}^-$.
\smallskip

{\noindent\bf Case III ($n-s=p=s$).} $E$ has multiplicity $0$ along $D_{r}^+$ and $D_{r}^-$ for $\overline{\mathcal V_{(p-r+1,r-1)}^-}$ and $\overline{\mathcal V_{(p-1,1)}^+}$ are divisors of $G(p,n)$.
\smallskip

We conclude Lemma \ref{excep}.\,\,\,
$\endpf$
\medskip
\begin{remark}\label{ivdw}
By Corollary \ref{itsmoo}, one can show that the center of each intermediate blow-up in  (\ref{cblow1}) or (\ref{csblow}) is a scheme-theoretic union of a smooth submanifold and a divisor (which can be dropped). 
\end{remark}

Denote by  $\mathcal O_{G(p,n)}(1)$  the hyperplane line bundle on $G(p,n)$. We have that

\begin{lemma}\label{picb}
The Picard group Pic$(\mathcal T_{s,p,n})$ of $\mathcal T_{s,p,n}$ is a torsion free abelian group over $\mathbb Z$. When $p<s$ and $p\neq n-s$, Pic$(\mathcal T_{s,p,n})$ has a $\mathbb Z$-basis
\begin{equation}\label{basis1}
    \left\{(R_{s,p,n})^*\left(\mathcal O_{ G(p,n)}(1)\right),D^+_1,\cdots,D^+_{r},D^-_1,\cdots,D^-_{r} \right\};
\end{equation}
when $p=n-s<s$, Pic$(\mathcal T_{s,p,n})$ has a $\mathbb Z$-basis
\begin{equation}
    \left\{(R_{s,p,n})^*\left(\mathcal O_{ G(p,n)}(1)\right),D^+_1,\cdots,D^+_{r},D^-_1,\cdots,D^-_{r-1} \right\}\,;
\end{equation}
when $p=s=n-s$, Pic$(\mathcal T_{s,p,n})$ has a $\mathbb Z$-basis
\begin{equation}
    \left\{(R_{s,p,n})^*\left(\mathcal O_{ G(p,n)}(1)\right),D^+_1,\cdots,D^+_{r-1},D^-_1,\cdots,D^-_{r-1} \right\}\,.
\end{equation}
\end{lemma}
{\noindent\bf Proof of  Lemma \ref{picb}.} 
First assume that $p<s$ and $p\neq n-s$.
Let $Z\subset \mathcal T_{s,p,n}$ be an irreducible divisor.  If $Z$ is contained in the exceptional divisor of $R_{s,p,n}$, then $Z$ is a $\mathbb Z$-linear combination of $D_j^{\pm}$, $1\leq j\leq r$, by Lemma \ref{excep}. Otherwise, the image  $R_{s,p,n}(Z)$ is an irreducible divisor of $G(p,n)$; hence
$R_{s,p,n}(Z)=h\cdot \mathcal O_{G(p,n)}(1)$
for a certain positive integer $h$. Pulling pack $R_{s,p,n}(Z)$, we can conclude by the smoothness of $\mathcal T_{s,p,n}$ and $G(p,n)$ that
\begin{equation}
  R_{s,p,n}^{-1}\left(R_{s,p,n}(Z)\right)=Z+a_r^+\cdot D^+_{r}+\sum_{i=1}^{r-1}\left(a^+_i\cdot D^+_i+a^-_{r-i}D^-_{r-i}\right)+a_r^-\cdot D^-_{r}
\end{equation}
where $a_1^+,a_2^+,\cdots,a_r^+$, $a_1^-,a_2^-,\cdots,a_r^-$ are positive integers.

As a conclusion, the divisors in (\ref{basis1}) generate Pic$(\mathcal T_{s,p,n})$ over $\mathbb Z$. 

Next assume  the following relation in Pic$(\mathcal T_{s,p,n})$.
\begin{equation}\label{chow}
    0=h\cdot (R_{s,p,n})^*\left(\mathcal O_{ G(p,n)}(1)\right)+a_r^+\cdot D^+_{r}+\sum_{i=1}^{r-1}\left(a^+_i\cdot D^+_i+a^-_{r-i}D^-_{r-i}\right)+a_r^-\cdot D^-_{r}\,.
\end{equation}
It is clear that $h=0$. Then, there is a rational function $f$ on  $\mathcal T_{s,p,n}$ such that 
\begin{equation}
    (f)=a_r^+\cdot D^+_{r}+\sum_{i=1}^{r-1}\left(a^+_i\cdot D^+_i+a^-_{r-i}D^-_{r-i}\right)+a_r^-\cdot D^-_{r}\,,
\end{equation}
where $(f)$ is the associated principle divisor of $f$. Since $\mathcal T_{s,p,n}$ is birational to $G(p,n)$,  $f$ induces a rational function $\widetilde f$ on $G(p,n)$. Notice that $\widetilde f$ is regular outside the center of the blow-up $R_{s,p,n}$ which is of codimension at lest two. Therefore, $\widetilde f$ extends to a holomorphic function on $G(p,n)$ for $G(p,n)$ is smooth. Then $\widetilde f$ must be a constant function, and hence $a_i^{\pm}=0$ for $1\leq i\leq r$.

We can complete the proof for other cases in the same way.\,\,\,$\endpf$
\medskip

Denote by $K_{\mathcal T_{s,p,n}}$ be the canonical bundle of $\mathcal T_{s,p,n}$.
We express $K_{\mathcal T_{s,p,n}}$ in terms of the exceptional divisor and  $(R_{s,p,n})^*(\mathcal O_{G(p,n)}(1))$ as follows.
\begin{lemma}\label{kan}
When $r=p$ $(p\leq n-s)$,
    \begin{equation}\label{bcp}
    \begin{split}
     K_{\mathcal T_{s,p,n}}=-n&\cdot (R_{s,p,n})^*(\mathcal O_{G(p,n)}(1))+\sum_{i=1}^{r}\big((p-i+1)(n-s-i+1)-1\big)\cdot D_{i}^-\\
     & +\sum_{i=1}^{r}\big((p-i+1)(s-i+1)-1\big)\cdot D_{i}^+\,\,.
    \end{split}
    \end{equation}
When $r=n-s$ $(n-s\leq p)$,
    \begin{equation}\label{bcr}
    \begin{split}
     K_{\mathcal T_{s,p,n}}=-n&\cdot (R_{s,p,n})^*(\mathcal O_{G(p,n)}(1))+\sum_{i=1}^{r}\big((p-i+1)(n-s-i+1)-1\big)\cdot D_{i}^-\\
     & +\sum_{i=1}^{r}\big((n-p-i+1)(n-s-i+1)-1\big)\cdot D_{i}^+\,\,.
    \end{split}
    \end{equation}
\end{lemma}
{\noindent\bf Proof of Lemma \ref{kan}.} For simplicity, we only prove for the case $p\leq n-s$. 

Recall that the canonical bundle of $G(p,n)$ is \begin{equation}
  K_{G(p,n)}=-n\cdot\mathcal O_{G(p,n)}(1).  
\end{equation} We can take a rational section $\alpha$ of $K_{G(p,n)}$ as follows.  $\alpha$ is defined in on $U_0$ by
\begin{equation}\label{se1}
\small
    \begin{split}
        \alpha&=1\cdot\left(dy_{11}\wedge\cdots\wedge dy_{1(s-p)}\right)\wedge\cdots\wedge\left(dy_{i1}\wedge\cdots\wedge dy_{i(s-p)}\right)\wedge\cdots\wedge\left(dy_{p1}\wedge\cdots\wedge dy_{p(s-p)}\right)\,\\
        &\wedge\left(dw_{1(s+1)}\wedge\cdots\wedge dw_{1n}\right)\wedge\cdots\wedge\left(dw_{i(s+1)}\wedge\cdots\wedge dw_{in}\right)\wedge\cdots\wedge\left(dw_{p(s+1)}\wedge\cdots\wedge dw_{pn}\right)
    \end{split}
\end{equation}
where the holomorphic coordinates ${(y_{ij})}$ and ${(w_{ij})}$ are defined by (\ref{ulx}).
$\alpha$ is defined on $U_r$ by
\begin{equation}\label{se2}
    \begin{split}
        \alpha=\frac{-1}{H^n}\cdot&\left(dz_{11}\wedge\cdots\wedge dz_{1s}\right)\wedge\cdots\wedge\left(dz_{i1}\wedge\cdots\wedge dz_{is}\right)\wedge\cdots\wedge\left(dz_{p1}\wedge\cdots\wedge dz_{ps}\right)\,\\
        &\wedge\left(dx_{1(s+p+1)}\wedge\cdots\wedge dx_{1n}\right)\wedge\cdots\wedge\left(dx_{i(s+p+1)}\wedge\cdots\wedge dx_{in}\right)\wedge\cdots\\
        &\,\,\,\,\,\,\,\,\,\,\,\,\,\,\,\,\wedge\left(dx_{p(s+p+1)}\wedge\cdots\wedge dx_{pn}\right)
    \end{split}
\end{equation}
where the holomorphic coordinates ${(z_{ij})}$ and ${(x_{ij})}$ are defined by (\ref{ulx}), and $H$ is a polynomial defined by a determinant  
\begin{equation}
    H=\left|\begin{matrix}
     z_{1(s-p+1)}&z_{1(s-p+1)}&\cdots&z_{1(s-p+1)}\\
     z_{2(s-p+1)}&z_{2(s-p+1)}&\cdots&z_{2(s-p+1)} \\
    \vdots&\vdots&\ddots&\vdots\\
    z_{p(s-p+1)}&z_{p(s-p+1)}&\cdots&z_{p(s-p+1)} \\
    \end{matrix}\right|\,.
\end{equation}
We note  that $H$ is just a holomorphic section of $O_{G(p,n)}(1)$.

To prove Lemma \ref{kan}, it suffices to compute the pull back of $\alpha$ under the canonical blow-up $R_{s,p,n}$ in $U_0$ and $U_r$, for $\mathcal T_{s,p,n}\backslash\left( R_{s,p,n}^{-1}(U_0)\bigcup R_{s,p,n}^{-1}(U_r)\right)$ is of codimension at least two in $\mathcal T_{s,p,n}$. 

Recall that $R_{s,p,n}$ is given by the maps $\Gamma_{0}^{\tau}$ and $\Gamma_{r}^{\tau}$ defined in (\ref{ws}), in terms of the local holomorphic coordinates associated with the Van der Waerden representation of $R_{s,p,n}^{-1}(U_0)$ and $R_{s,p,n}^{-1}(U_r)$.  Substitute (\ref{ws}) into (\ref{se1}) and (\ref{se2}), we can conclude  Lemma \ref{kan}.\,\,\,$\endpf$
\begin{remark}
We can prove Lemma \ref{kan} alternatively as follows. Recall that for a blow-up $Y$ of a smooth complex manifold $X$ along a smooth submanifold $Z\subset X$, we have that
\begin{equation}\label{ctra}
K_Y=f^*(K_X)+(c-1)E
\end{equation}
where $E$ is the exceptional divisor and $c$ is the codimension of $Z$ in $X$.  Consider the blow-ups (\ref{cblow1}) and (\ref{csblow}) for $U_0$ and $U_r$ respectively. By Remark \ref{ivdw},  we can conclude Lemma (\ref{kan}) by applying (\ref{ctra}) repeatedly.
\end{remark}

We would like to give an interesting characterization of $\overline{\mathcal V_{(p-k,k)}^-}$, $1\leq k\leq p$.  This is inspired by the beautiful work of Hwang (see \cite{Hw}).  For simplicity, let $s=n-p$ and $\mathfrak a$ a pint of $U_p$ such that $z_{ij}(\mathfrak a)=0$ for $1\leq i\leq p$ and $1\leq j\leq s$, where $z_{ij}$ are the holomorphic coordinates defined by (\ref{u0}).
\begin{lemma}\label{key}
Let $1\leq k\leq p-1 $. Denote by $\mathcal C_{\mathfrak a}^k$ the set of smooth rational curves of $G(p,n)$ which pass through $\mathfrak a$ and are of degree $k$ with respect to $\mathcal O_{G(p,n)}(1)$. Then, 
\begin{equation}
    \overline{\bigcup_{\gamma\in\mathcal C_\mathfrak a^k}\gamma}=\overline{\mathcal V_{(p-k,k)}^-}.
\end{equation}
\end{lemma}

{\noindent\bf Proof of Lemma \ref{key}.} 
By taking special smooth rational curves of degree $k$ (for instance, certian integral curves associated with $\mathbb C^*$-action $\psi_{s,p,n}$), it is easy to verify that \begin{equation}
    \overline{\bigcup_{\gamma\in\mathcal C_\mathfrak a^k}\gamma}\supset\overline{\mathcal V_{(k,p-k)}^+}.
\end{equation}

We prove the other direction in the following. Take a smooth rational curve $\gamma$ passing through $\mathfrak a$ with degree $k$.  We can parametrize it by $\gamma(t):=\left(z_{ij}(t)\right)_{1\leq i\leq p,\,1\leq j\leq s}=t\cdot W(t)$ where
$W (t):=\left(\frac{ f_{ij}(t)}{g_{ij}(t)}\right)_{1\leq i\leq p,\,1\leq j\leq s}$
such that $f_{ij}$ and $g_{ij}$ are polynomials in $t$ and that $t\cdot f_{ij}$ and $g_{ij}$ are coprime. 

For each $I=(i_1,\cdots,i_{p-l},i_{p-l+1},\cdots, i_p)\in\mathbb I^{p-l}_{s,p,n}$, $0\leq l\leq p$, we can conclude that 
\begin{equation}\label{curv}
P_{I}(\gamma(t))=t^l\cdot
\frac{F_{I}(t)}{G_{I^*}(t)}\,\,,    
\end{equation}
where ${F_{I}(t)}$ and ${G_{I}(t)}$ are  polynomials in $t$ such that $t^l\cdot{F_{I}(t)}$ and ${G_{I}(t)}$ are coprime.

To prove Lemma \ref{key}, it suffices to show that $P_{I}(\gamma(t))\equiv0$ for each  $I\in\mathbb I^{p-l}_{s,p,n}$ and  $k+1\leq l\leq p$. Suppose that $P_{I}(\gamma(t))\not\equiv 0$ for a certain  and $I^*\in\mathbb I^{p-l^*}_{s,p,n}$ where $k+1\leq l^*\leq p$. Choose a generic hyperplane $H\in\mathbb{CP}^{N_{p,n}}$ defined by
\begin{equation}\label{hypp}
\sum_{I\in\mathbb I_{p,n}}\beta_I\cdot z_I=0\,.
\end{equation} 
Substituting (\ref{curv}) into  (\ref{hypp}), we have that
\begin{equation}\label{solu}
  0=\sum_{l=0}^pt^l\left(\sum_{I\in\mathbb I_{s,p,n}^l}\beta_I\cdot  \frac{F_{I}(t)}{G_{I}(t)}\right)\,.  
\end{equation}
Under the above assumption we can arrange $\beta_I$, $I\in\mathbb I_{p,n}$, such that  (\ref{solu}) is reduced to a polynomial equation of degree at least $k+1$.  It contradicts the fact that $\gamma$ has degree $k$ with respect to $\mathcal O_{G(p,n)}(1)$. \,\,\,$\endpf$

\section{Basic properties of \texorpdfstring{$\mathcal M_{s,p,n}$ }{jj}} \label{basicm}

In this section, we define a special subvariety $\mathcal M_{s,p,n}$ of $\mathcal T_{s,p,n}$, and then study its geometric properties such as the existence of the moduli maps, isomorphisms, fibration structures, invariant divisors, etc. 

We assume that   $2p\leq n\leq  2s$  and let $r$ be the rank of $\mathcal T_{s,p,n}$.

\subsection{The definition of \texorpdfstring{$\mathcal M_{s,p,n}$ }{jj} and the existence of the flat map \texorpdfstring{$\mathcal P_{s,p,n}:\mathcal T_{s,p,n}\rightarrow\mathcal M_{s,p,n}$ }{jj}}
\label{basicp}
\begin{definition}\label{mspn}
Define  $\mathcal M_{s,p,n}$ to be the source of $\mathcal T_{s,p,n}$ associated with the $\mathbb C^*$-action $\Psi_{s,p,n}$.
\end{definition}
\begin{remark}
$\mathcal M_{s,p,n}=\mathcal D_{(p,0)}=D^-_1$. In the following, we will use $\mathcal M_{s,p,n}$ and $D^-_1$ interchangeably to indicate different aspects.
\end{remark}
{\bf\noindent Proof of Proposition \ref{msmooth}.} This follows from Lemma \ref{11}.
\,\,\,\,$\endpf$
\medskip

Next, we will construct a natural flat map  $\mathcal P_{s,p,n}:\mathcal T_{s,p,n}\rightarrow\mathcal M_{s,p,n}$. Define a projection map
$\mathfrak {P}:\mathcal T_{s,p,n}\rightarrow \mathbb {CP}^{N^0_{s,p,n}}\times\cdots\times\mathbb {CP}^{N^p_{s,p,n}}$ to be the restriction to $\mathcal T_{s,p,n}$ of the following natural projection of the ambient space.
\begin{equation}
\begin{split}
\widetilde {\mathfrak P}&:\,\,\mathbb {CP}^{N_{p,n}}\times\mathbb {CP}^{N^0_{s,p,n}}\times\cdots\times\mathbb {CP}^{N^p_{s,p,n}}\longrightarrow\mathbb {CP}^{N^0_{s,p,n}}\times\cdots\times\mathbb {CP}^{N^p_{s,p,n}}
\end{split}.
\end{equation}
\begin{lemma}\label{emb}
$\mathfrak P$ induces an embedding of $\mathcal M_{s,p,n}$ into $\mathbb {CP}^{N^0_{s,p,n}}\times\cdots\times\mathbb {CP}^{N^p_{s,p,n}}$.
\end{lemma}
{\bf\noindent Proof of Lemma \ref{emb}.} First notice that the rational map $f^0_s$ is well-defined on the source $\mathcal V_{(p,0)}$ and $f^0_s\big|_{\mathcal V_{(p,0)}}:\mathcal V_{(p,0)}\rightarrow \mathbb {CP}^{N^0_{s,p,n}}$ is an embedding. This is due to the fact that $\mathcal V_{(p,0)}$ is isomorphic to a sub-grassmannian $G(p,s)$ of $G(p,n)$, and, moreover,  the sutset of the Pl\"ucker coordinate functions  $\{P_I\}_{I\in\mathbb I^0_{s,p,n}}$ of $G(p,n)$ is the (full) set of the Pl\"ucker coordinate functions $\{P_I\}_{I\in\mathbb I_{p,s}}$ of  $G(p,s)$.

Then the projection of 
$\mathcal M_{s,p,n}$ to $\mathbb {CP}^{N^0_{s,p,n}}\times\cdots\times\mathbb {CP}^{N^p_{s,p,n}}$ is injective. Since $\mathcal M_{s,p,n}$ is smooth, Lemma \ref{emb} follows.
\,\,\,\,$\endpf$.
\medskip

\begin{lemma}\label{iemb}
The image of $\mathcal T_{s,p,n}$ under $\mathfrak P$ is the image of $\mathcal M_{s,p,n}$ under $\mathfrak P$.
\end{lemma}
{\bf\noindent Proof of Lemma \ref{iemb}.}  It suffices to show that there is an dense open subset $U$ of $\mathcal T_{s,p,n}$ such that $\mathfrak P(U)\subset \mathfrak P\left(\mathcal M_{s,p,n}\right)$.

Consider the open subset $U_0$ of $G(p,n)$ and the Van der Waerden representation of $R_{s,p,n}^{-1}(U_0)$. Let $\left(A^{\tau},(J_0^{\tau})^{-1}\right)$ be a holomorphic coordinate chart in the Van der Waerden representation where $\tau=\left(
	\begin{matrix}
	1&2&\cdots&r\\
	s+1&s+2&\cdots&s+r\\
	\end{matrix}\right)\in\mathbb J_0$.
Let $\mathfrak a$ be an arbitrary point of $A^{\tau}$ with  the holomorphic coordinates $\left(\widetilde Y(\mathfrak a),\overrightarrow B^1(\mathfrak a),\cdots,\overrightarrow B^r(\mathfrak a)\right)$. 
By writing down the P\"ucker coordinate functions explicitly, we can show that the image of $\mathfrak a$ under the projection $\mathfrak P$ is the same as $\mathfrak P\left(\check {\mathfrak a}\right)$ where $\check {\mathfrak a}$ is a  point of $\mathcal M_{s,p,n}$ such that
\begin{equation}
\begin{split}
&a_{1(s+1)}(\check {\mathfrak a})=0\\
&\xi^{(1)}_{1j}(\check {\mathfrak a})=\xi^{(1)}_{1j}({\mathfrak a})\,,\,\,\,s+2\leq j\leq n,\,\,\,\,\,\,\xi^{(1)}_{i(s+1)}(\check {\mathfrak a})=\xi^{(1)}_{i(s+1)}({\mathfrak a})\,,\,\,\,2\leq i\leq p,\\
&\widetilde Y(\check {\mathfrak a})=\widetilde Y({\mathfrak a}),\,\,\,\,\,\,\overrightarrow B^l(\check {\mathfrak a})=\overrightarrow B^{l}({\mathfrak a})\,,\,\,\,2\leq l\leq r.
\end{split}
\end{equation}

We thus conclude Lemma \ref{iemb}.\,\,\,\,$\endpf$.
\medskip

Denote by 
\begin{equation}
    \check {\mathfrak P}:\mathcal M_{s,p,n}\rightarrow \mathfrak P\left(\mathcal M_{s,p,n}\right)
\end{equation}
the induced isomorphism of $\mathcal M_{s,p,n}$ to its image $\mathfrak P\left(\mathcal M_{s,p,n}\right)$ under $\mathfrak P$. Denote by $\left(\check {\mathfrak P}\right)^{-1}$ the inverse map.
\begin{definition}
Define a map $\mathcal P_{s,p,n}:\mathcal T_{s,p,n}\rightarrow\mathcal M_{s,p,n}$ by
\begin{equation}
    \mathcal P_{s,p,n}:=\left(\check {\mathfrak P}\right)^{-1}\circ \mathfrak P.
\end{equation}
We call $\mathcal P_{s,p,n}$ the {\it moduli map}. 
\end{definition}
It follows from definition that
\begin{corollary}\label{retra}
$\mathcal P_{s,p,n}:\mathcal T_{s,p,n}\rightarrow\mathcal M_{s,p,n}$ is a retraction, that is, the restriction of $\mathcal P_{s,p,n}$ to $\mathcal M_{s,p,n}$ is the identity map.
\end{corollary}
\begin{remark}\label{gfi}
There is a geometric way to visualize the map $\mathcal P_{s,p,n}$ as follows. For each point $\mathfrak a\in\mathcal T_{s,p,n}$ we can flow it back  to the source by $\Lambda^+$ (defined by (\ref{pl})) through  a sequence of adjacent integral curves associated with the $\mathbb C^*$-action $\Psi_{s,p,n}$. It is well-defined since each $\mathfrak a\in\mathcal T_{s,p,n}$ has a unique negative direction unless it is in the source. 
\end{remark}
\begin{remark}
We give an explicit formula of the map $\mathcal P_{s,p,n}$ in local coordinate charts as follows. For each index $\tau=\left(\begin{matrix} i_1&i_2&\cdots&i_r\\ j_1&j_2&\cdots&j_r\\
\end{matrix}\right)\in\mathbb J_l$, $0\leq l\leq r$, consider the holomorphic chart $\left(A^{\tau},(J^{\tau}_l)^{-1}\right)$ .  Define a holomorphic map $p_l:\mathbb C^{p(n-p)}\rightarrow \mathcal V_{(p,0)}$ by
\begin{equation}
p_l\left(\widetilde X,\widetilde Y,\overrightarrow B^1,\cdots,\overrightarrow B^r\right)=\left(
\begin{matrix}
\sum_{k=r-l+1}^{r}\theta_k^T\cdot\Omega_k &0_{l\times(p-l)}&0_{l\times l}&0_{l\times (n-s-l)}\\\widetilde Y&I_{(p-l)\times(p-l)}&0_{(p-l)\times l}&0_{(p-l)\times (n-s-l)}\\
\end{matrix}\right)\,,
\end{equation}
where for $r-l+1\leq k\leq r$,
\begin{equation}
\begin{split}
&\theta_k=\big(0,\cdots,0,\underset{\substack{\uparrow\\(k-r+l)^{th}}}{1,}\,0,\cdots,0\big)\,\\
&\Omega_{k}=
\left(\xi^{(k)}_{i_k1},\xi^{(k)}_{i_k2},\cdots,{\xi^{(k)}_{i_k(j_1-1)}},0,{\xi^{(k)}_{i_k(j_1+1)}},\cdots,{\xi^{(k)}_{i_k(j_2-1)}},0,{\xi^{(k)}_{i_k(j_2+1)}},\cdots,\right.\\
&\,\,\,\,\,\,\,\,\,\,\,\,\,\,\,\,\,\,\,\,\,\,\,\,\,\left.{\xi^{(k)}_{i_k(j_{k-1}-1)}},0,{\xi^{(k)}_{i_k(j_{k-1}+1)}},\cdots,\xi^{(k)}_{i_k(j_k-1)},1,\xi^{(k)}_{i_k(j_k+1)},\cdots,\xi^{(k)}_{i_ks}\right).\\
\end{split}
\end{equation}
Then $\mathcal P_{s,p,n}\big|_{A^{\tau}}=(p_l,f_s^0\circ\Gamma_{l}^{\tau},\cdots,f_s^p\circ\Gamma_{l}^{\tau})$.
\end{remark}

\begin{lemma}\label{flat}
    $\mathcal P_{s,p,n}:\mathcal T_{s,p,n}\rightarrow\mathcal M_{s,p,n}$ is flat. 
\end{lemma}
{\bf \noindent Proof of Lemma \ref{flat}.}
Since $\mathcal T_{s,p,n}$ and $\mathcal M_{s,p,n}$ are smooth, the following are equivalent by the holomoprhic version of the miracle flatness theorem (see Corollary in Page 158 of \cite{Fi}).
\begin{enumerate}[label=(\alph*).]
    \item $\mathcal P_{s,p,n}$  is flat.
    \item $\mathcal P_{s,p,n}$ is open.
    \item Each fiber of $\mathcal P_{s,p,n}:\mathcal T_{s,p,n}\rightarrow\mathcal M_{s,p,n}$ is of pure dimension $\dim \mathcal T_{s,p,n}-\dim \mathcal M_{s,p,n}=1$. 
\end{enumerate}

Since $\mathcal M_{s,p,n}$ is smooth of dimension $p(n-p)-1$, then each fiber can be defined by $p(n-p)-1$ equations locally. Therefore, for each point $\mathfrak a\in\mathcal T_{s,p,n}$ the fiber $\mathcal P_{s,p,n}^{-1}\left(\mathcal P_{s,p,n}\left(\mathfrak a\right)\right)$ is of dimension at least $1$ at $\mathfrak a$.
Next, we will show that each fiber of $\mathcal P_{s,p,n}$ is of dimension at most $1$.

Let $\check {\mathfrak a}$ be an arbitrary point of $\mathcal M_{s,p,n}$. It suffices to compute the preimage  $\mathcal P_{s,p,n}^{-1}\left(\check {\mathfrak a}\right)$ in local coordinate charts in the Van der Waerden representation of $R_{s,p,n}^{-1}(U_l)$, $0\leq l\leq r$. For each $\tau=\left(\begin{matrix} i_1&i_2&\cdots&i_r\\ j_1&j_2&\cdots&j_r\\
\end{matrix}\right)\in\mathbb J_l$,
let $\left(A^{\tau},(J_l^{\tau})^{-1}\right)$  be the associated local coordinate chart.  By computing explicitly the Pl\"ucker coordinate functions in terms of the holomorphic coordinates $\left(\widetilde X,\widetilde Y,\overrightarrow B^1,\cdots,\overrightarrow B^r\right)$, we can show that of each point $\mathfrak b\in\mathcal P_{s,p,n}^{-1}\left(\check {\mathfrak a}\right)\cap A^{\tau}$ the following holomorphic coordinates  take the same values respectively.

\begin{equation}
\begin{split}
&\xi^{(1)}_{i_1j}\left(\mathfrak b\right)=\xi^{(1)}_{i_1j}\left(\check {\mathfrak a}\right)\,,\,\,\,j\in\{s+l+1,s+l+2,\cdots,n\}\backslash\{j_1\},\\
&\xi^{(1)}_{ij_1}\left(\mathfrak b\right)=\xi^{(1)}_{ij_1}\left(\check {\mathfrak a}\right)\,,\,\,\,i\in\{l+1,l+2,\cdots,p\}\backslash\{i_1\},\\
&\xi^{(r-l+1)}_{i_{r-l+1}j}\left(\mathfrak b\right)=\xi^{(r-l+1)}_{i_{r-l+1}j}\left(\check {\mathfrak a}\right)\,,\,\,\,j\in\{1,2,\cdots,s-p+l\}\backslash\{j_{r-l+1}\},\\
&\xi^{(r-l+1)}_{ij_{r-l+1}}\left(\mathfrak b\right)=\xi^{(r-l+1)}_{ij_{r-l+1}}\left(\check {\mathfrak a}\right)\,,\,\,\,i\in\{1,2,\cdots,l\}\backslash\{i_{r-l+1}\},\\
&\widetilde X\left(\mathfrak b\right)=\widetilde X\left(\check {\mathfrak a}\right),\,\,\,\,\\
&\widetilde Y\left(\mathfrak b\right)=\widetilde Y\left(\check {\mathfrak a}\right),\,\\
&\overrightarrow B^l\left(\mathfrak b\right)=\overrightarrow B^l\left(\check {\mathfrak a}\right)\,,\,\,\,2\leq l\leq r.
\end{split}
\end{equation}
When $l=0$ (resp. $r$), the remaining variable is $a_{i_1j_1}$ (resp. $b_{i_rj_r}$); then  $\mathcal P_{s,p,n}^{-1}\left(\check {\mathfrak a}\right)\cap A^{\tau}$ is of dimension at most $1$. When $1\leq l\leq r-1$, the remaining variables are $a_{i_{r-l+1}j_{r-l+1}}$ and $b_{i_1j_1}$. Define an index $I^{**}\in\mathbb J_l$ by
\begin{equation}
    I^{**}:=(j_1,s+l,s+l-1,\cdots,\widehat {s+i_{r-l+1}},\cdots,\widehat {s-p+i_{1}},\cdots,s-p+l+2, s-p+l+1, j_{r-l+1})\,.
\end{equation}
Computing the Pl\"ucker coordinate function $P_{I^{**}}$, we can show that each point $\mathfrak b\in\mathcal P_{s,p,n}^{-1}\left(\check {\mathfrak a}\right)\cap A^{\tau}$ satisfies the following equation.
\begin{equation}
    a_{i_{r-l+1}j_{r-l+1}}\left(\mathfrak b\right)\cdot b_{i_1j_1}\left(\mathfrak b\right)=C\left(\check {\mathfrak a}\right)
\end{equation}
where $C(\check {\mathfrak a})$ is a certain function depending on $\check {\mathfrak a}$  only. We thus conclude that $\mathcal P_{s,p,n}^{-1}\left(\check {\mathfrak a}\right)\cap A^{\tau}$ is of dimension at most $1$.

We complete the proof of Lemma \ref{flat}.\,\,\,$\endpf$

\begin{remark}
$\mathcal P_{s,p,n}$ is a degeneration of smooth rational curves such that each fiber is a union of rational curves by the geometric interpretation in Remark \ref{gfi}. Fujiki  proved the existence of a flat family parametrizing the integral curves (see  \cite{Fu} or \cite{BS}). We construct Fujiki families for Grassmannian manifolds in an explicitly way such that it desingularizes the foliation.
\end{remark}
\begin{remark}
Figures 1 and 2 illustrate the degeneration and the foliation on $G(2,4)$ and $G(3,6)$ respectively. The red (resp. blue) arrows indicate the positive (resp. negative) directions. In Figure 1 a generic rational curve of degree $2$ breaks into two rational curves of degree $1$. In Figure 2 a generic rational curve of degree $3$ breaks into a rational curves of degree $1$ and a rational curve of degree  $2$ in two ways, depending on whether the degree $1$ curve is connected to the source or the sink. Finally, it breaks into three rational curves of degree $1$. 
\end{remark}
\newpage

\begin{figure}[htbp]
	\vspace{.3in}
	\centering
\hspace{-1.5in}	\includegraphics[angle=-90,origin=c,height=11cm,width=16cm]{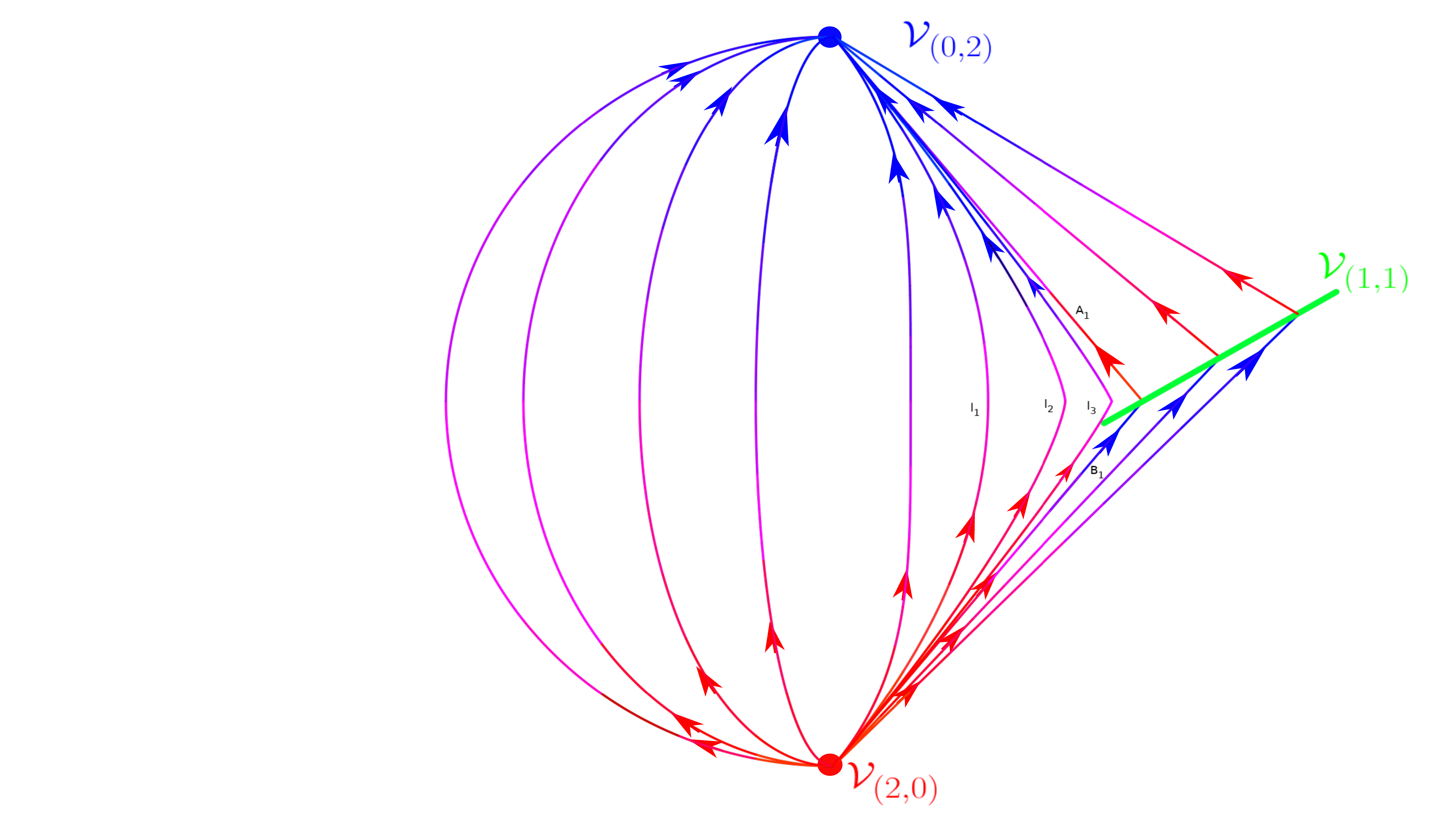}
		\vspace{-1in}
	\caption{The foliation on $G(2,4)$}
	\label{G2}
	\vspace{.5in}
\end{figure}

\begin{figure}[htbp]
	\centering
\hspace{-.7in}	\includegraphics[angle=-90,origin=c,height=11cm,width=16cm]{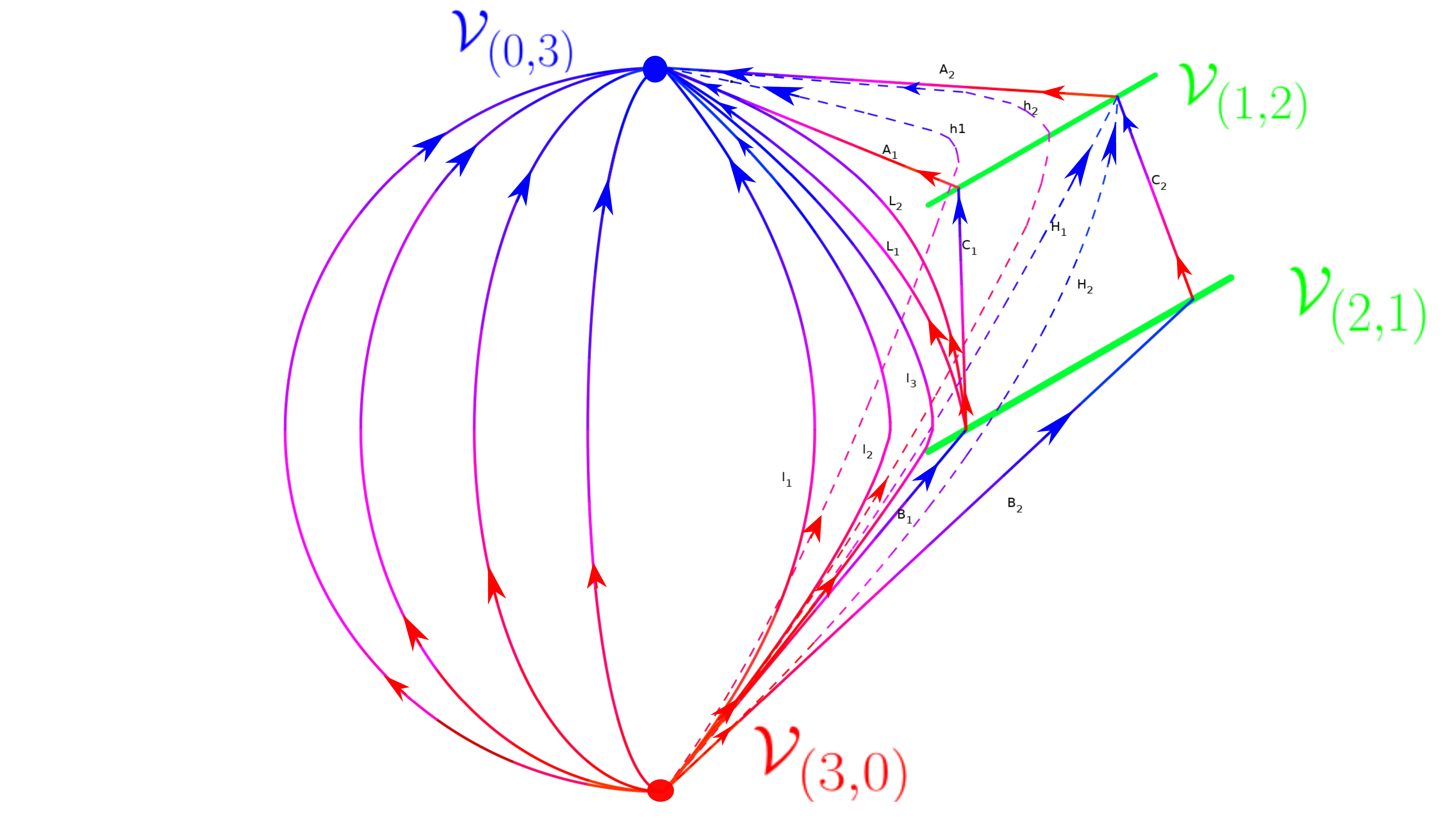}
	\vspace{-1in}
	\caption{The foliation on $G(3,6)$}
	\label{G3}
		\vspace{-2in}
\end{figure}

\newpage

\subsection{Isomporphisms and fibrations of \texorpdfstring{$\mathcal M_{s,p,n}$ }{jj}} \label{basicif}

In this subsection, we introduce an isomorphism $\mathcal L_{s,p,n}$ from the sink to the source.  As a consequence, we can define two types of isomorphisms of  $\mathcal M_{s,p,n}$ induced from the USD and DUAL isomorphisms of  $\mathcal T_{s,p,n}$. We note that the restriction of USD or DUAL  itself is not an isomorphism of $\mathcal M_{s,p,n}$. Another application of $\mathcal L_{s,p,n}$ is to give different fibration structures on generic $\mathcal M_{s,p,n}$.

\begin{lemma}\label{sti}
$\mathcal P_{s,p,n}$ induces an isomorphsim from $D_1^+$ to $D_1^-$.
\end{lemma}

{\noindent\bf Proof of  Lemma \ref{sti}.} 
Notice that $\mathcal P_{s,p,n}$ defines a holomorphic map $\mathcal P_{s,p,n}\big|_{D_1^+}$ from $D_1^+$ to $D_1^-$.  
Similarly to Lemma \ref{emb}, we can show that
the rational map $f^r_s$ is well-defined on the sink $\mathcal V_{(p-r,r)}$ and $f^r_s\big|_{\mathcal V_{(p-r,r)}}:\mathcal V_{(p-r,r)}\rightarrow G(p,n)$ is an embedding.
Then the projection of 
$D^+$ to $\mathbb {CP}^{N^0_{s,p,n}}\times\cdots\times\mathbb {CP}^{N^p_{s,p,n}}$ is injective.

It is then easy to show that $\mathcal P_{s,p,n}\big|_{D_1^+}$ is an isomorphism.  We conclude Lemma \ref{sti}.\,\,\,$\endpf$
\begin{remark}
Notice that $D_1^-$ is the source and $D_1^+$ is the sink, and that each point of $\mathcal D_{p-k,k}$, $1\leq k\leq r-1$, has one positive direction and one negative direction. Then geometrically it is clear that  $\mathcal P_{s,p,n}\big|_{D_1^+}$ is a bijection. 
\end{remark}
\begin{definition}\label{fi}
 Denote  $\mathcal P_{s,p,n}\big|_{D_1^+}$  by $\mathcal L_{s,p,n}:D_1^+\rightarrow D_1^-$ and call it the {\it flow isomorphism}.
\end{definition}

We can define isomorphisms as follows. (Recall that $D^-_1=\mathcal M_{s,p,n}$.)
\begin{definition}\label{musd}
Define an isomorphism ${\rm Usd}:\mathcal M_{s,p,n}\rightarrow \mathcal M_{n-s,p,n}$ by
\begin{equation}
{\rm Usd}:=\mathcal L_{n-s,p,n}\circ\left({\rm USD}\big|_{\mathcal M_{s,p,n}}\right)\,,
\end{equation}
where $\left({\rm USD}\big|_{\mathcal M_{s,p,n}}\right)$ is the restriction of the map $\rm USD$ defined by (\ref{USD}).
\end{definition}
\begin{remark}
    ${\rm Usd}$ is an automorphism of $\mathcal M_{s,p,2s}$. 
\end{remark}

By (\ref{dtrans}) we can show that the dual map from $G(p,n)$ to $G(n-p,n)$ maps the source $\mathcal V_{(p,0)}$ of $G(p,n)$ to the sink $\mathcal V_{(p-r,r)}$ of $G(n-p,n)$.
Hence the {\rm DUAL} isomorphism defined in (\ref{DUAL}) maps the source $D_1^-$ of $\mathcal T_{s,p,n}$ to the sink $D_1^+$ of $\mathcal T_{s,n-p,n}$. 

\begin{definition}\label{mdual}
Define an isomorphism ${\rm Dual}:\mathcal M_{s,p,n}\rightarrow \mathcal M_{s,n-p,n}$ by
\begin{equation}
{\rm Dual}:=\mathcal L_{s,n-p,n}\circ\left({\rm DUAL}\big|_{\mathcal M_{s,p,n}}\right)\,.
\end{equation}
\end{definition}

Next we turn to the fibration structures on $\mathcal M_{s,p,n}$. Recall the following realization of 
the variety of complete collineations (see (5) of Example 2.10 in \cite{Pez} when $p=q$).
\begin{example}\label{ccolli}
Let $M_{p\times q}$ be the set of all $p\times q$ complex matrices and $P(M_{p\times q})$ the projectivization of $M_{p\times q}$; that is,
\begin{equation}\label{ccolli1}
P(M_{p\times q}):=\bigslant{\left(M_{p\times q}\middle/\{0\}\right)}{(a_{ij})\sim (\lambda a_{ij})}\,,\,\,\,\lambda\in\mathbb C^*.
\end{equation}
It is clear that $P(M_{p\times q})$ is a $GL(p,\mathbb C)\times GL(q,\mathbb C)$-variety under the left multiplication by the group $GL(p,\mathbb C)$ and the right multiplication by the group $GL(q,\mathbb C)$. Define $r:=\min\left\{p,q\right\}$.
Then $P(M_{p\times q})$ has $r$ $GL(p,\mathbb C)\times GL(q,\mathbb C)$-orbits, whose closures $Z_{r}\supset Z_{r-1}\supset \cdots\supset Z_1$ are given by the condition
that $Z_i$ is the set of points corresponding to matrices of rank at most $i$.  Blowing up $P(M_{p\times q})$ successively along the (strict transform of)
$Z_1,\cdots,Z_{r-1}$, we obtain a smooth variety $\widetilde P(M_{p\times q})$.  By the Van der Waerden representation,  we can show that $\widetilde P(M_{p\times q})\cong\mathcal M_{p,p,p+q}$. In particular, when $p=q$, $\mathcal M_{p,p,2p}$ is the group compactification of $PGL(p,\mathbb C)$.
\end{example}

{\bf\noindent Proof of Proposition \ref{fiber}.}  It is clear that the image of $\mathcal M_{s,p,n}$ under $R_{s,p,n}$ is isomorphic to $G(p,s)$. Consider the Van der Waerden representation of $R_{s,p,n}^{-1}(U_0)$. By fixing   variables $\widetilde Y$, the holomorphic coordinate charts $\left\{\left(A^{\tau}, (J^{\tau}_0)^{-1}\right)\right\}_{\tau\in\mathbb J_0}$ induce a holomorphic atlas of each fiber of $R_{s,p,n}\big|_{\mathcal M_{s,p,n}}$.  Comparing the holomorphic atlas with that of the variety of complete collineations in Example \ref{ccolli},  we can conclude each fiber  of $R_{s,p,n}\big|_{\mathcal M_{s,p,n}}$ is isomorphic to the variety of complete collineations.

The other fibration structure can be derived by pulling back the fibration structures on $\mathcal M_{n-s,p,n}$ by Usd.

We complete the proof of Propostion \ref{fiber}.\,\,\,$\endpf$ 
\medskip

The fibration structures on $\mathcal M_{s,p,n}$ can be realized as blow-ups of projective bundles over sub-grassmannians as follows.  Denote by $N_{\mathcal V_{(p,0)}/G(p,n)}$  the normal bundle of the source $\mathcal V_{(p,0)}$ in $G(p,n)$, and $\mathbb P(N_{\mathcal V_{(p,0)}/G(p,n)})$ its projectivization. Recall   (\ref{cblow1}).
By Theorem \ref{bb} and  Lemma \ref{kl}, $(g^{\sigma_0}_0)^{-1}(\mathcal V_{(p,0)})\cong\mathbb P(N_{\mathcal V_{(p,0)}/G(p,n)})$; by a slight abuse of notation, we denote $(g^{\sigma_0}_0)^{-1}(\mathcal V_{(p,0)})$ by $\mathbb P(N_{\mathcal V_{(p,0)}/G(p,n)})$. Then the following diagram commutes, where $\tau^1_{s,p,n}:=\left(g^{\sigma_0}_1\circ\cdots\circ g^{\sigma}_{p}\right)\big|_{\mathcal M_{s,p,n}}$,  $\mathcal F^1_{s,p,n}:=R_{s,p,n}\big|_{\mathcal M_{s,p,n}}$,   and $\kappa^1_{s,p,n}:=\left(g^{\sigma_0}_{0}\right)\big|_{\mathbb P(N_{\mathcal V_{(p,0)}/G(p,n)})}$.
\vspace{-.05in}
\begin{equation}\label{checkr}
\begin{tikzcd}
&\mathcal M_{s,p,n}\arrow{d}{\hspace{-.6in} \mathcal F^1_{s,p,n}}\arrow{r}{\tau^1_{s,p,n}}&[4em]\mathbb P(N_{\mathcal V_{(p,0)}/G(p,n)}) \arrow{ld}{\kappa^1_{s,p,n}} \\
&\hspace{-.7in}G(p,s)\cong\mathcal V_{(p,0)}&\\
\end{tikzcd}\vspace{-20pt}\,.
\end{equation}

Similarly, we can realize another fibration structures by (\ref{csblow}).   Denote by $N_{\mathcal V_{(p-r,r)}/G(p,n)}$  the normal bundle of the sink $\mathcal V_{(p-r,r)}$ in $G(p,n)$, and $\mathbb P(N_{\mathcal V_{(p-r,r)}/G(p,n)})$ its projectivization. Define maps
$\tau^2_{s,p,n}:=\left(\left(g^{\sigma_p}_1\circ\cdots\circ g^{\sigma_p}_{p}\right)\big|_{\mathcal M_{s,p,n}}\right)\circ \mathcal L_{s,p,n}^{-1}$, 
$\mathcal F^2_{s,p,n}:=\left(R_{s,p,n}\big|_{\mathcal M_{s,p,n}}\right)\circ\mathcal L^{-1}_{s,p,n}$, and $\kappa^2_{s,p,n}:=\left(g^{\sigma_p}_{0}\right)\big|_{\mathbb P(N_{\mathcal V_{(p-r,r)}/G(p,n)})}$.  Then we have the following commutative diagram.
\vspace{-.05in}
\begin{equation}\label{checkr1}
\begin{tikzcd}
&\mathcal M_{s,p,n}\arrow{d}{\hspace{-.6in} \mathcal F^2_{s,p,n}}\arrow{r}{\tau^2_{s,p,n}}&[4em]\mathbb P(N_{\mathcal V_{(p-r,r)}/G(p,n)}) \arrow{ld}{\kappa^2_{s,p,n}} \\
&\hspace{-.1in}\mathcal V_{(p-r,r)}&\\
\end{tikzcd}\vspace{-25pt}\,.
\end{equation}
It is clear that fibers of $\mathcal F^i_{s,p,n}$, $i=1,2$, are isomorphic to the variety of complete collineations. 
\smallskip

To determine the automorphism groups of $\mathcal M_{s,p,n}$ in Section \ref{sym}, we further assign certain convenient local coordinate charts for the projective bundles. The readers are referred to Appendix \ref{section:projbc}.

\begin{remark}
    The group $GL(s,\mathbb C)\times GL(n-s,\mathbb C)$ acts naturally on  $\mathbb P(N_{\mathcal V_{(p,0)}/G(p,n)})$ and  $\mathbb P(N_{\mathcal V_{(p-r,r)}/G(p,n)})$.  We may write down the transition functions of local trivializations of the projective bundles explicitly; it is given by the conjugate  action of the transition functions of the local trivializations of the base sub-grassmannians.
\end{remark}

\begin{remark}
    The two fibration structures on $\mathcal M_{s,p,n}$ leads to the following interesting process:  blow up a projective bundle over $G(p,s)$; then blow down the exceptional divisors to derive a projective bundle over  $G(p,n-s)$ when $p>n-s$ (or $G(s+p-n,s)$ when $p<n-s$).
\end{remark}

\subsection{The invariant divisors of \texorpdfstring{$\mathcal M_{s,p,n}$ }{jj}}
\label{basicid}
For $2\leq i\leq r$, define \begin{equation}\label{checkd}
   \check D_i:= D_{1}^-\cap D_{i}^-=\mathcal M_{s,p,n}\cap D_{i}^-\,.
\end{equation} 
It is clear  that the above divisors are invariant under the natural  $GL(s,\mathbb C)\times GL(n-s,\mathbb C)$-action of $\mathcal M_{s,p,n}$ induced from that of $\mathcal T_{s,p,n}$.
For convenience, we denote by $\check D_1$ the restriction of the line bundle $D_1^-$ to itself.

\begin{remark}
It is easy to show that $\check D_i=\mathcal P_{s,p,n}\left(\mathcal D_{p-i+1,i-1}\right)$ for $2\leq i\leq r$  (see Lemma \ref{11} for the definition of $\mathcal D_{p-i+1,i-1}$). In fact, $\check D_i\cong\mathcal D_{p-i+1,i-1}$.
\end{remark}

Similarly to Lemmas \ref{snc} and \ref{excep}, we have that
\begin{lemma}\label{msnc}
 Then
$\check D_{2}$, $\check D_{3}$, $\cdots$, $\check D_{r}$ are distinct smooth divisors of $\mathcal M_{s,p,n}$, and the following divisor is simple normal crossing.
\begin{equation}
   \check D_{2}+\check D_{3}+\cdots+\check D_{r}\,.
\end{equation}
\end{lemma}

{\noindent\bf Proof of  Lemma \ref{msnc}.} 
It follows from a similar argument to Lemma \ref{snc} by applying the Van der Waerden representation. \,\,\,$\endpf$
\medskip

Let $E$ be the exceptional divisor of the blow-up $\tau^1_{s,p,n}:\mathcal M_{s,p,n}\rightarrow \mathbb P(N_{\mathcal V_{(p,0)}/G(p,n)})$ defined in (\ref{checkr}). Then,
\begin{lemma}\label{mexcep}
\begin{equation}
E=\left\{\begin{array}{ll}
   \check D_{2}+\check D_{3}+\cdots+\check D_{r}  & {\rm when}\,\,n-s\neq p  \\
   \check D_{2}+\check D_{3}+\cdots+\check D_{r-1}  & {\rm when}\,\,n-s=p  
\end{array}\right..
\end{equation}
\end{lemma}
{\noindent\bf Proof of  Lemma \ref{mexcep}.} 
The proof is similar to  Lemma \ref{excep}. We can apply the intermediate  Van der Waerden representation in Appendix \ref{section:coverl} to track the iterated blow-ups (\ref{cblow1}). $\endpf$
\medskip

{\noindent\bf Proof of  Proposition \ref{mwond}.} 
It is easy to verify that the  complement of the open $GL(s,\mathbb C)\times GL(n-s,\mathbb C)$-orbit of $\mathcal M_{s,p,n}$ is the union of 
$\check  D_2,\,\check D_3,\, \cdots,\, \check D_r$.

In the following,  we will  determine the  $GL(s,\mathbb C)\times GL(n-s,\mathbb C)$-orbits of $\mathcal M_{s,p,n}$. Let $\left(A^{\tau},(J^{\tau}_0)^{-1}\right)$ be a coordinate chart  in the  Van der Waerden representation of $R_{s,p,n}^{-1}(U_0)$ where $\tau=\left(\begin{matrix} i_1&i_2&\cdots&i_r\\ j_1&j_2&\cdots&j_r\\
\end{matrix}\right)\in\mathbb J_0$. Let $\check{\mathfrak a}$ be a point of $A^{\tau}\cap \mathcal M_{s,p,n}$ with the holomorphic coordinates $\left(\widetilde Y(\check{\mathfrak a}),\overrightarrow B^1(\check{\mathfrak a}),\cdots,\overrightarrow B^r(\check{\mathfrak a})\right)$. Then $\check{\mathfrak b}\in A^{\tau}\cap \mathcal M_{s,p,n}$ is in the same
$GL(s,\mathbb C)\times GL(n-s,\mathbb C)$-orbit of  $\check{\mathfrak a}$ if and only if for $1\leq m\leq r$ the following holds.
\begin{equation}
b_{i_mj_m}\left(\check {\mathfrak b}\right)=0\Longleftrightarrow b_{i_mj_m}(\check{\mathfrak a})=0\,.\\
\end{equation}
Applying the $GL(s,\mathbb C)\times GL(n-s,\mathbb C)$-action, we can verify Property (C) in Definition \ref{wond}.

We complete the proof.\,\,\,$\endpf$
\medskip

{\noindent\bf Proof of  Proposition \ref{gwond}.} 
Let $\mathfrak a\in\mathcal T_{s,p,n}$ be a point not contained in divisors $D_1^-,D^-_2,\cdots$, $D_r^1$, $D_1^+,D^+_2,\cdots$, $D_r^+$. Then   $R_{s,p,n}(\mathfrak a)\subset G(p,n)$ is contained in the stable manifold of the source $\mathcal V_{(p,0)}$ and the unstable manifold of the sink $\mathcal V_{(p-r,r)}$. Then it is clear that $\mathfrak a$ has an open dense $GL(s,\mathbb C)\times GL(n-s,\mathbb C)$-orbit in $\mathcal T_{s,p,n}$. By Lemma \ref{snc}, we can conclude that the complement of the open
$GL(s,\mathbb C)\times GL(n-s,\mathbb C)$-orbit in $\mathcal T_{s,p,n}$ is  a simple normal crossing divisor which  consists of $2r$ smooth, irreducible divisors 
$D^-_1, D^-_2,\cdots,D^-_r,D^+_1,D^+_2,\cdots,D^+_r$.

Let $\left(A^{\tau},(J^{\tau}_0)^{-1}\right)$ be a coordinate chart  in the  Van der Waerden representation of $R_{s,p,n}^{-1}(U_l)$ where $0\leq l\leq r$ and $\tau=\left(\begin{matrix} i_1&i_2&\cdots&i_r\\ j_1&j_2&\cdots&j_r\\
\end{matrix}\right)\in\mathbb J_l$. Similarly to Proposition \ref{mwond}, we can show that
$\mathfrak a, \mathfrak b\in A^{\tau}$
are in the same orbit if and only if for $1\leq m\leq r-l$,
\begin{equation}
b_{i_mj_m}\left(\mathfrak b\right)=0\Longleftrightarrow b_{i_mj_m}(\mathfrak a)=0\,,\\
\end{equation}
and for $r-l+1\leq m\leq r$,
\begin{equation}
a_{i_mj_m}\left(\mathfrak b\right)=0\Longleftrightarrow a_{i_mj_m}(\mathfrak a)=0\,.\\
\end{equation}
Then we can verify that
$GL(s,\mathbb C)\times GL(n-s,\mathbb C)$-orbit of $\mathcal T_{s,p,n}$ one to one corresponds to the  quasi-projective variety $X_{(I^-,I^+)}$ defined by (\ref{inrule}). Moreover, the closure of each $GL(s,\mathbb C)\times GL(n-s,\mathbb C)$-orbit  is smooth.

It is clear that $\mathcal P_{s,p,n}:\mathcal T_{s,p,n}\rightarrow \mathcal M_{s,p,n}=D_1^-$  is a $GL(s,\mathbb C)\times GL(n-s,\mathbb C)$-equivariant flat map by Lemma \ref{flat}. The restriction of $\mathcal P_{s,p,n}$ on $D^+_1$ is an isomorphism by Lemma \ref{sti}; $\mathcal P_{s,p,n}$ is a holomorphic retraction by Corollary \ref{retra}.

Notice that the restriction of the $\mathbb C^*$-action $\psi_{s,p,n}$ on $\mathcal M_{s,p,n}$ is the identity map. Then by the Bia{\l}ynicki-Birula decomposition for $\mathcal T_{s,p,n}$ in Lemma \ref{11}, we can conclude that $\mathcal P_{s,p,n}(D^-_i)=\mathcal P_{s,p,n}(D^+_{r+2-i})$ for $2\leq i\leq r$.

By Proposition \ref{mwond}, 
$D_1^-=\mathcal M_{s,p,n}$ is a wonderful $GL(s,\mathbb C)\times GL(n-s,\mathbb C)$-variety  with  $(r-1)$ $GL(s,\mathbb C)\times GL(n-s,\mathbb C)$-stable divisors $\check D_2,\cdots,\check D_r$.

For $\mathcal T_{s,p,2p}$ (resp. $\mathcal T_{s,p,2s}$), take   $\sigma$ to be USD (resp. DUAL). Property (D) thus follows from Lemmas \ref{dis}, \ref{dis2} (to be proved later in Section \ref{sym}).

We complete the proof of Proposition \ref{gwond}.\,\,\,$\endpf$
\medskip

For convenience, we introduce the following notation.
\begin{definition}
   Denote by $\check R_{s,p,n}$ the blow up map $R_{s,p,n}\big|_{\mathcal M_{s,p,n}}:\mathcal M_{s,p,n}\rightarrow\mathcal V_{(p,0)}$ in (\ref{checkr}). 
\end{definition}
Since the hyperplane bundle of $\mathcal V_{(p,0)}\cong G(p,s)$ is the restriction of the hyperplane bundle of $G(p,n)$, we  denote  by $(\check R_{s,p,n})^*(\mathcal O_{G(p,n)}(1))$ the restriction of the line bundle $(R_{s,p,n})^*(\mathcal O_{G(p,n)}(1))$ to $\mathcal M_{s,p,n}$.  
\begin{lemma}\label{mpicb}
The Picard group Pic$(\mathcal M_{s,p,n})$ of $\mathcal M_{s,p,n}$ is a torsion free abelian group over $\mathbb Z$. When $p<s$ and $p\neq n-s$, Pic$(\mathcal M_{s,p,n})$ has a $\mathbb Z$-basis 
\begin{equation}\label{mbasis1}
    \left\{(\check R_{s,p,n})^*\left(\mathcal O_{ G(p,n)}(1)\right),\check D_1,\check D_{2},\cdots,\check D_{r} \right\};
\end{equation}
when $p<s$ and $p=n-s$, Pic$(\mathcal M_{s,p,n})$ has a $\mathbb Z$-basis 
\begin{equation}
     \left\{(\check R_{s,p,n})^*\left(\mathcal O_{ G(p,n)}(1)\right),\check D_1,\check D_{2},\cdots,\check D_{r-1} \right\};
\end{equation}
when $p=s=n-s$, Pic$(\mathcal M_{s,p,n})$ has a $\mathbb Z$-basis 
\begin{equation}
   \left\{\check D_1,\check D_{2},\cdots,\check D_{r-1} \right\}.
\end{equation}
\end{lemma}
{\noindent\bf Proof of  Lemma \ref{mpicb}.} The proof is the same as in  Lemma \ref{picb}. 
\,\,\,$\endpf$
\medskip

Denote by $K_{\mathcal M_{s,p,n}}$ the canonical bundle of $\mathcal M_{s,p,n}$. We have the following canonical bundle formula.

\begin{lemma}\label{mkan}
Assume that $p\leq s$. Then,
\begin{equation}\label{mkan1}
\begin{split}
K_{\mathcal M_{s,p,n}}=&-n\cdot(\check R_{s,p,n})^*(\mathcal O_{G(p,n)}(1))+p(n-s)\cdot\check D_1\\
+\sum_{i=2}^{r}& \big((p-i+1)(n-s-i+1)-1\big)\cdot\check D_{i}\,.
\end{split}
\end{equation}
\end{lemma}

{\noindent\bf Proof of  Lemma \ref{mkan}.} 
Notice that $D^+_i\cap D^-_1=\emptyset$ for $1\leq i\leq r$. 
Lemma \ref{mkan} follows form  Lemma \ref{kan} by the adjunction formula. \,\,\,$\endpf$

\section{Invariant divisors and curves in \texorpdfstring{$\mathcal T_{s,p,n}$ }{jj} and \texorpdfstring{$\mathcal M_{s,p,n}$ }{jj}} \label{curanddiv}

In this section, we introduce invariant divisors and curves of  $\mathcal T_{s,p,n}$ and $\mathcal M_{s,p,n}$ with respect to certain group actions. By computing the intersection numbers, we prove the (semi-)positivity of the anti-canonical bundles of $\mathcal T_{s,p,n}$ and $\mathcal M_{s,p,n}$.

In the following, we assume that $2p\leq n\leq 2s$ and let $r$ be  the rank of $\mathcal T_{s,p,n}$.

\subsection{The \texorpdfstring{$G$-stable and $B$-stable divisors}{rr}} \label{curanddivg}
For convenience, we denote by
$G$ the group $GL(s,\mathbb C)\times GL(n-s,\mathbb C)$. Define a Borel subgroup $B$  of $G$ by \begin{equation}
B:=\left\{\left. \left(
\begin{matrix}
g_1&0\\
0&g_2\\
\end{matrix}\right)\right\vert_{} 
\begin{matrix}
g_1\in GL(s,\mathbb C) \,\, {\rm is\,\, a\,\, lower\,\, triangular\,\, matrix\,};\,\,\,\,\\
g_2\in GL(n-s,\mathbb C)\,\, {\rm is\,\, an \,\, upper \,\, triangular \,\, matrix}  
\end{matrix}
\right\}.
\end{equation} Let $T=T_1\times T_2$  be  a maximal torus of $B$ defined by \begin{equation}\label{maxt}
T_1\times T_2:=\left\{\left. \left(
\begin{matrix}
g_1&0\\
0&g_2\\
\end{matrix}\right)\right\vert_{} 
\begin{matrix}
g_1\in GL(s,\mathbb C) \,\, {\rm and\,\,}g_2\in GL(n-s,\mathbb C)
\,\, {\rm are\,\, diagonal \,\, matrices} \end{matrix}
\right\}.
\end{equation}
It is clear that $\mathcal T_{s,p,n}$ is a spherical $G$-variety for $B$ has an open orbit in $\mathcal T_{s,p,n}$.

\begin{definition}
We call  $D_{1}^-$, $D_{2}^-$, $\cdots$, $D_{r}^-$,  $D_{1}^+$, $D_{2}^+$, $\cdots$, $D_{r}^+$ the $G$-stable divisors of  $\mathcal T_{s,p,n}$.
\end{definition}

Let $I_j$  be the index defined in (\ref{I_k}), that is, $I_j=(s+j,s+j-1,\cdots,s-p+j+1)\in\mathbb I^j_{s,p,n}$. 
For 
$0\leq j\leq r$, define divisors $b_j$ of $G(p,n)$ by 
\begin{equation}\label{}
 b_{j}:=\left\{x\in G(p,n) \big|\,P_{I_j}(x)=0\,\right\}\,.
\end{equation}
Notice that $b_j$ is biholomporphic to the infinity hyperplane section of $G(p,n)$ which is the closure of a complex Euclidean space; hence, $b_j$ is irreducible for $0\leq j\leq r$.
\begin{definition}
Define divisors $B_{j}\subset\mathcal T_{s,p,n}$,
$0\leq j\leq r$,  to be the strict transformation of $b_j$ under the canonical blow-up $R_{s,p,n}$. We call $B_{0}$, $B_1,\cdots$,$B_r$ the $B$-stable divisors of $\mathcal T_{s,p,n}$.
\end{definition}

It is easy to verify that  $B_j$, $0\leq j\leq r$, is $B$-invariant and irreducible.

\begin{lemma}\label{gb}
Let $\mathfrak D$ be an irreducible divisor of $\mathcal T_{s,p,n}$. If $\mathfrak D$ is $G$-invariant, then
\begin{equation}\label{gi}
\mathfrak D\in \{D_{1}^-, D_{2}^-, \cdots, D_{r}^-, D_{1}^+, D_{2}^+, \cdots, D_{r}^+\}\,.
\end{equation}
If $\mathfrak D$ is $B$-invariant, then
\begin{equation}\label{bi}
\mathfrak D\in \{D_{1}^-, D_{2}^-, \cdots, D_{r}^-, D_{1}^+, D_{2}^+, \cdots, D_{r}^+, B_0, B_1, \cdots, B_r\}\,.
\end{equation}
\end{lemma}
{\noindent\bf Proof of Lemma \ref{gb}.} If $\mathfrak D$ is $G$-invariant, (\ref{gi}) holds by Proposition \ref{gwond}.

Assume that $\mathfrak D$ is $B$-invariant but $\mathfrak D\not\in\left\{D_{1}^-,\cdots, D_{r}^-, D_{1}^+, \cdots, D_{r}^+, B_0,\cdots,B_r\right\}$.  Denote by $\mathfrak d$ the image of $\mathfrak D$ under $R_{s,p,n}$. It is clear that $\mathfrak d$ is a $B$-invariant divisor of $G(p,n)$.  Recall that the exceptional divisor of $R_{s,p,n}$ is contained in the union of the $G$-stable divisors by Lemma \ref{excep}. Then  $P_{I_j}\not\equiv 0$ on $\mathfrak d$ for $0\leq j\leq r$. We can verify that $\mathfrak d$ contains a point $\mathfrak a$ with 
a matrix representative $\widetilde {\mathfrak a}$ defined by \begin{equation}\label{1pp}
\widetilde {\mathfrak a}:=\left(0_{p\times(s-p)}
\hspace{-.13in}\begin{matrix}
  &\hfill\tikzmark{g}\\
  &\hfill\tikzmark{h}
  \end{matrix}\,\,\,\,\begin{matrix}
  I_{p\times p} \hspace{-.1in}\begin{matrix}
  &\hfill\tikzmark{c}\\
  &\hfill\tikzmark{d}
  \end{matrix}\hspace{-.1in}\begin{matrix}
  &\hfill\tikzmark{e}\\
  &\hfill\tikzmark{f}\end{matrix}\,\,\,\,\,I_{p\times p}
  \end{matrix}\hspace{-.1in}
  \begin{matrix}
  &\hfill\tikzmark{a}\\
  &\hfill\tikzmark{b}\end{matrix}\,\,\,\,\,0_{p\times(n-s-p)}\right)\,\,{\rm when}\,\, r=p
  \tikz[remember picture,overlay]   \draw[dashed,dash pattern={on 4pt off 2pt}] ([xshift=0.5\tabcolsep,yshift=7pt]a.north) -- ([xshift=0.5\tabcolsep,yshift=-2pt]b.south);\tikz[remember picture,overlay]   \draw[dashed,dash pattern={on 4pt off 2pt}] ([xshift=0.5\tabcolsep,yshift=7pt]c.north) -- ([xshift=0.5\tabcolsep,yshift=-2pt]d.south);\tikz[remember picture,overlay]   \draw[dashed,dash pattern={on 4pt off 2pt}] ([xshift=0.5\tabcolsep,yshift=7pt]e.north) -- ([xshift=0.5\tabcolsep,yshift=-2pt]f.south);\tikz[remember picture,overlay]   \draw[dashed,dash pattern={on 4pt off 2pt}] ([xshift=0.5\tabcolsep,yshift=7pt]g.north) -- ([xshift=0.5\tabcolsep,yshift=-2pt]h.south);
\end{equation}
or
\begin{equation}\label{1pr}
\widetilde {\mathfrak a}:=\left(\begin{matrix}
  0_{r\times(s-p)}\\
   0_{(p-r)\times(s-p)}
 \end{matrix}
\hspace{-.13in}\begin{matrix}
  &\hfill\tikzmark{g}\\
  &\hfill\tikzmark{h}
  \end{matrix}\,\,\,\,\begin{matrix}
  I_{r\times r}\\
   0_{(p-r)\times r}
  \end{matrix}\hspace{-.1in}\begin{matrix}
  &\hfill\tikzmark{c}\\
  &\hfill\tikzmark{d}
  \end{matrix}\,\,\,\,\,\begin{matrix}
  0_{r\times(p-r)}\\
   I_{(p-r)\times(p-r)}
  \end{matrix}\hspace{-.1in}
\begin{matrix}
  &\hfill\tikzmark{e}\\
  &\hfill\tikzmark{f}\end{matrix}\hspace{-.1in}
  \begin{matrix}
  &\hfill\tikzmark{a}\\
  &\hfill\tikzmark{b}\end{matrix}\,\,\,\,\,
  \begin{matrix}
  I_{r\times r}\\
   0_{(p-r)\times r}
  \end{matrix}\right)\,\,{\rm when}\,\,r=n-s. 
  \tikz[remember picture,overlay]   \draw[dashed,dash pattern={on 4pt off 2pt}] ([xshift=0.5\tabcolsep,yshift=7pt]a.north) -- ([xshift=0.5\tabcolsep,yshift=-2pt]b.south);\tikz[remember picture,overlay]   \draw[dashed,dash pattern={on 4pt off 2pt}] ([xshift=0.5\tabcolsep,yshift=7pt]c.north) -- ([xshift=0.5\tabcolsep,yshift=-2pt]d.south);\tikz[remember picture,overlay]   \draw[dashed,dash pattern={on 4pt off 2pt}] ([xshift=0.5\tabcolsep,yshift=7pt]e.north) -- ([xshift=0.5\tabcolsep,yshift=-2pt]f.south);\tikz[remember picture,overlay]   \draw[dashed,dash pattern={on 4pt off 2pt}] ([xshift=0.5\tabcolsep,yshift=7pt]g.north) -- ([xshift=0.5\tabcolsep,yshift=-2pt]h.south);
\end{equation}
Then $\mathfrak a$ is in a dense open $B$-orbit of $G(p,n)$, which is a contradiction.

We complete the proof of Lemma \ref{gb}.\,\,\,$\endpf$.
\begin{remark}\label{coin}
When $p=n-s$, $B_r=D^-_r$; when  $p=s$, $B_0=D^+_r$. The notion of $B$-stable divisors in this paper differs slightly from that in the theory of spherical varieties. 
$B_0$ (resp. $B_r$) is not a $B$-stable divisor when $p=s$ (resp. $p=n-s$) therein. 
\end{remark}

We define the following convenient local coordinate charts to investigate the invariant divisors in a quantitative way.
\begin{definition}\label{taul}
For $0\leq l\leq r$, let $\left(A^{\tau_l},\left( J^{\tau_l}_l\right)^{-1}\right)$  be the local coordinate chart in the Van der Waerden representation of $R_{s,p,n}^{-1}(U_l)$ associated with the following index. 
\begin{equation}
  \tau_l:=\left(\begin{matrix}
		1+1&l+2&\cdots&r&l&l-1&\cdots&1\\
		s+l+1&s+l+2&\cdots&s+r&s-p+l&s-p+l-1&\cdots&s-p+1\\
	\end{matrix}\right)\in\mathbb J_l\,. 
\end{equation}
We call $\left(A^{\tau_l},\left( J^{\tau_l}_l\right)^{-1}\right)$  (or $A^{\tau_l}$ by a slight abuse of notation) the $l$-th main coordinate chart.
\end{definition}

\begin{lemma}\label{orbits}
If $B_j$ contains a non-empty $G$-orbit of  $\mathcal T_{s,p,n}$ for a certain $0\leq j\leq r$, then   $j=r$ and $p=n-s$, or $j=0$ and $p=s$.
\end{lemma}

{\noindent\bf Proof of Lemma \ref{orbits}.} By Proposition \ref{gwond},  for each $G$-invariant closed subvariety $\mathfrak C$ of  $\mathcal T_{s,p,n}$ there is an integer $0\leq l\leq r$ and a dense open subset $\mathfrak c$ of $\mathfrak C$, such that $\mathfrak c$ is contained in each local coordinate chart in the Van der Waerden representation of $R_{s,p,n}^{-1}(U_l)$.
Therefore, there is an integer  $0\leq l\leq r$ such that $B_j\cap A^{\tau_l}\neq\emptyset$  where $\left(A^{\tau_l},\left( J^{\tau_l}_l\right)^{-1}\right)$ is the $l$-th main coordinate chart.

In terms of the holomorphic coordinates of  $\left(A^{\tau_l},\left( J^{\tau_l}_l\right)^{-1}\right)$, the defining function of the total transformation of $b_j$ takes the following form. 
\begin{equation}\label{te}
    (R_{s,p,n})^*(P_{I_j})=\left\{ \begin{array}{cl}
    1\,\,\,\,\,\,\,\,\,\,\,\,\,\,&{\rm if}\,\,j=l\\
    \prod\limits_{t=j+1}^{l} a^{t-j}_{t(s-p+t)}&{\rm if}\,\,j\leq l-1\\
    \prod\limits_{t=l+1}^{j}b^{j+1-t}_{t(s+t)}\,\,\,\,\,&{\rm if}\,\,j\geq l+1\,\,\\
    \end{array}\right..
\end{equation}
It is clear that $B_j\cap A^{\tau_l}\neq\emptyset$ if and only if $B_j$ coincides with a certain $G$-invariant divisor.

We complete the proof of Lemma \ref{orbits}.\,\,\,$\endpf$
\medskip

Brion (\cite{Br5}) expressed the canonical bundles of spherical varieties in terms of $G$-invariant and $B$-invariant divisors. We compute the coefficients therein for $K_{\mathcal T_{s,p,n}}$ explicitly as follows. (When $n=2s=2p$ it follows from \cite{DP} and \cite{Ri2}.)
\begin{lemma}\label{wk}
When  $p<n-s\leq s$,
\begin{equation}
    K_{\mathcal T_{s,p,n}}=-(s-p+1)\cdot B_0 -2\sum_{j=1}^{p-1}B_j -(n-s-p+1)\cdot B_p -\sum_{i=1}^pD^-_i-\sum_{i=1}^pD^+_i\,;
\end{equation}
when $n-s=p<s$,
\begin{equation}
    K_{\mathcal T_{s,p,n}}=-(s-p+1)\cdot B_0 -2\sum_{j=1}^{p-1}B_j -\sum_{i=1}^pD^-_i-\sum_{i=1}^pD^+_i\,;
\end{equation}
when $n-s<p<s$ ($r=n-s$),
\begin{equation}
    K_{\mathcal T_{s,p,n}}=-(s-p+1)\cdot B_0-2\sum_{j=1}^{r-1}B_j-(p-r+1)\cdot B_r-\sum_{i=1}^rD^-_i-\sum_{i=1}^rD^+_i\,;
\end{equation}
when $n-s=p=s$,
\begin{equation}
    K_{\mathcal T_{s,p,n}}= -2\sum_{j=1}^{p-1}B_j  -\sum_{i=1}^pD^-_i-\sum_{i=1}^pD^+_i\,.
\end{equation}
\end{lemma}
{\noindent\bf Proof of  Lemma \ref{wk}.}  Compute the defining functions of $B_j$ in the $0$-th and $r$-th main  coordinate charts. Then by (\ref{te}) we can derive the following for $0\leq j\leq r$.
\begin{equation}\label{bst}
\begin{split}
 B_j&=(R_{s,p,n})^*\left(\mathcal O_{G(p,n)}(1)\right)-\sum_{i=1}^{r-j}(r-j+1-i)\cdot D^+_i-\sum_{i=1}^j(j+1-i)\cdot D^-_i.
\end{split}
\end{equation}
Notice that by Remark \ref{coin} if $p=s$ (resp. $p=n-s$) we should modify (\ref{bst}) for $B_0$ (resp. $B_r$) as follows. When $p=s$, 
\begin{equation}\label{bstt0}
\begin{split}
B_0=D^+_r&=(R_{s,p,n})^*\left(\mathcal O_{G(p,n)}(1)\right)-\sum_{i=1}^{r-1}(r+1-i)\cdot D^+_i\,;
\end{split}
\end{equation}
when $p=n-s$,
\begin{equation}\label{bsttr}
\begin{split}
  B_r=D^-_r&=(R_{s,p,n})^*\left(\mathcal O_{G(p,n)}(1)\right)-\sum_{i=1}^{r-1}(r+1-i)\cdot D^-_i\,.
\end{split}
\end{equation}

Combining with Lemma \ref{kan}, we can conclude Lemma \ref{wk}. \,\,\,\,$\endpf$
\medskip

In the following, we give a geometric interpretation of the $B$-stable divisors.
For $0\leq j\leq p$, denote by  $H_j$ the positive generator of the Picard group of $\mathbb {CP}^{N^j_{s,p,n}}$; by a slight abuse of notation, denote by $H_j$ its pullback to  $\mathbb {CP}^{N_{p,n}}\times\mathbb {CP}^{N^0_{s,p,n}}\times\cdots\times\mathbb {CP}^{N^p_{s,p,n}}$ under the corresponding projection.  \begin{lemma}\label{partialline} Let $H_j\big|_{\mathcal T_{s,p,n}}$ be the restriction of  $H_j$ to $\mathcal T_{s,p,n}$. Then for $0\leq j\leq r$,
\begin{equation}
      H_j\big|_{\mathcal T_{s,p,n}}=R^*_{s,p,n}(\mathcal O_{G(p,n)}(1))-\sum\limits_{i=1}^{r-j}(r-i+1-j)\cdot D^+_{i}-\sum\limits_{i=1}^{j}(j+1-i)\cdot D^-_{i}\,.\,\,
    \end{equation}
\end{lemma}
{\bf\noindent Proof of Lemma \ref{partialline}.} 
For $0\leq j\leq r$, compute the Pl\"ucker coordinate functions $P_{I}$, $I\in\mathbb I^j_{s,p,n}$, in the $0$-th and $r$-th main coordinate charts. Then Lemma \ref{partialline} follows from (\ref{a1}).\,\,\,\,\,$\endpf$
\medskip

Recall that a line bundle  $L$ is called basepoint-free if the intersection of the zero sets of all global sections of $L$ is empty. 
\begin{lemma}\label{basepoint}
The line bundle $B_j$ is basepoint-free in the following cases.
\begin{enumerate}
    \item $p\neq s$, $p\neq n-s$ and $0\leq j\leq r$.
    \item $n-s=p<s$ and $0\leq j\leq p-1$.
    \item $n-s=p=s$ and $1\leq j\leq p-1$.
\end{enumerate}
\end{lemma}
{\bf\noindent Proof of Lemma \ref{basepoint}.} Recall that $  H_j\big|_{\mathcal T_{s,p,n}}$ is a basepoint-free line bundle of $\mathcal T_{s,p,n}$ for $0\leq j\leq p$ provided it is nontrivial. Then Lemma \ref{basepoint} follows from (\ref{bst}), (\ref{bstt0}), (\ref{bsttr}), and Lemma \ref{partialline}.\,\,\,\,\,$\endpf$
\medskip

For a sequence of integers $a_j$, $0\leq j\leq r$, define a line bundle $L_{a_0,\cdots,a_r}$ by
\begin{equation}
   L_{a_0,\cdots,a_r}:=\sum_{j=0}^{r}a_j\cdot B_j+  R^*_{s,p,n}(\mathcal O_{G(p,n)}(1)).
\end{equation}

\begin{lemma}\label{ample}
 $L_{a_0,\cdots,a_r}$ is ample in the following cases.
\begin{enumerate}[label=(\alph*).]
    \item $p\neq n-s$, $p\neq s$, and $a_j>0$ for $0\leq j\leq r$.
    \item $p=n-s<s$, $a_p=0$, and $a_j>0$ for $0\leq j\leq p-1$.
    \item $p=n-s=s$, $a_0=a_p=0$, and $a_j>0$ for $1\leq j\leq p-1$.
\end{enumerate}
\end{lemma}
{\bf\noindent Proof of Lemma \ref{ample}.} Recall that $\sum_{j=0}^{r}H_j+R^*_{s,p,n}(\mathcal O_{G(p,n)}(1))$ is an ample line bundle on $\mathbb {CP}^{N_{p,n}}\times\mathbb {CP}^{N^0_{s,p,n}}\times\cdots\times\mathbb {CP}^{N^p_{s,p,n}}$. Notice that $H_0$ is trivial when $p=n-s<s$ and $H_0$ and $H_r$ are trivial when $p=n-s=s$. Then Lemma \ref{ample} follows from (\ref{bst}) and Lemma \ref{partialline}.\,\,\,\,\,$\endpf$
\begin{remark}\label{cls}
One can show that the complete linear series of $H_j\big|_{\mathcal T_{s,p,n}}$ on $\mathcal T_{s,p,n}$ is isomorphic to $\mathbb {CP}^{N^j_{s,p,n}}$.
This is due to the fact that for a grassmannian $G(p,n)$ the complete linear series of $\mathcal O_{G(p,n)}(1)$ is generated by the Pl\"ucker coordinate functions (see \cite{Da}).
\end{remark}

We can define invariant divisors of $\mathcal M_{s,p,n}$ in a similar way for  $\mathcal M_{s,p,n}$ is a spherical $G$-variety with respect to $B$ as well.
\begin{definition}
We call the divisors $\check D_{2}$, $\check D_{3}$, $\cdots$, $\check D_{r}$ defined by (\ref{checkd}) the $G$-stable divisors of  $\mathcal M_{s,p,n}$.
\end{definition}

\begin{definition}
For $0\leq i\leq r$, define divisors $\check B_i$ of $\mathcal M_{s,p,n}$ by
\begin{equation}
    \check B_i:=D_{1}^-\cap B_{i}=\mathcal M_{s,p,n}\cap B_{i}\,.
\end{equation}
We call  $\check B_i$, $0\leq i\leq r$,  the $B$-stable divisors of  $\mathcal M_{s,p,n}$.
\end{definition}
\begin{remark}\label{checkcoin}
$\check B_0=D_{1}^-\cap B_{0}=D_{1}^-\cap D^+_{r}=\emptyset$ is trivial when $p=s$; $\check B_r=\check D_r$ when $p=n-s$. By a slightly abuse of notation, we view $\check B_0$ as a $B$-stable divisor even if $p=s$.
\end{remark}

\begin{lemma}\label{checkbir}
$\check B_j$ is an irreducible divisor if $1\leq j\leq r$, or $j=0$ and $p<s$.
\end{lemma}
{\bf\noindent Proof of Lemma
\ref{checkbir}.} When $1\leq j\leq r$,  or $j=0$ and $p<s$, we can conclude that \begin{equation}
    \mathcal P_{s,p,n}^{-1}\left(\check B_j\right)=B_j\,.
\end{equation}
Then Lemma \ref{checkbir} follows from the fact that $B_j$ is irreducible.
\,\,\,\,$\endpf$
\medskip

Similarly to Lemmas \ref{gb} and \ref{orbits}, we have that
\begin{lemma}\label{checkgb}
Let $\check {\mathfrak D}$ be an irreducible divisor of $\mathcal M_{s,p,n}$. If $\check {\mathfrak D}$ is $G$-invariant, then
\begin{equation}\label{checkgi}
\check {\mathfrak D}\in \{\check D_{2}, \check D_{3}, \cdots, \check D_{r}\}\,.
\end{equation}
If $\check {\mathfrak D}$ is $B$-invariant, then
\begin{equation}\label{checkbi}
\check {\mathfrak D}\in \{\check D_{2}, \check D_{3}, \cdots, \check D_{r}, \check B_0, \check B_1, \cdots, \check B_r\}\,.
\end{equation}
\end{lemma}
\begin{lemma}\label{checkorbits}
If $\check B_j$ contains a non-empty $G$-orbit of  $\mathcal M_{s,p,n}$ for a certain $0\leq j\leq r$, then   $j=r$ and $p=n-s$.
\end{lemma}

Notice that $D^+_i\cap D^-_0=\emptyset$ for $1\leq i\leq r$. By restricting (\ref{bst}), (\ref{bsttr}) to $\mathcal M_{s,p,n}$, we then have
\begin{equation}\label{mb=b}
\begin{split}
&\check B_i=(\check R_{s,p,n})^*(\mathcal O_{G(p,n)}(1))-i\cdot \check D_1-\sum_{k=2}^i (i+1-k)\cdot\check D_k\,,\,\,\,\,0\leq i\leq r\,;
\end{split}
\end{equation}
when $p=n-s$, we  modify (\ref{mb=b}) for $\check B_r$ by 
\begin{equation}\label{mb=br}
 \check B_r=\check D_r=(\check R_{s,p,n})^*(\mathcal O_{G(p,n)}(1))-r\cdot \check D_1-\sum_{k=2}^{r-1}(r+1-k)\cdot\check D_k\,.
\end{equation}

Similarly to Lemma \ref{wk}, we have
\begin{lemma}\label{mkb}
When  $p<n-s$ and $p<s$,
\begin{equation}
    K_{\mathcal M_{s,p,n}}=-(s-p+1)\cdot \check B_0 -2\sum_{j=1}^{p-1}\check B_j -(n-s-p+1)\cdot \check B_p -\sum_{i=2}^p\check D_i\,;
\end{equation}
when $n-s=p<s$,
\begin{equation}
    K_{\mathcal M_{s,p,n}}=-(s-p+1)\cdot \check B_0 -2\sum_{j=1}^{p-1}\check B_j -\sum_{i=2}^p\check D_i\,;
\end{equation}
when $n-s<p$ and $p<s$ ($r=n-s$),
\begin{equation}\label{m3}
    K_{\mathcal M_{s,p,n}}=-(s-p+1)\cdot \check B_0-2\sum_{j=1}^{r-1}\check B_j-(p-r+1)\cdot \check B_r-\sum_{i=2}^r\check D_i\,;
\end{equation}
when $n-s=p=s$,
\begin{equation}
    K_{\mathcal M_{s,p,n}}= -2\sum_{j=1}^{p-1}\check B_j  -\sum_{i=2}^p\check D_i\,.
\end{equation}
\end{lemma}
{\noindent\bf Proof of  Lemma \ref{mkb}.} 
Lemma \ref{mkb} follows from  Lemma \ref{wk} by the adjunction formula.\,\,\,$\endpf$

\subsection{\texorpdfstring{$T$-invariant curves}{rr}} \label{curanddivi}

In this subsection, we introduce $T$-invariant curves  $\gamma_l$, $\zeta^l_j$, $\zeta_{u,v}^{l,k}$, $\delta_{m_1,m_2}^l$, $\Delta_{m_1,m_2}^l$ of $\mathcal T_{s,p,n}$ with various parameters as the closures of certain affine lines in the $l$-th main coordinate chart $A^{\tau_l}$, $0\leq l\leq r$.

By the geometric interpretation of the $B$-invariant divisors given in Section \ref{curanddivg}, we compute the degrees of the line bundles of  $\mathcal T_{s,p,n}$ on such curves.  Here we give a detailed proof for the most important case $\gamma_l$, and leave the remainder to  Appendix \ref{section:inters}.

As an immediate corollary, we can calculate the intersection numbers of such curves and the canonical bundle $K_{\mathcal T_{s,p,n}}$, which is used in the next subsection to establish the (semi-)positivity of the anti-canonical bundles. Since the proof is routine, we omit it.
\smallskip

{\bf\noindent (0). Curves $\gamma_l$, $0\leq l\leq r-1$.}
\smallskip

In terms of the local coordinates, we can define an affine line $\mathring{\gamma}_l:\mathbb C\rightarrow A^{\tau_l}$, $0\leq l\leq r-1$,  as follows. $\widetilde X\left(\mathring{\gamma}_l(t)\right)$ and $\widetilde Y\left(\mathring{\gamma}_l(t)\right)$ are zero matrices for each $t\in\mathbb C$; $\overrightarrow B^{k}\left(\mathring{\gamma}_l(t)\right)$ is a zero vector for each $t\in\mathbb C$ and $2\leq k\leq r$\,;
\begin{equation}
\begin{split}
\overrightarrow B^{1}\left(\mathring{\gamma}_l(t)\right)&=\left(b_{(l+1)(s+l+1)},\xi^{(1)}_{(l+1)(s+l+2)},\xi^{(1)}_{(l+1)(s+l+3)},\cdots,\xi^{(1)}_{(l+1)n},\right.\\ 
&\,\,\,\,\,\,\,\,\,\,\,\,\,\,\,\,\,\,\,\,\,\,\left.\xi^{(1)}_{(l+2)(s+l+1)},\xi^{(1)}_{(l+3)(s+l+1)},\cdots,\xi^{(1)}_{p(s+l+1)}\right)\left(\mathring{\gamma}_l(t)\right)=\left(t,0,0\cdots,0\right)\,.
\end{split}
\end{equation}

Define $\gamma_l$ to be the closure of ${\mathring{\gamma}_l}$ in $\mathcal T_{s,p,n}$.

\begin{lemma}\label{i1} For $0\leq l\leq r-1$,
\begin{equation}
\left\{
\begin{array}{ll}
   \gamma_l\cdot R^*_{s,p,n}((\mathcal O_{G(p,n)}(1))=1\,\,\,\, & \\[1em]
  \gamma_l\cdot D^-_{i}=\left\{ \begin{array}{cl}
       1\,\,\,\,&   i=l+1\\
       -1\,\,\,\,& i=l+2\,\,{\rm and\,\,}1\leq i\leq r \\
       0 \,\,\,\,   & {\rm otherwise\,\,} \\
  \end{array}\right.  &  \\[2.3em]
  \gamma_l\cdot D^+_{i}=\left\{ \begin{array}{cl}
       1\,\,\,\,&   i=r-l\\
       -1\,\,\,\,& i=r-l+1\,\, {\rm and\,\,} 1\leq i\leq r \\
       0 \,\,\,\,   & {\rm otherwise\,\,}  \\
  \end{array}\right.   &  \\
\end{array}\right..
\end{equation}
\end{lemma}
{\bf\noindent Proof of Lemma \ref{i1}.} We view ${\gamma}_l$ as a curve in $\mathbb {CP}^{N_{p,n}} \times\mathbb{CP}^{N^0_{s,p,n}} \times\cdots\times\mathbb {CP}^{N^p_{s,p,n}}$. 

By computing the corresponding Pl\"ucker coordinate functions, it is easy to verify that the projections of ${\gamma}_l$ to the components of the ambient space satisfy the following properties.
\begin{equation}
\begin{split}
&(1)\,\,\,\,{\rm The\,\, projection\,\, of\,\,}{\gamma}_l\,\,{\rm  to}\,\,\mathbb {CP}^{N_{p,n}}\,\,{\rm is\,\, a\,\, line};\\
&(2)\,\,\,\,{\rm The\,\, projection\,\, of\,\,}{\gamma}_l\,\,{\rm  to}\,\,\mathbb {CP}^{N^j_{s,p,n}}\,\,{\rm is\,\, a\,\, point\,\,for\,\,}0\leq j\leq p.\\
\end{split}
\end{equation}
Therefore,  we have 
\begin{equation}
\begin{split}
&(1)\,\,\,\,{\gamma}_l\cdot R^*_{s,p,n}((\mathcal O_{G(p,n)}(1))=1;\\
&(2)\,\,\,\,{\gamma}_l\cdot H_j\big|_{\mathcal T_{s,p,n}}=0,\,\,0\leq j\leq p,\,\,{\rm where\,\,}H_j\big|_{\mathcal T_{s,p,n}}\,\,{\rm is\,\,the\,\,restriction\,\, of\,\,the\,\, positive\,\,}\\
&\,\,\,\,\,\,\,\,\,\,\,\,\,\,{\rm  generator\,\,}
H_j\,\,{\rm of\,\, the\,\, Picard\,\, group\,\, of \,\,}\mathbb {CP}^{N^j_{s,p,n}}.\\
\end{split}
\end{equation}

Notice that $\gamma_l$ is disjoint with the divisors
\begin{equation}
    D_1^-,D^-_2,\cdots,D^-_l,D^+_1,D^+_2,\cdots,D^+_{r-l-1}\,.
\end{equation}
Combined with Lemma \ref{partialline}, we can conclude Lemma \ref{i1} by solving linear equations.\,\,\,\,$\endpf$
\medskip

Combining with Lemma \ref{kan}, we can derive the following by a direct computation.
\begin{lemma}\label{vki1} When $r=1$,
\begin{equation}
\small
     -K_{\mathcal T_{s,p,n}}\cdot\gamma_0=2\,.
\end{equation}
When $r\geq 2$,
\begin{equation}
\small
-K_{\mathcal T_{s,p,n}}\cdot\gamma_l=\left\{
\begin{array}{ll}
   1&\,\,\,{\rm when}\,\,\, l=0\\[1em]
   0&\,\,\,{\rm when}\,\,\, 1\leq l\leq r-2\\[1em]
  1&\,\,\,{\rm when}\,\,\, l=r-1\\
\end{array}\right..
\end{equation}
\end{lemma}

{\bf\noindent (1). Curves $\zeta^0_j$, $\zeta_{u,v}^{0,k}$, $\delta_{m_1,m_2}^0$ in $R_{s,p,n}^{-1}\left(\mathcal V_{(p,0)}\right)$.}
\smallskip

For $2\leq j\leq r$, define an affine line $\mathring{\zeta}^0_j\subset A^{\tau_0}$ in terms of the local coordinates as follows. $\widetilde Y\left(\mathring{\zeta}^0_j(t)\right)\equiv \overrightarrow 0$\,;
$\overrightarrow B^{k}\left(\mathring{\zeta}^0_j(t)\right)\equiv \overrightarrow 0$ for $1\leq k\leq r$ and $k\neq j$\,;
\begin{equation}
\small
\begin{split}
\overrightarrow B^{j}\left(\mathring{\zeta}^0_j(t)\right)&=\left(b_{j(s+j)},\xi^{(j)}_{j(s+j+1)},\xi^{(j)}_{j(s+j+2)},\cdots,\xi^{(j)}_{jn},\xi^{(j)}_{(j+1)(s+j)},\xi^{(j)}_{(j+2)(s+j)},\cdots,\xi^{(j)}_{p(s+j)}\right)\left(\mathring{\zeta}_j(t)\right)\\ 
&=\left(t,0,0\cdots,0\right)\,.
\end{split}
\end{equation}

Define  $\zeta^0_j$,  $2\leq j\leq r$, to be the closure of $\mathring{\zeta}^0_j$.

\begin{lemma}\label{i2} 
For $2\leq j\leq r$,
\begin{equation}
\left\{
\begin{array}{ll}
   \zeta^0_j\cdot R^*_{s,p,n}((\mathcal O_{G(p,n)}(1))=0\,\,\,\, & \\[1em]
  \zeta_j^0\cdot D^-_{i}=\left\{ \begin{array}{cl}
       -1\,\,\,\,&   i=j-1\\
       2\,\,\,\,& i=j\,\,\\
       -1\,\,\,\,&   i=j+1\,\,{\rm and\,\,}1\leq i\leq r \\
       0 \,\,\,\,   & {\rm otherwise\,\,}  \\
  \end{array}\right.  &  \\[2.3em]
  \zeta_j^0\cdot D^+_{i}=0  & \hspace{-1.9in} 1\leq i\leq r \\
\end{array}\right..
\end{equation}
\end{lemma}

\begin{lemma}\label{vki2} 
For $2\leq j\leq r$,
\begin{equation}
\small
-K_{\mathcal T_{s,p,n}}\cdot\zeta^0_j=\left\{
\begin{array}{ll}
   3&\,\,\,{\rm when}\,\,\,  j=r\geq 2\\[1em]
   2&\,\,\,{\rm when}\,\,\, 2\leq j\leq r-1\\
\end{array}\right..
\end{equation}
\end{lemma}

Assume that $u=k$ and $s+k+1\leq v\leq n$, or $v=s+k$  and $k+1\leq u\leq p$ where  $1\leq k\leq r$.  Define an affine line $\mathring{\zeta}_{u,v}^{0,k}\subset A^{\tau_0}$ in terms of the local coordinates as follows. $\widetilde Y\left(\mathring{\zeta}_{u,v}^{0,k}(t)\right)\equiv \overrightarrow 0$\,; $\overrightarrow B^{m}\left(\mathring{\zeta}_{u,v}^{0,k}(t)\right)\equiv \overrightarrow 0$ for $1\leq m\leq r$ and $m\neq k$\,;
\begin{equation}
\small
\begin{split}
\overrightarrow B^{k}\left(\mathring{\zeta}_{u,v}^{0,k}(t)\right)&=\left(b_{k(s+k)},\xi^{(k)}_{k(s+k+1)},\xi^{(k)}_{k(s+k+2)},\cdots,\xi^{(k)}_{kn},\xi^{(k)}_{(k+1)(s+k)},\cdots,\xi^{(k)}_{p(s+k)}\right)\left(\mathring{\zeta}_{u,v}^{0,k}(t)\right)\\
&=\left(0,\right.0,\cdots,0,\underset{\substack{\uparrow\\\xi^{(k)}_{uv}}}{t,}0,\cdots,0\left.\right)\,.
\end{split}
\end{equation}

Let $\zeta_{u,v}^{0,k}$ be the closure of  $\mathring{\zeta}_{u,v}^{0,k}$.

\begin{lemma}\label{i3} 
Assume that  $1\leq k\leq r$. When $u=k$ and $v=s+k+1\leq n$, or $v=s+k$ and $u=k+1\leq p$,
\begin{equation}
\left\{
\begin{array}{ll}
   \zeta^{0,k}_{u,v}\cdot R^*_{s,p,n}((\mathcal O_{G(p,n)}(1))=0\,\,\,\, & \\[1em]
  \zeta^{0,k}_{u,v}\cdot D^-_{i}=\left\{ \begin{array}{cl}
       -1\,\,\,\,&   i=k\\
       2\,\,\,\,& i=k+1\,\,{\rm and\,\,}1\leq i\leq r \\
       -1\,\,\,\,&   i=k+2\,\,{\rm and\,\,}1\leq i\leq r \\
       0 \,\,\,\,   &{\rm otherwise\,\,}\\
  \end{array}\right.  &  \\[2.3em]
  \zeta^{0,k}_{u,v}\cdot D^+_{i}=0  & \hspace{-1.9in} 1\leq i\leq r \\
\end{array}\right..
\end{equation}
When $u=k$ and  $s+k+2\leq v\leq n$, 
\begin{equation}
\left\{
\begin{array}{ll}
   \zeta^{0,k}_{u,v}\cdot R^*_{s,p,n}((\mathcal O_{G(p,n)}(1))=0\,\,\,\, & \\[1em]
  \zeta^{0,k}_{u,v}\cdot D^-_{i}=\left\{ \begin{array}{cl}
       -1\,\,\,\,&   i=k\\
       1\,\,\,\,& i=k+1\,\,{\rm and\,\,}1\leq i\leq r \\
       1\,\,\,\,&   i=v-s\,\,{\rm and\,\,}1\leq i\leq r \\
       -1\,\,\,\,&   i=v-s+1\,\,{\rm and\,\,}1\leq i\leq r \\
       0 \,\,\,\,   & {\rm otherwise\,\,}\\
  \end{array}\right.  &  \\[2.3em]
  \zeta^{0,k}_{u,v}\cdot D^+_{i}=0  & \hspace{-2.2in} 1\leq i\leq r \\
\end{array}\right..
\end{equation}
When $v=s+k$ and  $k+2\leq u\leq p$, 
\begin{equation}
\left\{
\begin{array}{ll}
   \zeta^{0,k}_{u,v}\cdot R^*_{s,p,n}((\mathcal O_{G(p,n)}(1))=0\,\,\,\, & \\[1em]
  \zeta^{0,k}_{u,v}\cdot D^-_{i}=\left\{ \begin{array}{cl}
       -1\,\,\,\,&   i=k\\
       1\,\,\,\,& i=k+1\,\,{\rm and\,\,}1\leq i\leq r \\
       1\,\,\,\,&   i=u\,\,\,\,\,{\rm and\,\,}\,\,\,1\leq i\leq r \\
       -1\,\,\,\,&   i=u+1\,\,{\rm and\,\,}1\leq i\leq r \\
       0 \,\,\,\,   & {\rm otherwise\,\,} \\
  \end{array}\right.  &  \\[2.3em]
  \zeta^{0,k}_{u,v}\cdot D^+_{i}=0  & \hspace{-1.9in} 1\leq i\leq r \\
\end{array}\right..
\end{equation}
\end{lemma}

\begin{lemma}\label{vki3}
Assume that  $1\leq k\leq r$. Then,
\begin{equation}
\small
-K_{\mathcal T_{s,p,n}}\cdot\zeta^{0,k}_{k,v}=\left\{
\begin{array}{ll}
   2(v-s-k)&\,\,{\rm when}\,\, s+k+1\leq v\leq r+s-1\\[1em]
   2(r-k)+1&\,\,{\rm when}\,\, v=r+s\,\,{\rm and}\,\,1\leq k\leq r-1\\[1em]
   n-s+p-2k+1&\,\,{\rm when}\,\, r+s+1\leq v\leq n\,\,{\rm and}\,\,1\leq k\leq r-1\\[1em]
   n-s+p-2r&\,\,{\rm when}\,\, r+s+1\leq v\leq n\,\,{\rm and}\,\,k=r\\
\end{array}\right.
\end{equation}
and
\begin{equation}
\small
-K_{\mathcal T_{s,p,n}}\cdot\zeta^{0,k}_{u,s+k}=\left\{
\begin{array}{ll}
   2(u-k)&\,\,{\rm when}\,\, k+1\leq u\leq r-1\,\,{\rm and}\\[1em]
   2(r-k)+1&\,\,{\rm when}\,\, u=r\,\,{\rm and}\,\,1\leq k\leq r-1\\[1em]
   n-s+p-2k+1&\,\,{\rm when}\,\, r+1\leq u\leq p\,\,{\rm and}\,\,1\leq k\leq r-1\\[1em]
   n-s+p-2r&\,\,{\rm when}\,\, r+1\leq u\leq p\,\,{\rm and}\,\,k=r\\
\end{array}\right..
\end{equation}
\end{lemma}

When $1\leq m_1\leq p$ and $1\leq m_2\leq s-p$, define an affine line $\mathring{\delta}_{m_1,m_2}^0\subset  A^{\tau_0}$ in terms of the local coordinates as follows.  $y_{m_1(s-p+1-m_2)}\left(\mathring{\delta}_{m_1,m_2}^0(t)\right)=t$ and the other variables are constantly zero.

Let $\delta_{m_1,m_2}^0$  be the closure of   $\mathring{\delta}_{m_1,m_2}^0$.

\begin{lemma}\label{i4}
For $1\leq m_1\leq p$ and $1\leq m_2\leq s-p$, 
\begin{equation}
\left\{
\begin{array}{ll}
   \delta_{m_1,m_2}^0\cdot R^*_{s,p,n}((\mathcal O_{G(p,n)}(1))=1\,\,\,\, & \\[1em]
  \delta_{m_1,m_2}^0\cdot D^-_{i}=\left\{ \begin{array}{cl}
       1\,\,\,\,&   i=m_1\,\, {\rm and}\,\,1\leq i\leq r \\
       -1\,\,\,\,&   i=m_1+1\,\,{\rm and}\,\,1\leq i\leq r\\
       0 \,\,\,\,   &{\rm otherwise}\,\,\\
  \end{array}\right.  &  \\[2.3em]
  \delta_{m_1,m_2}^0\cdot D^+_{i}=0  & \hspace{-1.6in} 1\leq i\leq r \\
\end{array}\right..
\end{equation}
\end{lemma}

\begin{lemma}\label{vki4}
For $1\leq m_1\leq p$ and $1\leq m_2\leq s-p$, 
\begin{equation}
\small
-K_{\mathcal T_{s,p,n}}\cdot\delta_{m_1,m_2}^0=\left\{
\begin{array}{ll}
   2m_1+s-p-1&\,\,{\rm when}\,\,1\leq m_1\leq r-1\\[1em]
   2r+s-p&\,\,{\rm when}\,\,m_1=r\\[1em]
   n&\,\,{\rm when}\,\,r+1\leq m_1\leq p\\
\end{array}\right..
\end{equation}
\end{lemma}

{\bf\noindent (2). Curves  $\zeta^l_j$, $\zeta_{u,v}^{l,k}$, $\delta_{m_1,m_2}^l$, $\Delta_{m_1,m_2}^l$  in $R_{s,p,n}^{-1}\left(\mathcal V_{(p-l,l)}\right)$, $1\leq l\leq r-1$.}
\smallskip

Assume that $1\leq l\leq r-1$.

For $2\leq j\leq r-l$, define an affine line $\mathring{\zeta}^l_j\subset A^{\tau_l}$  in terms of the local coordinates as follows.  $\widetilde X\left(\mathring{\zeta}^l_j(t)\right)\equiv \overrightarrow 0$\,; $\widetilde Y\left(\mathring{\zeta}^l_j(t)\right)\equiv \overrightarrow 0$\,; $\overrightarrow B^{k}\left(\mathring{\zeta}^l_j(t)\right)\equiv \overrightarrow 0$ 
for $1\leq k\leq r$ and $k\neq j$;
\begin{equation}
\begin{split}
\overrightarrow B^{j}\left(\mathring{\zeta}^l_j(t)\right)&=\left(b_{(l+j)(s+l+j)},\xi^{(j)}_{(l+j)(s+l+j+1)},\xi^{(j)}_{(l+j)(s+l+j+2)},\cdots,\xi^{(j)}_{(l+j)n},\right.\,\,\,\\
&\,\,\,\,\,\,\,\,\,\,\,\,\,\,\,\,\,\,\,\,\,\,\,\,\,\,\,\,\,\,\left.\xi^{(j)}_{(l+j+1)(s+l+j)},\xi^{(j)}_{(l+j+2)(s+l+j)},\cdots,\xi^{(j)}_{p(s+l+j)}\right)\left(\mathring{\zeta}^l_j(t)\right)\\ 
&=\left(t,0,0\cdots,0\right)\,.
\end{split}
\end{equation}

For $r-l+2\leq j\leq r$, define an affine line $\mathring{\zeta}^l_j\subset A^{\tau_l}$ as follows.  $\widetilde X\left(\mathring{\zeta}^l_j(t)\right)\equiv \overrightarrow 0$\,; $\widetilde Y\left(\mathring{\zeta}^l_j(t)\right)\equiv \overrightarrow 0$\,;
$\overrightarrow B^{k}\left(\mathring{\zeta}^l_j(t)\right)\equiv \overrightarrow 0$ for $1\leq k\leq r$ and $k\neq j$;
\begin{equation}
\begin{split}
\overrightarrow B^{j}\left(\mathring{\zeta}^l_j(t)\right)&=\left(a_{(r+1-j)(s-p+r+1-j)},\xi^{(j)}_{(r+1-j)1},\xi^{(j)}_{(r+1-j)2},\cdots,\xi^{(j)}_{(r+1-j)(s-p+r-j)},\right.\,\,\,\\
&\,\,\,\,\,\,\,\,\,\,\,\,\,\,\,\,\,\,\,\,\,\,\,\,\,\,\,\,\,\,\left.\xi^{(j)}_{1(s-p+r+1-j)},\xi^{(j)}_{2(s-p+r+1-j)},\cdots,\xi^{(j)}_{(r-j)(s-p+r+1-j)}\right)\left(\mathring{\zeta}^l_j(t)\right)\\ 
&=\left(t,0,0\cdots,0\right)\,.
\end{split}
\end{equation}

Let $\zeta^l_j$ be the closure of  $\mathring{\zeta}^l_j$ for $2\leq j\leq r-l$ or $r-l+2\leq j\leq r$.

\begin{lemma}\label{i5} 
For $2\leq j\leq r-l$, 
\begin{equation}\label{zei5}
\left\{
\begin{array}{ll}
   \zeta^l_j\cdot R^*_{s,p,n}((\mathcal O_{G(p,n)}(1))=0\,\,\,\, & \\[1em]
  \zeta^l_j\cdot D^-_{i}=\left\{ \begin{array}{cl}
       -1\,\,\,\,&   i=l+j-1\\
       2\,\,\,\,& i=l+j\,\,\\
       -1\,\,\,\,&   i=l+j+1\,\,{\rm and\,\,}1\leq i\leq r \\
       0 \,\,\,\,   & {\rm otherwise\,\,} \\
  \end{array}\right.  &  \\[2.3em]
  \zeta^l_j\cdot D^+_{i}=0  & \hspace{-2in} 1\leq i\leq r \\
\end{array}\right..
\end{equation}
For $r-l+2\leq j\leq r$,
\begin{equation}\label{zer5}
\left\{
\begin{array}{ll}
   \zeta^l_j\cdot R^*_{s,p,n}((\mathcal O_{G(p,n)}(1))=0\,\,\,\, & \\[1em]
   \zeta^l_j\cdot D^-_{i}=0  & \hspace{-1.88in} 1\leq i\leq r \\[1em]
  \zeta^l_j\cdot D^+_{i}=\left\{ \begin{array}{cl}
       -1\,\,\,\,&   i=j-1\\
       2\,\,\,\,& i=j\,\,\\
       -1\,\,\,\,&   i=j+1\,\,{\rm and\,\,}1\leq i\leq r \\
       0 \,\,\,\,   & {\rm otherwise\,\,}  \\
  \end{array}\right.  &  \\
\end{array}\right..
\end{equation}
\end{lemma}

\begin{lemma}\label{vki5} 
For $1\leq l\leq r-1$, 
\begin{equation}
\small
-K_{\mathcal T_{s,p,n}}\cdot\zeta^l_j=\left\{
\begin{array}{ll}
   2&\,\,{\rm when}\,\,2\leq j\leq r-l-1\\[1em]
   3&\,\,{\rm when}\,\,2\leq j=r-l\\[1em]
   2&\,\,{\rm when}\,\,r-l+2\leq j\leq r-1\\[1em]
   3&\,\,{\rm when}\,\,r-l+2\leq j=r\\
\end{array}\right..
\end{equation}
\end{lemma}

Assume that one of the following holds for integers $k,u,v$.
\begin{equation}\label{l1kuv}
\begin{split}
&(a).\,\,\,\,1\leq k\leq r-l\,; \,u=l+k\,\,{\rm and\,\,}s+l+k+1\leq v\leq n,\,\,{\rm or\,\,}v=s+l+k\,\,{\rm and}\\
&\,\,\,\,\,\,\,\,\,\,\,\,\,\,\,\,\,\,    l+k+1\leq u\leq p.\\
&(b).\,\,\,\,r-l+1\leq k\leq r\,; u=r-k+1\,\,{\rm and\,\,}1\leq v\leq s-p+r-k,\,\,{\rm or\,\,}\\
&\,\,\,\,\,\,\,\,\,\,\,\,\,\,\,\,\,\,  v=s-p+r-k+1\,\,{\rm and\,\,}1\leq u\leq r-k.
    \end{split}
\end{equation}
Then we can define an affine line  $\mathring{\zeta}_{u,v}^{l,k}\subset A^{\tau_l}$ in terms of the local coordinates as follows. $\xi^{(k)}_{uv}\left(\mathring{\zeta}_{u,v}^{l,k}(t)\right)=t$ and the other variables are constantly zero. 

Let $\zeta_{u,v}^{l,k}$ be the closure of $\mathring{\zeta}_{u,v}^{l,k}$.

\begin{lemma}\label{li1} 
Assume that $1\leq k\leq r-l$. When $u=l+k$ and $v=s+l+k+1\leq n$, or $v=s+l+k$ and $u=l+k+1\leq p$,
\begin{equation}
\left\{
\begin{array}{ll}
   \zeta^{l,k}_{u,v}\cdot R^*_{s,p,n}((\mathcal O_{G(p,n)}(1))=0\,\,\,\, & \\[1em]
  \zeta^{l,k}_{u,v}\cdot D^-_{i}=\left\{ \begin{array}{cl}
       -1\,\,\,\,&   i=l+k\\
       2\,\,\,\,& i=l+k+1\,\,{\rm and\,\,}1\leq i\leq r \\
       -1\,\,\,\,&   i=l+k+2\,\,{\rm and\,\,}1\leq i\leq r \\
       0 \,\,\,\,   & {\rm otherwise\,\,} \\
  \end{array}\right.  &  \\[2.3em]
  \zeta^{l,k}_{u,v}\cdot D^+_{i}=0  & \hspace{-2in} 1\leq i\leq r \\
\end{array}\right..
\end{equation}
When $u=l+k$ and  $s+l+k+2\leq v\leq n$, 
\begin{equation}
\left\{
\begin{array}{ll}
   \zeta^{l,k}_{u,v}\cdot R^*_{s,p,n}((\mathcal O_{G(p,n)}(1))=0\,\,\,\, & \\[1em]
  \zeta^{l,k}_{u,v}\cdot D^-_{i}=\left\{ \begin{array}{cl}
       -1\,\,\,\,&   i=l+k\\
       1\,\,\,\,& i=l+k+1\,\,{\rm and\,\,}1\leq i\leq r \\
       1\,\,\,\,&   i=v-s\,\,\,\,{\rm and\,\,\,\,}1\leq i\leq r \\
       -1\,\,\,\,&   i=v-s+1\,\,{\rm and\,\,}1\leq i\leq r \\
       0 \,\,\,\,   & {\rm otherwise\,\,}  \\
  \end{array}\right.  &  \\[3em]
  \zeta^{l,k}_{u,v}\cdot D^+_{i}=0  & \hspace{-2in} 1\leq i\leq r \\
\end{array}\right..
\end{equation}
When $v=s+l+k$ and  $l+k+2\leq u\leq p$, 
\begin{equation}
\left\{
\begin{array}{ll}
   \zeta^{l,k}_{u,v}\cdot R^*_{s,p,n}((\mathcal O_{G(p,n)}(1))=0\,\,\,\, & \\[1em]
  \zeta^{l,k}_{u,v}\cdot D^-_{i}=\left\{ \begin{array}{cl}
       -1\,\,\,\,&   i=l+k\\
       1\,\,\,\,& i=l+k+1\,\,{\rm and\,\,}1\leq i\leq r \\
       1\,\,\,\,&   i=u\,\,\,\,{\rm and\,\,\,\,}1\leq i\leq r \\
       -1\,\,\,\,&   i=u+1\,\,\,\,{\rm and\,\,\,\,}1\leq i\leq r \\
       0 \,\,\,\,   &{\rm otherwise\,\,}  \\
  \end{array}\right.  &  \\[3em]
  \zeta^{l,k}_{u,v}\cdot D^+_{i}=0  & \hspace{-2in} 1\leq i\leq r \\
\end{array}\right..
\end{equation}
\end{lemma}

\begin{lemma}\label{li2} 
Assume that $r-l+1\leq k\leq r$. When $u=r-k+1$ and $v=1\leq s-p+r-k$, or $v=s-p+r-k+1$ and $u=r-k\geq 1$, 
\begin{equation}
\left\{
\begin{array}{ll}
   \zeta^{l,k}_{u,v}\cdot R^*_{s,p,n}((\mathcal O_{G(p,n)}(1))=0\,\,\,\, & \\[1em]
   \zeta^{l,k}_{u,v}\cdot D^-_{i}=0  & \hspace{-2in} 1\leq i\leq r \\[1em]
  \zeta^{l,k}_{u,v}\cdot D^+_{i}=\left\{ \begin{array}{cl}
       -1\,\,\,\,&   i=k\\
       2\,\,\,\,& i=k+1\,\,{\rm and\,\,}1\leq i\leq r \\
       -1\,\,\,\,&   i=k+2\,\,{\rm and\,\,}1\leq i\leq r \\
       0 \,\,\,\,   &{\rm otherwise\,\,}  \\
  \end{array}\right.  &  \\
\end{array}\right..
\end{equation}
When $u=r-k+1$ and  $1\leq v\leq s-p+r-k-1$, 
\begin{equation}
\left\{
\begin{array}{ll}
   \zeta^{l,k}_{u,v}\cdot R^*_{s,p,n}((\mathcal O_{G(p,n)}(1))=0\,\,\,\, & \\[1em]
    \zeta^{l,k}_{u,v}\cdot D^-_{i}=0  & \hspace{-2.4in} 1\leq i\leq r \\[1em]
  \zeta^{l,k}_{u,v}\cdot D^+_{i}=\left\{ \begin{array}{cl}
       -1\,\,\,\,&   i=k\\
       1\,\,\,\,& i=k+1\,\,{\rm and\,\,}1\leq i\leq r \\
       1\,\,\,\,&   i=s-p+r+1-v\,\,{\rm and\,\,}1\leq i\leq r \\
       -1\,\,\,\,&   i=s-p+r+2-v\,\,{\rm and\,\,}1\leq i\leq r \\
       0 \,\,\,\,   & {\rm otherwise\,\,} \\
  \end{array}\right.  &  \\
\end{array}\right..
\end{equation}
When $v=s-p+r-k+1$ and  $1\leq u\leq r-k-1$, 
\begin{equation}
\left\{
\begin{array}{ll}
   \zeta^{l,k}_{u,v}\cdot R^*_{s,p,n}((\mathcal O_{G(p,n)}(1))=0\,\,\,\, & \\[1em]
  \zeta^{l,k}_{u,v}\cdot D^-_{i}=0  & \hspace{-2.in} 1\leq i\leq r \\[1em]
  \zeta^{l,k}_{u,v}\cdot D^+_{i}=\left\{ \begin{array}{cl}
       -1\,\,\,\,&   i=k\\
       1\,\,\,\,& i=k+1\,\,{\rm and\,\,}1\leq i\leq r \\
       1\,\,\,\,&   i=r+1-u\,\,{\rm and\,\,}1\leq i\leq r \\
       -1\,\,\,\,&   i=r+2-u\,\,{\rm and\,\,}1\leq i\leq r \\
       0 \,\,\,\,   & {\rm otherwise\,\,}  \\
  \end{array}\right.  &  \\
\end{array}\right..
\end{equation}
\end{lemma}

\begin{lemma}\label{vkli1}
Assume that $1\leq k\leq r-l$. Then, 
\begin{equation}
\small
-K_{\mathcal T_{s,p,n}}\cdot\zeta^{l,k}_{l+k,v}=\left\{
\begin{array}{ll}
   2(v-s-l-k)&\,\,{\rm when}\,\, s+l+k+1\leq v\leq r+s-1\\[1em]
   2(r-l-k)+1&\,\,{\rm when}\,\, v=r+s\,\,{\rm and}\,\, k\leq r-l-1\\[1em]
   n-s+p-2(l+k)+1&\,\,{\rm when}\,\, r+s+1\leq v\leq n\,\,{\rm and}\,\, k\leq r-l-1\\[1em]
   n-s+p-2r&\,\,{\rm when}\,\, r+s+1\leq v\leq n\,\,{\rm and}\,\,k=r-l\\
\end{array}\right.
\end{equation}
and
\begin{equation}
\small
-K_{\mathcal T_{s,p,n}}\cdot\zeta^{l,k}_{u,s+l+k}=\left\{
\begin{array}{ll}
   2(u-l-k)&\,\,{\rm when}\,\, l+k+1\leq u\leq r-1\\[1em]
   2(r-l-k)+1&\,\,{\rm when}\,\, u=r\,\,{\rm and}\,\,k \leq r-l-1\\[1em]
   n-s+p-2(l+k)+1&\,\,{\rm when}\,\, r+1\leq u\leq p\,\,{\rm and}\,\,k\leq r-l-1\\[1em]
   n-s+p-2r&\,\,{\rm when}\,\, r+1\leq u\leq p\,\,{\rm and}\,\,k=r-l\\
\end{array}\right..
\end{equation}
\end{lemma}

\begin{lemma}\label{vkli2}
Assume that $r-l+1\leq k\leq r$. Then,
\begin{equation}
\small
-K_{\mathcal T_{s,p,n}}\cdot\zeta^{l,k}_{r-k+1,v}=\left\{
\begin{array}{ll}
   2(s-p+r-k+1-v)&\,\,{\rm when}\,\, s-p+2\leq v\leq s-p+r-k\,\,\\[1em]
   2(r-k)+1&\,\,{\rm when}\,\, v=s-p+1\,\,{\rm and}\,\, k\leq r-1\\[1em]
   2(r-k)+s-p+1&\,\,{\rm when}\,\, 1\leq v\leq s-p\,\,{\rm and}\,\, k\leq r-1\\[1em]
   s-p&\,\,{\rm when}\,\, 1\leq v\leq s-p\,\,{\rm and}\,\,k=r\\
\end{array}\right.
\end{equation}
and
\begin{equation}
\small
-K_{\mathcal T_{s,p,n}}\cdot\zeta^{l,k}_{u,s-k+r-p+1}=\left\{
\begin{array}{ll}
   2(r-k+1-u)&\,\,{\rm when}\,\, 2\leq u\leq r-k\\[1em]
   2(r-k)+1&\,\,{\rm when}\,\, u=1\,\,{\rm and}\,\, k\leq r-1\\
\end{array}\right..
\end{equation}
\end{lemma}

\bigskip

For $1\leq m_1\leq p-l$ and $1\leq m_2\leq s-p+l$,  define an affine line 
$\mathring{\delta}_{m_1,m_2}^l\subset A^{\tau_l}$  as follows. $y_{(l+m_1)(s-p+l+1-m_2)}\left(\mathring{\delta}_{m_1,m_2}^l(t)\right)=t$ and the other variables are constantly zero.  

Let $\delta_{m_1,m_2}^l$ be the closure of 
$\mathring{\delta}_{m_1,m_2}^l$.

\begin{lemma}\label{li4}
For $1\leq m_1\leq p-l$ and $1\leq m_2\leq s-p+l$, 
\begin{equation}
\left\{
\begin{array}{ll}
   \delta_{m_1,m_2}^l\cdot R^*_{s,p,n}((\mathcal O_{G(p,n)}(1))=1\,\,\,\, & \\[1em]
  \delta_{m_1,m_2}^l\cdot D^-_{i}=\left\{ \begin{array}{cl}
       1\,\,\,\,&   i=l+m_1\,\,{\rm and}\,\,1\leq i\leq r\\
       -1\,\,\,\,&   i=l+m_1+1\,\,{\rm and}\,\,1\leq i\leq r\\
       0 \,\,\,\,   &{\rm otherwise\,\,}  \\
  \end{array}\right.  &  \\[2.3em]
  \delta_{m_1,m_2}^l\cdot D^+_{i}=\left\{ \begin{array}{cl}
       1\,\,\,\,&   i=r-l+m_2\,\,{\rm and}\,\,1\leq i\leq r\\
       -1\,\,\,\,&   i=r-l+m_2+1\,\,{\rm and}\,\,1\leq i\leq r\\
       0 \,\,\,\,   &{\rm otherwise\,\,} \\
  \end{array}\right.  &  \\
\end{array}\right..
\end{equation}
\end{lemma}

\begin{lemma}\label{vkli4}
For $1\leq m_1\leq p-l$ and $1\leq m_2\leq s-p+l$, 
\begin{equation}
\footnotesize
-K_{\mathcal T_{s,p,n}}\cdot\delta_{m_1,m_2}^l=\left\{
\begin{array}{ll}
2m_1+2m_2-2&\,\,{\rm when}\,\,1\leq m_1\leq r-l-1\,\,{\rm and}\,\,1\leq m_2\leq l-1\\[1em]
2m_1+2l-1&\,\,{\rm when}\,\,1\leq m_1\leq r-l-1\,\,{\rm and}\,\, m_2=l\\[1em]
2m_1+2l-1+s-p&\,\,{\rm when}\,\,1\leq m_1\leq r-l-1\,\,{\rm and}\,\,m_2\geq l+1\\[1em]
2(r-l+m_2)-1&\,\,{\rm when}\,\,m_1=r-l\,\,{\rm and}\,\,1\leq m_2\leq l-1\\[1em]
2r&\,\,{\rm when}\,\,m_1=r-l\,\,{\rm and}\,\, m_2=l\\[1em]
2r+s-p&\,\,{\rm when}\,\,m_1=r-l\,\,{\rm and}\,\,m_2\geq l+1\\[1em]
2(r-l+m_2)-1+s+p-n&\,\,{\rm when}\,\,m_1\geq r-l+1\,\,{\rm and}\,\,1\leq m_2\leq l-1\\[1em]
n-s+p&\,\,{\rm when}\,\,m_1\geq r-l+1\,\,{\rm and}\,\, m_2=l\\[1em]
n&\,\,{\rm when}\,\,m_1\geq r-l+1\,\,{\rm and}\,\,m_2\geq l+1\\
\end{array}\right..
\end{equation}
\end{lemma}

\bigskip

For $1\leq m_1\leq n-s-l$ and $1\leq m_2\leq l$,  define  an affine line $\mathring{\Delta}_{m_1,m_2}^l\subset A^{\tau_l}$ as follows. $x_{(l+1-m_2)(s+l+m_1)}\left(\mathring{\Delta}_{m_1,m_2}^l(t)\right)=t$ and the other variables are constantly zero.

Let  $\Delta_{m_1,m_2}^l$  be the closure of $\mathring{\Delta}_{m_1,m_2}^l$.

\begin{lemma}\label{li5}
For $1\leq m_1\leq n-s-l$ and $1\leq m_2\leq l$, 
\begin{equation}
\left\{
\begin{array}{ll}
   \Delta_{m_1,m_2}^l\cdot R^*_{s,p,n}((\mathcal O_{G(p,n)}(1))=1\,\,\,\, & \\[1em]
  \Delta_{m_1,m_2}^l\cdot D^-_{i}=\left\{ \begin{array}{cl}
       1\,\,\,\,&   i=l+m_1\\
       -1\,\,\,\,&   i=l+m_1+1\,\,{\rm and}\,\,1\leq i\leq r\\
       0 \,\,\,\,   &{\rm otherwise\,\,}  \\
  \end{array}\right.  &  \\[2.3em]
  \Delta_{m_1,m_2}^l\cdot D^+_{i}=\left\{ \begin{array}{cl}
       1\,\,\,\,&   i=r-l+m_2\\
       -1\,\,\,\,&   i=r-l+m_2+1\,\,{\rm and}\,\,1\leq i\leq r\\
       0 \,\,\,\,   &{\rm otherwise\,\,} \\
  \end{array}\right.  &  \\
\end{array}\right..
\end{equation}
\end{lemma}

\begin{lemma}\label{vkli5}
For $1\leq m_1\leq n-s-l$ and $1\leq m_2\leq l$, 
\begin{equation}
\footnotesize
-K_{\mathcal T_{s,p,n}}\cdot\Delta_{m_1,m_2}^l=\left\{
\begin{array}{ll}
2m_1+2m_2-2&\,\,{\rm when}\,\,1\leq m_1\leq r-l-1\,\,{\rm and}\,\,1\leq m_2\leq l-1\\[1em]
2m_1+2l-1&\,\,{\rm when}\,\,1\leq m_1\leq r-l-1\,\,{\rm and}\,\, m_2=l\\[1em]
2(r-l+m_2)-1&\,\,{\rm when}\,\,m_1=r-l\,\,{\rm and}\,\,1\leq m_2\leq l-1\\[1em]
2r&\,\,{\rm when}\,\,m_1=r-l\,\,{\rm and}\,\, m_2=l\\[1em]
2(r-l+m_2)-1+n-s-p&\,\,{\rm when}\,\,m_1\geq r-l+1\,\,{\rm and}\,\,1\leq m_2\leq l-1\\[1em]
n-s+p&\,\,{\rm when}\,\,m_1\geq r-l+1\,\,{\rm and}\,\, m_2=l\\
\end{array}\right..
\end{equation}
\end{lemma}

{\bf\noindent (3). Curves  $\zeta^r_j$, $\zeta_{u,v}^{r,k}$, $\delta_{m_1,m_2}^r$, $\Delta_{m_1,m_2}^r$ in $R_{s,p,n}^{-1}\left(\mathcal V_{(p-r,r)}\right)$.}
\smallskip

For $2\leq j\leq r$, define an affine line  $\mathring{\zeta}^r_j\subset A^{\tau_r}$  as follows.  $a_{(r+1-j)(s-p+r+1-j)}\left(\mathring{\zeta}^r_j(t)\right)=t$ and the other variables are constantly zero.
 
Let $\zeta^r_j$ be the closure of    $\mathring{\zeta}^r_j$.

\begin{lemma}\label{ri1} 
For $2\leq j\leq r$,
\begin{equation}
\left\{
\begin{array}{ll}
   \zeta^r_j\cdot R^*_{s,p,n}((\mathcal O_{G(p,n)}(1))=0\,\,\,\, & \\[1em]
  \zeta^r_j\cdot D^-_{i}=0  & \hspace{-2.4in} 1\leq i\leq r \\ [1em]
  \zeta^r_j\cdot D^+_{i}=\left\{ \begin{array}{cl}
       -1\,\,\,\,&   i=j-1\\
       2\,\,\,\,& i=j\,\,{\rm and\,\,}1\leq i\leq r \\
       -1\,\,\,\,&   i=j+1\,\,{\rm and\,\,}1\leq i\leq r \\
       0 \,\,\,\,   & {\rm otherwise\,\,}  \\
  \end{array}\right.  &  \\
\end{array}\right..
\end{equation}
\end{lemma}

\begin{lemma}\label{rki2} 
For $2\leq j\leq r$,
\begin{equation}
\small
-K_{\mathcal T_{s,p,n}}\cdot\zeta^r_j=\left\{
\begin{array}{ll}
   3&\,\,\,{\rm when}\,\,\, 2\leq j=r\\[1em]
   2&\,\,\,{\rm when}\,\,\, 2\leq j\leq r-1\\
\end{array}\right..
\end{equation}
\end{lemma}

Assume that $u=r-k+1$ and $1\leq v\leq s-p+r-k$, or $v=s-p+r-k+1$ and $1\leq u\leq r-k$, where $1\leq k\leq r$. Define an affine line $\mathring{\zeta}_{u,v}^{r,k}\subset  A^{\tau_r}$ as follows.  $\xi^{(k)}_{uv}\left(\mathring{\zeta}_{u,v}^{r,k}(t)\right)=t$ and other variables are constantly zero.

Let $\zeta_{u,v}^{r,k}$ be the closure of $\mathring{\zeta}_{u,v}^{r,k}$.

\begin{lemma}\label{ri2} 
Assume that $1\leq k\leq r$. When $u=r-k+1$ and $v=s-p+r-k\geq 1$, or $v=s-p+r-k+1$ and $u=r-k\geq 1$, \begin{equation}
\left\{
\begin{array}{ll}
   \zeta^{r,k}_{u,v}\cdot R^*_{s,p,n}((\mathcal O_{G(p,n)}(1))=0\,\,\,\, & \\[1em]
   \zeta^{r,k}_{u,v}\cdot D^-_{i}=0  & \hspace{-2in} 1\leq i\leq r \\[1em]
  \zeta^{r,k}_{u,v}\cdot D^+_{i}=\left\{ \begin{array}{cl}
       -1\,\,\,\,&   i=k\\
       2\,\,\,\,& i=k+1\,\,{\rm and\,\,}1\leq i\leq r \\
       -1\,\,\,\,&   i=k+2\,\,{\rm and\,\,}1\leq i\leq r \\
       0 \,\,\,\,   & {\rm otherwise\,\,}  \\
  \end{array}\right.  &  \\
\end{array}\right..
\end{equation}
When $u=r-k+1$ and  $1\leq v\leq s-p+r-k-1$, 
\begin{equation}
\left\{
\begin{array}{ll}
   \zeta^{r,k}_{u,v}\cdot R^*_{s,p,n}((\mathcal O_{G(p,n)}(1))=0\,\,\,\, & \\[1em]
    \zeta^{r,k}_{u,v}\cdot D^-_{i}=0  & \hspace{-2.4in} 1\leq i\leq r \\[1em]
  \zeta^{r,k}_{u,v}\cdot D^+_{i}=\left\{ \begin{array}{cl}
       -1\,\,\,\,&   i=k\\
       1\,\,\,\,& i=k+1\,\,{\rm and\,\,}1\leq i\leq r \\
       1\,\,\,\,&   i=s-p+r+1-v\,\,{\rm and\,\,}1\leq i\leq r \\
       -1\,\,\,\,&   i=s-p+r+2-v\,\,{\rm and\,\,}1\leq i\leq r \\
       0 \,\,\,\,   & {\rm otherwise\,\,} \\
  \end{array}\right.  &  \\
\end{array}\right..
\end{equation}
When $v=s-p+r-k+1$ and  $1\leq u\leq r-k-1$, 
\begin{equation}
\left\{
\begin{array}{ll}
   \zeta^{r,k}_{u,v}\cdot R^*_{s,p,n}((\mathcal O_{G(p,n)}(1))=0\,\,\,\, & \\[1em]
  \zeta^{r,k}_{u,v}\cdot D^-_{i}=0  & \hspace{-2in} 1\leq i\leq r \\[1em]
  \zeta^{r,k}_{u,v}\cdot D^+_{i}=\left\{ \begin{array}{cl}
       -1\,\,\,\,&   i=k\\
       1\,\,\,\,& i=k+1\,\,{\rm and\,\,}1\leq i\leq r \\
       1\,\,\,\,&   i=r+1-u\,\,{\rm and\,\,}1\leq i\leq r \\
       -1\,\,\,\,&   i=r+2-u\,\,{\rm and\,\,}1\leq i\leq r \\
       0 \,\,\,\,   &{\rm otherwise\,\,}  \\
  \end{array}\right.  &  \\
\end{array}\right..
\end{equation}
\end{lemma}

\begin{lemma}\label{rri2} 
Assume that $1\leq k\leq r$. Then,
\begin{equation}
\small
-K_{\mathcal T_{s,p,n}}\cdot\zeta^{r,k}_{r-k+1,v}=\left\{
\begin{array}{ll}
   2(s-p+r-k+1-v)&\,\,{\rm when}\,\, s-p+2\leq v\leq s-p+r-k\\[1em]
   2r-2k+1&\,\,{\rm when}\,\, v=s-p+1\,\,{\rm and}\,\,1\leq k\leq r-1\\[1em]
   2r-2k+s-p+1&\,\,{\rm when}\,\, 1\leq v\leq s-p\,\,{\rm and}\,\,1\leq k\leq r-1\\[1em]
   s-p&\,\,{\rm when}\,\, 1\leq v\leq s-p\,\,{\rm and}\,\,k=r\\
\end{array}\right.
\end{equation}
and
\begin{equation}
\small
-K_{\mathcal T_{s,p,n}}\cdot\zeta^{r,k}_{u,s-p+r-k+1}=\left\{
\begin{array}{ll}
   2(r-k+1-u)&\,\,{\rm when}\,\, 2\leq u\leq r-k\\[1em]
   2r-2k+1&\,\,{\rm when}\,\, u=1\,\,{\rm and}\,\,1\leq k\leq r-1\\
\end{array}\right..
\end{equation}
\end{lemma}

When $n-s<p$,   $1\leq m_1\leq s+p-n$, and $1\leq m_2\leq n-p$,  define an affine line $\mathring{\delta}_{m_1,m_2}^r\subset A^{\tau_{r}}$  as follows. $y_{(r+m_1)(s-p+r+1-m_2)}\left(\mathring{\delta}_{m_1,m_2}^r(t)\right)=t$ and other variables are constantly zero.

Let $\delta_{m_1,m_2}^r$ be the closure of  $\mathring{\delta}_{m_1,m_2}^r$.

\begin{lemma}\label{ri3}
Assume that $n-s<p$,   $1\leq m_1\leq s+p-n$, and $1\leq m_2\leq n-p$. Then,
\begin{equation}
\left\{
\begin{array}{ll}
   \delta_{m_1,m_2}^r\cdot R^*_{s,p,n}((\mathcal O_{G(p,n)}(1))=1\,\,\,\, & \\[1em]
   \delta_{m_1,m_2}^r\cdot D^-_{i}=0  & \hspace{-1.3in} 1\leq i\leq r \\[1em]
  \delta_{m_1,m_2}^r\cdot D^+_{i}=\left\{ \begin{array}{cl}
       1\,\,\,\,&   i=m_2\,\,{\rm and}\,\,1\leq i\leq r\\
       -1\,\,\,\,&   i=m_2+1\,\,{\rm and}\,\,1\leq i\leq r\\
       0 \,\,\,\,   &{\rm otherwise}\\
       \end{array}\right.
\end{array}\right..
\end{equation}
\end{lemma}

\begin{lemma}\label{rki3}
Assume that $n-s<p$. Then,
\begin{equation}
-K_{\mathcal T_{s,p,n}}\cdot {\delta}_{m_1,m_2}^r=\left\{
\begin{array}{ll}
   2m_2+s+p-n-1&\,\,{\rm when}\,\,1\leq m_2\leq r-1\\[1em]
   n-s+p&\,\,{\rm when}\,\,m_2=r\\[1em]
   n&\,\,{\rm when}\,\,r+1\leq m_2\leq n-p\\
\end{array}\right..
\end{equation}
\end{lemma}

When  $p<n-s$, $1\leq m_1\leq n-s-r$, and $1\leq m_2\leq r$,    define an affine line $\mathring{\Delta}_{m_1,m_2}^r\subset A^{\tau_l}$  as follows. $x_{(r+1-m_2)(s+r+m_1)}\left(\mathring{\Delta}_{m_1,m_2}^r(t)\right)=t$ and the other variables are constantly zero.

Let $\Delta_{m_1,m_2}^r$ be the closure of  $\mathring{\Delta}_{m_1,m_2}^r$.

\begin{lemma}\label{ri4}
Assume that $p<n-s$, $1\leq m_1\leq n-s-r$ and $1\leq m_2\leq r$. Then,
\begin{equation}
\left\{
\begin{array}{ll}
   {\Delta}_{m_1,m_2}^r\cdot R^*_{s,p,n}((\mathcal O_{G(p,n)}(1))=1\,\,\,\, & \\[1em]
   {\Delta}_{m_1,m_2}^r\cdot D^-_{i}=0  & \hspace{-1.1in} 1\leq i\leq r \\[1em]
 {\Delta}_{m_1,m_2}^r\cdot D^+_{i}=\left\{ \begin{array}{cl}
       1\,\,\,\,&   i=m_2\,\,{\rm and}\,\,1\leq i\leq r\\
       -1\,\,\,\,&   i=m_2+1\,\,{\rm and}\,\,1\leq i\leq r\\
       0 \,\,\,\,   &{\rm otherwise\,\,}  \\
       \end{array}\right.
\end{array}\right..
\end{equation}
\end{lemma}

\begin{lemma}\label{rki4}
Assume that $p<n-s$. Then,
\begin{equation}
-K_{\mathcal T_{s,p,n}}\cdot {\Delta}_{m_1,m_2}^r=\left\{
\begin{array}{ll}
   2m_2+n-s-p-1&\,\,{\rm when}\,\,1\leq m_2\leq r-1\\[1em]
   n-s+p&\,\,{\rm when}\,\,m_2=r\\
\end{array}\right..
\end{equation}
\end{lemma}

\subsection{(Semi-)positivity of the anti-canonical bundles} \label{curanddivp}

In this subsection, we will study the (semi)-positivity of the anti-canonical bundles of $\mathcal T_{s,p,n}$  and $\mathcal M_{s,p,n}$.

Recall the following notion of Chow groups. Let $X$ be a smooth projective manifold over $\mathbb C$. An algebraic cycle on $X$ means a finite linear combination of closed subvarieties of $X$ with integer coefficients. For a natural number $i$, the group $Z_{i}(X)$ of $i$-dimensional cycles (or $i$-cycles, for short) on $X$ is the free abelian group on the set of $i$-dimensional subvarieties of $X$. For a variety $W$ of dimension $i+1$ and any rational function $f$ on $W$ which is not identically zero, the divisor of $f$ is the $i$-cycle
\begin{equation}
    (f)=\sum _{Z}\operatorname {ord} _{Z}(f)Z
\end{equation}
where the sum runs over all $i$-dimensional subvarieties $Z$ of $W$ and the integer $ \operatorname {ord} _{Z}(f)$ denotes the order of vanishing of $f$ along $Z$. The group of $i$-cycles rationally equivalent to zero is the subgroup of $ Z_{i}(X)$ generated by the cycles $(f)$ for all $(i+1)$-dimensional subvarieties $W$ of $X$ and all nonzero rational functions $f$ on $W$. The Chow group $A_{i}(X)$ of $i$-dimensional cycles on $X$ is the quotient group of $Z_{i}(X)$ by the subgroup of cycles rationally equivalent to zero.

Brion characterized the Chow groups of spherical varieties  as follows.
\begin{lemma}[\cite{Br3}]\label{cone}
Let $X$ be an irreducible, complete spherical variety of complex dimension $n$.
The cone of effective cycles in $A_{i}(X)\otimes_{\mathbb Z}\mathbb Q$ is a polyhedral convex cone generated by the classes of the closures of the $B$-orbits.
\end{lemma}

{\bf\noindent Proof of Theorem \ref{fano}.} Recall that by Kleiman's criterion $-K_{\mathcal T_{s,p,n}}$ is ample if and only if its has positive degree on every nonzero element of the closure of the cone of effective curves in $A_{1}(X)\otimes_{\mathbb Z}\mathbb Q$. By Lemma \ref{cone}, the closure of the cone of effective curves in $A_{1}(X)\otimes_{\mathbb Z}\mathbb Q$ is generated by finitely many irreducible curves which are the closures of the $B$-orbits.  Therefore, to prove Theorem \ref{fano}, it suffices to show that the following holds. 
\begin{enumerate}[label=(\alph*)]
    \item $-K_{\mathcal T_{s,p,n}}$ has non-negative degrees on all the irreducible curves  of $\mathcal T_{s,p,n}$ which are the closures of the $B$-orbits of the dimension $1$. Moreover, it has strictly positive degrees on such curves if and only if $r\leq 2$.
    \item $-K_{\mathcal T_{s,p,n}}$ is big.
\end{enumerate}

We assume that $2p\leq n\leq 2s$, and split the proof into two steps. In the first step, we show that the anti-canonical bundle $-K_{\mathcal T_{s,p,n}}$ is nef (or ample if $r\leq 2$); in the second step, we show that  $-K_{\mathcal T_{s,p,n}}$ is big. 
\smallskip

{\bf\noindent Step I (nef part).} Recall that $\mathcal T_{s,p,n}$ is covered by the images of the open sets $A^{\tau_0}$, $A^{\tau_1},\cdots$, $A^{\tau_r}$ under  the actions of the permutation matrices contained in $GL(s,\mathbb C)\times GL(n-s,\mathbb C)$. Here $\left(A^{\tau_l},\left(J_{l}^{\tau_l}\right)^{-1}\right)$ is the $l$-th main coordinate chart defined in Definition \ref{taul}; a permutation matrix is a square  matrix that has exactly one entry of $1$ in each row and each column and $0$s elsewhere.

Let $\gamma$ be an irreducible curve of $\mathcal T_{s,p,n}$ which is the closure of a $1$-dimensional $B$-orbit. Then it is easy to verify that there is a permutation matrix $g\in GL(s,\mathbb C)\times GL(n-s,\mathbb C)$ and an integer  $0\leq l^*\leq r$ such that $A^{\tau_l^*}$ contains an open subset of the image  curve $g(\gamma)$.  

Let $T$ be the maximal torus of $B$ defined by (\ref{maxt}). Let $\mathfrak a$ be a generic point of $g(\gamma)$.  It is clear that the image curve $g(\gamma)$ is invariant under the $T$-action. Consider the local coordinates of $\mathfrak a$ with respect to the $i^*$-th main coordinate charts   $\left(A^{\tau_l^*},\left(J^{\tau_l^*}\right)^{-1}\right)$. We can show that only one of the local coordinates of $\mathfrak a$ is nonzero, for otherwise $g(\gamma)$ is of dimension at least two. Therefore, $g(\gamma)$ coincides with one of the $T$-invariant curves $\gamma_l,\zeta_j^l,\zeta_{u,v}^{l,k},\delta_{m_1,m_2}^l,  \Delta_{m_1,m_2}^l$ defined in Section \ref{curanddivi}. 

Since $K_{\mathcal T_{s,p,n}}$ is invariant under the holomorphic automorphism of $\mathcal T_{s,p,n}$, we have
\begin{equation}
    -K_{\mathcal T_{s,p,n}}\cdot\gamma=-K_{\mathcal T_{s,p,n}}\cdot g(\gamma).
\end{equation}
Recall Lemmas \ref{vki1},    \ref{vki2},  \ref{vki3}, \ref{vki4}, \ref{vki5}, \ref{vkli1}, \ref{vkli2}  \ref{vkli4},  \ref{vkli5},  \ref{rki2}, \ref{rri2}, \ref{rki3}, \ref{rki4}. Then, we can conclude that $-K_{\mathcal T_{s,p,n}}\cdot\gamma\geq 0$, and  $-K_{\mathcal T_{s,p,n}}\cdot\gamma=0$ if and only if $g(\gamma)$ coincides with $\gamma_l$ for a certain $1\leq l\leq r-2$.

Therefore, we can conclude that $-K_{\mathcal T_{s,p,n}}$ is a numerical effective line bundle of $\mathcal T_{s,p,n}$, and $-K_{\mathcal T_{s,p,n}}$ is ample if and only if $r\leq 2$.
\smallskip

{\bf\noindent Step II (big part).} We assume that $r\geq 3$. By Kodaira's lemma, $-K_{\mathcal T_{s,p,n}}$ is a big line bundle if it has a decomposition $A+E$ (as $\mathbb Q$-divisors), with $A$ ample and $E$ effective. 

First assume that  $p<n-s\leq s$. By (\ref{bst}) and Lemma \ref{wk}, we have 
\begin{equation}
\begin{split}
-K_{\mathcal T_{s,p,n}}&=(s-p+1)\cdot B_0 +2\sum_{j=1}^{p-1}B_j +(n-s-p+1)\cdot B_p +\sum_{i=1}^pD^-_i+\sum_{i=1}^pD^+_i\\
&\hspace{-.56in}=\left(s-p+1-\frac{1}{p}\right)\cdot B_0 +2\sum_{j=1}^{p-1}B_j +\left(n-s-p+1\right)\cdot B_p+ \sum_{i=1}^pD^-_i\\
&+\sum_{i=1}^p\left(1-\frac{p+1-i}{p}\right)\cdot D^+_i+\frac{1}{p}\left(B_0+\sum_{i=1}^p\left({p+1-i}\right)\cdot D^+_i\right)\\
\end{split}
\end{equation}
\begin{equation*}
\begin{split}
&=\left(s-p+1-\frac{1}{p}\right)\cdot B_0 +2\sum_{j=1}^{p-1}B_j +\left(n-s-p+1\right)\cdot B_p+\frac{1}{p}\cdot(R_{s,p,n})^*\left(\mathcal O_{G(p,n)}(1)\right)\\
&\,\,\,\,\,\,\,\,+\sum_{i=1}^pD^-_i+\sum_{i=1}^p\left(1-\frac{p+1-i}{p}\right)\cdot D^+_i\,.\\
\end{split}
\end{equation*}
By Lemma \ref{ample},  we conclude that $-K_{\mathcal T_{s,p,n}}$ is big.

Similarly, when $n-s=p<s$,
\begin{equation}
\begin{split}
    -K_{\mathcal T_{s,p,n}}&=\left(s-p+1-\frac{1}{p}\right)\cdot B_0 +2\sum_{j=1}^{p-1}B_j +\frac{1}{p}\cdot(R_{s,p,n})^*\left(\mathcal O_{G(p,n)}(1)\right)\\
&\,\,\,\,\,\,\,\,\,+\sum_{i=1}^pD^-_i+\sum_{i=1}^p\left(1-\frac{p+1-i}{p}\right)\cdot D^+_i\,;\\
\end{split}
\end{equation}
when $n-s<p<s$ ($r=n-s$),
\begin{equation}
\begin{split}
-K_{\mathcal T_{s,p,n}}&=\left(s-p+1-\frac{1}{r}\right)\cdot B_0+2\sum_{j=1}^{r-1}B_j+(p-r+1)\cdot B_r+\frac{1}{r}\cdot(R_{s,p,n})^*\left(\mathcal O_{G(p,n)}(1)\right)\\
&\,\,\,\,\,\,\,\,+\sum_{i=1}^rD^-_i+\sum_{i=1}^r\left(1-\frac{r+1-i}{r}\right)\cdot D^+_i\,;\\
\end{split}
\end{equation}
when $n-s=p=s$,
\begin{equation}
    -K_{\mathcal T_{s,p,n}}= 2\sum_{j=1}^{p-1}B_j  +\frac{1}{p}\cdot(R_{s,p,n})^*\left(\mathcal O_{G(p,n)}(1)\right)+\sum_{i=1}^pD^-_i+\sum_{i=1}^{p-1}\left(1-\frac{p+1-i}{p}\right)\cdot D^+_i+D^+_p\,.
\end{equation}
By Lemma \ref{ample},  we conclude that $-K_{\mathcal T_{s,p,n}}$ is big as well.

We complete the proof of Theorem \ref{fano}.\,\,\,\,\,$\endpf$

\medskip

{\bf\noindent Proof of Theorem \ref{mfano}.} The proof is the same as in Theorem \ref{fano}.
 
Recall that the effective cone of curves is generated by  $B$-invariant curves by Lemma \ref{cone}. Let $\check\gamma\subset\mathcal M_{s,p,n}$ be an irreducible curve  which is the closure of a $1$-dimensional $B$-orbit. Then there is a permutation matrix $g\in GL(s,\mathbb C)\times GL(n-s,\mathbb C)$ such that $A^{\tau_0}$ contains an open subset of the image  curve $g(\check\gamma)$.  We can further show that, $g(\check\gamma)$ coincides with one of the $T$-invariant curves $\zeta_j^0,\zeta_{u,v}^{0,k},\delta_{m_1,m_2}^0$ defined in Section \ref{curanddivi}.
 
Notice that for each curve $\gamma\subset\mathcal M_{s,p,n}$,
\begin{equation}
    -K_{\mathcal M_{s,p,n}}\cdot \gamma=\left(-K_{\mathcal T_{s,p,n}}-D_1^-\right)\cdot\gamma
\end{equation}
Then Theorem \ref{mfano} follows from Lemmas    \ref{vki2},  \ref{vki3}, \ref{vki4} by Kleiman's criterion. \,\,\,\,\,$\endpf$
\begin{remark}\label{mcls}
Similarly, we can show that  the line bundle $-K_{\mathcal T_{s,p,n}}+H_j\big|_{\mathcal T_{s,p,n}}-D_1^-$ on $\mathcal T_{s,p,n}$ is nef and big for  $0\leq j\leq r$. Then  $
H^1\left(\mathcal T_{s,p,n},\,\,H_j\big|_{\mathcal T_{s,p,n}}-D_1^-\right)=0$ by Kawamata-Viehweg vanishing theorem; hence the following natural map is surjective.
\begin{equation}
\Gamma\left(\mathcal T_{s,p,n},\,\,H_j\big|_{\mathcal T_{s,p,n}}\right)\rightarrow \Gamma\left(\mathcal M_{s,p,n},\,\,H_j\big|_{\mathcal M_{s,p,n}}\right)\,.
\end{equation}
By Remark \ref{cls}, we can conclude that the complete linear series of $H_j\big|_{\mathcal M_{s,p,n}}$ on $\mathcal M_{s,p,n}$ is isomorphic to $\mathbb {CP}^{N^j_{s,p,n}}$ as well.
\end{remark}

\section{Symmetries of \texorpdfstring{$\mathcal T_{s,p,n}$ }{jj} and \texorpdfstring{$\mathcal M_{s,p,n}$ }{jj} } \label{sym}
In this section, we will determine the holomorphic automorphisms of  $\mathcal T_{s,p,n}$ and $\mathcal M_{s,p,n}$. The proof consists two  technical parts. Firstly, we show that the induced automorphisms of the Picard group consist of at most two elements; secondly, we complete the proof under the extra assumption that the induced automorphism of the Picard group is the identity map.  We note that the former  is easier for the case $\mathcal M_{s,p,n}$ and the latter is easier for the case $\mathcal T_{s,p,n}$.

We will assume in this section that $2p\leq n\leq 2s$.

\subsection{Discrete symmetries} \label{dsym}

We compute in this subsection the  isomorphisms of the Picard groups induced by the following discrete automorphisms.
\begin{itemize}
    \item The {\rm USD} isomorphisms of $\mathcal T_{s,p,2s}$; the {\rm DUAL} isomorphisms of $\mathcal T_{s,p,2p}$.
    \item The {\rm Usd} isomorphisms of $\mathcal M_{s,p,2s}$; the {\rm Dual} isomorphisms of $\mathcal M_{s,p,2p}$.    
\end{itemize}

\begin{lemma}\label{dis}
Assume that  $1\leq p\leq s$.  The  automorphism  {\rm USD} induces an automorphism of the Picard group of $\mathcal T_{s,p,2s}$ as follows.
\begin{equation}
\begin{split}
&({\rm USD})^*(D^+_i)=D^-_i,\,\,\, 1\leq i\leq r\,;\\
&({\rm USD})^*(D^-_i)=D^+_i,\,\,\, 1\leq i\leq r\,;\\
&({\rm USD})^*\left((R_{s,p,2s})^*(\mathcal O_{G(p,2s)}(1))\right)=(R_{s,p,2s})^*(\mathcal O_{G(p,2s)}(1))\,;\\
&({\rm USD})^*\left(B_i\right)=B_{r-i}\,,\,\,0\leq i\leq r.
\end{split}
\end{equation}
\end{lemma}
{\noindent\bf Proof of Lemma \ref{dis}.} We have $({\rm USD})^*\left((R_{s,p,2s})^*(\mathcal O_{G(p,2s)}(1))\right)=(R_{s,p,2s})^*(\mathcal O_{G(p,2s)}(1))$ for USD induces an automorphism of $G(p,2s)$. By (\ref{dtrans}) and the Van der Waerden representation, we can conclude that {\rm USD} interchanges $D^+_i$ and $D^-_i$ for $1\leq i\leq r$\,. \,\,\,$\endpf$

\begin{lemma}\label{dis2}Assume that  $1\leq p\leq s$. The {\rm DUAL} automorphism  induces an automorphism of the Picard group of 
$\mathcal T_{s,p,2p}$ as follows.
\begin{equation}
\begin{split}
&({\rm DUAL})^*(D^+_i)=D^-_i,\,\,\, 1\leq i\leq r\,;\\
&({\rm DUAL})^*(D^-_i)=D^+_i,\,\,\, 1\leq i\leq r\,;\\
&({\rm DUAL})^*\left((R_{s,p,2p})^*(\mathcal O_{G(p,2p)}(1))\right)=(R_{s,p,2p})^*(\mathcal O_{G(p,2p)}(1))\,;\\
&({\rm DUAL})^*\left(B_i\right)=B_{r-i}\,,\,\,0\leq i\leq r.
\end{split}
\end{equation}
\end{lemma}
{\noindent\bf Proof of Lemma \ref{dis2}.} The proof is the same as in Lemma \ref{dis}. \,\,\,$\endpf$
\medskip

We next turn to the corresponding actions on the Picard groups of $\mathcal M_{s,p,n}$.
\begin{lemma}\label{mdis}Assume that  $1\leq p<s$. Then the  automorphism  {\rm Usd} induces an automorphism of the Picard group of 
$\mathcal M_{s,p,2s}$ as follows.
\begin{equation}
\begin{split}
&({\rm Usd})^*(\check D_i)=\check D_{r+2-i},\,\,\, 2\leq i\leq r\,;\\
&({\rm Usd})^*(\check D_1)=-\sum_{i=1}^r\check D_i\,;\\
&({\rm Usd})^*\left((\check R_{s,p,2s})^*(\mathcal O_{G(p,2s)}(1))\right)=(\check R_{s,p,2s})^*(\mathcal O_{G(p,2s)}(1))-\sum_{i=1}^r(r+1-i)\cdot\check D_i\,;\\
&({\rm Usd})^*\left(\check B_i\right)=\check B_{r-i}\,,\,\,0\leq i\leq r.
\end{split}
\end{equation}
\end{lemma}
{\noindent\bf Proof of Lemma \ref{mdis}.} Since ${\rm USD}(D_1^-\cap D_i^-)=D_1^+\cap D_i^+$ and $\mathcal P_{s,p,2s}(D^+_i)=\mathcal P_{s,p,2s}(D^-_{r+2-i})$ for $2\leq i\leq r$, we conclude that
\begin{equation}
    \mathcal L_{s,p,2s}(D_1^+\cap D_i^+)=D_1^-\cap D_{r+2-i}^-\,,\,\,2\leq i\leq r\,.
\end{equation}
Therefore, $({\rm Usd})^*(\check D_i)=\check D_{r+2-i}$ for $2\leq i\leq r$.

By (\ref{bst}) we have that
\begin{equation}
\begin{split}
&B_0=(R_{s,p,2s})^*(\mathcal O_{G(p,2s)}(1))-\sum_{i=1}^{r}(r+1-i)\cdot D^+_i\,;\\
&B_p=(R_{s,p,2s})^*(\mathcal O_{G(p,2s)}(1))-\sum_{i=1}^{r}(r+1-i)\cdot D^-_i\,;\\
&B_1=(R_{s,p,2s})^*(\mathcal O_{G(p,2s)}(1))-\sum_{i=1}^{r-1}(r-i)\cdot D^+_i-D^-_1\,;\\
&B_{p-1}=(R_{s,p,2s})^*(\mathcal O_{G(p,2s)}(1))-D^+_1-\sum_{i=1}^{r-1}(r-i)\cdot D^-_i\,.
\end{split}
\end{equation}

Notice that $(R_{s,p,2s})^*(\mathcal O_{G(p,2s)}(1))\big|_{D^-_1}=B_0\cap  D_1^-$; $\mathcal L_{s,p,2s} (B_0\cap  D_1^+)=B_0\cap  D_1^-$;  by Lemma \ref{dis} ${\rm USD}^*(B_0)=B_p$ and $(R_{s,p,2s})^*(\mathcal O_{G(p,2s)}(1))$ is invariant under USD\,. We thus conclude that
\begin{equation}
\begin{split}
&({\rm Usd})^*\left((\check R_{s,p,2s})^*(\mathcal O_{G(p,2s)}(1))\right)=({\rm Usd})^*\left(B_0\big|_{D^-_1}\right)=\left({\rm USD}\big|_{\mathcal M_{s,p,2s}}\right)^*\left(\mathcal L_{s,p,2s}^*\left(B_0\big|_{D^-_1}\right)\right)\\
&=\left({\rm USD}\big|_{\mathcal M_{s,p,2s}}\right)^*\left(B_0\big|_{D^+_1}\right)=B_p\big|_{D^-_1}=(\check R_{s,p,2s})^*(\mathcal O_{G(p,2s)}(1))-\sum_{i=1}^{r}(r+1-i)\cdot\check D_i\,.
\end{split}
\end{equation}
Similarly, we have that
\begin{equation}
\begin{split}
&({\rm Usd})^*\left(\check B_1\right)=\left({\rm USD}\big|_{\mathcal M_{s,p,2s}}\right)^*\left(\mathcal L_{s,p,2s}^*\left(B_1\big|_{D^-_1}\right)\right)=\left({\rm USD}\big|_{\mathcal M_{s,p,2s}}\right)^*\left(B_1\big|_{D^+_1}\right)\\
&=B_{p-1}\big|_{D^-_1}=(\check R_{s,p,2s})^*(\mathcal O_{G(p,2s)}(1))-\sum_{i=1}^{r-1}(r-i)\cdot\check D_i\,.
\end{split}
\end{equation}
Therefore,
\begin{equation}
\begin{split}
&\,({\rm Usd})^*\left(\check D_1\right)=({\rm Usd})^*\left((\check R_{s,p,2s})^*(\mathcal O_{G(p,2s)}(1))-\check B_1\right)\\
&=(\check R_{s,p,2s})^*(\mathcal O_{G(p,2s)}(1))-\sum_{i=1}^{r}(r+1-i)\check D_i-\left((\check R_{s,p,2s})^*(\mathcal O_{G(p,2s)}(1))-\sum_{i=1}^{r-1}(r-i)\check D_i\right)\\
&=-\sum_{i=1}^{r}\check D_i\,\,.
\end{split}
\end{equation}
Substituting the above into (\ref{bst}), we can conclude that $({\rm Usd})^*\left(\check B_i\right)=(\check B_{r-i})$ for $0\leq i\leq r$.
 \,\,\,$\endpf$

\begin{lemma}\label{mdis2}Assume that  $1\leq p<s$. Then the  automorphism  {\rm Dual} induces an automorphism of the Picard group of 
$\mathcal M_{s,p,2p}$ as follows.
\begin{equation}
\begin{split}
&({\rm Dual})^*(\check D_i)=\check D_{r+2-i},\,\,\, 2\leq i\leq r\,;\\
&({\rm Dual})^*(\check D_1)=-\sum_{i=1}^r\check D_i\,;\\
&({\rm Dual})^*\left((\check R_{s,p,2p})^*(\mathcal O_{G(p,2p)}(1))\right)=(\check R_{s,p,2p})^*(\mathcal O_{G(p,2p)}(1))-\sum_{i=1}^r(r+1-i)\cdot\check D_i\,;\\
&({\rm Dual})^*\left(\check B_i\right)=\check B_{r-i}\,,\,\,0\leq i\leq r.
\end{split}
\end{equation}
\end{lemma}
{\noindent\bf Proof of Lemma \ref{mdis2}.} For $0\leq j\leq r$ define a divisor $\widetilde b_{j}:=\left\{x\in G(p,2p) \big|\,P_{I_j}(x)=0\right\}$ where $\widetilde I_{j}:=(2p,2p-1,\cdots,2p-j+1,2p-s-j,2p-s-j-1,\cdots,1)\in\mathbb I^j_{s,p,2p}$;
let $\widetilde B_{j}\subset\mathcal T_{s,p,2p}$ be the strict transformation of $\widetilde b_{j}$ under the canonical blow-up $R_{s,p,2p}$. 

It is easy to show by (\ref{dtrans}) that the image of $B_{r-j}$ under ${\rm DUAL}$ is $\widetilde B_{j}$ for $0\leq j\leq r$. Moreover, similarly to (\ref{bst}) we can derive the following for $0\leq j\leq r$.
\begin{equation}
 \widetilde   B_j=(R_{s,p,2p})^*(\mathcal O_{G(p,2p)}(1))-\sum_{i=1}^{r-j}(r-j+1-i)\cdot D^+_i-\sum_{i=1}^j(j+1-i)\cdot D^-_i=B_j\,.
\end{equation}

Then the remainder of the proof is the same as in Lemma \ref{mdis}. \,\,\,$\endpf$
\medskip

Recall that when $n=2p=2s$, $(\check R_{p,p,2p})^*(\mathcal O_{G(p,2p)}(1))$ and hence $\sum_{i=1}^r(r+1-i)\cdot\check D_i$ are trivial. We can derive that
\begin{lemma}\label{mdis3}Assume that  $p\geq1$.  The  automorphisms  {\rm Usd} and {\rm Dual} induce the following automorphism of the Picard group of $\mathcal M_{p,p,2p}$.
\begin{equation}
\begin{split}
&({\rm Usd})^*(\check D_i)=({\rm Dual})^*(\check D_i)=\check D_{r+2-i},\,\,\, 2\leq i\leq r\,;\\
&({\rm Usd})^*(\check D_1)=({\rm Dual})^*(\check D_1)=-\sum_{i=2}^r\frac{i-1}{r}\cdot\check D_i\,;\\
&({\rm Usd})^*\left(\check B_i\right)=({\rm Dual})^*\left(\check B_i\right)=\check B_{r-i}\,,\,\,1\leq i\leq r-1.
\end{split}
\end{equation}
\end{lemma}

\subsection{The automorphism groups of \texorpdfstring{$\mathcal T_{s,p,n}$}{ff}}  \label{symt} In this subsection, we determine the automorphism groups of  $\mathcal T_{s,p,n}$. We first prove the following.

\begin{lemma}\label{des}
Let $\sigma\in{\rm Aut}(\mathcal T_{s,p,n})$. Denote by $\Sigma$ the birational self-map of $G(p,n)$ induced by $\sigma$. If the induced automorphism $\sigma^*$ of the Picard group of $\mathcal T_{s,p,n}$ is the identity map, then  $\Sigma$ extends to an automorphism of $G(p,n)$.
\end{lemma}
{\noindent\bf Proof of Lemma \ref{des}.} For convenience, we use the following notation to distinguish the source and the target.
\begin{equation}\label{idesent}
\begin{array}{ccc}
\vspace{.03in}
\mathcal T^1_{s,p,n}&\xrightarrow{\,\,\,\,\,\,\sigma\,\,\,\,\,} &\mathcal T^2_{s,p,n}\\
\vspace{-.04in}
\Big\downarrow\llap{$\scriptstyle R^1_{s,p,n}\,\,\,\,\,$}&&\Big\downarrow\rlap{$\scriptstyle R^2_{s,p,n}\,\,\,\,\,$}\\
G^1(p,n)&\xdashrightarrow{\,\,\,\Sigma\,\,\,}&G^2(p,n)\\
\end{array}\,\,\,\,\,\,\,\,\,\,\,\,\,\,\,\,.
\end{equation}

We claim that any connected curve $E\subset\mathcal  T_{s,p,n}^1$ which is contracted by $R_{s,p,n}^1$ is also contracted by $R_{s,p,n}^2$, that is,  if the image $R_{s,p,n}^1(E)$ is a point then so is the image $R_{s,p,n}^2(\sigma(E))$. Assume the contrary. Calculating the intersection numbers, we have a contradiction as follows.
\begin{equation}
\begin{split}
     0&=\mathcal O_{G^1(p,n)}(1)\cdot R_{s,p,n}^1(E)\\
     &=(R_{s,p,n}^1)^*(\mathcal O_{G^1(p,n)}(1))\cdot E\\
     &=(\sigma^{-1})^*\left((R_{s,p,n}^1)^*\left(\mathcal O_{G^1(p,n)}(1) \right)\right)\cdot \sigma (E) \\
     &=(R_{s,p,n}^2)^*(\mathcal O_{G^2(p,n)}(1))\cdot \sigma(E) \\
     &=\mathcal O_{G^2(p,n)}(1)\cdot R_{s,p,n}^2\left(\sigma(E)\right)\\
     &\geq 1\,.
\end{split}
\end{equation} 
Here the second and the fifth equalities are by the push forward formula (see \cite{Ful}); the third equality is because $\sigma$ is biholomorphic; the fourth equality is due to the assumption that $\sigma^*$ is the identity map; the first equality and the last inequality are by the assumption that $R_{s,p,n}^1(E)$ is a point while $R_{s,p,n}^2(\sigma(E))$ is a curve.

We next prove that $\Sigma$ extends to a holomorphic map 
from $G^1(p,n)$ to $ G^2(p,n)$. It suffices to prove that the map $\Sigma$ is bounded near each point of $G^1(p,n)$. Take an arbitrary point $x\in G^1(p,n)$.  Recall that  $\mathcal T_{s,p,n}$ can be locally realized as a sequence of blow-ups of $G^1(p,n)$ along smooth submanifolds (see Remark \ref{ivdw}), and hence each intermediate blow-up has connected fibers.  We conclude that the fiber $(R_{s,p,n}^1)^{-1}(x)$ is connected.

If $(R_{s,p,n}^2\circ\sigma)\left((R_{s,p,n}^1)^{-1}(x)\right)$ consists of more than one point, it must be of dimension at least $1$ for fibers are connected. We can show that there is a curve $E\subset (R_{s,p,n}^1)^{-1}(x)$ such that  $R_{s,p,n}^2\left(\sigma(E)\right)$ is a curve. This is a contradiction.

If the image $(R_{s,p,n}^2\circ\sigma)\left((R_{s,p,n}^1)^{-1}(x)\right)$  consists of a single point, there is a neighborhood $U_x$ of $x$ in $G^1(p,n)$ such that $\Sigma(U_x)$ is contained in a small neighborhood of $(R_{s,p,n}^2\circ\sigma)\left((R_{s,p,n}^1)^{-1}(x)\right)$. Then in certain local coordinate charts, we can write $\Sigma$ as  a bounded vector-valued rational function which is bounded near $x$. We thus conclude $\Sigma$ extend to a holomorphic map near $x$.

Hence, $\Sigma$ extends to a holomorphic map 
from $G^1(p,n)$ to $ G^2(p,n)$. Similarly, we can show that the inverse $\Sigma^{-1}$ has a holomorphic extension. As a conclusion, we conclude that   ${\Sigma}$ extends to an automorphism of $G(p,n)$.
\,\,\,$\endpf$
\medskip

The following lemma classifies the induced automorphisms of the Picard groups.
\begin{lemma}\label{sigma2}
Let $\sigma$ be an automorphism of $\mathcal T_{s,p,n}$ and $\sigma^*$ the induced automorphism of the Picard group. When $n\neq 2s$ and $n\neq 2p$, $\sigma^*$ is the identity map. When $n=2s$, $\sigma^*$ is the identity map or $({\rm USD})^*$. When  $n=2p$, $\sigma^*$ is the identity map or $({\rm DUAL})^*$.
\end{lemma}
{\bf\noindent Proof of Lemma \ref{sigma2}.} See Appendix \ref{section:rigidmc}.
\,\,\,$\endpf$
\medskip

{\noindent\bf Proof of Proposition \ref{auto1}.}  Without loss of generality, we may assume that $2p\leq n\leq 2s$. We prove Proposition \ref{auto1} based on a case by case argument.
\smallskip

{\bf\noindent Case 1 ($p=s=n-s=1$).} $\mathcal T_{1,1,2}\cong\mathbb {CP}^1$ and hence the automorphism group is $PGL(2,\mathbb C)$.
\smallskip

{\bf\noindent Case 2 ($p=n-s=1$ and $s\geq 2$).} Let $\sigma$ be an automorphism of $\mathcal T_{s,p,n}$. By Lemma \ref{sigma2}, $\sigma^*$ is the identity map. By Lemma \ref{des}, $\sigma$ induces an automorphism  $\Sigma$ of $G(p,n)$. Notice that each irreducible component of the exceptional divisor is invarinat under $\sigma$.   Therefore,  $\Sigma$ maps $\overline{\mathcal V_{(0,1)}^{+}}$ to $\overline{\mathcal V_{(0,1)}^{+}}$.

Write $\Sigma$ as a $n\times n$ matrix
\begin{equation}\label{block}
    \Sigma=\left(\begin{matrix}
    A&B\\
    C&D\\
    \end{matrix}\right)\,,
\end{equation}
where $A$, $B$, $C$, $D$ are submatrices of sizes $s\times s$, $s\times (n-s)$, $(n-s)\times s$, $(n-s)\times (n-s)$ respectively.
If  $C$ is not a zero matrix, there is a point $x\in \overline{V_{(0,1)}^{+}}$ such that $\Sigma(x)\notin\overline{V_{(0,1)}^{+}}$, which is a contradiction. Therefore,  $C$ is a zero matrix and $\Sigma\in P$ where $P$ is defined by (\ref{parabolic}).

\smallskip

{\bf\noindent Case 3 ($p=1$ and $2\leq n-s<s$).}  Let $\sigma$ be an automorphism of $\mathcal T_{s,p,n}$.  By Lemma \ref{sigma2}, $\sigma^*$ is the identity map. By Lemma \ref{des}, $\sigma$ induces an automorphism  $\Sigma$ of $G(p,n)$. Similarly to Case 2, $\Sigma$ maps $\overline{\mathcal V_{(0,1)}^{+}}$ to $\overline{\mathcal V_{(0,1)}^{+}}$ and  $\overline{\mathcal V_{(1,0)}^{-}}$ to $\overline{\mathcal V_{(1,0)}^{-}}$.

Recall that  the automorphism group of $G(p,n)$ is generated by $PGL(n,\mathbb C)$ when $n\neq 2p$. Write $\Sigma$ as (\ref{block}). 
If $\Sigma\not\in GL(s,\mathbb C)\times GL(n-s,\mathbb C)$, either $C$ or $B$ is not a zero matrix. If  $C$ is not a zero matrix, there is a point $x\in \overline{V_{(0,1)}^{+}}$ such that $\Sigma(x)\notin\overline{V_{(0,1)}^{+}}$; if  $B$ is not a zero matrix, there is a point $x\in \overline{V_{(1,0)}^{-}}$ such that $\Sigma(x)\notin\overline{V_{(1,0)}^{-}}$. 

Therefore,  $\Sigma\in GL(s,\mathbb C)\times GL(n-s,\mathbb C)$.

\smallskip

{\bf\noindent Case 4 ($2\leq p$ and $2p<n<2s$).} The proof is the same as in Case 3. Let $\sigma$ be an automorphism of $\mathcal T_{s,p,n}$. By Lemma \ref{sigma2}, $\sigma^*$ is the identity map. By Lemma \ref{des}, $\sigma$ induces an automorphism  $\Sigma$ of $G(p,n)$.

Since $r\geq 2$, we can conclude that $\Sigma$ maps $\overline{\mathcal V_{(p-r+1,r-1)}^{+}}$ to $\overline{\mathcal V_{(p-r+1,r-1)}^{+}}$ and  $\overline{\mathcal V_{(p-1,1)}^{-}}$ to $\overline{\mathcal V_{(p-1,1)}^{-}}$. Write $\Sigma$ as (\ref{block}). If  $C$ is not a zero matrix, there is a point $x\in \overline{V_{(p-r+1,r-1)}^{+}}$ such that $\Sigma(x)\notin\overline{V_{(p-r+1,r-1)}^{+}}$; if  $B$ is not a zero matrix, there is a point $x\in \overline{V_{(p-1,1)}^{-}}$ such that $\Sigma(x)\notin\overline{V_{(p-1,1)}^{-}}$.

Therefore, 
$\sigma\in GL(s,\mathbb C)\times GL(n-s,\mathbb C)$.

\smallskip

{\bf\noindent Case 5 ($1\leq p<s=n-s$).} Let $\sigma$ be an automorphism of $\mathcal T_{s,p,n}$. By Lemma \ref{sigma2}, $\sigma^*$ is either the identity map or the automrophism ${\rm USD}^*$.  By Lemma \ref{des}, either $\sigma$ or ${\rm USD}\circ\sigma$ induces an automorphism  $\Sigma$ of $G(p,n)$. Then
$\Sigma\in GL(s,\mathbb C)\times GL(n-s,\mathbb C)$  similarly to Case 4.

\smallskip

{\bf\noindent Case 6 ($2\leq p<s$ and $n=2p$).} Let $\sigma$ be an automorphism of $\mathcal T_{s,p,n}$. Recall that the automorphism group of $G(p,2p)$ is generated by $PGL(2p,\mathbb C)$ and the dual automorphism ${\empty}^*$ (see Section \ref{isomt}). By Lemma \ref{sigma2}, $\sigma^*$ is either the identity map or the automrophism ${\rm DUAL}^*$.   By Lemma \ref{des}, either $\sigma$ or ${\rm DUAL}\circ\sigma$ induces an automorphism  $\Sigma$ of $G(p,n)$. 

If $\Sigma\in PGL(n,\mathbb C)$,  we can conclude that $\Sigma\in GL(s,\mathbb C)\times GL(n-s,\mathbb C)$ similarly to Case 4.

If $\Sigma\not\in PGL(n,\mathbb C)$, we have that  $\widetilde\Sigma:=\Sigma\circ {\empty}^*\in PGL(n,\mathbb C)$ where ${\empty}^*$ takes the following form (see (\ref{dtrans})). 
\begin{equation}
\left({ \left(I_{p\times p} \hspace{-0.13in}\begin{matrix}
  &\hfill\tikzmark{c1}\\
  &\hfill\tikzmark{d1}
  \end{matrix}\hspace{-0.1in}\begin{matrix}
  &\hfill\tikzmark{c3}\\
  &\hfill\tikzmark{d3}
  \end{matrix}\,\,\,X\hspace{-0.13in}\begin{matrix}
  &\hfill\tikzmark{c4}\\
  &\hfill\tikzmark{d4}
  \end{matrix}\,\,\, A \right)}\right)^*=\left(\begin{matrix} -X^T\\
      -A^T
  \end{matrix} \hspace{-0.13in}\begin{matrix}
  &\hfill\tikzmark{c2}\\
  \\
  &\hfill\tikzmark{d2}
  \end{matrix}\hspace{-0.1in}\begin{matrix}
  &\hfill\tikzmark{c5}\\
  \\
  &\hfill\tikzmark{d5}
  \end{matrix}\,\,\, \begin{matrix} I_{(s-p)\times (s-p)}\\
      0
  \end{matrix} \hspace{-0.13in}\begin{matrix}
  &\hfill\tikzmark{c6}\\
  \\
  &\hfill\tikzmark{d6}
  \end{matrix}\,\,\,\begin{matrix} 0\\
      I_{(2p-s)\times (2p-s)}
  \end{matrix}\right).
  \tikz[remember picture,overlay]   \draw[dashed,dash pattern={on 4pt off 2pt}] ([xshift=0.5\tabcolsep,yshift=7pt]c1.north) -- ([xshift=0.5\tabcolsep,yshift=-2pt]d1.south);
  \tikz[remember picture,overlay]   \draw[dashed,dash pattern={on 4pt off 2pt}] ([xshift=0.5\tabcolsep,yshift=7pt]c2.north) -- ([xshift=0.5\tabcolsep,yshift=-2pt]d2.south);\tikz[remember picture,overlay]   \draw[dashed,dash pattern={on 4pt off 2pt}] ([xshift=0.5\tabcolsep,yshift=7pt]c3.north) -- ([xshift=0.5\tabcolsep,yshift=-2pt]d3.south);\tikz[remember picture,overlay]   \draw[dashed,dash pattern={on 4pt off 2pt}] ([xshift=0.5\tabcolsep,yshift=7pt]c4.north) -- ([xshift=0.5\tabcolsep,yshift=-2pt]d4.south);\tikz[remember picture,overlay]   \draw[dashed,dash pattern={on 4pt off 2pt}] ([xshift=0.5\tabcolsep,yshift=7pt]c5.north) -- ([xshift=0.5\tabcolsep,yshift=-2pt]d5.south);\tikz[remember picture,overlay]   \draw[dashed,dash pattern={on 4pt off 2pt}] ([xshift=0.5\tabcolsep,yshift=7pt]c6.north) -- ([xshift=0.5\tabcolsep,yshift=-2pt]d6.south);
\end{equation}
where $X$ is an $p\times(s-p)$ matrix.
Write $\widetilde\Sigma$ as
\begin{equation}
    \widetilde\Sigma=\left(\begin{matrix}
    A&B\\
    C&D\\
    \end{matrix}\right)\,
\end{equation}
where $A$, $B$, $C$, $D$ are matrices of sizes $s\times s$, $s\times (2p-s)$, $(2p-s)\times s$, $(2p-s)\times (2p-s)$ respectively. It is easy to verify that $\overline{\mathcal V_{(p,0)}^{-}}$ is not invariant under  $\Sigma=\widetilde\Sigma\circ{\empty}^*$. This is a contradiction since $D_1^-$ is invariant under $\sigma$.

\smallskip

{\bf\noindent Case 7 ($2\leq p=s=n-s$).} Let $\sigma$ be an automorphism of $\mathcal T_{s,p,n}$. By Lemma \ref{sigma2}, $\sigma^*$ is either the identity map or the automrophism ${\rm USD}^*$. By Lemma \ref{des}, either $\sigma$ or ${\rm USD}\circ\sigma$ induces an automorphism of $G(p,2p)$.  If  $\Sigma\in PGL(n,\mathbb C)$, then $\Sigma\in GL(s,\mathbb C)\times GL(n-s,\mathbb C)$ in the same way as Case 4. Otherwise, 
${\rm USD}\circ{\rm DUAL}\circ\sigma\in GL(s,\mathbb C)\times GL(n-s,\mathbb C)$ induces an automorphism  $\widetilde \Sigma$ of $G(p,2p)$ which is contained in $PGL(n,\mathbb C)$. Hence $\widetilde \Sigma\in GL(s,\mathbb C)\times GL(n-s,\mathbb C)$ in the same way as Case 4.

\smallskip

We complete the proof of Proposition \ref{auto1}.  \,\,\,$\endpf$

\subsection{The automorphism groups of \texorpdfstring{$\mathcal M_{s,p,n}$}{ff}} \label{symm}
In this subsection, we will determine the automorphism groups of $\mathcal M_{s,p,n}$. 

\begin{lemma}\label{pre}
Assume that $p\neq q$ are positive integers. Let $\sigma$ be an automorphism of $\mathcal M_{p,p,p+q}$ such that $\sigma^*$ is the identity map on the Picard group of $\mathcal M_{p,p,(p+q)}$.  Then 
\begin{equation}
    \sigma\in PGL(p,\mathbb C)\times PGL(q,\mathbb C).
\end{equation} 
\end{lemma}

{\noindent\bf Proof of Lemma \ref{pre}.} 
When $p=1$ or $q=1$, this is trivial for $\mathcal M_{p,p,p+q}$ is isomorphic to a projective space.

Recall that $\mathcal M_{p,p,p+q}\cong\mathcal M_{q,q,p+q}\cong P(M_{p\times q})$, where $\widetilde P(M_{p\times q})$ is the variety of complete collineations in   Example \ref{ccolli}. Without loss of generality, we can assume that $2\leq p<q$.

Since $\sigma^*$ is the identity map on the Picard group of $\mathcal M_{p,p,p+q}$, similarly to Lemma \ref{des} we can show that $\sigma$ descends to an automorphism of $P(M_{p\times q})$ as follows.
\begin{equation}
\begin{array}{ccc}
\vspace{.03in}
\mathcal M_{p,p,p+q}&\xrightarrow{\,\,\,\,\,\sigma\,\,\,\,\,} &\mathcal M_{p,p,p+q}\\
\vspace{-.03in}
\Big\downarrow\llap{$\scriptstyle R_{p,p,p+q}|_{\mathcal M_{p,p,p+q}}\,\,\,\,\,$}&&\Big\downarrow\rlap{$\scriptstyle R_{p,p,p+q}|_{\mathcal M_{p,p,p+q}}\,\,\,\,\,$}\\
P(M_{p\times q})&\xrightarrow{\,\,\,\,\,\Sigma\,\,\,\,\,}&P(M_{p\times q})\\
\end{array}\,\,\,\,\,\,\,\,\,\,\,\,\,\,\,\,\,\,\,\,\,\,\,\,\,\,\,\,\,\,\,\,.
\end{equation}
Here $P(M_{p\times q})\cong\mathbb {CP}^{pq-1}$ is
the projectivization of the matrix group (see (\ref{ccolli1})). Moreover, we can conclude that $\Sigma$ preserves the ranks of the matrices in $P(M_{p\times q})$ for the exceptional divisors are $\sigma$-invariant.

It is clear that $\Sigma\in PGL(pq,\mathbb C)$. Hence, to prove Lemma \ref{msigma} it suffices to show that $\Sigma\in GL(p,\mathbb C)\times GL(q,\mathbb C)$ by viewing  $\Sigma$ as an element of $GL(pq,\mathbb C)$.

For $1\leq i\leq p$ and $1\leq j\leq q$, denote by $E_{ij}$ the $p\times q$ matrix such that the $(i,j)^{th}$ entry is $1$ and zero elsewhere. Write $\Sigma$  as a linear transformation
\begin{equation}
\Sigma(E_{uv})=\sum_{i=1}^{p}\sum_{j=1}^{q}a_{uv}^{ij}\cdot E_{ij},\,\,\,\,\,\,a_{uv}^{ij}\in\mathbb C\,,\,\,1\leq u\leq p\,,\,\,1\leq v\leq q.
\end{equation}

Since  $\Sigma$ maps a matrix of rank $1$ to another matrix of rank $1$, there are integers  $1\leq i\leq p$ and $1\leq j\leq q$ such that $a_{11}^{ij}\neq 0$. Composing $\Sigma$  with an element of $GL(p,\mathbb C)\times GL(q,\mathbb C)$, we can assume that $a_{11}^{11}\neq 0$; indeed there is an element  $\tau\in GL(p,\mathbb C)\times GL(q,\mathbb C)$ such that
$(\tau\circ\Sigma)(E_{11})=E_{11}$. By a slight abuse of notation, we use $\Sigma$ for any composition $\tau\circ\Sigma$, $\tau\in GL(p,\mathbb C)\times GL(q,\mathbb C)$, in the following.

We claim that there are integers $2\leq i\leq p$ and $2\leq j\leq q$ such that $a_{22}^{ij}\neq 0$.  Assume the contrary. Then there are integers  $2\leq j\leq q$ and $2\leq i\leq p$ such that $a_{22}^{i1}\neq 0$ and $a_{22}^{1j}\neq 0$, for  $\Sigma(E_{11}+E_{22})$ has rank $2$. However, this contradicts the fact that $\Sigma(E_{22})$ has rank $1$. 

Composing $\Sigma$ with a certain element of $GL(p,\mathbb C)\times GL(q,\mathbb C)$  we can derive that $\Sigma(E_{11})=E_{11}$ and $\Sigma(E_{22})=E_{22}$. Similarly, we may assume that $\Sigma(E_{ii})=E_{ii}$  for $1\leq i\leq p$.

Since $E_{11}+\lambda E_{12}$ (resp. $E_{22}+\lambda E_{12}$) is of rank $1$ for $\lambda\in\mathbb C$, we can show that
\begin{equation}\label{150}
\Sigma(E_{12})=\sum_{i=1}^{p}a_{12}^{i1}\cdot E_{i1}\,\,\,{\rm or}\,\,\,\sum_{j=1}^{q}a_{12}^{1j}\cdot E_{1j}\,\,\left({\rm resp}.\,\, \Sigma(E_{12})=\sum_{i=1}^{p}a_{12}^{i2}\cdot E_{i1}\,\,\,{\rm or}\,\,\,\sum_{j=1}^{q}a_{12}^{2j}\cdot E_{1j}\right).
\end{equation}
Then $\Sigma(E_{12})=a_{12}^{12}\cdot E_{12}$  or $a_{12}^{21}\cdot E_{21}$.
We claim that $\Sigma(E_{12})=a_{12}^{12}\cdot E_{12}$. Assume the contrary. Then $\Sigma(E_{1i})=a_{1i}^{i1}\cdot E_{i1}$ for $2\leq i\leq p$. Then regardless of the value of $\Sigma(E_{1(p+1)})$, it violates the assumption that the ranks are preserved under $\Sigma$.

By the same argument, we can assume that
\begin{enumerate}
    \item $\Sigma(E_{ij})=a^{ij}_{ij}\cdot E_{ij}$ for $1\leq i\leq p$ and $1\leq j\leq q$.
    \item $a^{ij}_{ij}=1$ when $i\equiv j\,\, ({\rm mod}\, p)$.
\end{enumerate}
Composing an element of  $PGL(p,\mathbb C)$ by the left action, we can fix $a^{21}_{21}=\cdots a^{p1}_{p1}=1$. Then  $a^{ij}_{ij}=1$, $1\leq i\leq p$, $1\leq j\leq q$, for  $E_{11}+E_{1j}+E_{i1}+E_{ij}$ has rank $1$. 

We conclude that $\Sigma\in GL(p,\mathbb C)\times GL(q,\mathbb C)$, and hence complete the proof of Lemma \ref{pre}.
\,\,\,$\endpf$
\medskip

We generalize Lemma \ref{pre} as follows.

\begin{lemma}\label{fibpre}
Assume that $2p\leq n\leq 2s$ and $2\leq p<s$. Let $\sigma$ be an automorphism of $\mathcal M_{s,p,n}$ such that $\sigma^*$ is the identity map on the Picard group of $\mathcal M_{s,p,n}$.  Then 
$\sigma\in PGL(s,\mathbb C)\times PGL(n-s,\mathbb C)$.
\end{lemma}
{\bf\noindent Proof of Lemma \ref{fibpre}.} See Appendix \ref{section:rigidtc} (or Appendix \ref{section:rigidtc2} for an alternative proof).\,\,\,$\endpf$.
\medskip

Similarly to Lemma \ref{sigma2}, we have 
\begin{lemma}\label{msigma}
Assume that $2p\leq n\leq 2s$ and $2\leq p<s$. Let $\sigma$ be an automorphism of $\mathcal M_{s,p,n}$ and $\sigma^*$ the induced automorphism of the Picard group. When $n\neq 2s$ and $n\neq 2p$, $\sigma^*$ is the identity map. When $n=2s$, $\sigma^*$ is the identity map or $({\rm Usd})^*$. When  $n=2p$, $\sigma^*$ is the identity map or $({\rm Dual})^*$.
\end{lemma}
{\noindent\bf Proof of Lemma \ref{msigma}.} Denote by  Con$(\mathcal M_{s,p,n})\subset A_{n-1}(\mathcal M_{s,p,n})\otimes_{\mathbb Z}\mathbb Q$ the cone  of effective divisors of $\mathcal M_{s,p,n}$. By Lemmas \ref{checkgb} and \ref{cone}, Con$(\mathcal M_{s,p,n})$ is generated by 
\begin{equation}
 \{\check D_{2}, \check D_{3}, \cdots, \check D_{r}, \check B_0, \check B_1, \cdots, \check B_r\}\,.
\end{equation}
Recalling (\ref{mb=b}), (\ref{mb=br}),   we can determine $\mathfrak G$ the set of extremal rays of Con$(\mathcal M_{s,p,n})$ as follows.
When $p\neq n-s$ and $p<s$, 
\begin{equation}
    \mathfrak G=\{\check D_{2}, \check D_{3}, \cdots, \check D_{r}, \check B_0, \check B_r\};
\end{equation}
when $n-s=p<s$,  
\begin{equation}
    \mathfrak G=\{\check D_{2}, \check D_{3}, \cdots, \check D_{r-1}, \check B_0,  \check B_{r}\}.
\end{equation}
Since $\sigma^*$ preserves Con$(\mathcal M_{s,p,n})$, it induces a permutation of $\mathfrak G$.

We next prove Lemma \ref{msigma} based on a case by case argument. For convenience, we denote by $H$ the divisor $(\check R_{s,p,n})^*\left(\mathcal O_{ G(p,n)}(1)\right)$; also recall that $2p\leq n\leq 2s$.
\smallskip

{\bf\noindent Case 1 ($n-s=1$).} 
Con$(\mathcal M_{s,p,n})$ is generated by $\check B_0,\check B_1$. Recall that the canonical bundle $K_{\mathcal M_{s,p,n}}$ is invariant under $\sigma$. Moreover, by (\ref{m3}) we have
\begin{equation}
\begin{split}
K_{\mathcal M_{s,p,n}}&=-\left(n-p\right)\check B_0-p\check B_1\,.
\end{split}
\end{equation}
If $\sigma^*$ is not the identity map, then $n=2p$ and $\sigma^*=({\rm Dual})^*$.
\smallskip

{\bf\noindent Case 2 ($p<s$, $p\neq n-s$, and $n-s\geq 2$).} It is clear that $r\geq 2$, and
Con$(\mathcal M_{s,p,n})$ is generated by  $\check D_2,\cdots, \check D_r,\check B_0,\check B_r$. 
By (\ref{mkan1}) we have that \begin{equation}\label{mkan3}
\small
\begin{split}
&K_{\mathcal M_{s,p,n}}=-n\check B_0+\frac{p(n-s)}{r}\left(\check B_0-\check B_r-\sum_{k=2}^r (r+1-k)\check D_k\right)+\sum_{i=2}^{r}\big((p-i+1)(n-s-i+1)-1\big) \check D_{i}\\
&=\sum_{i=2}^{r}\left((p-i+1)(n-s-i+1)-1-\frac{p(n-s)(r+1-i)}{r}\right)\check D_{i}-\left(n-\frac{p(n-s)}{r}\right)\check B_0-\frac{p(n-s)}{r}\check B_r\,.
\end{split}
\end{equation}

We first assume that $r=2$.
If $p=2$, then $K_{\mathcal M_{s,p,n}}=-2\check D_{2}-s\check B_0-(n-s)\check B_r$; if $n-s=2$, then $K_{\mathcal M_{s,p,n}}=-2\check D_{2}-(n-p)\check B_0-p\check B_r$. Then $\sigma^*$ is  the identity map.

Let $r\geq 3$ in the following. By (\ref{mb=b}) we have 
\begin{equation}\label{iden1}
\sigma^*(\check D_1)=\frac{1}{r}\sigma^*(\check B_0)-\frac{1}{r}\sigma^*(\check B_r)-\sum_{i=2}^r\frac{r+1-i}{r} \sigma^*(\check D_i)\,.
\end{equation}
Notice that, with respect to the basis
$\left\{H, \check D_2,\cdots,\check D_r\right\}$, (\ref{iden1}) has integer coefficients by Lemma \ref{mpicb}. Checking the coefficient of $H$, we can conclude the following possibilities.
\begin{enumerate}[label=(\alph*).]
\item $\sigma^*$ is the identity map.
\item $\sigma^*(\check B_0)=\check B_r$,  $\sigma^*(\check B_r)=\check B_0$.
\item There is an integer $2\leq l\leq r$ such that $\sigma^*(\check D_l)=\check B_0$ and $\sigma^*(\check D_{r+2-l})=\check B_r$.
\item $\sigma^*(\check B_0)=\check B_0$ and $\sigma^*(\check D_{r})=\check B_r$, or $\sigma^*(\check B_0)=\check B_r$ and $\sigma^*(\check D_{r})=\check B_0$.
\item $\sigma^*(\check B_r)=\check B_0$ and $\sigma^*(\check D_{2})=\check B_r$, or $\sigma^*(\check B_r)=\check B_r$ and $\sigma^*(\check D_2)=\check B_0$.
\end{enumerate}

For Case (b), by checking the coefficients of $\check D_2,\cdots, \check D_r$  in (\ref{iden1})  we can conclude that $\sigma^*(\check D_i)=\check D_{r+2-i}$ for $2\leq i\leq r$. By (\ref{mkan3}) we have 
$n-\frac{p(n-s)}{r}=\frac{p(n-s)}{r}$;
hence $n=2s$ and $\sigma^*=({\rm Usd})^*$,  or $n=2p$ and $\sigma^*=({\rm Dual})^*$.

For Case (c) checking the coefficient of $\check D_2$ in (\ref{iden1}), we can conclude that $l=2$. By (\ref{mkan3}), we can conclude that $s=p$ which is a contradiction.

For Cases (d) and (e), we have  $p=n-s$ or $p=s$ by (\ref{mkan3}) which contradicts the assumption. 

\smallskip

{\bf\noindent Case 3 ($p=n-s<s$).} 
If $r=p=n-s=2$, Con$(\mathcal M_{s,p,n})$ is generated by $\check B_0,\check D_2$. Then,
\begin{equation}
\begin{split}
K_{\mathcal M_{s,p,n}}=-s\cdot \check B_0-2\cdot\check D_r\,.
\end{split}
\end{equation}
If $\sigma^*$ is not the identity map, $s=2=p$. This contradicts our assumption that $p<s$.

Let $r=p=n-s\geq 3$.
Con$(\mathcal M_{s,p,n})$ is generated by $\check B_0,\check B_r,\check D_2,\cdots, \check D_{r-1}$. By (\ref{mkan1}) we have that
\begin{equation}\label{mkan4}
\small
\begin{split}
K_{\mathcal M_{s,p,n}}&=-n \check B_0+p\left(\check B_0-\check B_r-\sum_{k=2}^{r-1} (r+1-k)\check D_k\right)+\sum_{i=2}^{r-1}\big((p-i+1)(p-i+1)-1\big)\check D_{i}\\
&=-\sum_{i=2}^{r-1}\big((p-i+1)(i-1)+1\big)\check D_{i}-\left(n-p\right)\check B_0-p\check B_r\,.
\end{split}
\end{equation}
By (\ref{mb=b}) and (\ref{mb=br})  we have 
\begin{equation}\label{iden2}
\sigma^*(\check D_1)=\frac{1}{r}\sigma^*(\check B_0)-\frac{1}{r}\sigma^*(\check B_r)-\sum_{i=2}^{r-1}\frac{r+1-i}{r} \sigma^*(\check D_i)\,.
\end{equation}
 Checking the coefficient of $H$, we can conclude the following possibilities.
\begin{enumerate}[label=(\alph*).]
\item $\sigma^*$ is the identity map.
\item $\sigma^*(\check B_0)=\check B_r$,  $\sigma^*(\check B_r)=\check B_0$.
\item There is an integer $3\leq i\leq r-1$ such that $\sigma^*(\check D_i)=\check B_0$ and $\sigma^*(\check D_{r+2-i})=\check B_r$.
\item $\sigma^*(\check B_r)=\check B_0$ and $\sigma^*(\check D_{2})=\check B_r$.
\item $\sigma^*(\check B_r)=\check B_r$ and $\sigma^*(\check D_2)=\check B_0$.
\end{enumerate}

For Case (b), by  (\ref{mkan4}) we have 
$n=2p$ which is a contradiction.

For Case (c) checking the coefficient of $\check D_2$ in (\ref{iden2}), we can derive a contradiction.

For Cases (d) and (e), we have $p=n-p$ by (\ref{mkan4}) which contradicts the assumption. 
\smallskip

We thus complete the proof of Lemma \ref{msigma}.
\,\,\,$\endpf$

\medskip

{\noindent\bf Proof of Proposition \ref{mauto}.}   
We  prove Proposition \ref{mauto} based on a case by case argument. Note that $2p\leq n\leq 2s$.

\smallskip

{\bf\noindent Case 1 ($p=1$).} 

If $p=s=1$ and $n=2$, $\mathcal M_{1,1,2}$ is trivially a point.

If $p=n-s=1$ and $n\geq 3$, $\mathcal M_{s,1,s+1}\cong\mathbb {CP}^{n-1}$ and  ${\rm Aut}(\mathcal M_{s,1,s+1})=PGL(n-1,\mathbb C)$.

If $p=1$, $2\leq n-s<s$, $\mathcal M_{s,1,2s}\cong\mathbb {CP}^{s-1}\times \mathbb {CP}^{n-s-1}$ and  
\begin{equation}
    {\rm Aut}(\mathcal M_{s,1,n})=PGL(s,\mathbb C)\times PGL(n-s,\mathbb C).
\end{equation}

If $p=1$, $2\leq n-s=s$, $\mathcal M_{s,1,n}\cong\mathbb {CP}^{s-1}\times \mathbb {CP}^{s-1}$ and  
\begin{equation}
    {\rm Aut}(\mathcal M_{s,1,2s})=PGL(s,\mathbb C)\times PGL(s,\mathbb C)\rtimes \mathbb Z/2\mathbb Z.
\end{equation}

{\bf\noindent Case 2 ($2\leq p$, $n\neq 2s$ and $n\neq 2p$).} 
Let $\sigma\in{\rm Aut}(\mathcal M_{s,p,n})$.
By Lemma \ref{msigma}, $\sigma^*$ is the identity map on the Picard group of $\mathcal M_{s,p,n}$. By By Lemma \ref{fibpre}, we have that
\begin{equation}
    \sigma\in PGL(s,\mathbb C)\times PGL(n-s,\mathbb C)\,.
\end{equation}

\smallskip

{\bf\noindent Case 3 ($2\leq p<s$ and $n=2s$).}
Let $\sigma\in{\rm Aut}(\mathcal M_{s,p,n})$.  By Lemma \ref{msigma}, $\sigma^*$ is either the identity map or the automrophism $({\rm Usd})^*$.  By Lemma \ref{fibpre}, either $\sigma$ or ${\rm Usd}\circ\sigma$ is in $PGL(s,\mathbb C)\times PGL(n-s,\mathbb C)$. Then, \begin{equation}
    {\rm Aut}(\mathcal M_{s,p,2s})=PGL(s,\mathbb C)\times PGL(n-s,\mathbb C)\rtimes \mathbb Z/2\mathbb Z.
\end{equation}

{\bf\noindent Case 4 ($2\leq p<s$ and $n=2p$).}
Let $\sigma\in{\rm Aut}(\mathcal M_{s,p,n})$. By Lemma \ref{msigma}, $\sigma^*$ is either the identity map or the automrophism $({\rm Dual})^*$.   By Lemma \ref{fibpre}, either $\sigma$ or ${\rm Dual}\circ\sigma$ is in $PGL(s,\mathbb C)\times PGL(n-s,\mathbb C)$. Then \begin{equation}
    {\rm Aut}(\mathcal M_{s,p,2p})=PGL(s,\mathbb C)\times PGL(n-s,\mathbb C)\rtimes \mathbb Z/2\mathbb Z.
\end{equation}

{\bf\noindent Case 5a ($4=2p=2s=n$).} $\mathcal M_{2,2,4}\cong\mathbb {CP}^3$, and hence  ${\rm Aut}(\mathcal M_{2,2,4})=PGL(4,\mathbb C)$.

\smallskip

{\bf\noindent Case 5b ($6\leq 2p=2s=n$).} Brion \cite{Br6} showed that 
\begin{equation}
  {\rm Aut}(\mathcal M_{p,p,2p})=\left(PGL(p,\mathbb C)\times PGL(p,\mathbb C)\right) \rtimes\mathbb Z/2\mathbb Z\rtimes\mathbb Z/2\mathbb Z.  
\end{equation}

\smallskip

We complete the proof of Proposition \ref{mauto}.\,\,\,\,$\endpf$

\section{The existence of K\"ahler-Einstein metrics} \label{ke}
In this section, we will study the existence of K\"ahler-Einstein metrics on $\mathcal T_{s,p,n}$ ($r\leq 2$), and $\mathcal M_{s,p,n}$. The main tool is Delcroix's criterion for spherical varieties \cite{De1,De2}. Readers are referred to the foundational work of \cite{LV}, \cite{Kn}, \cite{Br1, Br5, Br6}, \cite{BPa},  etc., on the general theory of spherical varieties.

\subsection{The existence of K\"ahler-Einstein metrics on \texorpdfstring{$\mathcal T_{s,p,n}$ }{jj}  of rank \texorpdfstring{$r\leq 2$ }{jj}}\label{ttest}
According to Theorem \ref{fano}, $\mathcal T_{s,p,n}$ is Fano if and only if the  rank $r\leq 2$. In this subsection, we will investigate the existence of K\"ahler-Einstein metrics on $\mathcal T_{s,p,n}$ with $r\leq 2$. 

For convenience, we denote by $G$ the group $GL(s,\mathbb C)\times GL(n-s,\mathbb C)$. Let $T=T_1\times T_2$  be  the maximal torus of $B$ defined by (\ref{maxt}).
Define coordinates of $T$ by
\begin{equation}
   g_1=\left(\begin{matrix}
    u_1&&&\\
    &u_2&&\\
    &&\ddots&\\
    &&&u_s\\
    \end{matrix}\right)\,\,{\rm and}\,\,
g_2=\left(\begin{matrix}
    v_1&&&\\
    &v_2&&\\
    &&\ddots&\\
    &&&v_{n-s}\\
    \end{matrix}\right)\,,
\end{equation}
where $g_1,g_2$ are elements of $T_1$, $T_2$ respectively,  and $u_1,\cdots,u_s,v_1,\cdots,v_{n-s}\in\mathbb C^*$;  we denote the elements of $T$ by $\big((u_1,\cdots,u_s),(v_1,\cdots,v_{n-s})\big)$ for simplicity.

Let $\mathfrak X(T)$ be the group of algebraic characters of $T$ and  $\mathfrak Y(T)$  the group of algebraic one parameter subgroups of $T$. Then $\mathfrak X(T)\cong\mathbb Z^n$ under the isomorphism
\begin{equation}\label{xt}
\begin{split}
\nu\,\,\,\,\,\,\,\,\,:\,\,\,\,\,\,\,\,\,\mathbb Z^n&\longrightarrow\mathfrak X(T)\\
\big((\alpha_1,\cdots,\alpha_s),(\beta_1,\cdots,\beta_{n-s})\big)&\mapsto \nu_{\alpha_1,\,\cdots,\,\alpha_s,\,\beta_1,\,\cdots,\,\beta_{n-s}}
\end{split}
\end{equation}
where $\nu_{\alpha_1,\,\cdots,\,\alpha_s,\,\beta_1,\,\cdots,\,\beta_{n-s}}$ is a holomorphic map defined by
\begin{equation}\label{nu}
\begin{split}
\nu_{\alpha_1,\,\cdots,\,\alpha_s,\,\beta_1,\,\cdots,\,\beta_{n-s}}:T&\longrightarrow \mathbb C^*\\
\big((x_1,\cdots,x_s),(y_1,\cdots,y_{n-s})\big)&\mapsto\left(\prod_{i=1}^s x_i^{\alpha_i}\right)\cdot\left(\prod_{j=1}^{n-s}y_j^{\beta_j}\right) \,\,\,.
\end{split}
\end{equation}
Similarly, $\mathfrak Y(T)\cong\mathbb Z^n$ under the isomorphism
\begin{equation}
\begin{split}
\mu\,\,\,\,\,\,\,\,\,:\,\,\,\,\,\,\,\,\,\mathbb Z^n&\longrightarrow\mathfrak Y(T)\\
\big((m_1,\cdots,m_s),(l_1,\cdots,l_{n-s})\big)&\mapsto \mu_{m_1,\,\cdots,\,m_s,\,l_1,\,\cdots,\,l_{n-s}}
\end{split}
\end{equation}
where $\mu_{m_1,\,\cdots,\,m_s,\,l_1,\,\cdots,\,l_{n-s}}$ is a holomorphic map defined by
\begin{equation}
\begin{split}
\mu_{m_1,\,\cdots,\,m_s,\,l_1,\,\cdots,\,l_{n-s}}:\mathbb C^*&\longrightarrow T\\
t&\mapsto \big((t^{m_1},\cdots,t^{m_s}),(t^{l_1},\cdots,t^{l_{n-s}})\big)\,\,.
\end{split}
\end{equation}
There is a natural pairing $\langle\cdot,\cdot\rangle$ 
\begin{equation}
\begin{split}
\langle\cdot,\cdot\rangle\,\,\,\,\,\,\,\,\,:\,\,\,\,\,\,\,\,\,\mathfrak X(T)\times \mathfrak Y(T)&\longrightarrow \mathbb Z\\
\big((\alpha_1,\cdots,\alpha_s),(\beta_1,\cdots,\beta_{n-s})\big)\times\big((m_1,\cdots,m_s),(l_1,\cdots,l_{n-s})\big)&\mapsto \sum_{i=1}^s\alpha_im_i+\sum_{j=1}^{n-s}\beta_jl_j\,\,.
\end{split}
\end{equation}

Let $\Phi\subset\mathfrak X(T)$ be the root system of $(G, T)$ and $\Phi^+$ the positive
roots determined by $B$. Let $\mathcal M\subset\mathfrak X(T)$
be the set of characters of $B$-semi-invariant functions in the function field $K(O)$   where  $O$ is the open and dense orbit of  $\mathcal T_{s,p,n}$ under the action of
$G$; denote by $\mathcal N:={\rm Hom}(\mathcal M,\mathbb Z)$  its $\mathbb Z$-dual.  We can determine the structure of $\mathcal M$ and $\mathcal N$ as follows.  Define a quasi projective variety $Z\subset G(p,n)$ in matrix representatives as follows. When  $p\leq n-s$ ($r=p$),
\begin{equation}
 Z:=\left\{ \left. \widetilde b:= \underbracedmatrixl{
  0 & \cdots&0 \\
  0 & \cdots&0 \\
  &\ddots& \\
  0 & \cdots&0 \\}{(s-p)\,\rm columns}\hspace{-.25in}\begin{matrix}
  &\hfill\tikzmark{c}\\
  \\
  \\
  \\
  &\hfill\tikzmark{d}
  \end{matrix}\,\,\,\, \begin{matrix}
  1 & & & \\
   & 1& &\\
  &&\ddots& \\
   & &&1 \\
  \end{matrix}
  \hspace{-.13in}\begin{matrix}
  &\hfill\tikzmark{a}\\
  \\
  \\
  \\
  &\hfill\tikzmark{b}  
  \end{matrix} \hspace{-.1in}
\begin{matrix}
  &\hfill\tikzmark{e}\\
  \\
  \\
  \\
  &\hfill\tikzmark{f}\end{matrix}\,\,\,\begin{matrix}
   \widetilde b_{11} & & & \\
   &\widetilde b_{22}& &\\
   & &\ddots &\\
   & &&\widetilde b_{rr} \\
  \end{matrix}
   \hspace{-.13in}\begin{matrix}
  &\hfill\tikzmark{g}\\
  \\
  \\
  \\
  &\hfill\tikzmark{h}
  \end{matrix}\hspace{-.12in}
  \underbracedmatrixr{
  0 & \cdots&0 \\
  0 & \cdots&0 \\
  &\ddots& \\
  0 & \cdots&0 \\}{(n-s-p) \,\rm columns}\right \vert_{} \begin{matrix}
  \widetilde b_{ii}\in\mathbb C^*,\\
  1\leq i\leq r
  \end{matrix}\right\};
  \tikz[remember picture,overlay]   \draw[dashed,dash pattern={on 4pt off 2pt}] ([xshift=0.5\tabcolsep,yshift=7pt]a.north) -- ([xshift=0.5\tabcolsep,yshift=-2pt]b.south);\tikz[remember picture,overlay]   \draw[dashed,dash pattern={on 4pt off 2pt}] ([xshift=0.5\tabcolsep,yshift=7pt]c.north) -- ([xshift=0.5\tabcolsep,yshift=-2pt]d.south);\tikz[remember picture,overlay]   \draw[dashed,dash pattern={on 4pt off 2pt}] ([xshift=0.5\tabcolsep,yshift=7pt]e.north) -- ([xshift=0.5\tabcolsep,yshift=-2pt]f.south);\tikz[remember picture,overlay]   \draw[dashed,dash pattern={on 4pt off 2pt}] ([xshift=0.5\tabcolsep,yshift=7pt]g.north) -- ([xshift=0.5\tabcolsep,yshift=-2pt]h.south);
\end{equation}
when $n-s<p$ ($r=n-s$), 
\begin{equation}
 Z:=\left\{ \left. \widetilde b:=
 \underbracedmatrixl{
  0  & \cdots&0 \\
  0 & \cdots&0  \\
    \vdots&\ddots&\vdots \\
  0 & \cdots&0  \\
0 & \cdots&0  \\
  \vdots&\ddots&\vdots \\
  0 & \cdots&0 \\}{(s-p)\,\rm columns}\hspace{-.23in}\begin{matrix}
  &\hfill\tikzmark{g}\\
  \\
  \\
  \\
  \\
  \\
  \\
  \\
  &\hfill\tikzmark{h}
  \end{matrix}\,\,
 \underbracedmatrix{
    1 & & & & & & \\
   & 1& && & & \\
  &&\ddots&& & &  \\
   & &&1& & &  \\
  & & &  & 1& &\\
  & & & &&\ddots& \\
  & & &  & &&1 \\}{p\,\rm columns}\hspace{-.17in}
\begin{matrix}
  &\hfill\tikzmark{e}\\
  \\
  \\
  \\
  \\
  \\
  \\
  \\
  &\hfill\tikzmark{f}\end{matrix}\hspace{-.1in}
  \begin{matrix}
  &\hfill\tikzmark{a}\\
  \\
  \\
  \\
  \\
  \\
  \\
  \\
  &\hfill\tikzmark{b}\end{matrix}\,
  \underbracedmatrixr{
   \widetilde b_{11} & & & \\
  & \widetilde b_{22}&&\\
  &&\ddots& \\
   & &&\widetilde b_{rr} \\
  0 &0 & \cdots&0 \\
  \vdots&\vdots&\ddots&\vdots \\
  0 &0 & \cdots&0 \\}{(n-s) \,\rm columns}\right \vert_{} \begin{matrix}
  \widetilde b_{ii}\in\mathbb C^*,\\
  1\leq i\leq r
  \end{matrix}\right\}.
  \tikz[remember picture,overlay]   \draw[dashed,dash pattern={on 4pt off 2pt}] ([xshift=0.5\tabcolsep,yshift=7pt]a.north) -- ([xshift=0.5\tabcolsep,yshift=-2pt]b.south);\tikz[remember picture,overlay]   \draw[dashed,dash pattern={on 4pt off 2pt}] ([xshift=0.5\tabcolsep,yshift=7pt]e.north) -- ([xshift=0.5\tabcolsep,yshift=-2pt]f.south);\tikz[remember picture,overlay]   \draw[dashed,dash pattern={on 4pt off 2pt}] ([xshift=0.5\tabcolsep,yshift=7pt]g.north) -- ([xshift=0.5\tabcolsep,yshift=-2pt]h.south);
\end{equation}
It is clear that $\mathcal K_{s,p,n}$ is well-defined on $Z$. By a slight abuse of notation, we denote a point $\widetilde a$ of $Z$ by its coordinates $(\widetilde a_{11},\cdots,\widetilde a_{rr})$, and identify $Z$ with its image $\mathcal K_{s,p,n}(Z)$ in $\mathcal T_{s,p,n}$.  Let $U$  be the unipotent radical of $B$; then the following is an open dense set in $\mathcal T_{s,p,n}$.
\begin{equation}
\big\{\mathfrak U\cdot\big(\mathcal K_{s,p,n}(z)\big)\big|\,\mathfrak U\in U\,\,{\rm and}\,\,z\in Z\big\}.    
\end{equation} 
Notice that each $B$-semi-invariant function $f\in K(O)$ is $U$-invariant, and 
\begin{equation}
    f\left(\,\widetilde b\,\right) = c \cdot \prod_{i=1}^{r} \left(\,\widetilde b_{ii}\,\right)^{\lambda_i},\,\,\,\,\,\widetilde b=(\widetilde b_{11},\cdots,\widetilde b_{rr})\in Z\,,
\end{equation}
where $\lambda_1,\cdots,\lambda_r$ are integers. Let $f_{1},f_{2},\cdots,f_{r}$ be $B$-semi-invariant functions such that 
\begin{equation}\label{bsemi}
\begin{split}
 &f_{1}\left(\,\widetilde b\,\right) = \left(\widetilde b_{11}\right)\,,\,\,\,\,f_{i}\left(\,\widetilde b\,\right) = \frac{\widetilde b_{ii}}{\widetilde b_{(i-1)(i-1)}}  \,,\, \,\,2\leq i\leq r\,\,\,\,\,{\rm for}\,\,\widetilde b=(\widetilde b_{11},\cdots,\widetilde b_{rr})\in Z.
\end{split}
\end{equation} 
According to (\ref{xt}), we can define a $\mathbb Z$-basis  $\{\chi_1,\cdots,\chi_r\}$ of $\mathcal M\subset\mathfrak X(T)$   by
\begin{equation}\label{chi1}
\begin{split}
&\chi_1:=\big((0,\cdots,0,\underset{\substack{\uparrow\\(s-p+1)^{th}}}{1,}\,0,\cdots,0),(\underset{\substack{\uparrow\\1^{st}}}{-1,}\,0,\cdots,0)\big)\,,\\
&\chi_i:=\big((0,\cdots,0,\underset{\substack{\uparrow\\(s-p+i-1)^{th}}}{-1,}\,1,0,\cdots,0),(0,\cdots,0,\underset{\substack{\uparrow\\{(i-1)}^{th}}}{1,}\,-1,0,\cdots,0)\big)\,\,{\rm for}\,\,\,2\leq i\leq r\,.
\end{split}
\end{equation}
It is clear that $f_{i}$, $1\leq i\leq r$, is a $B$-semi-invariant function of weight $\chi_i$ in the sense that $f(b^{-1}(x))=\chi_i(b)\cdot f(x)$, $b\in B$, $x\in O$.

We extend the basis $\{\chi_1,\cdots,\chi_r\}$ of $\mathcal M$ to a basis $\mathfrak B$ of $\mathfrak X(T)$ as follows. \begin{equation}
   \mathfrak B:=\{\chi_1,\cdots,\chi_r,\epsilon_1,\cdots, \epsilon_{r}, \tau_1,\cdots, \tau_{s-p}, \kappa_{1},\cdots,\kappa_{|n-s-p|}\}\,.
\end{equation} 
Here, when $p\leq n-s\leq s$,
\begin{equation}
\begin{split}
&\epsilon_{i}:=\big((0,\cdots,0,\underset{\substack{\uparrow\\(s-p+i)^{th}}}{1,}\,0,\cdots,0),(0,\cdots,0,\underset{\substack{\uparrow\\i^{th}}}{1,}\,0,\cdots,0)\big)\,,\,\,1\leq i\leq r\,,\\
&\tau_{i}:=\big((0,\cdots,0,\underset{\substack{\uparrow\\i^{th}}}{1,}\,0,\cdots,0),(0,\cdots,0)\big)\,,\,\,\,\,\,\,\,\,\,\,\,\,\,\,\,\,\,\,\,\,1\leq i\leq s-p\,,\\
&\kappa_{i}=\big((0,\cdots,0),(0,\cdots,0,\underset{\substack{\uparrow\\{(p+i)}^{th}}}{1,}\,0,\cdots,0)\big) \,,\,\,1\leq i\leq n-s-p\,;    
\end{split}  
\end{equation}
when $n-s<p<s$, $\epsilon_i$ and $\tau_i$ are defined the same as above, and
\begin{equation}
\begin{split}
&\kappa_{i}:=\big((0,\cdots,0,\underset{\substack{\uparrow\\(n-p+i)^{th}}}{1,}\,0,\cdots,0),(0,\cdots,0)\big) \,,\,\,1\leq i\leq s+p-n\,. 
\end{split}  
\end{equation}

Define a dual basis $ \{\gamma_1,\cdots,\gamma_r, \theta_1,\cdots, \theta_{r},\xi_1,\cdots, \xi_{s-p}, \delta_{1},\cdots,\delta_{|n-s-p|}\}$ of 
$\mathfrak Y(T)$ 
as follows. For $1\leq i\leq r$,
\begin{equation}
   \begin{split}
      \gamma_i\,\,\,:\,\,\,\, &\mathfrak X(T)\rightarrow\mathbb Z\,,\,\,\, \\
       & \chi_i\mapsto 1\\
      &\{\chi_1,\cdots,\chi_r,\epsilon_1,\cdots, \epsilon_{r}, \tau_1,\cdots, \tau_{s-p}, \kappa_{1},\cdots,\kappa_{|n-s-p|}\}\backslash\{\chi_i\}\mapsto 0\\
   \end{split}\,;
\end{equation}
for $1\leq i\leq r$, 
\begin{equation}
   \begin{split}
 \theta_i\,\,\,:\,\,\,\, &\mathfrak X(T)\rightarrow\mathbb Z\,,\,\,\, \\
      & \epsilon_i\mapsto 1\\
      & \{\chi_1,\cdots,\chi_r,\epsilon_1,\cdots, \epsilon_{r}, \tau_1,\cdots, \tau_{s-p}, \kappa_{1},\cdots,\kappa_{|n-s-p|}\}\backslash\{\epsilon_i\}\mapsto 0\\
   \end{split}\,;
\end{equation}
for $1\leq i\leq s-p$, 
\begin{equation}
   \begin{split}
 \xi_i\,\,\,:\,\,\,\, &\mathfrak X(T)\rightarrow\mathbb Z\,,\,\,\, \\
      &\tau_i\mapsto 1\\
      & \{\chi_1,\cdots,\chi_r,\epsilon_1,\cdots, \epsilon_{r}, \tau_1,\cdots, \tau_{s-p}, \kappa_{1},\cdots,\kappa_{|n-s-p|}\}\backslash\{\tau_i\}\mapsto 0\\
   \end{split}\,;
\end{equation}
for $1\leq i\leq |n-s-p|$,
\begin{equation}
   \begin{split}
 \delta_{i}\,\,\,:\,\,\,\, &\mathfrak X(T)\rightarrow\mathbb Z\,,\,\,\, \\
      &\kappa_{i}\mapsto 1\\
      & \{\chi_1,\cdots,\chi_r,\epsilon_1,\cdots, \epsilon_{r}, \tau_1,\cdots, \tau_{s-p}, \kappa_{1},\cdots,\kappa_{|n-s-p|}\}\backslash\{\kappa_{i}\}\mapsto 0\\
   \end{split}\,.
\end{equation}
In what follows, by a slight abuse notation, we also denote by $\gamma_i$, $1\leq i\leq r$, its image in $\mathcal N$ under the natural quotient map $\pi:\mathfrak Y(T)\otimes_{\mathbb Z}\mathbb Q\rightarrow \mathcal N\otimes_{\mathbb Z}\mathbb Q$. We note that $\left\{\gamma_1,\cdots,\gamma_r\right\}$ is the dual basis of $\mathcal N$ to $\{\chi_1,\cdots,\chi_r\}$ of $\mathcal M$.

The valuation cone
$\mathcal V$  with respect to $B$ can be defined as the set of the  elements of the vector space $\mathcal N\otimes_{\mathbb Z}\mathbb Q$ induced by $G$-invariant valuations on the function field $K(O)$. 
Define a map $\iota: \mathcal V\rightarrow \mathcal N\otimes_{\mathbb Z}\mathbb Q$ such that $\langle\iota(\nu), \chi_i\rangle:=\nu(f_{i})$, $1\leq i\leq r$; denote the valuations associated with $D_{i}^-$ and $D_{i}^+$ by $v_{D_{i}^-}$ and $v_{D_{i}^+}$  respectively.
Similarly, we define a  map $\rho: \left\{B_0,\cdots,B_r\right\}\rightarrow \mathcal N$ such that \begin{equation}
    \langle\rho(B_j),\chi_i\rangle=v_{B_j}(f_{i})\,,\,\,\,1\leq i\leq r\,,\,\,0\leq j\leq r\,,
\end{equation}
where $v_{B_j}$ is the discrete valuation associated to $B_j$.

\begin{lemma}\label{bfor} For $r\geq 2$,
\begin{equation}\label{bdn}
    \rho(B_j)=\left\{\begin{array}{ll}
    -\gamma_{1}+\gamma_{2}  & j=0\\
    \gamma_{j}-2\gamma_{j+1}+\gamma_{j+2} & 1\leq j\leq r-2\,\,\,\,\,\,\, \\
  
    \gamma_{r-1}-2\gamma_{r}  & j=r-1\\
    \gamma_{r}  & j=r\\
    \end{array}\right.;
\end{equation}
for $r=1$,
\begin{equation}
    \rho(B_j)=\left\{\begin{array}{ll}
    -\gamma_{1}& j=0\\
    \gamma_{1}  & j=1\\
    \end{array}\right..
\end{equation}
\end{lemma}
{\noindent\bf Proof of Lemma \ref{bfor}.} Checking  on a certain open set near $Z$, we can conclude that
\begin{equation}
   f_{1}=\frac{R_{s,p,n}^*\left(P_{I_{1}}\right)}{R_{s,p,n}^*\left(P_{I_{0}}\right)}\,\,\,\,{\rm and}\,\,\,\,  f_{k}=\frac{R_{s,p,n}^*\left(P_{I_{k}}\right)\cdot R_{s,p,n}^*\left(P_{I_{k-2}}\right)}{R_{s,p,n}^*\left(P_{I_{k-1}}\right)\cdot R_{s,p,n}^*\left(P_{I_{k-1}}\right)}\,\,{\rm for}\,\,2\leq k\leq r,
\end{equation}
where $P_{I_{k}}$ is the Pl\"ucker coordinate function associated with the index $I_{k}$ defined by (\ref{I_k}).

Then by (\ref{bst}), (\ref{bstt0}) and (\ref{bsttr}), we can conclude that the principal divisors $(f_k)$ takes the following form. (Notice that here the complicity is due to the fact that   $B_0=D^+_r$  when $p=s$ and  $B_r=D^-_r$ when $p=n-s$.)  When  $p<n-s\leq s$,
\begin{equation}\label{vo}
\begin{split}
&(f_{1})=B_{1}-B_{0}+D^-_1-\sum_{i=1}^{r}D^+_i\,;\\
&(f_{k})=B_{k}-2B_{k-1}+B_{k-2}+D^-_k+D^+_{r+2-k}\,\,{\rm for}\,\,2\leq k\leq r\,.\\
\end{split}
\end{equation}
When $1=n-s=p<s$,
\begin{equation}\label{vo21}
\begin{split}
&(f_{1})=B_{1}-B_{0}-D^+_1=-B_{0}+D^-_1-D^+_1\,;\\
\end{split}
\end{equation}
when $2\leq n-s=p<s$,
\begin{equation}\label{vo22}
\begin{split}
&(f_{1})=B_{1}-B_{0}+D^-_1-\sum_{i=1}^{r}D^+_i\,,\\
&(f_{k})=B_{k}-2B_{k-1}+B_{k-2}+D^-_k+D^+_{r+2-k}\,\,{\rm for}\,\,2\leq k\leq r-1\,,\\
&(f_{r})=B_{r}-2B_{r-1}+B_{r-2}+D^+_{2}=D^-_r-2B_{r-1}+B_{r-2}+D^+_{2}\,.\\
\end{split}
\end{equation}
When $1=n-s=p=s$,
\begin{equation}\label{vo31}
\begin{split}
&(f_{1})=B_{1}-B_{0}=D^-_1-D^+_1\,;\\
\end{split}
\end{equation}
when $2=n-s=p=s$,
\begin{equation}\label{vo32}
\begin{split}
&(f_{1})=B_{1}-B_{0}+D^-_1-D^+_1=B_{1}-D^+_2+D^-_1-D^+_1\,,\\
&(f_{2})=B_{2}-2B_{1}+B_{0}=D_{2}^--2B_{1}+D_2^+\,;\\
\end{split}
\end{equation}
when $3\leq n-s=p=s$,
\begin{equation}\label{vo3}
\begin{split}
&(f_{1})=B_{1}-B_{0}+D^-_1-\sum_{i=1}^{r-1}D^+_i=B_{1}-D^+_r+D^-_1-\sum_{i=1}^{r-1}D^+_i\,,\\
&(f_{2})=B_{2}-2B_{1}+B_0+D^-_2=B_{2}-2B_{1}+D^+_r+D^-_2\,,\\
&(f_{k})=B_{k}-2B_{k-1}+B_{k-2}+D^-_k+D^+_{r+2-k}\,\,{\rm for}\,\,3\leq k\leq r-1\,,\\
&(f_{r})=B_{r}-2B_{r-1}+B_{r-2}+D^+_{2}=D^-_r-2B_{r-1}+B_{r-2}+D^+_{2}\,.\\
\end{split}
\end{equation}

We complete the proof of Lemma \ref{bfor}.\,\,\,$\endpf$
\medskip

We describe  the valuation cone $\mathcal V$ by
\begin{lemma}\label{val}
$v_{D_{i}^+}$ and $v_{D_{i}^-}$ take the following from over the  basis $\left\{\gamma_1,\cdots,\gamma_{r}\right\}$ of $\mathcal N$.
\begin{equation}\label{dualc}
\begin{split}
&v_{D_{1}^+}=-\gamma_1\,\,\,{\rm and} \,\,\,v_{D_{i}^+}=-\gamma_1+\gamma_{r+2-i}\,\,\,{\rm for}\,\,\,2\leq i\leq  r\,\,;\,\,\,\,\,\\
&v_{D_{i}^-}=\gamma_i\,,\,\,\,\,1\leq i\leq r\,.
\end{split}
\end{equation}
$\mathcal V\subset \mathcal N\otimes_{\mathbb Z} \mathbb Q$ is the cone generated by $v_{D_{1}^-}$, $v_{D_{2}^-}$, $\cdots$, $v_{D_{r}^-}$, $v_{D_{1}^+}$ over $\mathbb Z$. 
\end{lemma}
{\noindent\bf Proof of  Lemma \ref{val}.} 
Let $v\in\mathcal V$ be an $G$-invariant discrete valuation. Then its center $\mathcal Z_{v}$ is a $G$-stable closed subvariety of  $\mathcal T_{s,p,n}$ (see \cite{Kn}). By Proposition \ref{gwond},  each $G$-stable closed subvariety of  $\mathcal T_{s,p,n}$ is a certain intersection of  $D_{1}^-$, $D_{2}^-$, $\cdots$, $D_{r}^-$,  $D_{1}^+$, $D_{2}^+$, $\cdots$, $D_{r}^+$. We assume that 
\begin{equation} \mathcal Z_{v}:=\left(\,\bigcap_{i=1}^{m_{-}}D_{k_i}^-\right)\scaleobj{1.1}{\bigcap} \left(\,\bigcap_{i=1}^{m_{+}}D_{j_i}^+\right)\,, \end{equation}
where $1\leq k_1< k_2<\cdots<k_{m_{-}}\leq r$ and  $1\leq j_1< j_2<\cdots<j_{m_{+}}\leq r$.
Recalling Lemma \ref{orbits} and the fact that $D_{1}^-+\cdots+D_{r}^-+D_{1}^++\cdots+D_{r}^+$ is a simple normal crossing divisor, we can conclude by (\ref{vo}) that
\begin{equation}
    v=\sum_{i=1}^{m_{-}}v_{D_{k_i}^-}+\sum_{i=1}^{m_{+}}v_{D_{j_i}^+}\,.
\end{equation}
(\ref{dualc}) follows from (\ref{vo}) as well. 

Noticing that $v_{D_{i}^+}=v_{D_{1}^+}+v_{D_{r+2-i}^-}$ for $2\leq i\leq r$, we complete the proof of Lemma \ref{val}.\,\,\,$\endpf$
\medskip

Following \cite{AB} we define a polytope  $Q_{\mathcal T_{s,p,n}}\subset\mathcal N\otimes\mathbb Q$ as a convex hull of certain points as follows. When $p<n-s\leq s$ ($r=p$),
\begin{equation}
  Q_{\mathcal T_{s,p,n}}:={\rm conv} \left\{\frac{\rho(B_0)}{s-p+1}, \frac{\rho(B_1)}{2},\cdots,\frac{\rho(B_{p-1})}{2}, \frac{\rho(B_p)}{n-s-p+1}, v_{D^-_1}, \cdots,v_{D^-_r},v_{D^+_1},\cdots,v_{D^+_r}\right\}\,;
\end{equation}
when $n-s=p<s$ ($r=p$),
\begin{equation}
  Q_{\mathcal T_{s,p,n}}:={\rm conv} \left\{\frac{\rho(B_0)}{s-p+1}, \frac{\rho(B_1)}{2},\cdots,\frac{\rho(B_{p-1})}{2}, v_{D^-_1}, \cdots,v_{D^-_r},v_{D^+_1},\cdots,v_{D^+_r}\right\}\,;
\end{equation}
when $n-s<p<s$ ($r=n-s$),
\begin{equation}
  Q_{\mathcal T_{s,p,n}}:={\rm conv} \left\{\frac{\rho(B_0)}{s-p+1}, \frac{\rho(B_1)}{2},\cdots,\frac{\rho(B_{r-1})}{2}, \frac{\rho(B_r)}{p-r+1}, v_{D^-_1}, \cdots,v_{D^-_r},v_{D^+_1},\cdots,v_{D^+_r}\right\}\,;
\end{equation}
when $n-s=p=s$ ($r=p$),
\begin{equation}
  Q_{\mathcal T_{s,p,n}}:={\rm conv} \left\{\frac{\rho(B_1)}{2},\cdots,\frac{\rho(B_{r-1})}{2}, v_{D^-_1}, \cdots,v_{D^-_r},v_{D^+_1},\cdots,v_{D^+_r}\right\}\,.
\end{equation}
Notice that the denominators of $B_i$ in the above come from the coefficients of 
$B_i$ in the canonical bundle formulas in Lemma \ref{wk}.

Denote by $Q^*_{\mathcal T_{s,p,n}}$ the dual polytope to $Q_{\mathcal T_{s,p,n}}\subset\mathcal N\otimes_{\mathbb Z}\mathbb Q$, that is,
\begin{equation}
Q^*_{\mathcal T_{s,p,n}}:=\left\{u\in\left(\mathcal N\otimes_{\mathbb Z}\mathbb Q\right)^*=\mathcal M\otimes_{\mathbb Z}\mathbb Q : \,\langle u, v\rangle\leq 1\,\,{\rm\,\, for\,\, every\,\,}
v\in Q_{\mathcal T_{s,p,n}}\right\}\,.   
\end{equation}
Denote by $\Delta_{\mathcal T_{s,p,n}}^+\subset\mathfrak X(T)\otimes_{\mathbb Z}\mathbb R$ the moment polytope of ${\mathcal T_{s,p,n}}$  with respect to $B$. By Proposition 3.3 in \cite{Br2}, we have that
\begin{equation}\label{tmp}
    \Delta_{\mathcal T_{s,p,n}}^+=2\rho_P+Q^*_{\mathcal T_{s,p,n}}\,.
\end{equation}
Here $P$ is the stabilizer of the open orbit of $B$ in $\mathcal T_{s,p,n}$, and $2\rho_P$ is the weight of a $B$-semi-invariant section of  $-K_{\mathcal T_{s,p,n}}$.

In what follows, we will compute the weight  $2\rho_P$ explicitly in terms of the roots associated with the unipotent radical of $P$.
Recall that, when $p\leq n-s\leq s$ $(r=p)$,
\begin{equation}\label{paro1}
P=\left\{\left. \left(
\begin{matrix}
V_{1}&0&0&0\\
W_{1}&U_{1}&0&0\\
0&0&U_2&W_2\\
0&0&0&V_2\\
\end{matrix}\right)\right\vert_{} \footnotesize 
\begin{matrix}
V_1\in GL(s-p,\mathbb C) \, ,V_2\in   GL(n-s-p,\mathbb C)\,;\,\,\,\\
U_1\in GL(p,\mathbb C)\,\,{\rm is\,\, lower\,\, triangular \,\,};\\
U_2 \in GL(p,\mathbb C)\,\,{\rm is\,\, upper\,\, triangular \,\,};\\
W_{1}\,\,{\rm \,\,is\,\,a}\,\,p\times (s-p)\,\,{\rm matrix\,\,};\\
W_{2}\,\,{\rm \,\,is\,\,a}\,\,p\times (n-s-p)\,\,{\rm matrix\,\,}\\
\end{matrix}
\right\}\,;
\end{equation}
when $n-s<p<s$ $(r=n-s)$,
\begin{equation}\label{paro2}
P=\left\{\left. \left(
\begin{matrix}
V_{11}&0&0&0\\
W_{11}&U_{1}&0&0\\
W_{12}&W_{13}&V_{12}&0\\
0&0&0&U_2\\
\end{matrix}\right)\right\vert_{} \footnotesize 
\begin{matrix}
V_{11}\in GL(s-p,\mathbb C) \, ,V_{12}\in   GL(p-r,\mathbb C)\,;\,\,\,\\
U_1\in GL(r,\mathbb C)\,\,{\rm is\,\, lower\,\, triangular \,\,};\\
U_2 \in GL(r,\mathbb C)\,\,{\rm is\,\, upper\,\, triangular \,\,};\\
W_{11}, W_{12}, W_{13}\,\,{\rm \,\,are\,\,matrices\,\,of\,\,sizes}\,\,r\times (s-p), \\
(p-r)\times (s-p), (p-r)\times r,\,\,{\rm respectively}\\
\end{matrix}
\right\}\,.
\end{equation}
The  unipotent radical $P^{u}$  of $P$ takes the following form. When $p\leq n-s$,
\begin{equation}\label{paro3}
P^u=\left\{\left. \left(
\begin{matrix}
I_{(s-p)\times(s-p)}&0&0&0\\
W_{1}&U^u_{1}&0&0\\
0&0&U^u_2&W_2\\
0&0&0&I_{(n-s-p)\times(n-s-p)}\\
\end{matrix}\right)\right\vert_{} \footnotesize 
\begin{matrix}
U^u_1 \,\,({\rm resp.}\,\,U^u_2)\in GL(p,\mathbb C)\,\,\\
{\rm is\,\, lower}\,\,({\rm resp,\,\, upper})\,\, {\rm triangular }\\
{\rm with\,\,diagonal\,\, entries\,\, }1\\
\end{matrix}
\right\}\,;
\end{equation}
when $n-s<p$,
\begin{equation}\label{paro4}
P^u=\left\{\left. \left(
\begin{matrix}
I_{(s-p)\times(s-p)}&0&0&0\\
W_{11}&U^u_{1}&0&0\\
W_{12}&W_{13}&I_{(p-r)\times(p-r)}&0\\
0&0&0&U^u_2\\
\end{matrix}\right)\right\vert_{} \footnotesize 
\begin{matrix}
U^u_1 \,\,({\rm resp.}\,\,U^u_2)\in GL(r,\mathbb C)\,\,\\
({\rm is\,\, lower}\,\,{\rm resp,\,\, upper})\,\, {\rm triangular }\\
{\rm with\,\,diagonal\,\, entries\,\, }1\\
\end{matrix}
\right\}\,.
\end{equation}
Then, 
\begin{equation}\label{canwe}
    2\rho_P=\sum_{\alpha\in\Phi_{P^u}}\alpha\,,
\end{equation}
where $\Phi_{P^u}$ is the set of roots of $P^u$.

\begin{lemma}\label{rhop}
When $p\leq n-s$,
\begin{equation}
\begin{split}
    2\rho_P&=\sum_{i=1}^{p}\left(s-\frac{n}{2}+i-1\right)(p+1-i)\cdot \chi_i+ \left(\frac{n}{2}-p\right)\sum_{i=1}^{p} \epsilon_i-p\sum_{i=1}^{s-p} \tau_i- p\sum_{i=1}^{n-s-p} \kappa_i\,;
\end{split}
\end{equation}
when $n-s\leq p$,
\begin{equation}
\begin{split}
    2\rho_P&=\sum_{i=1}^{n-s}\left(\frac{n}{2}-p+i-1\right)(n-s+1-i)\cdot \chi_i+ \left(\frac{n}{2}-p\right)\sum_{i=1}^{n-s}\epsilon_i-p\sum_{i=1}^{s-p} \tau_i+(n-p)\sum_{i=1}^{s+p-n} \kappa_i\,.
\end{split}
\end{equation}
\end{lemma}
{\bf\noindent Proof of Lemma \ref{rhop}. }  According to (\ref{xt}) the  simple roots of $\Phi^+$ are given by
\begin{equation}\label{sr}
\begin{split}
&e_{i}:=\big((0,\cdots,0,\underset{\substack{\uparrow\\(s-i)^{th}}}{-1,}\,1,0,\cdots,0),(0,0,\cdots,0)\big)\,,\,\,1\leq i\leq s-1\,,\\
&\widetilde e_{i}:=\big((0,0,\cdots,0),(0,\cdots,\underset{\substack{\uparrow\\i^{th}}}{1,}-1,0,\cdots,0)\big)\,,\,\,\,\,\,\,\,\,\,\,\,\,1\leq i\leq n-s-1\,.\\
\end{split}  
\end{equation}
When $r=p$,  $\Phi_{P^u}$ consists of the following roots.
\begin{equation}
\begin{split}
    &e_{ij}:=\sum_{k=i}^{j-1}e_k\,\,{\rm\,\,for\,\,all\,}\,1\leq i\leq p\,\,{\rm and\,\,}i+1\leq j\leq s\,;\\
    &\widetilde e_{ij}:=\sum_{k=i}^{j-1}\widetilde e_k\,\,{\rm\,\,for\,\,all\,}\,1\leq i\leq p\,\,{\rm and\,\,}i+1\leq j\leq n-s\,.\\
\end{split}
\end{equation}
When $r=n-s$,  $\Phi_{P^u}$ consists of the following roots.
\begin{equation}
\begin{split}
    &e_{ij}:=\sum_{k=i}^{j-1}e_k\,\,{\rm\,\,for\,\,all\,}\,1\leq i\leq p\,\,{\rm and\,\,}\max\{i,s+p-n\}+1\leq j\leq s\,;\\
    &\widetilde e_{ij}:=\sum_{k=i}^{j-1}\widetilde e_k\,\,{\rm\,\,for\,\,all\,}\,1\leq i\leq r-1\,\,{\rm and\,\,}i+1\leq j\leq r\,.\\
\end{split}
\end{equation}

Then direct computation yields Lemma \ref{rhop}. \,\,\,\,$\endpf$
\medskip

Recall that the Duistermaat–Heckman measure can be given by \begin{equation}\label{duiheck}
    \prod_{\alpha\in\Phi_{P^u}}\kappa(\alpha,p)\,dp
\end{equation}
where $\kappa$ is the Killing form and $dp$ is the Lebesgue measure on $\mathfrak X(T)\otimes_{\mathbb Z}\mathbb R$. %
We denote by ${\rm bar}_{DH}(\Delta_{\mathcal T_{s,p,n}}^+)$ the barycenter of the polytope
$\Delta_{\mathcal T_{s,p,n}}^+$ with respect to the Duistermaat–Heckman measure. 


\begin{lemma}\label{dhm}
Write
\begin{equation}
  p=\sum_{i=1}^{r}x_i\cdot \chi_i + \sum_{i=1}^{r}y_i\cdot \epsilon_i+\sum_{i=1}^{s-p}z_i\cdot \tau_i + \sum_{i=1}^{|n-s-p|}w_i\cdot \kappa_i\,,\,\, x_i,y_i,z_i,w_i\in\mathbb C\,.
\end{equation} Make the convention that $x_{r+1}=0$.  Then up to a constant  depending only on the dimension $n$, the Duistermaat–Heckman measure takes the following form. When $p\leq n-s$,
\begin{equation}
\begin{split}
&\prod_{\alpha\in\Phi_{P^u}}\kappa(\alpha,p)\,dp=\left(\prod_{1\leq i<j\leq r}\big((x_{i+1}-x_{i})+y_{i}-(x_{j+1}-x_{j})-y_{j}\big)\right)\\
&\cdot\left(\prod_{1\leq i<j\leq r}\big((x_{i+1}-x_{i})+y_{j}-(x_{j+1}-x_{j})-y_{i}\big)\right)\cdot\left(\prod_{i=1}^r\prod_{j=1}^{n-s-p}\big((x_{i+1}-x_{i})+y_i-w_j\big)\right)\\
&\cdot\left(\prod_{i=1}^r\prod_{j=1}^{s-p}\big(y_{i}-(x_{i+1}-x_{i})-z_j\big)\right)\cdot\left(\prod_{i=1}^rdx_i\right)\cdot\left(\prod_{i=1}^rdy_i\right)\cdot\left(\prod_{i=1}^{s-p}dz_i\right)\cdot\left(\prod_{i=1}^{n-s-p}dw_i\right)\,.\\
\end{split}
\end{equation}
When $p\geq n-s$,
\begin{equation}
\begin{split}
&\prod_{\alpha\in\Phi_{P^u}}\kappa(\alpha,p)\,dp=\left(\prod_{1\leq i<j\leq r}\big((x_{i+1}-x_{i})+y_{i}-(x_{j+1}-x_{j})-y_{j}\big)\right)\\
&\cdot\left(\prod_{1\leq i<j\leq r}\big((x_{i+1}-x_{i})+y_{j}-(x_{j+1}-x_{j})-y_{i}\big)\right)\cdot\left(\prod_{i=1}^r\prod_{j=1}^{s+p-n}\big((x_{i+1}-x_{i})-y_i+w_j\big)\right)\\
&\,\,\,\,\,\,\cdot\left(\prod_{i=1}^r\prod_{j=1}^{s-p}\big(y_{i}-(x_{i+1}-x_{i})-z_j\big)\right)\cdot\left(\prod_{i=1}^{s-p}\prod_{j=1}^{s+p-n}(w_{j}-z_i)\right)\\
&\,\,\,\,\,\,\,\,\,\,\,\cdot\left(\prod_{i=1}^rdx_i\right)\cdot\left(\prod_{i=1}^rdy_i\right)\cdot\left(\prod_{i=1}^{s-p}dz_i\right)\cdot\left(\prod_{i=1}^{s+p-n}dw_i\right)\,.\\
\end{split}
\end{equation}
\end{lemma}
{\bf\noindent Proof of Lemma \ref{dhm}. }
Notice that up to a constant, the Killing form on $\mathfrak X(T)$ is the Euclidean metric with respect to the standard basis given by (\ref{xt}).  Computation
yields Lemma \ref{dhm}.\,\,\,\,$\endpf$
\medskip

Denote by $\pi^{-1}(-\mathcal V)\subset \mathfrak Y(T)\otimes\mathbb R$ the closure of the inverse image by $\pi$ of the opposite of the valuation cone $\mathcal V$. Then $\pi^{-1}(-\mathcal V)$ is the cone generated by 
\begin{equation}
     \left \{-v_{D_{1}^-}, \cdots, -v_{D_{r}^-}, -v_{D_{1}^+}, \cdots, -v_{D_{r}^+}, \pm\theta_1,\cdots, \pm\theta_{r},\pm\xi_1,\cdots, \pm\xi_{s-p}, \pm\delta_{1},\cdots,\pm\delta_{|n-s-p|}\right\}\,.
\end{equation}
Denote by
$\Xi_{\mathcal T_{s,p,n}}\subset\mathfrak X(T)\otimes \mathbb R$ the dual cone to  $\pi^{-1}(-\mathcal V)$.

Recall the following criterion for the existence of K\"ahler-Einstein metrics.
\begin{theorem}[\cite{De1, De2}]\label{dcr}
    Let $X$ be a complete, spherical Fano  manifold. The following are equivalent.
\begin{enumerate}
    \item There exists a K\"ahler-Einstein metric on $X$.
    \item The barycenter ${\rm bar}_{DH}(X)$ is in the relative interior of the cone $2\rho_P+\Xi_{X}$.
\end{enumerate} 
\end{theorem}
{\bf\noindent Proof of Proposition \ref{tke}. }
The case when $p=1$ has been studied in \cite{De3}. There are K\"ahler-Einstein metrics on $\mathcal T_{s,1,n}^+$ if and only if $n=2s$.

Assume that $n-s=1$. It is easy to verify that $\Xi_{\mathcal T_{n-1,p,n}}$ consists of the origin;
\begin{equation}
Q^*_{\mathcal T_{n-1,p,n}}=\left\{x_1\cdot\chi_1\big|\,-1\leq x_1\leq 1\,\right\}\,.  \end{equation}
By Lemma \ref{rhop}, we have that
\begin{equation}
\begin{split}
    2\rho_P&=\left(\frac{n}{2}-p\right)\chi_1+ \left(\frac{n}{2}-p\right)\epsilon_1-p\sum_{i=1}^{n-p-1} \tau_i+(n-p)\sum_{i=1}^{p-1} \kappa_i\,;
\end{split}
\end{equation}
hence by (\ref{tmp}),
\begin{equation}
    \Delta_{\mathcal T_{n-1,p,n}}^+=\left\{\left.\left(x_1+\frac{n}{2}-p\right)\chi_1+ \left(\frac{n}{2}-p\right)\epsilon_1-p\sum_{i=1}^{n-p-1} \tau_i+(n-p)\sum_{i=1}^{p-1} \kappa_i\right\vert_{}\,-1\leq x_1\leq 1\,\right\}\,.
\end{equation}
Let $\widetilde x_1=x_1+\frac{n}{2}-p$. Then by Lemma \ref{dhm}, $\rho\left(\widetilde x_1\right)d\widetilde x_1$ the restriction of the Duistermaat–Heckman measure on $\Delta_{\mathcal T_{n-1,p,n}}^+$ takes the following form  up to a constant.
\begin{equation}
\begin{split}
&\rho\left(\widetilde x_1\right)d\widetilde x_1=\left(\frac{n}{2}-\widetilde x_1\right)^{p-1}\left(\frac{n}{2}+\widetilde x_1\right)^{n-p-1} d\widetilde x_1\,.\\
\end{split}
\end{equation}

By Theorem \ref{dcr}, there are K\"ahler-Einstein metrics on $\mathcal T_{n-1,p,n}$ if and only if barycenter ${\rm bar}_{DH}(\Delta_{\mathcal T_{n-1,p,n}}^+)=2\rho_P$, that is, 
\begin{equation}
\begin{split}
0=\displaystyle\int_{\frac{n}{2}-p-1}^{\frac{n}{2}-p+1}\left(\widetilde x_1+p-\frac{n}{2}\right)\cdot\rho\left(\widetilde x_1\right)d\widetilde x_1=\displaystyle\int_{\frac{n}{2}-p-1}^{\frac{n}{2}-p+1}\left(\widetilde x_1+p-\frac{n}{2}\right)\left(\frac{n}{2}-\widetilde x_1\right)^{p-1}\left(\frac{n}{2}+\widetilde x_1\right)^{n-p-1}\,.\\
\end{split}
\end{equation}
Notice that the above integral is exactly the one appearing in \cite{De3}. Therefore, similarly to \cite{De3} we can conclude that 
there exists a K\"ahler-Einstein metric on $\mathcal T_{n-1,p,n}$ if and only if $n=2p$.

Next we assume that $2=p\leq n-s$. Similarly, we can show that
\begin{equation}
   \Xi_{\mathcal T_{s,2,n}}=\left\{(x_1,x_2)\in\mathbb{R}^{2}\,\big|\,x_1=0\,,\,\,x_2\geq 0\, \right\}\,,
\end{equation}
and
\begin{equation}
\begin{split}
    2\rho_P&=\sum_{i=1}^{2}\left(s-\frac{n}{2}+i-1\right)(3-i)\cdot \chi_i+ \left(\frac{n}{2}-2\right)\sum_{i=1}^{2}\epsilon_i-2\sum_{i=1}^{s-2} \tau_i-2\sum_{i=1}^{n-s-2} \kappa_i\,.
\end{split}
\end{equation}
The moment polytope $\Delta_{\mathcal T_{s,2,n}}^+$ consists of the points $A(x_1,x_2)$ with the following form
\begin{equation}
\begin{split}
    A(x_1,x_2):=&\left(x_1+2\left(s-\frac{n}{2}\right)\right)\chi_1+\left(x_2+\left(s-\frac{n}{2}+1\right)\right)\chi_2\\
    &\hspace{-.3in}+ \left(\frac{n}{2}-2\right)(\epsilon_1+\epsilon_2)-2\sum_{i=1}^{s-2} \tau_i-2\sum_{i=1}^{n-s-2} \kappa_i\,,
\end{split}
\end{equation}
where $(x_1,x_2)$ is in the region $\Delta(\cong Q^*_{\mathcal T_{s,2,n}})$ defined as follows.
\begin{equation}
    \Delta:=\big\{(x_1,x_2)\in\mathbb R^2\,\big|\,-1\leq x_1\leq 1\,,\,\,x_2\leq 1\,,\,\,-x_1+x_2\leq 1\,,\,\,x_1-2x_2\leq 2\,\big\}\,.
\end{equation}
$\rho\left(x_1,x_2\right)dx_1dx_2$ the restriction of the Duistermaat–Heckman measure on $\Delta_{\mathcal T_{s,2,n}}^+$ takes the following form  up to a constant. 
\begin{equation}\label{crt0}
\begin{split}
&\rho\left(x_1,x_2\right)dx_1dx_2=(2x_2-x_1+2)^2(x_2-x_1+n-s+1)^{n-s-2}(-x_2+n-s-1)^{n-s-2}\\
&\hspace{.5in}\cdot(x_1-x_2+s-1)^{s-2}(x_2+s+1)^{s-2}dx_1dx_2\,.\\
\end{split}
\end{equation}

By Theorem \ref{dcr}, there exists a K\"ahler-Einstein metric on $\mathcal T_{s,2,n}$ if and only if barycenter ${\rm bar}_{DH}(\Delta_{\mathcal T_{s,2,n}}^+)=2\rho_P$, that is, 
\begin{equation}
\begin{split}
0=\displaystyle\int_{\Delta}x_1\cdot\rho\left(x_1,x_2\right)dx_1dx_2\,\,\,\,\,{\rm and}\,\,\,\,\,0<\displaystyle\int_{\Delta}x_2\cdot\rho\left(x_1,x_2\right)dx_1dx_2\,.\\
\end{split}
\end{equation}
By Lemma \ref{numcal} in Appendix \ref{section:nc}, this is equivalent to that $n=2s$.

Similarly, when $2=n-s\leq p$ there exists a K\"ahler-Einstein metric on $\mathcal T_{n-2,p,n}$ if and only if 
\begin{equation}
\begin{split}
0=\displaystyle\int_{\Delta}x_1\cdot\widetilde\rho\left(x_1,x_2\right)dx_1dx_2\,\,\,\,\,{\rm and}\,\,\,\,\,0<\displaystyle\int_{\Delta}x_2\cdot\widetilde\rho\left(x_1,x_2\right)dx_1dx_2\,.\\
\end{split}
\end{equation}
where
\begin{equation}
\begin{split}
&\widetilde\rho\left(x_1,x_2\right)dx_1dx_2=(2x_2-x_1+2)^2(x_2-x_1+p+1)^{p-2}(-x_2+p-1)^{p-2}\\
&\hspace{.5in}\cdot(x_1-x_2+n-p-1)^{n-p-2}(x_2+n-p+1)^{n-p-2}dx_1dx_2\,.\\
\end{split}
\end{equation}
Setting $s=n-p$ and $p=n-s$ in Lemma \ref{numcal}, we can conclude that this is equivalent to $n=2p$.

We complete the proof of Proposition \ref{tke}.\,\,\,\,$\endpf$


\subsection{Numerical test for \texorpdfstring{$\mathcal M_{s,p,n}$ }{jj}}

In this subsection, we will derive a numerical test for $\mathcal M_{s,p,n}$ in a similar way to Section \ref{ttest}.

Let $\check {\mathcal M}\subset\mathfrak X(T)$
be the set of characters of $B$-semi-invariant functions in the function field $K(\check O)$   where  $\check O$ is the open $G$-orbit of  $\mathcal M_{s,p,n}$; denote by $\check{\mathcal N}:={\rm Hom}(\check{\mathcal M},\mathbb Z)$ its $\mathbb Z$-dual.
Then  $\{ \chi_2,\cdots, \chi_{r}\}$ (see (\ref{chi1})) is a $\mathbb Z$-basis  of $\check {\mathcal M}\subset\mathfrak X(T)$. By a slight abuse of notation, we denote by $\gamma_i$, $2\leq i\leq r$, its image in $\check {\mathcal N}$; then $\left\{ \gamma_2,\cdots, \gamma_{r}\right\}$ is the basis of $\check {\mathcal N}$ dual to $\{ \chi_2,\cdots, \chi_{r}\}$. 

Let $\check {\mathcal V}\subset\check { \mathcal N}\otimes_{\mathbb Z} \mathbb Q$ be the valuation cone with respect to $B$ on $K(\check O)$. Denote  by $v_{\check D_{i}}$, $2\leq i\leq r$, the valuations associated with $\check D_{i}$. Notice that for $2\leq i\leq r$ the B-semi-invariant  function $f_i$ defined by (\ref{bsemi}) is regular at generic points of $\mathcal M_{s,p,n}$; let $\check f_i:=f_i\big|_{\mathcal M_{s,p,n}}$ be the restriction of $f_i$; $\check f_{i}$ is a $B$-semi-invariant of $K(\check O)$ associated with the weight $\chi_i$.  Define a map $\check \rho: \left\{\check B_0,\cdots,\check B_r\right\}\rightarrow \check{\mathcal N}$ such that $\langle\check\rho(\check B_j), \chi_i\rangle=v_{\check B_j}(\check f_{i})$, $0\leq j\leq r$, $2\leq i\leq r$,
where $v_{\check B_j}$ is the discrete valuation associated to $\check B_j$.

Similarly to Lemmas \ref{bfor}, \ref{val}, we have
\begin{lemma}\label{mval}
$\check {\mathcal V}$ is the cone generated by $v_{\check D_{2}}$, $v_{\check D_{3}}$, $\cdots$, $v_{\check D_{r}}$  over $\mathbb Z$. Over the basis $\left\{ \gamma_2,\cdots, \gamma_{r}\right\}$,
\begin{equation}\label{checkvd}
v_{\check D_{i}} = \gamma_i\,\,,\,\,\,\,\,\, 2\leq i\leq r\,.
\end{equation} 
Moreover, when $p<s$,
\begin{equation}\label{mbdn}
    \check \rho(\check B_j)=
 \gamma_{j}-2 \gamma_{j+1}+ \gamma_{j+2}\,,\,\,\,\, 0\leq j\leq r\,;
\end{equation}
when $p=s=n-s$,
\begin{equation}
    \check \rho(\check B_j)=
 \gamma_{j}-2 \gamma_{j+1}+ \gamma_{j+2}\,,\,\,\,\, 1\leq j\leq r\,.
\end{equation}
Here we use the convention that $\gamma_i=0$ when $i>r$ or $i<2$.
\end{lemma}

Define a polytope  $Q_{\mathcal M_{s,p,n}}\subset\check {\mathcal N}\otimes_{\mathbb Z}\mathbb Q$ based on the canonical bundle formulas in  Lemma \ref{mkb} as follows. When $p<n-s$ and $p<s$ ($r=p$),
\begin{equation}
  Q_{\mathcal M_{s,p,n}}:={\rm conv} \left\{\frac{\check\rho(\check B_0)}{s-p+1},\frac{\check \rho(\check B_1)}{2},\cdots,\frac{\check \rho(\check B_{p-1})}{2}, \frac{\check \rho(\check B_p)}{n-s-p+1}, v_{\check D_2}, \cdots,v_{\check D_r}\right\}\,;
\end{equation}
when $n-s=p<s$ ($r=p$),
\begin{equation}
  Q_{\mathcal M_{s,p,n}}:={\rm conv} \left\{ \frac{\check\rho(\check B_0)}{s-p+1},\frac{\check \rho(\check B_1)}{2},\cdots,\frac{\check\rho(\check B_{p-1})}{2}, v_{\check D_2}, \cdots,v_{\check D_r}\right\}\,;
\end{equation}
when $n-s<p$ and $p<s$ ($r=n-s$),
\begin{equation}
  Q_{\mathcal M_{s,p,n}}:={\rm conv} \left\{\frac{\check\rho(\check B_0)}{s-p+1}, \frac{\check\rho(\check B_1)}{2},\cdots,\frac{\check\rho(\check B_{r-1})}{2}, \frac{\check\rho(\check B_r)}{p-r+1},v_{\check D_2}, \cdots,v_{\check D_r}\right\}\,;
\end{equation}
when $n-s=p=s$ ($r=p$),
\begin{equation}
  Q_{\mathcal M_{s,p,n}}:={\rm conv} \left\{\frac{\check \rho(\check B_1)}{2},\cdots,\frac{\check\rho(\check B_{p-1})}{2}, v_{\check D_2}, \cdots,v_{\check D_r}\right\}\,.
\end{equation}
\medskip

Denote by $Q^*_{\mathcal M_{s,p,n}}$ the dual polytope to $Q_{\mathcal M_{s,p,n}}\subset\check{\mathcal N}\otimes_{\mathbb Z}\mathbb Q$. Denote by $\Delta_{\mathcal M_{s,p,n}}^+\subset\mathfrak X(T) \otimes_{\mathbb Z}\mathbb R$ the moment polytope of ${\mathcal M_{s,p,n}}$ with respect to $B$.
Then,
\begin{equation}
    \Delta_{\mathcal M_{s,p,n}}^+=2\rho_{P}+Q^*_{\mathcal M_{s,p,n}}\,,
\end{equation}
where $2\rho_{P}$ is the the weight of a $B$-semi-invariant section of the anti-canonical bundle of $\mathcal M_{s,p,n}$ which is also given by (\ref{canwe}). 


We denote by ${\rm bar}_{DH}(\Delta^+)$ the barycenter of the polytope
$\Delta^+$ with respect to the Duistermaat–Heckman measure defined by (\ref{duiheck}).

Let $\check\pi:\mathfrak Y(T)\otimes_{\mathbb Z}\mathbb Q\rightarrow \check{\mathcal N}\otimes_{\mathbb Z}\mathbb Q$ be the natural quotient map. $\check\pi^{-1}(-\check {\mathcal V})$ is the cone generated by 
\begin{equation}\label{sproo}
     \left \{-v_{\check  D_{2}}, \cdots, -v_{\check  D_{r}}, \pm\gamma_1, \pm\theta_1,\cdots, \pm\theta_{r},\pm\xi_1,\cdots, \pm\xi_{s-p}, \pm\delta_{1},\cdots,\pm\delta_{|n-s-p|} \right\}\,.
\end{equation}
Let
$\Xi_{\mathcal M_{s,p,n}}\subset{\mathfrak X(T)}\otimes_{\mathbb Z} \mathbb R$ be the dual cone to  $\check\pi^{-1}(-\check {\mathcal V})$.
\medskip

{\bf\noindent Proof of Proposition \ref{mr13}.} When $r=1$, it is clear that each $\mathcal M_{s,p,n}$ is homogeneous, and hence there are K\"ahler-Einstein metrics by  \cite{Be}.  Without loss of generality, we may assume that $2p\leq n\leq 2s$ and $r\geq 2$. Recall that $\mathcal M_{p,s,n}\cong\mathcal M_{n-p,n-s,n}$. Hence,  we may further assume that $r=n-s\leq p$.

It is clear that 
\begin{equation}
    \Xi_{\mathcal M_{s,p,n}}=\left\{(x_2,\cdots,x_r)\in\mathbb{R}^{r-1}\,\big|\,x_i\geq 0\,,\,\,2\leq i\leq r\,\right\}\,;
\end{equation}
\begin{equation}
\small
\begin{split}
    2\rho_P&=\sum_{i=1}^{n-s}\left(\frac{n}{2}-p+i-1\right)(n-s+1-i)\cdot \chi_i+ \left(\frac{n}{2}-p\right)\sum_{i=1}^{n-s}\epsilon_i-p\sum_{i=1}^{s-p} \tau_i+(n-p)\sum_{i=1}^{s+p-n} \kappa_i\,.
\end{split}
\end{equation}
Similarly to Proposition \ref{tke}, we can show that
the moment polytope $\Delta_{\mathcal M_{s,p,n}}^+$ consists of the points $A(x_2,x_3,\cdots,x_r)$ with the following form
\begin{equation}
\begin{split}
    A(x_2,x_3,\cdots,x_r):=&\sum_{i=2}^{n-s}\left(x_i+\left(\frac{n}{2}-p+i-1\right)(n-s+1-i)\right)\cdot\chi_i\\
    &\hspace{-1in}+\left(\frac{n}{2}-p\right)(n-s)\cdot\chi_1+ \left(\frac{n}{2}-p\right)\sum_{i=1}^{n-s}\epsilon_i-p\sum_{i=1}^{s-p} \tau_i+(n-p)\sum_{i=1}^{s+p-n} \kappa_i\,,
\end{split}
\end{equation}
where $(x_2,x_3,\cdots,x_r)$ is in the region $\Delta(\cong Q^*_{\mathcal M_{s,p,n}})$ defined as follows.
\begin{equation}\label{region}
\Delta:=\left\{(x_2,\cdots,x_r)\in\mathbb{R}^{r-1}\,\rule[-.2in]{0.01in}{.5in}\small{\begin{matrix}
\,\,x_i\leq 1\,\,{\rm for}\,\,2\leq i\leq r\,;\\
\,\,x_{j-1}-2x_{j}+x_{j+1}\leq 2\,\,{\rm for}\,\,2\leq j\leq r\,\,{\rm with}\,\\
\,{\rm the\,\,convention\,\,that\,\,}x_1=x_{r+1}=0
\end{matrix}}\right\}\,.
\end{equation} 
By Lemma \ref{dhm}, $\rho(x_2,\cdots,x_r)dx_2\cdots dx_r$ the restriction of the Duistermaat–Heckman measure on $\Delta_{\mathcal M_{s,p,n}}^+$ takes the following form up to a constant. 
\begin{equation}
\begin{split}
&\rho(x_2,\cdots,x_r)dx_2\cdots dx_r=dx_2\cdots dx_r\left(\prod_{1\leq i<j\leq r}\big((x_{i+1}-x_{i})-(x_{j+1}-x_{j})-2i+2j\big)\right)^2\\
&\cdot\left(\prod_{i=1}^r\left((x_{i+1}-x_{i})+n-s+p-2i+1\right)\right)^{s+p-n}\left(\prod_{i=1}^r\left(s-p+2i-1-(x_{i+1}-x_{i})\right)\right)^{s-p}\\
\end{split}
\end{equation}

By Theorem \ref{dcr}, there are K\"ahler-Einstein metrics on $\mathcal M_{s,p,n}$ if and only if barycenter ${\rm bar}_{DH}(\Delta_{\mathcal M_{s,p,n}}^+)$ is in the relative interior of $2\rho_P+ \Xi_{\mathcal M_{s,p,n}}$, that is,
\begin{equation}\label{criterion}
\begin{split}
\displaystyle\int_{\Delta}x_k\cdot\rho(x_2,\cdots,x_r)dx_2\cdots dx_r>0\,\,{\rm for}\,\, 2\leq k\leq r.\\
\end{split}
\end{equation}

Similarly we have 
\begin{equation}
    \Xi_{\mathcal M_{n-p,n-s,n}}=\left\{(x_2,\cdots,x_r)\in\mathbb{R}^{r-1}\,\big|\,x_i\geq 0\,,\,\,2\leq i\leq r\,\right\}\,.
\end{equation}
The corresponding $2\rho_P$ takes the form
\begin{equation}
\footnotesize
\begin{split}
    2\rho_P&=\sum_{i=1}^{n-s}\left(\frac{n}{2}-p+i-1\right)(n-s+1-i)\cdot \chi_i+ \left(s-\frac{n}{2}\right)\sum_{i=1}^{n-s} \epsilon_i-(n-s)\sum_{i=1}^{s-p} \tau_i- (n-s)\sum_{i=1}^{s+p-n} \kappa_i\,;
\end{split}
\end{equation}
The moment polytope $\Delta_{\mathcal M_{n-p,n-s,n}}^+$ consists of the points $A(x_2,x_3,\cdots,x_r)$ with the following form
\begin{equation}
\begin{split}
    A(x_2,x_3,\cdots,x_r):=&\sum_{i=2}^{n-s}\left(x_i+\left(\frac{n}{2}-p+i-1\right)(n-s+1-i)\right)\cdot\chi_i\\
    &\hspace{-1in}+\left(\frac{n}{2}-p\right)(n-s)\cdot\chi_1+ \left(s-\frac{n}{2}\right)\sum_{i=1}^{n-s}\epsilon_i-(n-s)\sum_{i=1}^{s-p} \tau_i-(n-s)\sum_{i=1}^{s+p-n} \kappa_i\,,
\end{split}
\end{equation}
where $(x_2,x_3,\cdots,x_r)$ is in the region $\Delta(\cong Q^*_{\mathcal M_{n-p,n-s,n}})$ defined by (\ref{region}).
By Lemma \ref{dhm}, $\widetilde \rho(x_2,\cdots,x_r)dx_2\cdots dx_r$ the restriction of the Duistermaat–Heckman measure on $\Delta_{\mathcal M_{n-p,n-s,n}}^+$ takes the following form up to a constant. 
\begin{equation}
\small
\begin{split}
&\widetilde\rho(x_2,\cdots,x_r)dx_2\cdots dx_r=dx_2\cdots dx_r\left(\prod_{1\leq i<j\leq r}\big((x_{i+1}-x_{i})-(x_{j+1}-x_{j})-2i+2j\big)\right)^2\\
&\cdot\left(\prod_{i=1}^r\left((x_{i+1}-x_{i})+n-s+p-2i+1\right)\right)^{s+p-n}\left(\prod_{i=1}^r\left(s-p+2i-1-(x_{i+1}-x_{i})\right)\right)^{s-p}\\
\end{split}
\end{equation}

By Theorem \ref{dcr}, there are K\"ahler-Einstein metrics on $\mathcal M_{n-p,n-s,n}$ if and only if barycenter ${\rm bar}_{DH}(\Delta_{\mathcal M_{n-p,n-s,n}}^+)$ is in the relative interior of $2\rho_P+ \Xi_{\mathcal M_{n-p,n-s,n}}$, that is,
\begin{equation}
\begin{split}
\displaystyle\int_{\Delta}x_k\cdot\widetilde\rho(x_2,\cdots,x_r)dx_2\cdots dx_r>0\,\,{\rm for}\,\, 2\leq k\leq r.\\
\end{split}
\end{equation}

Noticing that  $\widetilde\rho(x_2,\cdots,x_r)=\rho(x_2,\cdots,x_r)$, we complete the proof of Proposition \ref{mr13}.\,\,\,\, $\endpf$
\medskip

{\bf\noindent Proof of Proposition \ref{mr14}.} Without loss of generality, we may assume that $n-s\leq p$. Then for fixed $n-s$, we have
\begin{equation}
\begin{split}
&\lim_{s-p\rightarrow\infty,s+p-n\rightarrow\infty}\frac{\displaystyle\int_{\Delta}x_k\cdot\rho(x_2,\cdots,x_r)dx_2\cdots dx_r}{\prod_{i=1}^r(n-s+p-2i+1)^{s+p-n}(s-p+2i-1)^{s-p}}\\
&=\lim_{s-p\rightarrow\infty,s+p-n\rightarrow\infty}\displaystyle\int_{\Delta}dx_2\cdots dx_r\,x_k\left(\prod_{1\leq i<j\leq r}\big((x_{i+1}-x_{i})-(x_{j+1}-x_{j})-2i+2j\big)\right)^2\\
&\hspace{.3in}\cdot\left(\prod_{i=1}^r\left(1+\frac{x_{i+1}-x_{i}}{n-s+p-2i+1}\right)\right)^{s+p-n}\left(\prod_{i=1}^r\left(1-\frac{x_{i+1}-x_{i}}{s-p+2i-1}\right)\right)^{s-p}\\
\end{split}
\end{equation}
\begin{equation*}
\begin{split}
&=\displaystyle\int_{\Delta}dx_2\cdots dx_r\,x_k\left(\prod_{1\leq i<j\leq r}\big((x_{i+1}-x_{i})-(x_{j+1}-x_{j})-2i+2j\big)\right)^2\\
&\hspace{.3in}\cdot \lim_{s-p\rightarrow\infty,s+p-n\rightarrow\infty}\exp\left\{\sum_{i=1}^r\left(\frac{s+p-n}{n-s+p-2i+1}-\frac{s-p}{s-p+2i-1}\right)(x_{i+1}-x_{i})\right\}\,\,\,\,\,\,\,\,\,\,\,\,\,\,\,\,\\
&=\displaystyle\int_{\Delta}dx_2\cdots dx_r\,x_k\left(\prod_{1\leq i<j\leq r}\big((x_{i+1}-x_{i})-(x_{j+1}-x_{j})-2i+2j\big)\right)^2\,.\\
\end{split}
\end{equation*}

We conclude Proposition \ref{mr14}.\,\,\,\, $\endpf$
\medskip

{\bf\noindent Proof of Corollary \ref{mr15}.} There are K\"ahler-Einstein metrics on  $\mathcal M_{3,3,6}$ by \cite{De2}, and on $\mathcal M_{4,4,8}$ and $\mathcal M_{5,5,10}$ by Lemma \ref{critm1} in Appendix \ref{section:nc}.  

The cases $r=1$ and $r=2$ follow from  \cite{Be} and  Lemma \ref{numcal2} respectively. \,\,\,\, $\endpf$

\section{Blow-ups with  vector-valued parameters} \label{vps}

Recall the rational map $\mathcal K_{\overline s,p,n}$ defined by (\ref{r}). Without loss of generality, we may assume that $s_1\geq \cdots\geq s_t\geq 1$. Let $r_1,\cdots,r_m$ be positive integers such that $\sum_{i=1}^mr_m=t$ and that
\begin{equation}
s_1=s_2=\cdots=s_{r_1}>s_{r_1+1}=\cdots= s_{r_1+r_2}>\cdots> s_{r_1+\cdots r_{m-1}+1}=\cdots= s_{r_1+\cdots r_m}\,.
\end{equation} 
It is clear that the following group $G$ acts  on $\mathcal T_{\overline s,p,n}$.
\begin{equation}
    G:=\left(\bigslant{(GL(s_1,\mathbb C)\times GL(s_2,\mathbb C)\times\cdots\times GL(s_t,\mathbb C))}{Z(n,\mathbb C)}\right)\rtimes(\mathbb Z/r_1\mathbb Z\times\cdots\times\mathbb Z/r_m\mathbb Z)\,.
\end{equation}
We can introduce a $(\mathbb C^*)^t$-action on $\mathcal T_{\overline s,p,n}$ similarly to the $\mathbb C^*$-action $\Psi_{s,p,n}$ on $\mathcal T_{s,p,n}$. Here $(\mathbb C^*)^t=\mathbb C^*\times\cdots\times\mathbb C^*$ is a high-dimensional torus.

It is clear that $\mathcal T_{\overline s,p,n}$ can be constructed as iterated blow-ups in the same way as in Section \ref{iterated}.  
\begin{definitionlemma} $\left(\mathcal K_{\overline s,p,n}\right)^{-1}$ has a regular extension to $\mathcal T_{\overline s,p,n}$. Denote the blow-up by  $R_{\overline s,p,n}:\mathcal T_{\overline s,p,n}\rightarrow G(p,n)$. \end{definitionlemma}

\begin{lemma}\label{nonreg}
$\mathcal T_{\overline s,3,9}$  is not smooth when $\overline s=(3,3,3)$.
\end{lemma}
{\bf\noindent Proof of Lemma \ref{nonreg}.} Denote by  $S_{(k_1,k_2,k_3)}$ the ideal sheaf of $G(3,9)$ associated with  the index set $\mathbb I_{\overline{s},3,9}^{(k_1,k_2,k_3)}$ where $k_1+k_2+k_3=3$. Define an open set of $G(3,9)$ by
\begin{equation}
  U_0:=  \left\{\left. \left( I_{3\times 3}\hspace{-0.13in}\begin{matrix}
  &\hfill\tikzmark{c1}\\
  &\hfill\tikzmark{d1}
  \end{matrix}\,\,\,\,Y\,\hspace{-0.15in}\begin{matrix}
  &\hfill\tikzmark{c2}\\
  &\hfill\tikzmark{d2}
  \end{matrix}\,\,\,\, W\right)\right\vert_{}\footnotesize\begin{matrix}
  \,Y\,\,{\rm  is\,\, a\,\,} 3\times 3\,\,{\rm matrix}\,;\\
  W\,\,{\rm  is\,\, a\,\,} 3\times 3\,\,{\rm matrix}\,\\
  \end{matrix}\right\}\,.
  \tikz[remember picture,overlay]   \draw[dashed,dash pattern={on 4pt off 2pt}] ([xshift=0.5\tabcolsep,yshift=7pt]c1.north) -- ([xshift=0.5\tabcolsep,yshift=-2pt]d1.south);\tikz[remember picture,overlay]   \draw[dashed,dash pattern={on 4pt off 2pt}] ([xshift=0.5\tabcolsep,yshift=7pt]c2.north) -- ([xshift=0.5\tabcolsep,yshift=-2pt]d2.south);
\end{equation}
We can blow up the ideals $S_{(2,1,0)}$, $S_{(1,2,0)}$, $S_{(2,0,1)}$, $S_{(1,0,2)}$, and derive the following local coordinate chart similarly to the Van der Waerden representation.
\begin{equation}
\left( I_{3\times 3} \hspace{-0.1in}\begin{matrix}
  &\hfill\tikzmark{c1}\\
  \\
  &\hfill\tikzmark{d1}
  \end{matrix}\,\,\, \sum_{k=1}^3\,\Xi_k^T\cdot\Omega_k\cdot\prod_{l=1}^{k}a_{l}\hspace{-0.13in}\begin{matrix}
  &\hfill\tikzmark{c2}\\
  \\
  &\hfill\tikzmark{d2}
  \end{matrix}\,\,\,\,\sum_{k=1}^3\,\widetilde\Xi_k^T\cdot\widetilde \Omega_k\cdot\prod_{l=1}^{k}b_{l} \right)\,,
\tikz[remember picture,overlay]   \draw[dashed,dash pattern={on 4pt off 2pt}] ([xshift=0.5\tabcolsep,yshift=7pt]c2.north) -- ([xshift=0.5\tabcolsep,yshift=-2pt]d2.south);
\tikz[remember picture,overlay]   \draw[dashed,dash pattern={on 4pt off 2pt}] ([xshift=0.5\tabcolsep,yshift=7pt]c1.north) -- ([xshift=0.5\tabcolsep,yshift=-2pt]d1.south);
\end{equation} 
where  \begin{equation}
\begin{split}
&\Xi_1=\left(1,\xi^{(1)}_{24},\xi^{(1)}_{34}\right)\,,\,\,\Xi_2=\left(0,1,\xi^{(2)}_{35}\right)\,,\,\,\Xi_3=\left(0,0,1\right)\,,\\
&\Omega_1=\left(1,\xi^{(1)}_{15},\xi^{(1)}_{16}\right)\,,\,\,\Omega_2=\left(0,1,\xi^{(2)}_{26}\right)\,,\,\,\Omega_3=\left(0,0,1\right)\,,\\
&\widetilde\Xi_1=\left(\xi^{(1)}_{18},1,\xi^{(1)}_{38}\right)\,,\,\,\Xi_2=\left(1,0,\xi^{(2)}_{37}\right)\,,\,\,\Xi_3=\left(0,0,1\right)\,,\\
&\widetilde\Omega_1=\left(\xi^{(1)}_{27},1,\xi^{(1)}_{29}\right)\,,\,\,\Omega_2=\left(1,0,\xi^{(2)}_{19}\right)\,,\,\,\Omega_3=\left(0,0,1\right)\,.\\
\end{split}
\end{equation}	
By a slight abuse of notation,  we call the set of such local coordinate charts the partial Van der Waerden representation. Denote the partial blow-up by $\tau:Y\rightarrow G(3,9)$.

Next we will compute the pullbacks of the ideal sheaves $S_{(1,1,1)}$,$S_{(0,2,1)}$, $S_{(0,1,2)}$ under $\tau$. Notice that, in terms of the local coordinates, they are given by the Pl\"ucker coordinate functions associated with the indices sets  $\mathbb I_{\overline{s},p,n}^{(1,1,1)}$, $\mathbb I_{\overline{s},p,n}^{(0,2,1)}$, $\mathbb I_{\overline{s},p,n}^{(0,1,2)}$ of the following matrix.
\begin{equation}
\left( \begin{matrix}
  1&0&0\\
  0&1&0\\
  0&0&1\\
\end{matrix} \hspace{-0.1in}\begin{matrix}
  &\hfill\tikzmark{c1}\\
  \\
  &\hfill\tikzmark{d1}
  \end{matrix}\,\,\,\, \begin{matrix}
  a_1&0&0\\
  0&a_1a_2&0\\
  0&0&a_1a_2a_3\\
\end{matrix}\hspace{-0.13in}\begin{matrix}
  &\hfill\tikzmark{c2}\\
  \\
  &\hfill\tikzmark{d2}
  \end{matrix}\,\,\,\,\begin{matrix}
  b_1b_2&0&0\\
  0&b_1&0\\
  0&0&b_1b_2b_3\\
\end{matrix}\right)\,,
\tikz[remember picture,overlay]   \draw[dashed,dash pattern={on 4pt off 2pt}] ([xshift=0.5\tabcolsep,yshift=7pt]c2.north) -- ([xshift=0.5\tabcolsep,yshift=-2pt]d2.south);
\tikz[remember picture,overlay]   \draw[dashed,dash pattern={on 4pt off 2pt}] ([xshift=0.5\tabcolsep,yshift=7pt]c1.north) -- ([xshift=0.5\tabcolsep,yshift=-2pt]d1.south);
\end{equation} 
or equivalently,
\begin{equation}
\left( \begin{matrix}
  1&0&0\\
  0&1&0\\
  0&0&1\\
\end{matrix} \hspace{-0.1in}\begin{matrix}
  &\hfill\tikzmark{c1}\\
  \\
  &\hfill\tikzmark{d1}
  \end{matrix}\,\,\,\, \begin{matrix}
  1&0&0\\
  0&a_2&0\\
  0&0&a_2a_3\\
\end{matrix}\hspace{-0.13in}\begin{matrix}
  &\hfill\tikzmark{c2}\\
  \\
  &\hfill\tikzmark{d2}
  \end{matrix}\,\,\,\,\begin{matrix}
  b_2&0&0\\
  0&1&0\\
  0&0&b_2b_3\\
\end{matrix}\right)\,.
\tikz[remember picture,overlay]   \draw[dashed,dash pattern={on 4pt off 2pt}] ([xshift=0.5\tabcolsep,yshift=7pt]c2.north) -- ([xshift=0.5\tabcolsep,yshift=-2pt]d2.south);
\tikz[remember picture,overlay]   \draw[dashed,dash pattern={on 4pt off 2pt}] ([xshift=0.5\tabcolsep,yshift=7pt]c1.north) -- ([xshift=0.5\tabcolsep,yshift=-2pt]d1.south);
\end{equation}
Computation yields that in a neighborhood of the origin where all the local coordinates are close to zero, the ideal sheaf $S_{(1,1,1)}$ is trivial,   $S_{(0,2,1)}=(a_2 a_3,a_2 b_2b_3)$, and $S_{(0,1,2)}= (b_2 b_3,b_2 a_2a_3)$.
Then $\mathcal T_{\overline s,3,9}$ is locally isomorphic to 
\begin{equation}
    \begin{split}
    &{\rm Spec}\,\mathbb C[a_1,a_2,\cdots,\xi^{(k)}_{ij},\cdots]\times{\rm Spec}\,\left(\frac{\mathbb C[a_2,a_3,b_2,b_3,t,\tilde t\,]}{ \big(a_3\tilde t-b_2b_3,\,b_3t-a_2a_3\big)}\right)\,\,\,.\\
     \end{split}
\end{equation}

Then $\mathcal T_{\overline s,3,9}$ is not smooth by the Jacobi criterion.
\,\,\,$\endpf$
\begin{remark}\label{rmon}
As shown in Lemma \ref{nonreg}, there is a partial parameterization similar to the Van der Waerden representation such that the singularities are defined  explicitly by the equations of the following form. 
\begin{equation}
    \prod_{i=1}^vx_i^{\alpha_i}=\prod_{j=1}^uy_j^{\beta_j}\,
\end{equation}
where $x_i$, $y_j$ are certain local coordinates and $\alpha_i$, $\beta_j$ are positive integers. 
\end{remark}

\begin{question}
What is the structure of the singularities of $\mathcal T_{\overline s,p,n}$.
\end{question}

Similarly to Definition \ref{mspn}, we can define $\mathcal M_{{\overline s},p,n}$. For simplicity, assume that $s_i\geq p$ for each $1\leq i\leq t$. Then we can define $\mathcal M_{{\overline s},p,n}$ by the preimage  of the following sub-grassmannian $Z_1$ of $G(p,n)$ under the blow-up $R_{\overline s,p,n}$.
\begin{equation}
Z_1:=\left\{\left.\left(\,\,\begin{matrix}
  Z\\
\end{matrix}
  \hspace{-.1in}\begin{matrix}
  &\hfill\tikzmark{a}\\
  &\hfill\tikzmark{b}  
  \end{matrix} \,\,\,\,\,
  \begin{matrix}
  0\\
\end{matrix}\hspace{-.07in}
\begin{matrix}
  &\hfill\tikzmark{c}\\
  &\hfill\tikzmark{d}
  \end{matrix}\,\,\,\,
\begin{matrix}
\cdots\\
\end{matrix}\hspace{-.08in}
\begin{matrix}
  &\hfill\tikzmark{e}\\
  &\hfill\tikzmark{f}\end{matrix}\,\,\,\,\,\,
  \begin{matrix} 0\\ \end{matrix}\,\,\,\right)\right\vert_{}\,{Z\,\, {\rm is\,\,a}\,\,p\times s_1{\rm\,\,nondegenerate\,\,matrix}}\,\right\}.
  \tikz[remember picture,overlay]   \draw[dashed,dash pattern={on 4pt off 2pt}] ([xshift=0.5\tabcolsep,yshift=7pt]a.north) -- ([xshift=0.5\tabcolsep,yshift=-2pt]b.south);\tikz[remember picture,overlay]   \draw[dashed,dash pattern={on 4pt off 2pt}] ([xshift=0.5\tabcolsep,yshift=7pt]c.north) -- ([xshift=0.5\tabcolsep,yshift=-2pt]d.south);\tikz[remember picture,overlay]   \draw[dashed,dash pattern={on 4pt off 2pt}] ([xshift=0.5\tabcolsep,yshift=7pt]e.north) -- ([xshift=0.5\tabcolsep,yshift=-2pt]f.south);
\end{equation}
Similarly to Lemma \ref{emb}, we can derive
\begin{definitionlemma}\label{flatm}
By projecting  $\mathcal T_{\overline s,p,n}$  to $\mathbb {CP}^{N_{p,n}}\times\mathbb {CP}^{N_{\overline{s},p,n}^{K_1}}\times\cdots\times\mathbb {CP}^{N_{\overline{s},p,n}^{K_{L_{\overline s}}}}$, we have a regular map
\begin{equation}
   \mathcal P_{\overline s,p,n}:\mathcal T_{\overline s,p,n}\rightarrow \mathcal M_{\overline s,p,n}\,. 
\end{equation}
\end{definitionlemma}

There are also various fibration structures on $\mathcal M_{\overline s,p,n}$ by the following lemma.

\begin{lemma}\label{flatm2}
Assume that $s_i\geq p$ for each $1\leq i\leq t$. Define a sub-grassmannian $Z_l$, $1\leq l\leq t$, of $G(p,n)$ by
\begin{equation}
Z_l:=\left\{\left.\left(\,\,
  \begin{matrix}
  0\\
\end{matrix}\hspace{-.07in}
\begin{matrix}
  &\hfill\tikzmark{c}\\
  &\hfill\tikzmark{d}
  \end{matrix}\,\,\,\,
\begin{matrix}
\cdots\\
\end{matrix}\hspace{-.08in}
\begin{matrix}
  &\hfill\tikzmark{e}\\
  &\hfill\tikzmark{f}\end{matrix}\,\,\,\,\,\,
  \begin{matrix} 0\\ \end{matrix}\hspace{-.1in}\begin{matrix}
  &\hfill\tikzmark{a1}\\
  &\hfill\tikzmark{b1}  
  \end{matrix}\,\,\,\, \begin{matrix}
  Z\\
\end{matrix}
  \hspace{-.1in}\begin{matrix}
  &\hfill\tikzmark{a}\\
  &\hfill\tikzmark{b}  
  \end{matrix} \,\,\,\,\,\begin{matrix}
  0\\
\end{matrix}\hspace{-.07in}
\begin{matrix}
  &\hfill\tikzmark{c1}\\
  &\hfill\tikzmark{d1}
  \end{matrix}\,\,\,\,
\begin{matrix}
\cdots\\
\end{matrix}\hspace{-.08in}
\begin{matrix}
  &\hfill\tikzmark{e1}\\
  &\hfill\tikzmark{f1}\end{matrix}\,\,\,\,\,\,
  \begin{matrix} 0\\ \end{matrix}\,\,\right)\right\vert_{}
  \begin{matrix}
    \,{Z\,\, {\rm is\,\,a}\,\,p\times s_l{\rm\,\,nondegenerate\,\,matrix}}\,\\
    {\rm consisting\,\, of\,\,the\,\,}(\alpha_l+1)^{th}\,\,\cdots,(\alpha_l+s_l)^{th}\\
     {\rm rows\,\, where\,\,}\alpha_l=s_1+\cdots+s_{l-1}
  \end{matrix}\right\}.
  \tikz[remember picture,overlay]   \draw[dashed,dash pattern={on 4pt off 2pt}] ([xshift=0.5\tabcolsep,yshift=7pt]a.north) -- ([xshift=0.5\tabcolsep,yshift=-2pt]b.south);\tikz[remember picture,overlay]   \draw[dashed,dash pattern={on 4pt off 2pt}] ([xshift=0.5\tabcolsep,yshift=7pt]c.north) -- ([xshift=0.5\tabcolsep,yshift=-2pt]d.south);\tikz[remember picture,overlay]   \draw[dashed,dash pattern={on 4pt off 2pt}] ([xshift=0.5\tabcolsep,yshift=7pt]e.north) -- ([xshift=0.5\tabcolsep,yshift=-2pt]f.south);\tikz[remember picture,overlay]   \draw[dashed,dash pattern={on 4pt off 2pt}] ([xshift=0.5\tabcolsep,yshift=7pt]a1.north) -- ([xshift=0.5\tabcolsep,yshift=-2pt]b1.south);\tikz[remember picture,overlay]   \draw[dashed,dash pattern={on 4pt off 2pt}] ([xshift=0.5\tabcolsep,yshift=7pt]c1.north) -- ([xshift=0.5\tabcolsep,yshift=-2pt]d1.south);\tikz[remember picture,overlay]   \draw[dashed,dash pattern={on 4pt off 2pt}] ([xshift=0.5\tabcolsep,yshift=7pt]e1.north) -- ([xshift=0.5\tabcolsep,yshift=-2pt]f1.south);
\end{equation}
Denote $\mathcal Z_l$ the preimage  of $Z_l$ under the blow-up $R_{\overline s,p,n}$. Then
\begin{equation}
   \mathcal Z_l\cong \mathcal M_{\overline s,p,n}\,.
\end{equation}
\end{lemma}
{\bf\noindent Proof of Lemma \ref{flatm2}.} Notice that $\mathcal P_{\overline s,p,n}$ is an embedding restricted to $\mathcal Z_l$. Moreover,
\begin{equation}
\mathcal P_{\overline s,p,n}\left(\mathcal T_{\overline s,p,n}\right)=\mathcal P_{\overline s,p,n}\left(\mathcal Z_l\right)\,. 
\end{equation}

Then Lemma \ref{flatm2} follows.  \,\,\,$\endpf$

\medskip
\begin{example}
The partial  Van der Waerden representation in Lemma \ref{nonreg} gives local coordinate charts for $\mathcal T_{\overline s,1,n}$. We can verify that $\mathcal P_{\overline s,1,n}$ is smooth and flat.  In fact, $\mathcal T_{\overline s,1,n}$ is spherical.
\end{example}

Notice that the  $(\mathbb C^*)^t$-action induces a foliation on $\mathcal T_{\overline s,p,n}$. The generic leaves are the generic fibers of $\mathcal P_{\overline s,p,n}$  which are isomorphic to $\mathbb {CP}^{t-1}$.  The connected components of the set of the fixed points under the  $(\mathbb C^*)^t$-action  one to one correspond to the elements of $\mathbb K^{\overline s}$; each  component is a Cartesion product of certain grassmannians. On the other hand,  the foliation is more complicated due to the lack of the  sink-source interpretation. For a sub-action by one $\mathbb C^*$ factor, there is a partial order on $\mathbb K^{\overline s}$ according to Definition \ref{order} as follows.
\begin{definition}
Let $A=(i_1,\cdots,i_t)$ and  $B=(i^*_1,\cdots,i^*_t)$ be two indices in $\mathbb K^{\overline s}$. We say $A<B$ if and only if $A\neq B$ and the following holds,
	\begin{equation}
	\begin{split}
	&i_1\geq i^*_1,\\
	&i_1+i_2\geq i^*_1+i^*_2,\\
	&\cdots\\
	&i_1+i_2+\cdots+i_{t-1}\geq i^*_1+i^*_2+\cdots+i^*_{t-1}.\\
	\end{split}
	\end{equation}
\end{definition}

By the partial  Van der Waerden representation, we can show that locally $\mathcal T_{\overline s,p,n}$ is a vector bundle over $\mathcal Z_l$, $1\leq l\leq t$.  We would like to know
\begin{question}
Is there an analogue of the Bia{\l}ynicki-Birula decomposition associated with  algebraic $(\mathbb C^*)^t$-actions.
\end{question}

\section{Functoriality}\label{sat}
In this section, we study the functoriality of the canonical blow-ups.

For $0\leq p^{\prime}-p\leq n^{\prime}-n$, define an embedding between grassmannians $G(p,n)$, $G(p^{\prime},n^{\prime})$ by 
\begin{equation}
\begin{split}
e^{(p,n)}_{(p^{\prime},\,n^{\prime})}:\,\,\,\,\,\,\,\,G(p,n)&\rightarrow G(p^{\prime},n^{\prime})\\
U&\mapsto\left(\begin{matrix}
0_{p\times(p^{\prime}-p)}&U&0_{p\times(n^{\prime}-n+p-p^{\prime})}\\
I_{(p^{\prime}-p)\times(p^{\prime}-p)}&0_{(p^{\prime}-p)\times n}&0_{(p^{\prime}-p)\times(n^{\prime}-n+p-p^{\prime})}\\
\end{matrix}\right),
\end{split}
\end{equation}
where $U$ is a matrix representative of a point of $G(p,n)$; define an embedding  between the general linear groups $GL(n,\mathbb C)$, $GL(n^{\prime},\mathbb C)$ by
\begin{equation}
\begin{split}
E^{(p,n)}_{(p^{\prime},\,n^{\prime})}:\,\,\,\,\,\,\,\,&GL(n,\mathbb C)\rightarrow GL(n^{\prime},\mathbb C)\\
g&\mapsto\left(\begin{matrix}
I_{(p^{\prime}-p)\times(p^{\prime}-p)}&0_{(p^{\prime}-p)\times n} &0_{(p^{\prime}-p)\times(n^{\prime}-n+p-p^{\prime})}\\
0_{n\times(p^{\prime}-p)}&g&0_{n\times(n^{\prime}-n+p-p^{\prime})}\\
0_{(n^{\prime}-n+p-p^{\prime})\times(p^{\prime}-p)}&0_{(n^{\prime}-n+p-p^{\prime})\times n}&I_{(n^{\prime}-n+p-p^{\prime})\times(n^{\prime}-n+p-p^{\prime})}\\
\end{matrix}\right).
\end{split}
\end{equation}

It is clear that there is an embedding $\mathfrak e^{s,p,n}_{s+p^{\prime}-p,\,p^{\prime},\,n^{\prime}}:\mathcal T_{s,p,n}\rightarrow\mathcal T_{s+p^{\prime}-p,\,p^{\prime},\,n^{\prime}}$, $0\leq p^{\prime}-p\leq n^{\prime}-n$, such that the following diagram commutes and is compatible with the group actions.
\begin{equation}
\begin{tikzcd} &\mathcal T_{s,p,n}  \arrow{r}{\mathfrak e^{s,p,n}_{s+p^{\prime}-p,\,p^{\prime},\,n^{\prime}}}&[2em]\mathcal T_{s+p^{\prime}-p,\,p^{\prime},\,n^{\prime}}\arrow{d}{R_{s+p^{\prime}-p,\,p^{\prime},\,n^{\prime}}} \\ &G(p,n)\arrow[leftarrow]{u}{R_{s,p,n}} \arrow{r}{e^{(p,n)}_{(p^{\prime},\,n^{\prime})}}&G(p^{\prime},n^{\prime})\\ \end{tikzcd}\,\,.\vspace{-20pt} \end{equation}

It is easy to verify that
\begin{lemma}\label{func1}
For $0\leq p^{\prime}-p\leq n^{\prime}-n$ and $0\leq p^{\prime\prime}-p^{\prime}\leq n^{\prime\prime}-n^{\prime}$, 
\begin{equation}
    \mathfrak e^{s,p,n}_{s+p^{\prime}-p,\,p^{\prime},\,n^{\prime}}\circ \mathfrak e_{s+p^{\prime\prime}-p,\,p^{\prime\prime},\,n^{\prime\prime}}^{s+p^{\prime}-p,\,p^{\prime},\,n^{\prime}}=\mathfrak e^{s,p,n}_{s+p^{\prime\prime}-p,\,p^{\prime\prime},\,n^{\prime\prime}}\,,
\end{equation}
\end{lemma}
{\bf\noindent Proof of Lemma \ref{func1}. } We can prove Lemma \ref{func1} by considering the corresponding rational maps $\mathcal K_{s,p,n}$,  $\mathcal K_{s+p^{\prime}-p,\,p^{\prime},\,n^{\prime}}$, and
$\mathcal K_{s+p^{\prime\prime}-p,\,p^{\prime\prime},\,n^{\prime\prime}}$.\,\,\,$\endpf$
\begin{remark}
The above construction can be view as a blow-up of the Sato Grassmannian (\cite{S1,S2}).
\end{remark}

Next we show that $\mathcal T_{s,p,n}$ can be embedded into $\mathcal M_{s+1,p+1,n+2}$ (see \cite{MT} for the case $\mathcal T_{p,p,2p}$). 

Define an embedding between grassmannians $G(p,n)$ and $G(p+1,n+2)$ by
\begin{equation}
\begin{split}
h^{(p,n)}_{(p+1,n+2)}:\,\,\,\,\,\,\,\,G(p,n)&\rightarrow G(p+1,n+2)\\
U&\mapsto\left(\begin{matrix}
0_{p\times 1} &U&0_{p\times 1}\\
1&0_{1\times n} &1\\
\end{matrix}\right).
\end{split}
\end{equation}
There is an embedding $\mathfrak h^{s,p,n}_{s+1,p+1,n+2}:\mathcal T_{s,p,n}\rightarrow\mathcal T_{s+1,p+1,n+2}$ such that the following diagram commutes.
\begin{equation}
\begin{tikzcd} &\mathcal T_{s,p,n}  \arrow{r}{\mathfrak h^{s,p,n}_{s+1,p+1,n+2}}&[4em]\mathcal T_{s+1,p+1,n+2} \arrow{d}{R_{s+1,p+1,n+2}} \\ &G(p,n)\arrow[leftarrow]{u}{R_{s,p,n}} \arrow{r}{h^{(p,n)}_{(p+1,n+2)}}&G(p+1,n+2)\\ 
\end{tikzcd}\,\,.\vspace{-20pt} 
\end{equation}
Define $\mathfrak g^{s,p,n}_{s+1,p+1,n+2}:\mathcal T_{s,p,n}\rightarrow \mathcal M_{s+1,p+1,n+2}$ by
\begin{equation}
    \mathfrak g^{s,p,n}_{s+1,p+1,n+2}:=\mathcal P_{s+1,p+1,n+2}\circ\mathfrak h^{s,p,n}_{s+1,p+1,n+2}\,.
\end{equation}

It is easy to verify that
\begin{lemma}\label{func2}
$\mathfrak g_{s,p,n}^{s+1,p+1,n+2}$ is an embedding.
\end{lemma}
{\bf\noindent Proof of Lemma \ref{func2}. } It follows from the properties of $\mathcal P_{s+1,p+1,n+2}$. \,\,\,$\endpf$

\appendix

\section{\texorpdfstring{Intermediate Van der Waerden representation }{rr}}\label{section:cover}

In this appendix, we consider the iterated blow-ups of certain orders and construct locally the Van der Waerden type representation for each intermediate blow-up. As an application, we prove Lemmas \ref{coor2}, \ref{qcoor2}. 
\subsection{\texorpdfstring{A simple case}{rr}}\label{section:coverp}
To illustrate the idea, in what follows, we construct the local Van der Waerden representation over the open set $U_p\subset G(p,n)$ when $p=n-s$. 

Define a permutation  $\sigma_p$ by
\begin{equation}
  \sigma_p(k)=p-k,\,\,\,0\leq k\leq p.
\end{equation} Restricting $(\ref{sblow})$ to $U_p$, we have that
\vspace{-.05in}
\begin{equation}\label{blow1}
\footnotesize
\begin{tikzcd}
&Y^p_p\ar{r}{g^{\sigma_p}_p}&Y^{p}_{p-1}\ar{r}{g^{\sigma_p}_{p-1}}&\cdots\ar{r}{g^{\sigma_p}_1}&Y^{p}_0\ar{r}{g^{\sigma_p}_0}&U_p \\
&&(g^{\sigma_p}_0\circ\cdots\circ g^{\sigma_p}_{p-1})^{-1}(S_0\cap U_p)\ar[hook]{u}&\cdots&(g^{\sigma_p}_0)^{-1}(S_{p-1}\cap U_p)\ar[hook]{u}&S_{p}\cap U_p\ar[hook]{u}\\
\end{tikzcd}\vspace{-20pt}
\end{equation}
where for $0\leq k\leq p$,
\begin{equation}
   Y^p_{k}:=Y_k^{\sigma}\cap (g^{\sigma_p}_0\circ\cdots\circ g^{\sigma_p}_{k})^{-1}(U_p)\subset  G(p,n)\times\mathbb {CP}^{N^p_{s,p,n}}\times\cdots\times\mathbb {CP}^{N^{p-k}_{s,p,n}}\,\,\,\,. 
\end{equation}
It is clear that $Y^p_p\cong R_{s,p,n}^{-1}(U_p)$.

Next, we will define a system of  holomorphic coordinate charts for $Y^p_{k}$, $0\leq k\leq p$.

Define  index sets $\mathbb J_{p,k}$, $0\leq k\leq p$, by 
\begin{equation}
\mathbb J_{p,k}:=\left\{\left. \left(
\begin{matrix}
i_1&i_2&\cdots&i_k\\
j_1&j_2&\cdots&j_k\\
\end{matrix}\right)\right\vert_{}{\footnotesize\begin{matrix}
1\leq i_t\leq p\,,\,\, 1\leq t\leq k,\,\,{\rm and\,\,}i_{t_1}\neq i_{t_2}\,\,{\rm if\,\,}t_1\neq t_2\,;\\
1\leq j_t\leq s\,,\,\, 1\leq t\leq k,\,\,{\rm and\,\,}j_{t_1}\neq j_{t_2}\,\,{\rm if\,\,}t_1\neq t_2\\
\end{matrix}}
\right\}.	    
\end{equation}
For each $\tau=\left(\begin{matrix}
i_1&i_2&\cdots&i_k\\
j_1&j_2&\cdots&j_k\\
\end{matrix}\right)\in\mathbb J_{p,k}$, $0\leq k\leq p$, we have unique integers   $i^*_1,i^*_2,\cdots, i^*_{p-k}$ and $j^*_1,j^*_2$, $\cdots,$ $j^*_{s-k}$ such that the following holds.
\begin{enumerate}[label=(\alph*).]
\item $i^*_1<i^*_2<\cdots<i^*_{p-k}$\,\,\,and\,\,\,$j^*_1<j^*_2<\cdots<j^*_{s-k}\,.$
\item $\left\{i_1,\cdots,i_k,i^*_1,\cdots,i_{p-k}^*\right\}=\left\{1,2,\cdots,p\right\}$\,\,\,and\,$\left\{j_1,\cdots,j_k,j^*_1,\cdots,j_{s-k}^*\right\}=\left\{1,2,\cdots,s\right\}$.
\end{enumerate}
Associate $\tau$ with a complex Euclidean space $\mathbb {C}^{p(n-p)}$ equipped with holomorphic coordinates  $\left(\overrightarrow B^1,\cdots,\overrightarrow B^k,\overrightarrow B^*_k\right)$ where
$\overrightarrow B^j$ is defined by (\ref{bp}) for $1\leq j\leq k$, and
\begin{equation}\label{b8}
\overrightarrow B^{*}_k:=\left(x^{(k+1)}_{i^*_1j^*_1},x^{(k+1)}_{i^*_1j^*_2},\cdots, x^{(k+1)}_{i^*_1j^*_{s-k}},x^{(k+1)}_{i^*_2j^*_1},\cdots,x^{(k+1)}_{i^*_2j^*_{s-k}},\cdots,x^{(k+1)}_{i^*_{p-k}j^*_{1}}\cdots,x^{(k+1)}_{i^*_{p-k}j^*_{s-k}}\right)\in\mathbb C^{(p-k)(s-k)}.
\end{equation}

Define a holomorphic map $\Gamma_{p,k}^{\tau}:\mathbb C^{p(n-p)}\rightarrow U_p$ by
\begin{equation}\label{gk}
\Gamma_{p,k}^{\tau}\left(\overrightarrow B^1,\cdots,\overrightarrow B^k,\overrightarrow B^{*}_k\right):=\left(C^*_k+ \sum_{m=1}^k\,\Xi_m^T\cdot\Omega_m\cdot\prod_{t=1}^{m}a_{i_tj_t}\hspace{-0.11in}\begin{matrix}
  &\hfill\tikzmark{c2}\\
  \\
  &\hfill\tikzmark{d2}
  \end{matrix}\,\,\,\,I_{p\times p} \right)\,.
\tikz[remember picture,overlay]   \draw[dashed,dash pattern={on 4pt off 2pt}] ([xshift=0.5\tabcolsep,yshift=7pt]c2.north) -- ([xshift=0.5\tabcolsep,yshift=-2pt]d2.south);
\end{equation} 
Here  $\Xi_m$, $\Omega_m$ are defined by (\ref{nxi}), (\ref{nome}) respectively; $C^*_p$ is a null matrix and $C^*_k$, $0\leq k\leq p-1$, is defined by
\begin{equation}\label{w7}
C^*_k=\prod_{t=1}^{k} a_{i_tj_t}\cdot\left(\begin{matrix}
&\cdots&0&0&\cdots&0&0&\cdots&0&\cdots\\
&\cdots&0&x^{(k+1)}_{i^*_1j^*_1}&\cdots&x^{(k+1)}_{i^*_1j^*_2}&0&\cdots& x^{(k+1)}_{i^*_1j^*_{s-k}}&\cdots\\
&\cdots&0&0&\cdots&0&0&\cdots&0&\cdots\\
&\cdots&0&x^{(k+1)}_{i^*_2j^*_1}&\cdots&x^{(k+1)}_{i^*_2j^*_2}&0&\cdots& x^{(k+1)}_{i^*_2j^*_{s-k}}&\cdots\\
&\cdots&0&0&\cdots&0&0&\cdots&0&\cdots\\
&\ddots&\vdots&\vdots&\ddots&\vdots&\vdots&\ddots&\vdots&\ddots\\
&\cdots&0&0&\cdots&0&0&\cdots&0&\cdots\\
&\cdots&0&x^{(k+1)}_{i^*_{p-k}j^*_1}&\cdots&x^{(k+1)}_{i^*_{p-k}j^*_2}&0&\cdots& x^{(k+1)}_{i^*_{p-k}j^*_{s-k}}&\cdots\\
&\cdots&0&0&\cdots&0&0&\cdots&0&\cdots\\
\end{matrix}\right).\\
\end{equation}
Define a rational map $J_{p,k}^{\tau}:\mathbb C^{p(n-p)}\dashrightarrow\mathbb {CP}^{N_{p,n}}
\times\mathbb {CP}^{N^p_{s,p,n}}\times\cdots\times\mathbb {CP}^{N^{p-k}_{s,p,n}}$ by
\begin{equation}
    J^{\tau}_{p,k}:=\left(e,f_s^{p},f_s^{p-1},\cdots,f_s^{p-k}\right)\circ \Gamma_{p,k}^{\tau}\,,\,\,0\leq k\leq p.
\end{equation}
It is clear that $J^{\tau}_{p,p}$ is isomorphic to the holomorphic embedding $J^{\tau}$ defined in (\ref{jjtp}) after reordering the factors $\mathbb {CP}^{N^k_{s,p,n}}$ of the ambient space.

We can show that the rational map $J^{\tau}_{p,k}$ extends to a holomorphic embedding of $\mathbb C^{p(n-p)}$ for each $0\leq k\leq p$ and $\tau\in\mathbb J_{p,k}$, similarly Lemma \ref{em}. Denote by $A_{k}^{\tau}$ the image of $\mathbb C^{p(n-p)}$ in $\mathbb {CP}^{N_{p,n}}
\times\mathbb {CP}^{N^{p}_{s,p,n}}\times\cdots\times\mathbb {CP}^{N^{p-k}_{s,p,n}}$  under the holomorphic map $J_{p,k}^{\tau}$. By a slight abuse of notation, let $\left(J_{p,k}^{\tau}\right)^{-1}:A_{k}^{\tau}\rightarrow\mathbb C^{p(n-p)}$ be the inverse of $J_{p,k}^{\tau}$.

\begin{definition}\label{lvdwp}
For $0\leq k\leq p$, we call the above defined system of holomorphic  coordinate charts $\left\{\left(A^{\tau}_k,\left(J_{p,k}^{\tau}\right)^{-1}\right)\right\}_{\tau\in\mathbb J_{p,k}}$ the {\it intermediate Van der Waerden representation} of $Y_k^p$.
\end{definition}
\medskip

{\bf\noindent Proof of Lemma \ref{coor2}.} It suffices to prove the following by induction on $k$.

\begin{lemma}\label{kccord}
For $0\leq k\leq p$, $\left\{\left(A_{k}^{\tau},\left(J_{p,k}^{\tau}\right)^{-1}\right)\right\}_{\tau\in\mathbb J_{p,k}}$ is a holomorphic atlas for $Y^p_k$\,, or equivalently,    $\bigcup_{\tau\in\mathbb J_{p,k}}A_{k}^{\tau}=Y^p_k$.
\end{lemma}

{\noindent\bf Proof of  Lemma \ref{kccord}.} When $k=0$, $g^{\sigma_p}_0$ is a biholomorpshim between $Y^p_0$ and $U_p$, $\mathbb{CP}^{N^p_{s,p,n}}$ consists a single point and $f_s^{p}$ is a constant map.
\smallskip

{\noindent\bf {Step 1 ($k=1$).}} Take an arbitrary point $v_0\in Y_1^p\subset \mathbb {CP}^{N_{p,n}}\times\mathbb{CP}^{N^p_{s,p,n}}\times\mathbb{CP}^{N^{p-1}_{s,p,n}}$.  We will show that $v_0$ is in a certain coordinate chart $A_{1}^{\tau}$ where $\tau\in\mathbb J_{p,1}$

We can choose a sequence of points $\{y_i\}_{i=1}^{\infty}\subset U_p\subset G(p,n)$ such that the following holds.
   \begin{enumerate}[label=(\alph*).]
    \item The inverse rational map $(\Gamma_{p,1}^{\tau})^{-1}:U_p\rightarrow \mathbb C^{p(n-p)}$ is holomorphic on $\{y_i\}_{i=1}^{\infty}$ for each $\tau\in\mathbb J_{p,1}$\,.
    \item The rational map $f_s^{p-1}$ is well-defined on $\{y_i\}_{i=1}^{\infty}$\,.
    \item $\big(e(y_i), f_s^p(y_i),f_s^{p-1}(y_i)\big)\rightarrow v_0\,,\,\,i\rightarrow\infty$\,.
 \end{enumerate}

Write $f_s^{p-1}(y_i)$, $i\geq 1$, in terms of the homogeneous coordinates for $\mathbb{CP}^{N^{p-1}_{s,p,n}}$ as
\begin{equation}
    f_s^{p-1}(y_i)=\big[\cdots,\, P_I(y_i),\, \cdots\, \big]_{I\in \mathbb I_{s,p,n}^{p-1}}\,\,.
\end{equation}
After taking a subsequence, we can assume that there is an index $I^*\in \mathbb I_{s,p,n}^{p-1}$ such that 
\begin{equation}
\big|P_{I^*}(y_i)\big|={\rm max}\left\{\left|P_{I}(y_i)\right|\big|I\in \mathbb I_{s,p,n}^{p-1}\right\}\,\,\,\forall i\geq 1.
\end{equation}
Then $I^*=(n,n-1,\cdots,\widehat {i_1},\cdots,s+1,j_1)$ where $s+1\leq i_1\leq n$ and $1\leq j_1\leq s$.  Take an index $\tau=\left(\begin{matrix}
i_1-s\\
j_1\\
\end{matrix}\right)\in\mathbb J_{p,1}$.

In the following, we will prove that $v_0\in A_{1}^{\tau}$. Since $J_{p,1}^{\tau}$ is defined on $\mathbb C^{p(n-p)}$, it suffices to show that the  set  $  \left\{\left (\Gamma_{p,1}^{\tau}\right)^{-1}(y_i)\big|\,i\geq1\right\}$ is bounded in $\mathbb C^{p(n-p)}$.

In terms of the coordinates $\left(\overrightarrow B^1,\cdots,\overrightarrow B^k,\overrightarrow B^*_k\right)$, write $(\Gamma_{p,1}^{\tau})^{-1}(y_i)$ as
\begin{equation}\label{a1t}
\small
(\Gamma_{p,1}^{\tau})^{-1}(y_i)=\left(a_{i_1j_1}\left((\Gamma_{p,1}^{\tau})^{-1}(y_i)\right),\,\xi^{(1)}_{i_11}\left((\Gamma_{p,1}^{\tau})^{-1}(y_i)\right),\,\cdots\,,\,x^{(2)}_{i^*_1j^*_1}\left((\Gamma_{p,1}^{\tau})^{-1}(y_i)\right)\,,\,\cdots\,\right),\,\,i\geq 1.
\end{equation}

Similarly to Claims I, III, III$^{\prime}$ in Lemma \ref{em}, we can show that for each variable $\xi^{(1)}_{\alpha\beta}$ in (\ref{a1t}) there is an index $I\in \mathbb I_{s,p,n}^{p-1}$  such that 
\begin{equation}
    \left|\xi^{(1)}_{\alpha\beta}\left((\Gamma_{p,1}^{\tau})^{-1}(y_i)\right)\right|=\frac{\big|P_I(y_i)\big|}{\big|P_{I^*}(y_i)\big|}\,,\,\,\,\,\,\,i\geq 1.
\end{equation}
Therefore, the following coordinates of $(\Gamma_{p,1}^{\tau})^{-1}(y_i)$  are uniformly bounded for $i\geq 1$ thanks to the choice of $I^*$.
\begin{equation}\label{xb}
\left\{\begin{aligned}
&\xi^{(1)}_{i_11}\left((\Gamma_{p,1}^{\tau})^{-1}(y_i)\right),\xi^{(1)}_{i_12}\left((\Gamma_{p,1}^{\tau})^{-1}(y_i)\right),\cdots,\xi^{(1)}_{i_1(j_1-1)}\left((\Gamma_{p,1}^{\tau})^{-1}(y_i)\right),\\
&\xi^{(1)}_{i_1(j_1+1)}\left((\Gamma_{p,1}^{\tau})^{-1}(y_i)\right),\cdots,\xi^{(1)}_{i_1s}\left((\Gamma_{p,1}^{\tau})^{-1}(y_i)\right),\xi^{(1)}_{1j_1}\left((\Gamma_{p,1}^{\tau})^{-1}(y_i)\right),\cdots,\\
&\xi^{(1)}_{(i_1-1)j_1}\left((\Gamma_{p,1}^{\tau})^{-1}(y_i)\right),\xi^{(1)}_{(i_1+1)j_1}\left((\Gamma_{p,1}^{\tau})^{-1}(y_i)\right),\cdots,\xi^{(1)}_{pj_1}\left((\Gamma_{p,1}^{\tau})^{-1}(y_i)\right).
\end{aligned}\right\}_{i=1}^{\infty}.
\end{equation}

Notice that $(g^{\sigma_p}_0\circ g^{\sigma_p}_{1})(v_0)\in U_p$ and $\big(e(y_i), f_s^p(y_i),f_s^{p-1}(y_i)\big)\rightarrow v_0\,,\,\,i\rightarrow\infty$. Then the following set is bounded
\begin{equation}
    \left\{z_{\alpha\beta}(y_i)\big|1\leq \alpha\leq p,\, 1\leq\beta\leq s,\,\,{\rm and}\,\,i\geq 1\right\},
\end{equation}
where $z_{\alpha\beta}$ is the coordinate for $U_p$ defined by (\ref{u0}). Since $a_{i_1j_1}\left((\Gamma_{p,1}^{\tau})^{-1}(y_i)\right)=z_{i_1j_1}(y_i)$, we can conclude that the set $\left\{a_{i_1j_1}\left((\Gamma_{p,1}^{\tau})^{-1}(y_i)\right)\right\}_{i=1}^{\infty}$ is bounded.

By formula (\ref{gk}) we have that
\begin{equation}
     x^{(2)}_{i^*_aj^*_b}\left((\Gamma_{p,1}^{\tau})^{-1}(y_i)\right)=z_{i^*_aj^*_b}(y_i)-a_{i_1j_1}\left((\Gamma_{p,1}^{\tau})^{-1}(y_i)\right)\cdot\xi^{(1)}_{i_1j^*_b}\left((\Gamma_{p,1}^{\tau})^{-1}(y_i)\right)\cdot\xi^{(1)}_{i^*_aj_1}\left((\Gamma_{p,1}^{\tau})^{-1}(y_i)\right)
\end{equation}
for $1\leq a\leq p-1$ and $1\leq b\leq s-1$.
Then the following set is bounded \begin{equation}
    \left\{\left.x^{(2)}_{i^*_aj^*_b}\left((\Gamma_{p,1}^{\tau})^{-1}(y_i)\right)\right\vert_{}1\leq a\leq p-1\,,\,\,1\leq b\leq s-1,\,\,{\rm and}\,\,i\geq1\,\right\}.
\end{equation}

As a conclusion, $\left\{ (\Gamma_{p,1}^{\tau})^{-1}(y_i)\right\}_{i=1}^{\infty}$ is a bounded subset of $\mathbb C^{p(n-p)}$, and hence we complete the proof of Lemma \ref{kccord} when $k=1$.
\smallskip

{\noindent\bf {Step 2 ($k=m+1$).}}  Suppose that Lemma \ref{kccord} holds for all $1\leq k\leq m<p$ and we will prove it also holds for $k=m+1$. If $p=s$ and $m=p-1$, then $Y_p^p=Y_{p-1}^p$ and Lemma \ref{kccord} holds trivially. We assume that $p<s$, or $p=s$ and $m<p-1$ in the following. 

Take an arbitrary point $v_0\in Y_{m+1}^p\subset \mathbb {CP}^{N_{p,n}}
\times\mathbb {CP}^{N^{p}_{s,p,n}}\times\cdots\times\mathbb {CP}^{N^{p-m-1}_{s,p,n}}$.  We will show that $v_0$ is in a certain coordinate chart $A_{m+1}^{\kappa}$ where $\kappa\in\mathbb J_{p,m+1}$. Recalling (\ref{blow1}), we can derive the following commutative diagram.
\begin{equation}
\begin{tikzcd}
&Y_{m+1}^p \arrow[hookrightarrow]{r}&[2em]Y_m^p \times\mathbb{CP}^{N^{p-m-1}_{s,p,n}}\arrow{d}{} \\ &&Y_m^p\arrow[leftarrow]{lu}{g^{\sigma_p}_{m+1}}\\ \end{tikzcd}.\vspace{-20pt} \end{equation}
Denote by $\Pi$ the projection of $Y_{m}^p$ to its first factor $\mathbb {CP}^{N_{p,n}}$ (or equivalently $G(p,n)$). 

By induction hypothesis for $k=m$, there is an index $\tau=\left(\begin{matrix}
i_1&i_2&\cdots&i_m\\
j_1&j_2&\cdots&j_m\\
\end{matrix}\right)\in\mathbb J_{p,m}$ such that $g^{\sigma_p}_{m+1}(v_0)\in A^{\tau}_{m}$. We can choose a sequence of points $\{y_i\}_{i=1}^{\infty}\subset A^{\tau}_{m}$ such that the following  holds.
\begin{enumerate}[label=(\alph*).]
    \item The inverse rational map $(g^{\sigma_p}_{m+1})^{-1}$ is well-defined (hence holomorphic) on $\{y_i\}_{i=1}^{\infty}$.
    \item $(g^{\sigma_p}_{m+1})^{-1}(y_i)\rightarrow v_0\,,\,\,i\rightarrow\infty.$
    \item $\left\{(g^{\sigma_p}_{m+1})^{-1}(y_i)\right\}_{i=1}^{\infty}\subset A_{m+1}^{\kappa}$ for each $\kappa\in\mathbb J_{p,m+1}$ (hence the holomorphic map $(J^{\kappa}_{p,m+1})^{-1}$ is well-defined on $\left\{(g^{\sigma_p}_{m+1})^{-1}(y_i)\right\}_{i=1}^{\infty}$).
\end{enumerate}

Since $\lim_{i\rightarrow\infty}y_i=g^{\sigma_p}_{m+1}(v_0)$ and $v_0\in A^{\tau}_{m}$, it is clear that $\{(J^{\tau}_{p,m})^{-1}(y_i)\}_{i=1}^{\infty}$ is  a bounded subset of $A^{\tau}_{m}$. Therefore,
the following set is bounded in $\mathbb C^{p(n-p)}$.
\begin{equation}
    \left\{\overrightarrow B^1\left((J^{\tau}_{p,m})^{-1}(y_i)\right),\cdots,\overrightarrow B^l\left((J^{\tau}_{p,m})^{-1}(y_i)\right),\overrightarrow B^{*}_{m}\left((J^{\tau}_{p,m})^{-1}(y_i)\right)\right\}_{i=1}^{\infty}\,,
\end{equation}
where $\left(\overrightarrow B^1,\cdots,\overrightarrow B^m\right)$ and $\overrightarrow B^{*}_m$ are the holomorphic coordinates defined by (\ref{bp}) and (\ref{b8}) respectively.

Notice that for each pair of integers  $(i^*, j^*)$  such that \begin{equation}\label{dd}
i^*\in\left\{1,2,\cdots,p\}\backslash\{ {i_1},\cdots, {i_m}\right\}\,\, {\rm and}\,\,j^*\in\{1,\cdots,s\}\backslash\left\{ {j_1},\cdots, {j_m}\right\},    
\end{equation}
we can associate  it with a unique index ${I_{i^*j^*}}=(i^*_1,\cdots,i^*_{p-m-1},j^*_1,\cdots,j^*_{m+1})\in \mathbb I^{p-m-1}_{s,p,n}$ such that  \begin{equation}\label{new}
\begin{split}
    &\left \{i^*_1-s,\cdots,i^*_{p-m-1}-s,i^*,{i_1},\cdots, {i_m}\right\}=\left\{1,\cdots,p\right\}\,,\\
  &\left\{j^*_1,\cdots,j^*_{m+1}\right\}=\left\{ {j_1},\cdots, {j_m}, j^*\right\}.
\end{split}
\end{equation}
Similarly to Claims I, III and III$^{\prime}$ in Lemma \ref{em}, we can show that
\begin{equation}\label{hh}
\left|P_{I_{i^*j^*}}(\Pi (y_i))\right|=\left|\prod_{t=1}^{m}a^{m+2-t}_{i_tj_t}\left((J^{\tau}_{p,m})^{-1}(y_i)\right)\right|\cdot\left| x^{(m+1)}_{i^*j^*}\left((J^{\tau}_{p,m})^{-1}(y_i)\right)\right|.
\end{equation}

By taking a subsequence of $\{y_i\}_{i=1}^{\infty}$ we can  assume that there is a pair of integers $(i^*,j^*)$ satisfying (\ref{dd}) such that
\begin{equation}\label{yy}
\left|P_{I_{i^* j^*}}(\Pi(y_i))\right|={\rm max}\left\{|P_{I_{i^*j^*}}(\Pi (y_i))|\big|(i^*,j^*)\,\,{\rm satisfies}\,\, (\ref{dd}) \right\}\,\,\,\forall i\geq 1. 
\end{equation}
It is clear that $\left|P_{I_{i^* j^*}}(\Pi(y_i))\right|>0$ for $(g^{\sigma_p}_{m+1})^{-1}$ is well-defined on $\{y_i\}_{i=1}^{\infty}$.

Take $\kappa=\left(\begin{matrix}
i_1&i_2&\cdots&i_m&i^*\\
j_1&j_2&\cdots&j_m&j^*\\
\end{matrix}\right)\in\mathbb J_{p,m+1}$ and consider the holomorphic coordinate chart $\left(A_{m+1}^{\kappa},(J^{\kappa}_{p,m+1})^{-1}\right)$.
In the following, we will show that $v_0\in A_{m+1}^{\kappa}$. It suffices to prove that the following set is bounded.
\begin{equation}
\begin{split}
&\left\{\overrightarrow B^1\left((g^{\sigma_p}_{m+1}\circ J^{\kappa}_{p,m+1})^{-1}(y_i)\right),\cdots,\overrightarrow B^{m+1}\left((g^{\sigma_p}_{m+1}\circ J^{\kappa}_{p,m+1})^{-1}(y_i)\right),\right.\\
&\,\,\,\,\,\,\,\,\,\,\,\,\,\,\,\,\,\,\,\,\,\left.\overrightarrow B^{*}_{m+1}\left((g^{\sigma_p}_{m+1}\circ J^{\kappa}_{p,m+1})^{-1}(y_i)\right)\right\}_{i=1}^{\infty},
\end{split}
\end{equation}
where $\left(\overrightarrow B^1,\cdots,\overrightarrow B^{m+1},\overrightarrow B^{*}_{m+1}\right)$ are  defined by (\ref{bp}) and (\ref{b8}).

We compare the coordinates $\left(\overrightarrow B^1,\cdots,\overrightarrow B^{m+1},\overrightarrow B^{*}_{m+1}\right)$ to $\left(\overrightarrow B^1,\cdots,\overrightarrow B^m,\overrightarrow B^{*}_m\right)$. It is clear that
\begin{equation}
a_{i^*j^*}\left((g^{\sigma_p}_{m+1}\circ J^{\kappa}_{p,m+1})^{-1}(y_i)\right)=x^{(m+1)}_{i^*j^*}\left((J^{\tau}_{p,m+1})^{-1}(y_i)\right);
\end{equation}
for each pair of integers $(i^{\prime},j^{\prime})$ such that \begin{equation}
  i^{\prime} =i^*\,\,{\rm and}\,\, j^{\prime}\in\{1,\cdots,s\}\backslash\{j_1,\cdots,j_m,j^*\}  ,
\end{equation} or \begin{equation}
    i^{\prime}\in\{1,\cdots,p\}\backslash\{i_1,\cdots,i_m,i^*\}\,\,{\rm and}\,\,j^{\prime}=j^*,
\end{equation}  
we have
\begin{equation}
    \left|\xi^{(m+1)}_{i^{\prime}j^{\prime}}\left((g^{\sigma_p}_{m+1}\circ J^{\kappa}_{p,m+1})^{-1}(y_i)\right)\right|=\frac{\big|P_{I_{i^{\prime}j^{\prime}}}(\Pi(y_i))\big|}{\big|P_{I_{i^*j^*}}(\Pi(y_i))\big|}\,,
\end{equation}
where $I_{i^{\prime}j^{\prime}}$ is defined by (\ref{new}). 
Therefore, thanks to (\ref{yy}) and the fact that the coordinates $\left(\overrightarrow B^1,\cdots,\overrightarrow B^m,\overrightarrow B^{*}_m\right)$ of the points $(J^{\tau}_{p,m+1})^{-1}(y_i)$ are uniformly bounded, we can conclude that  $\left\{\overrightarrow B^{m+1}\left((g^{\sigma_p}_{m+1}\circ J^{\kappa}_{p,m+1})^{-1}(y_i) \right)\right\}_{i=1}^{\infty}$ is bounded.

Moreover, by  (\ref{gk}) we have that
\begin{equation}
  C^{*}_{m}= C^{*}_{m+1}+\left(\prod_{r=1}^{m}a_{i_rj_r}\right)\cdot a_{i^*j^*}\cdot\Xi^T_{m+1}\cdot\Omega_{m+1}\,.
\end{equation}
It is easy to verify that $\left\{\overrightarrow B^{*}_{m+1}\left((g^{\sigma_p}_{m+1}\circ J^{\kappa}_{p,m+1})^{-1}(y_i)\right)\right\}_{i=1}^{\infty}$ is bounded.

We conclude Lemma \ref{kccord} when $k=m+1$. 

By induction,  the proof of Lemma \ref{kccord} is complete. \,\,\,\,$\endpf$
\smallskip

We complete the proof of Lemma \ref{coor2}. \,\,\,\,$\endpf$

\begin{corollary}
 $Y_k^p$ is smooth for $0\leq k\leq p$.
\end{corollary}

\subsection{\texorpdfstring{General cases    }{rr}}\label{section:coverl}

In this subsection, we define the local Van der Waerden  representation in general.  By Remark \ref{dualusd}, we may assume that $2p\leq n\leq 2s$ (then either $p\leq n-s\leq s$ or $n-s<p<s$). 

Define a parameter (rank) $r$ by
\begin{equation}
    r:=\min\{p,n-s,n-p,s\}\,.
\end{equation}
In the following, we will focus on the case when $p\leq n-s$; the case $n-s<p$ can be easily derived by setting $r=n-s$.

\medskip

For $0\leq l\leq r$, define a permutation  $\sigma_l$ by
\begin{equation}
  \sigma_l(k)=\left\{  \begin{matrix}
        &l+1+k\,\,\,\,\,&{\rm when}\,\, &0\leq k\leq p-l-1\,\\
        & p-k\,\,\,\,\,&{\rm when}\,\, &p-l\leq k\leq p
    \end{matrix}\right..
\end{equation}
Restricting $(\ref{sblow})$ to $U_l\subset G(p,n)$, 
we have that
\vspace{-.05in}
\begin{equation}
\footnotesize
\begin{tikzcd}
&Y^l_p\ar{r}{g^{\sigma_l}_p}&Y^{l}_{p-1}\ar{r}{g^{\sigma_l}_{p-1}}&\cdots\ar{r}{g^{\sigma_l}_1}&Y^{l}_0\ar{r}{g^{\sigma_l}_0}&U_l \\
&&(g^{\sigma_l}_0\circ\cdots\circ g^{\sigma_l}_{p-1})^{-1}(S_{\sigma_l(p)}\cap U_l)\ar[hook]{u}&\cdots&(g^{\sigma_l}_0)^{-1}(S_{\sigma_l(1)}\cap U_l)\ar[hook]{u}&S_{\sigma_l(0)}\cap U_l\ar[hook]{u}\\
\end{tikzcd}\vspace{-20pt}
\end{equation}
where for $0\leq k\leq p$,
\begin{equation}
   Y^l_{k}:=Y_k^{\sigma_l}\cap (g^{\sigma_l}_0\circ\cdots\circ g^{\sigma_l}_{k})^{-1}(U_l)\subset  G(p,n)\times\mathbb {CP}^{N^{\sigma_l(0)}_{s,p,n}}\times\cdots\times\mathbb {CP}^{N^{\sigma_l(k)}_{s,p,n}}\,\,\,\,. 
\end{equation}
It is clear that $Y^l_p\cong R_{s,p,n}^{-1}(U_l)$. Moreover, $Y^l_{r-l}\cong Y^l_{r-l+1}\cong\cdots\cong Y^l_{p-l}$ for corresponding blow-up maps $g^{\sigma_l}_{k}$, $r-l+1\leq k\leq p-l$, are the identity maps.

Next, we will define a system of  holomorphic coordinate charts for $Y^l_{k}$, $0\leq k\leq p$.

For $0\leq k\leq r-l-1$, define an index set $\mathbb J_{l,k}$ by
\begin{equation}
\mathbb J_{l,k}:=\left\{\left. \left(
\begin{matrix}
i_1&i_2&\cdots&i_{k+1}\\
j_1&j_2&\cdots&j_{k+1}\\
\end{matrix}\right)\right\vert_{}\footnotesize\begin{matrix}
l+1\leq i_t\leq p\,\,\,{\rm and\,\,}i_{t_1}\neq i_{t_2}\,\,{\rm if\,\,}t_1\neq t_2\,;\,\,\\
s+l+1\leq j_t\leq n\,\,\,{\rm and\,\,}j_{t_1}\neq j_{t_2}\,\,{\rm if\,\,}t_1\neq t_2\,\,\\
\end{matrix}
\right\}.	    
\end{equation}
Each $\tau=\left(\begin{matrix}
i_1&i_2&\cdots&i_{k+1}\\
j_1&j_2&\cdots&j_{k+1}\\
\end{matrix}\right)\in\mathbb J_{l,k}$ determines uniquely integers   $i^*_1,i^*_2,\cdots, i^*_{p-l-k-1}$ and $j^*_1,j^*_2$, $\cdots,$ $j^*_{n-s-l-k-1}$ such that the following holds.
\begin{enumerate}[label=(\alph*).]
    \item $i^*_1<i^*_2<\cdots<i^*_{p-l-k-1}\,$,\,\,$j^*_1<j^*_2<\cdots<j^*_{n-s-l-k-1}\,;$
    \item $\left\{i_1,\cdots,i_{k+1},i^*_1,\cdots,i_{p-l-k-1}^*\right\}=\left\{l+1,l+2,\cdots,p\right\}\,;$
    \item $\left\{j_1,\cdots,j_{k+1},j^*_1,\cdots,j_{n-s-l-k-1}^*\right\}=\left\{s+l+1,s+l+2,\cdots,n\right\}$.
\end{enumerate}
Associate $\tau$ with a complex Euclidean space $\mathbb {C}^{p(n-p)}$ equipped with holomorphic coordinates 
$\left(\widetilde X, \widetilde Y, \widetilde M, \overrightarrow B^1,\cdots,\overrightarrow B^{k+1},\overrightarrow B^*_{k+1}\right)$ as follows. $\widetilde X, \widetilde Y$ are defined by (\ref{ulu}); $\overrightarrow B^j$ is defined by (\ref{qbp}) for $1\leq j\leq k+1$;
\begin{equation}\label{qb8}
\begin{split}
&\overrightarrow B^{*}_{k+1}:=\left(x^{(k+2)}_{i^*_1j^*_1},x^{(k+2)}_{i^*_1j^*_2},\cdots, x^{(k+2)}_{i^*_1j^*_{n-s-l-k-1}},\cdots,x^{(k+2)}_{i^*_2j^*_1},\cdots,x^{(k+2)}_{i^*_2j^*_{n-s-l-k-1}},\cdots,\right.\\
&\,\,\,\,\,\,\,\,\,\,\,\,\,\,\,\,\,\,\,\,\,\,\,\,\,\,\,\,\,\,\,\,\,\,\,\,\,\left.x^{(k+2)}_{i^*_{p-l-k-1}j^*_{1}},\cdots,x^{(k+2)}_{i^*_{p-l-k-1}j^*_{n-s-l-k-1}}\right)\,
\end{split}
\end{equation}
and
\begin{equation}
\,\,\,\,\,\widetilde M:=\left(\begin{matrix}
m^{(k+2)}_{11}&m^{(k+2)}_{12}&\cdots&m^{(k+2)}_{1(s-p+l)}\\
m^{(k+2)}_{21}&m^{(k+2)}_{22}&\cdots&m^{(k+2)}_{2(s-p+l)}\\
\vdots&\vdots&\ddots&\vdots\\
m^{(k+2)}_{l1}&m^{(k+2)}_{l2}&\cdots&m^{(k+2)}_{l(s-p+l)}\\
\end{matrix}\right).\,\,\,\,\,\,\,\,\,\,\,\,\,\,\,\,\,\,\,\,\,\,\,\,\,\,\,\,\,\,\,\,\,\,\,\,\,\,\,\,\,\,\,\,\,\,\,\,\,\,\,\,\,\,\,\,\,\,\,\,\,\,\,\,\,\,\,\,\,\,\,\,\,\,\,\,\,\,\,\,\,\,\,\,\,\,\,\,
\end{equation}	 
Define a holomorphic map $\Gamma_{l,k}^{\tau}:\mathbb C^{p(n-p)}\rightarrow U_l$, $0\leq k\leq r-l-1$,  by
\begin{equation}\label{qgk}
\footnotesize
\begin{split}
&\Gamma_{l,k}^{\tau}\left(\widetilde X, \widetilde Y,\widetilde M,\overrightarrow B^1,\cdots,\overrightarrow B^{k+1},\overrightarrow B^{*}_{k+1}\right):=\left(
\begin{matrix}
\widetilde M &0_{l\times(p-l)}&I_{l\times l}&\widetilde X\\ \widetilde Y&I_{(p-l)\times(p-l)}&0_{(p-l)\times l}&C^*_{k+1}+\sum\limits_{m=1}^{k+1}\left(\prod\limits_{t=1}^{m}b_{i_tj_t}\right)\cdot\Xi_m^T\cdot\Omega_m\\
\end{matrix}\right)\,.
\end{split}
\end{equation}
Here $\Xi_m$, $\Omega_m$ are defined by (\ref{w3}), (\ref{w4}) respectively;  $C^*_{k+1}$  is a $(p-l)\times (n-s-l)$ matrix defined by
\begin{equation}\label{qw}
 C^*_{k+1}=\left(\prod_{t=1}^{k+1} b_{i_tj_t}\right)\cdot \left(\begin{matrix}
&\cdots&0&0&\cdots&0&0&\cdots\\
&\cdots&0&x^{(k+2)}_{i^*_1j^*_1}&\cdots&x^{(k+2)}_{i^*_1j^*_{n-s-l-k-1}}&0&\cdots\\
&\cdots&0&0&\cdots&0&0&\cdots\\
&\ddots&\vdots&\vdots&\ddots&\vdots&\vdots&\ddots\\
&\cdots&0&0&\cdots&0&0&\cdots\\
&\cdots&0&x^{(k+2)}_{i^*_{p-l-k-1}j^*_1}&\cdots&x^{(k+2)}_{i^*_{p-l-k-1}j^*_{n-s-l-k-1}}&0&\cdots\\
&\cdots&0&0&\cdots&0&0&\cdots\\
\end{matrix}\right)\,.\,\,\,
\end{equation}

Define a rational map $J_{l,k}^{\tau}:\mathbb C^{p(n-p)}\dashrightarrow\mathbb {CP}^{N_{p,n}}
\times\mathbb {CP}^{N^{\sigma_l(0)}_{s,p,n}}\times\cdots\times\mathbb {CP}^{N^{\sigma_l(k)}_{s,p,n}}$ by
 \begin{equation}
    J_{l,k}^{\tau}:=\left(e,f_s^{\sigma_l(0)},f_s^{\sigma_l(1)},\cdots,f_s^{\sigma_l(k)}\right)\circ \Gamma_{l,k}^{\tau}\,,\,\,0\leq k\leq r-l-1.
\end{equation}
Similarly to Lemma \ref{qem}, we can show that $J_{l,k}^{\tau}$ is a holomorphic embedding of $\mathbb C^{p(n-p)}$ for each $\tau\in\mathbb J_{l,k}$ and $0\leq k\leq r-l-1$.
Denote by $A_{k}^{\tau}\subset \mathbb {CP}^{N_{p,n}}
\times\mathbb {CP}^{N^{\sigma_l(0)}_{s,p,n}}\times\cdots\times\mathbb {CP}^{N^{\sigma_l(k)}_{s,p,n}}$  be the image of $\mathbb C^{p(n-p)}$ under  $J_{l,k}^{\tau}$.

We skip the indices $r-l\leq k\leq p-l$ for $Y^l_{r-l}\cong Y^l_{r-l+1}\cong\cdots\cong Y^l_{p-l}$.

For $p-l+1\leq k\leq p$, define an index set $\mathbb J_{l,k}$ by
	\begin{equation}\label{jkp}
\mathbb J_{l,k}:=\left\{\left(
\begin{matrix}
i_1&i_2&\cdots&i_{r-l}&\cdots&i_{k-p+r}\\
j_1&j_2&\cdots&j_{r-l}&\cdots&j_{k-p+r}\\
\end{matrix}\right)\rule[-.37in]{0.01in}{.8in}\footnotesize\begin{matrix}
	l+1\leq\, i_t\,\leq p\,\,{\rm for}\,\,1\leq\,t\,\leq r-l\,,\,\,\\
	1\leq\, i_t\,\leq l\,\,{\rm for}\,\,r-l+1\leq\,t\,\leq k\,,\\
	s+l+1\leq\, j_t\,\leq n\,\,{\rm for}\,\,1\leq\,t\,\leq r-l\,,\,\,\\
	1\leq\, j_t\,\leq s-p+l\,\,{\rm for}\,\,r-l+1\leq\,t\,\leq k\,;\\
	i_{t_1}\neq i_{t_2}\,\,{\rm and\,\,} j_{t_1}\neq j_{t_2}\,\,{\rm for\,\,} t_1\neq t_2\\
	\end{matrix}\right\}.
	\end{equation}
Each $\tau=\left(
\begin{matrix}
i_1&i_2&\cdots&i_{r-l}&\cdots&i_k\\
j_1&j_2&\cdots&j_{r-l}&\cdots&j_k\\
\end{matrix}\right)\in\mathbb J_{l,k}$  determines uniquely integers $i^*_1$, $i^*_2$, $\cdots$, $i^*_{p-k}$ and $j^*_1$, $j^*_2$, $\cdots$, $j^*_{s-k}$ such that the following holds.
\begin{enumerate}[label=(\alph*).]
    \item $i^*_1<i^*_2<\cdots<i^*_{p-k}\,$,\,\,$j^*_1<j^*_2<\cdots<j^*_{s-k}\,$;
    \item $\left\{i_{r-l+1},\cdots,i_{k-p+r},i^*_1,\cdots,i_{p-k}^*\right\}$ $=\left\{1,2,\cdots,l\right\}\,$;
    \item $\left\{j_{r-l+1},\cdots,j_{k-p+r},j^*_1,\cdots,j_{s-k}^*\right\}=\left\{1,2,\cdots,s-p+l\right\}\,.$
\end{enumerate}
Associate $\tau$ with a complex Euclidean space $\mathbb {C}^{p(n-p)}$ equipped with holomorphic coordinates 
$\left(\widetilde X, \widetilde Y, \overrightarrow B^1,\cdots,\overrightarrow B^{k-p+r},\overrightarrow B^*_{k-p+r}\right)$  where $\widetilde X, \widetilde Y, \overrightarrow B^1,\cdots,\overrightarrow B^{k-p+r}$, and $\overrightarrow B^*_{k-p+r}$ are defined as follows.  $\widetilde X, \widetilde Y$ are defined by (\ref{ulu}); for $1\leq j\leq r-l$, $\overrightarrow B^{j}$ is defined by (\ref{qbp}); for $r-l+1\leq j\leq k-p+r$,
$\overrightarrow B^{j}$ is defined by (\ref{qbpb}) when $r=p$ and by (\ref{rqbpb}) when $r=n-s$;
\begin{equation}\label{qb7}
\begin{split}
\overrightarrow B^*_{k-p+r}:=&\left(x^{(k-p+r+1)}_{i^*_1j^*_1}, x^{(k-p+r+1)}_{i^*_1j^*_2},\cdots, x^{(k-p+r+1)}_{i^*_1j^*_{s-k}},\cdots,x^{(k-p+r+1)}_{i^*_2j^*_1},\cdots,x^{(k-p+r+1)}_{i^*_2j^*_{s-k}},\cdots,\right.\\
&\,\,\,\,\,\,\left.\cdots,x^{(k-p+r+1)}_{i^*_{p-k}j^*_1},\cdots,x^{(k-p+r+1)}_{i^*_{p-k}j^*_{n-p-k}}\right).
\end{split}
\end{equation}
Define a holomorphic map $\Gamma_{l,k}^{\tau}:\mathbb C^{p(n-p)}\rightarrow U_l$ by
\begin{equation}\label{qgk1}
\begin{split}
&\,\,\,\,\,\,\,\,\,\,\,\,\,\,\,\,\,\,\,\,\,\,\,\,\,\,\,\,\,\,\,\,\,\,\,\,\,\,\,\,\,\,\,\,\,\,\,\,\,\,\,\,\,\,\,\,\,\,\,\,\,\,\,\,\,\Gamma_{l,k}^{\tau}\left(\widetilde X, \widetilde Y,\overrightarrow B^1,\cdots,\overrightarrow B^{k-p+r},\overrightarrow B^{*}_{k-p+r}\right):=\\
&\left(
\begin{matrix}
C^*_{k-p+r}+\sum\limits_{m=r-l+1}^{k-p+r}\left(\prod\limits_{t=r-l+1}^{m}a_{i_{t}j_t}\right)\cdot\Xi_m^T\cdot\Omega_m &0_{l\times(p-l)}&I_{l\times l}&\widetilde X\\ \widetilde Y&I_{(p-l)\times(p-l)}&0_{(p-l)\times l}&\sum\limits_{m=1}^{r-l}\left(\prod\limits_{t=1}^{m}b_{i_tj_t}\right)\cdot\Xi_m^T\cdot\Omega_m\\
\end{matrix}\right)\,.      \end{split}
\end{equation}
Here $\Xi_m$, $\Omega_m$ are defined by (\ref{w3}), (\ref{w4}), (\ref{w5}), (\ref{w2}) when $r=p$ and (\ref{lw1}), (\ref{lw2}), (\ref{lw3}), (\ref{lw4}) when $r=n-s$;
\begin{equation}\label{qw2}
C^*_{k-p+r}=\left(\prod_{t=p-l+1}^{k-p+r} a_{i_tj_t}\right)\cdot \left(\begin{matrix}
&\cdots&0&0&\cdots&0&0&\cdots\\
&\cdots&0&x^{(k-p+r+1)}_{i^*_1j^*_1}&\cdots&x^{(k-p+r+1)}_{i^*_1j^*_{s-k}}&0&\cdots\\
&\ddots&\vdots&\vdots&\ddots&\vdots&\vdots&\ddots\\
&\cdots&0&x^{(k-p+r+1)}_{i^*_{p-k}j^*_1}&\cdots&x^{(k-p+r+1)}_{i^*_{p-k}j^*_{s-k}}&0&\cdots\\
&\cdots&0&0&\cdots&0&0&\cdots\\
\end{matrix}\right).
\end{equation}	 
Define a rational map $J_{l,k}^{\tau}:\mathbb C^{p(n-p)}\rightarrow\mathbb {CP}^{N_{p,n}}
\times\mathbb {CP}^{N^{\sigma_l(0)}_{s,p,n}}\times\cdots\times\mathbb {CP}^{N^{\sigma_l(k)}_{s,p,n}}$ by
\begin{equation}
J_{l,k}^{\tau}:=\left(e,f_s^{\sigma_l(0)},f_s^{\sigma_l(1)},\cdots,f_s^{\sigma_l(k)}\right)\circ \Gamma_{l,k}^{\tau}\,,\,\,p-l+1\leq k\leq p.
\end{equation}
Similarly to Lemma \ref{em}, we can show that $J_{l,k}^{\tau}$ is a holomorphic embedding for $\tau\in\mathbb J_{l,k}$ and $p-l+1\leq k\leq p$. Let $A_{k}^{\tau}\subset \mathbb {CP}^{N_{p,n}}\times\mathbb{CP}^{N^{\sigma_l(0)}_{s,p,n}}\times\cdots\times\mathbb{CP}^{N^{\sigma_l(k)}_{s,p,n}}$  be the image of $\mathbb C^{p(n-p)}$ under $J_{l,k}^{\tau}$.
By a slight abuse of notation, let $\left(J_{l,k}^{\tau}\right)^{-1}:A_{k}^{\tau}\rightarrow\mathbb C^{p(n-p)}$ be the inverse of $J_{l,k}^{\tau}$.

\begin{definition}
For $0\leq l\leq r$, and $0\leq k\leq r-l-1$ or $p-l+1\leq k\leq p$, we call the above defined system of holomorphic  coordinate charts $\left\{\left(A^{\tau}_k,\left(J_{l,k}^{\tau}\right)^{-1}\right)\right\}_{\tau\in\mathbb J_{l,k}}$ the {\it intermediate  Van der Waerden representation} of $Y_k^l$.
\end{definition}

{\bf\noindent Proof of Lemma \ref{qcoor2}.} It suffices to prove the following.

\begin{lemma}\label{lkccord}
For $0\leq k\leq r-l-1$ and $p-l+1\leq k\leq p$, $\left\{\left(A_{k}^{\tau},\left(J_{l,k}^{\tau}\right)^{-1}\right)\right\}_{\tau\in\mathbb J_{l,k}}$ is a holomorphic atlas for $Y^l_k$\,, or equivalently,    $\bigcup_{\tau\in\mathbb J_{l,k}}A_{k}^{\tau}=Y^l_k$.
\end{lemma}

{\noindent\bf Proof of  Lemma \ref{lkccord}.} Similarly to Lemma \ref{kccord}, we will prove by induction on $k$. 

{\noindent\bf {Step I ($k=0$).}} Take an arbitrary point $v_0\in Y_0^l\subset U_l\times\mathbb{CP}^{N_{s,p,n}^{\sigma_l(0)}}$.
We can choose a sequence of points $\{y_i\}_{i=1}^{\infty}\subset U_l$ such that the following holds.
\begin{enumerate}[label=(\alph*).]
    \item The inverse rational map $(\Gamma_{l,0}^{\tau})^{-1}$ is holomorphic on $\{y_i\}_{i=1}^{\infty}$ for each $\tau\in\mathbb J_{l,0}$.
     \item The rational map $f_s^{\sigma_l(0)}$ is well-defined on $\{y_i\}_{i=1}^{\infty}$\,.
    \item $\left(e(y_i), f_s^{\sigma_l(0)}(y_i)\right)\rightarrow v_0\,,\,\,i\rightarrow\infty.$
\end{enumerate}

Consider a subset $\mathcal I\subset \mathbb I_{s,p,n}^{\sigma_l(0)}$ which consists of the following indices, 
\begin{equation}
    I=(\alpha, s+l,s+l-1,\cdots,s+1,\cdots,\widehat{\beta},\cdots,s-p+l+1)
\end{equation} 
where $s+l+1\leq \alpha\leq n$ and $s-p+l+1\leq \beta\leq s$.
Since $ f_s^{\sigma_l(0)}$ is well-defined on $\{y_i\}_{i=1}^{\infty}$, after taking a subsequence we can assume that there is an index $I^*\in \mathcal I$ such that 
\begin{equation}
|P_{I^*}(y_i)|={\rm max}\big\{|P_{I}(y_i)|\big|I\in \mathcal I\big\}>0,\,\,i\geq 1.
\end{equation}
Write $I^*=(j_1, s+l,s+l-1,\cdots,s+1,\cdots,\widehat{i_1},\cdots,s-p+l+1)$; take $\tau=\left(\begin{matrix}
i_1\\
j_1\\
\end{matrix}\right)\in\mathbb J_{l,0}$.

In what follows, we will prove that $v_0\in A_{0}^{\tau}$. It suffices to prove that $\big\{ (\Gamma_{l,0}^{\tau})^{-1}(y_i)\big\}_{i=1}^{\infty}$ is contained in a bounded subset of $\mathbb C^{p(n-p)}$. 

Recall the following holomorphic coordinates of $(\Gamma_{l,0}^{\tau})^{-1}(y_i)$,
\begin{equation}\label{qxb}
\begin{split}
&\overrightarrow B^1\left((\Gamma_{l,0}^{\tau})^{-1}(y_i)\right)=\big(a_{i_1j_1}\left((\Gamma_{l,0}^{\tau})^{-1}(y_i)\right),\xi^{(1)}_{i_1(s+l+1)}\left((\Gamma_{l,0}^{\tau})^{-1}(y_i)\right),\xi^{(1)}_{i_1(s+l+2)}\left((\Gamma_{l,0}^{\tau})^{-1}(y_i)\right),\\
&\,\,\,\,\,\,\,\,\,\,\,\,\,\,\,\,\,\,\,\,\,\,\,\,\cdots,\xi^{(1)}_{i_1(j_1-1)}\left((\Gamma_{l,0}^{\tau})^{-1}(y_i)\right),\xi^{(1)}_{i_1(j_1+1)}\left((\Gamma_{l,0}^{\tau})^{-1}(y_i)\right),\cdots,\xi^{(1)}_{i_1n}\left((\Gamma_{l,0}^{\tau})^{-1}(y_i)\right),\\
&\,\,\,\,\,\,\,\,\,\,\,\,\,\,\,\,\,\,\,\,\,\,\,\,\,\,\,\,\,\,\,\,\,\,\,\xi^{(1)}_{(l+1)j_1}\left((\Gamma_{l,0}^{\tau})^{-1}(y_i)\right),\xi^{(1)}_{(l+2)j_1}\left((\Gamma_{l,0}^{\tau})^{-1}(y_i)\right),\cdots,\xi^{(1)}_{(i_1-1)j_1}\left((\Gamma_{l,0}^{\tau})^{-1}(y_i)\right),\\
&\,\,\,\,\,\,\,\,\,\,\,\,\,\,\,\,\,\,\,\,\,\,\,\,\,\,\,\,\,\,\,\,\,\,\,\,\,\,\,\,\,\,\,\left.\xi^{(1)}_{(i_1+1)j_1}\left((\Gamma_{l,0}^{\tau})^{-1}(y_i)\right),\cdots,\xi^{(1)}_{pj_1}\left((\Gamma_{l,0}^{\tau})(y_i)\right)\right).\\
\end{split}
\end{equation}
Then for each variable $\xi^{(1)}_{\alpha\beta}$ in (\ref{qxb}) there is an index $I\in \mathcal I$  such that 
\begin{equation}
    \left|\xi^{(1)}_{\alpha\beta}\left((\Gamma_{l,0}^{\tau})^{-1}(y_i)\right)\right|=\frac{\big|P_I(y_i)\big|}{\big|P_{I^*}(y_i)\big|}\,.
\end{equation}
By the choice of $I^*$, it is clear that $\left\{\xi^{(1)}_{\alpha\beta}\left((\Gamma_{l,0}^{\tau})^{-1}(y_i)\right)\right\}_{i=1}^{\infty}$ is bounded  for each  $\xi^{(1)}_{\alpha\beta}$ in (\ref{qxb}).

By (\ref{qgk}),  We can conclude that the remaining  coordinates of $(\Gamma_{l,0}^{\tau})^{-1}(y_i)$ in $\mathbb C^{p(n-p)}$ are uniformly bounded as well.
We complete  the proof of Lemma \ref{lkccord} when $k=0$.
\smallskip

{\noindent\bf Step II $(k=m+1\leq r-l-1)$.}  Suppose Lemma \ref{lkccord} holds for all $0\leq k\leq m<r-l-1$. We proceed to prove  it  for $k=m+1$.

Take an arbitrary point $v_0\in Y_{m+1}^l$.  Recalling (\ref{blow1}), we derive the following commutative diagram,
\begin{equation}
\begin{tikzcd}
&Y_{m+1}^l \arrow[hookrightarrow]{r}&[2em]  Y_m^l\times\mathbb{CP}^{N^{\sigma_l(m+1)}_{s,p,n}}\arrow{d}{} \\ &&Y_m^l\arrow[leftarrow]{lu}{g^{\sigma_l}_{m+1}}\\ \end{tikzcd}.\vspace{-20pt} \end{equation}
Denote by $\Pi$ the projection of $Y_{m+1}^l$ to its first factor $\mathbb {CP}^{N_{p,n}}$.

By induction hypothesis for $k=m$, we conclude that there is an index $\tau\in\mathbb J_{l,m}$ such that $g^{\sigma_l}_{m+1}(v_0)\in A^{\tau}_{m}$; write $\tau=\left(\begin{matrix}
i_1&i_2&\cdots&i_{m+1}\\
j_1&j_2&\cdots&j_{m+1}\\
\end{matrix}\right)$. 
We can choose a sequence of points $\{y_i\}_{i=1}^{\infty}\subset A^{\tau}_{m}$ such that the following holds.
\begin{enumerate}[label=(\alph*).]
\item The inverse rational maps $(g^{\sigma_l}_{m+1})^{-1}$ is well-defined holomorphic on $\{y_i\}_{i=1}^{\infty}$.
\item $(g^{\sigma_l}_{m+1})^{-1}(y_i)\rightarrow v_0\,,\,\,i\rightarrow\infty.$
\item The inverse rational map $(J^{\kappa}_{l,m+1})^{-1}$ is well-defined on $\{(g^{\sigma_l}_{m+1})^{-1}(y_i)\}_{i=1}^{\infty}$  for each $\kappa\in\mathbb J_{l,m+1}$.
\end{enumerate}
Since  $\lim_{i\rightarrow\infty}y_i=g^{\sigma_l}_{m+1}(v_0)$,  $\{y_i\}_{i=1}^{\infty}$ is contained in a bounded subset of $A^{\tau}_{m}$; in particular,
the coordinates $\widetilde X,\widetilde Y,\widetilde M, \overrightarrow B^1,\cdots$, $\overrightarrow B^{m+1}$, $\overrightarrow B^{*}_{m+1}$ of $(J^{\tau}_{l,m})^{-1}(y_i)$ are uniformly bounded for $i\geq 1$.

For each pair of integers  $(i^*, j^*)$  such that   \begin{equation}\label{ddl}
i^*\in\left\{l+1,l+2,\cdots,p\right\}\backslash\left\{ {i_1},\cdots, {i_{m+1}}\right\}\,\, {\rm and}\,\,j^*\in\left\{s+l+1,\cdots,n\right\}\backslash\left\{ {j_1},\cdots, {j_{m+1}}\right\},   
\end{equation}
we can associate it with a unique index \begin{equation}
    {I_{i^*j^*}}=(j^*_1,\cdots,j^*_{m+2},s+l,s+l-1,\cdots,s+1, i^*_1,\cdots,i^*_{p-l-m-2})\in \mathbb I^{\sigma(m+1)}_{s,p,n}
\end{equation} 
such that  
\begin{equation}\label{ttl}
\begin{split}
    &\left\{i^*_1+p-s,i^*_2+p-s,\cdots,i^*_{p-l-m-2}+p-s,i^*,i_1,\cdots,i_{m+1}\right\}=\left\{l+1,\cdots,p\right\}\,,\\ &\left\{ j^*_1,\cdots,j^*_{m+2}\right\}=\left\{ j^*,{j_1},\cdots, {j_{m+1}}\right\}.
\end{split}
\end{equation}
Moreover,  we have that
\begin{equation}\label{qhh}
\big|P_{I_{i^*j^*}}(\Pi(y_i))\big|=\left|\prod_{t=1}^{m+1}a^{m+3-t}_{i_tj_t}\left((J^{\tau}_{l,m})^{-1}(y_i)\right)\right|\cdot \left|x^{(m+2)}_{i^*j^*}\left((J^{\tau}_{l,m})^{-1}(y_i)\right)\right|.
\end{equation}
By taking a subsequence we can assume that there is a pair of integers $(i^*,j^*)$ satisfying (\ref{ddl}) such that
\begin{equation}\label{yyl}
\left|P_{I_{i^*j^*}}(\Pi(y_i))\right|={\rm max}\left\{|P_{I_{i^{\prime}j^{\prime}}}(\Pi (y_i))|\big|(i^{\prime},j^{\prime})\,\,{\rm satisfies}\,\, (\ref{ddl}) \right\}>0,\,\, i\geq 1.
\end{equation}
Take $\kappa=\left(\begin{matrix}
i_1&i_2&\cdots&i_{m+1}&i^*\\
j_1&j_2&\cdots&j_{m+1}&j^*\\
\end{matrix}\right)\in\mathbb J_{l,m+1}$ and consider the holomorphic map $J_{l,m+1}^{\kappa}$. 

We will prove that $v_0\in A_{m+1}^{\kappa}$. It suffices to show that $\left\{\left(g^{\sigma_l}_{m+1}\circ J^{\kappa}_{l,m+1}\right)^{-1}(y_i)\right\}_{i=1}^{\infty}$ is bounded in $\mathbb C^{p(n-p)}$. 
Compare the coordinates $\big(\overrightarrow B^1,\cdots,\overrightarrow B^{m+1},\overrightarrow B^{*}_{m+1}\big)$ with  $\big(\overrightarrow B^1,\cdots,\overrightarrow B^{m+2},\overrightarrow B^{*}_{m+2}\big)$. It is clear  that
\begin{equation}
    a_{i^*j^*}\left(\left(g^{\sigma_l}_{m+1}\circ J^{\kappa}_{l,m+1}\right)^{-1}(y_i)\right)=x^{(m+2)}_{i^*j^*}\left((y_i)\right);
\end{equation}  hence $\left\{ a_{i^*j^*}\left(\left(g^{\sigma_l}_{m+1}\circ J^{\kappa}_{l,m+1}\right)^{-1}(y_i)\right)\right\}_{i=1}^{\infty}$ is bounded.  Moreover,
\begin{equation}
    \left|\xi^{(m+2)}_{\alpha\beta}\left(\left(g^{\sigma_l}_{m+1}\circ J^{\kappa}_{l,m+1}\right)^{-1}(y_i)\right)\right|=\frac{\big|P_{I_{\alpha\beta}}(\Pi(y_i))\big|}{\big|P_{I_{i^*j^*}}(\Pi(y_i))\big|}\,;
\end{equation} hence, $\left\{\xi^{(m+2)}_{\alpha\beta}\left(\left(g^{\sigma_l}_{m+1}\circ J^{\kappa}_{l,m+1}\right)^{-1}(y_i)\right)\right\}_{i=1}^{\infty}$ is bounded by (\ref{yyl}).
Recalling  (\ref{qgk}),  (\ref{qw}) we have that
\begin{equation}
 {C}^{*}_{m+1}= C^{*}_{m+2}+\left(\prod_{r=1}^{m+1}a_{i_rj_r}\right)\cdot a_{i^*j^*}\cdot\Xi^T_{m+2}\cdot\Omega_{m+2}\,.
\end{equation}
Then it is easy to verify that $\left\{\overrightarrow B^{*}_{m+2}\left(\left(g^{\sigma_l}_{m+1}\circ J^{\kappa}_{l,m+1}\right)^{-1}(y_i)\right)\right\}_{i=1}^{\infty}$ is bounded.

Therefore, we proved that $\bigcup_{\tau\in\mathbb J_{l,k}}A_{k}^{\tau}=Y^l_k$ when $k=m+1$.

{\noindent\bf Step III $(k\geq p-l)$.} Suppose Lemma \ref{qcoor2} holds for $k\leq r-l-1$. Notice that $Y^l_{r-l}\cong Y^l_{r-l+1}\cong\cdots\cong Y^l_{p-l}$. We shall prove that  Lemma \ref{qcoor2} holds for $p-l+1\leq k\leq p$. The proof is the same as Steps I, II, we omit it here for simplicity.
\smallskip

By induction, we complete the proof of Lemma \ref{lkccord}. \,\,\,\,$\endpf$
\smallskip

We conclude Lemma \ref{qcoor2}. \,\,\,\,$\endpf$

\begin{corollary}\label{itsmoo}
$Y_k^l$ is smooth for $0\leq k\leq p$ and $0\leq l\leq r$.
\end{corollary}

\section{\texorpdfstring{Coordinate charts for the projective bundles }{rr}}\label{section:projbc}

In the following, we define certain convenient local coordinate charts for the projective bundles  $\mathbb P(N_{\mathcal V_{(p,0)}/G(p,n)})$ and $\mathbb P(N_{\mathcal V_{(p-r,r)}/G(p,n)})$. These local coordinate charts will be used in Appendix \ref{section:rigid} to prove the rigidity of the automorphism groups of $\mathcal M_{s,p,n}$ under the assumption that the Picard groups are fixed.

First assume $p<s$.
Define an open subset $U$ of $\mathcal V_{(p,0)}$ by

\begin{equation}
U:=\left\{\left.\underbracedmatrixl{
  x_{11} &x_{12}& \cdots&x_{1(s-p)}  \\
  x_{21} &x_{22}& \cdots&x_{2(s-p)}  \\
  \vdots&\vdots&\ddots&\vdots \\
  x_{p1} &x_{p2}& \cdots&x_{p(s-p)} \\}{(s-p)\,\rm columns}\hspace{-.22in}\begin{matrix}
  &\hfill\tikzmark{g}\\
  \\
  \\
  \\
  &\hfill\tikzmark{h}
  \end{matrix}\,\underbracedmatrix{ 1 & & & \\
   & 1& &\\
  &&\ddots& \\& & &1\\}
   {p\,\rm columns}\hspace{-.2in}
\begin{matrix}
  &\hfill\tikzmark{e}\\
  \\
  \\
  \\
  &\hfill\tikzmark{f}\end{matrix}\hspace{-.1in}
  \begin{matrix}
  &\hfill\tikzmark{a}\\
  \\
  \\
  \\
  &\hfill\tikzmark{b}\end{matrix}\,
  \underbracedmatrixr{ 0  &0& \cdots & 0\\
 0  &0 & \cdots&0\\
  \vdots&\vdots&\ddots&\vdots \\
 0  &0 & \cdots& 0\\} {(n-s) \,\rm columns}\right \vert_{} \begin{matrix}
  x_{uv}\in\mathbb C,\\
  1\leq u\leq p,\\
  1\leq v\leq s-p
  \end{matrix}\right\}.
  \tikz[remember picture,overlay]   \draw[dashed,dash pattern={on 4pt off 2pt}] ([xshift=0.5\tabcolsep,yshift=7pt]a.north) -- ([xshift=0.5\tabcolsep,yshift=-2pt]b.south);\tikz[remember picture,overlay]   \draw[dashed,dash pattern={on 4pt off 2pt}] ([xshift=0.5\tabcolsep,yshift=7pt]e.north) -- ([xshift=0.5\tabcolsep,yshift=-2pt]f.south);\tikz[remember picture,overlay]   \draw[dashed,dash pattern={on 4pt off 2pt}] ([xshift=0.5\tabcolsep,yshift=7pt]g.north) -- ([xshift=0.5\tabcolsep,yshift=-2pt]h.south);
\end{equation}
Locally, the normal bundle $N_{\mathcal V_{(p,0)}/G(p,n)}$ is a product of an open subset of $\mathcal V_{(p,0)}$  and $\mathbb {C}^{p(n-s)}$; equip the open subset $U^N\subset N_{\mathcal V_{(p,0)}/G(p,n)}$ over $U$ with the following local coordinates.
\begin{equation}\label{nbc1}
   \left((\cdots,x_{uv},\cdots)_{\substack{1\leq u\leq p\\1\leq v\leq s-p}}\,,\,\, (\cdots, a_{ij},\cdots)_{\substack{1\leq i\leq p\\s+1\leq j\leq n}}\right)=:(X,A)\,.
\end{equation}
We can identify $U^N$ with an open subset of $G(p,n)$ by
\begin{equation}\label{nin1}
    \begin{split}
        M\,\,\,\,:\,\,\,\,\,\,\,\,\,\,\,\,\,\,\, U^N \,\,\,\,\,&\rightarrow\,\,\,\,\,\,\,\, G(p,n)\\
        (X,A)\,\,\, &\mapsto \,\,\,\,\,\,\, \,M_{(X,A)}
    \end{split}\,\,\,\,\,\,
\end{equation}
where
\begin{equation}\label{mnin1}
 M_{(X,A)}:=\underbracedmatrixl{
  x_{11} &x_{12}& \cdots&x_{1(s-p)}  \\
  x_{21} &x_{22}& \cdots&x_{2(s-p)}  \\
  \vdots&\vdots&\ddots&\vdots \\
  x_{p1} &x_{p2}& \cdots&x_{p(s-p)} \\}{(s-p)\,\rm columns}\hspace{-.22in}\begin{matrix}
  &\hfill\tikzmark{g}\\
  \\
  \\
  \\
  &\hfill\tikzmark{h}
  \end{matrix}\,\underbracedmatrix{ 1 & & & \\
   & 1& &\\
  &&\ddots& \\& & &1\\}
   {p\,\rm columns}\hspace{-.23in}
\begin{matrix}
  &\hfill\tikzmark{e}\\
  \\
  \\
  \\
  &\hfill\tikzmark{f}\end{matrix}\hspace{-.1in}
  \begin{matrix}
  &\hfill\tikzmark{a}\\
  \\
  \\
  \\
  &\hfill\tikzmark{b}\end{matrix}\,
  \underbracedmatrixr{ a_{1(s+1)}  &a_{1(s+2)} & \cdots&a_{1(n-1)} & a_{1n}\\
 a_{2(s+1)}  &a_{2(s+2)} & \cdots&a_{2(n-1)} & a_{2n}\\
  \vdots&\vdots&\ddots&\vdots&\vdots \\
 a_{p(s+1)}  &a_{p(s+2)} & \cdots&a_{p(n-1)} & a_{pn}\\} {(n-s) \,\rm columns}.
  \tikz[remember picture,overlay]   \draw[dashed,dash pattern={on 4pt off 2pt}] ([xshift=0.5\tabcolsep,yshift=7pt]a.north) -- ([xshift=0.5\tabcolsep,yshift=-2pt]b.south);\tikz[remember picture,overlay]   \draw[dashed,dash pattern={on 4pt off 2pt}] ([xshift=0.5\tabcolsep,yshift=7pt]e.north) -- ([xshift=0.5\tabcolsep,yshift=-2pt]f.south);\tikz[remember picture,overlay]   \draw[dashed,dash pattern={on 4pt off 2pt}] ([xshift=0.5\tabcolsep,yshift=7pt]g.north) -- ([xshift=0.5\tabcolsep,yshift=-2pt]h.south);
\end{equation}
It is clear the the above identification gives a natural birational map between the normal bundle $N_{\mathcal V_{(p,0)}/G(p,n)}$ and the grassmannian $G(p,n)$. 
Let $U^P$ be the open subset of $\mathbb P(N_{\mathcal V_{(p,0)}/G(p,n)})$ over $U$. Then $U^P$ is equipped with the following local coordinates. 
\begin{equation}\label{trr1}
   \left((\cdots,x_{uv},\cdots)_{\substack{1\leq u\leq p\\1\leq v\leq s-p}}\,,\,\, [\cdots, a_{ij},\cdots]_{\substack{1\leq i\leq p\\s+1\leq j\leq n}}\right)=:(X,[A\,])\,
\end{equation}
where $[\cdots, a_{ij},\cdots]_{1\leq i\leq p,\,s+1\leq j\leq n}$ are the homogeneous coordinates for the fiber $\mathbb {CP}^{p(n-s)-1}$.

More generally, for a sequence of integers $1\leq i_1<i_2 <\cdots<i_p\leq s$, define an open subset $U_{i_1\cdots i_p}$ of $\mathcal V_{(p,0)}$ by
\begin{equation}
\small
U_{i_1\cdots i_p}:=\left\{\underbracedmatrixl{
 \cdots &\widetilde x_{1(i_1-1)}&1& \cdots&\widetilde x_{1(i_2-1)}&0& \cdots&\widetilde x_{1(i_p-1)}&0& \cdots  \\
  \cdots &\widetilde x_{2(i_1-1)}&0& \cdots&\widetilde x_{1(i_2-1)}&1& \cdots&\widetilde x_{2(i_p-1)}&0& \cdots  \\
  \ddots&\vdots&\vdots&\ddots&\vdots&\vdots&\ddots&\vdots&\vdots&\ddots\\
  \cdots &\widetilde x_{p(i_1-1)}&0& \cdots&\widetilde x_{1(i_2-1)}&0& \cdots&\widetilde x_{p(i_p-1)}&1& \cdots \\}{s\,\rm columns}\hspace{-.22in}
\begin{matrix}
  &\hfill\tikzmark{e}\\
  \\
  \\
  \\
  &\hfill\tikzmark{f}\end{matrix}\hspace{-.1in}
  \begin{matrix}
  &\hfill\tikzmark{a}\\
  \\
  \\
  \\
  &\hfill\tikzmark{b}\end{matrix}\,
  \underbracedmatrixr{ 0  &0& \cdots & 0\\
 0  &0 & \cdots&0\\
  \vdots&\vdots&\ddots&\vdots \\
 0  &0 & \cdots& 0\\} {(n-s) \,\rm columns}\right\}.
  \tikz[remember picture,overlay]   \draw[dashed,dash pattern={on 4pt off 2pt}] ([xshift=0.5\tabcolsep,yshift=7pt]a.north) -- ([xshift=0.5\tabcolsep,yshift=-2pt]b.south);\tikz[remember picture,overlay]   \draw[dashed,dash pattern={on 4pt off 2pt}] ([xshift=0.5\tabcolsep,yshift=7pt]e.north) -- ([xshift=0.5\tabcolsep,yshift=-2pt]f.south);
\end{equation}
Equip the open subset $U_{i_1\cdots i_p}^N\subset N_{\mathcal V_{(p,0)}/G(p,n)}$ over $U_{i_1\cdots i_p}$ with local coordinates
\begin{equation}
   \left(\left(\cdots,\widetilde x_{uv},\cdots\right)_{\substack{1\leq u\leq p,\,1\leq v\leq s\\v\neq i_1,\cdots,i_p}}\,,\,\, (\cdots, \widetilde a_{ij},\cdots)_{1\leq i\leq p,\,s+1\leq j\leq n}\right)=:(\widetilde X,\widetilde A)_{i_1\cdots i_p}\,.
\end{equation}
Identify  $U_{i_1\cdots i_p}^N$ with an open subset of $G(p,n)$ by $M^{i_1\cdots i_p}:U_{i_1\cdots i_p}^N\rightarrow G(p,n)$ where
\begin{equation}\label{tri1}
 M_{\left(\widetilde X,\widetilde A\right)_{i_1\cdots i_p}}:=\underbracedmatrixl{
 \cdots &\widetilde x_{1(i_1-1)}&1& \cdots&\widetilde x_{1(i_p-1)}&0& \cdots  \\
  \cdots &\widetilde x_{2(i_1-1)}&0& \cdots&\widetilde x_{2(i_p-1)}&0& \cdots  \\
  \ddots&\vdots&\vdots&\ddots&\vdots&\vdots&\ddots\\
  \cdots &\widetilde x_{p(i_1-1)}&0& \cdots&\widetilde x_{p(i_p-1)}&1& \cdots \\}{s\,\rm columns}\hspace{-.22in}
\begin{matrix}
  &\hfill\tikzmark{e}\\
  \\
  \\
  \\
  &\hfill\tikzmark{f}\end{matrix}\hspace{-.1in}
  \begin{matrix}
  &\hfill\tikzmark{a}\\
  \\
  \\
  \\
  &\hfill\tikzmark{b}\end{matrix}\,
  \underbracedmatrixr{\widetilde a_{1(s+1)}  & \cdots& \widetilde a_{1n}\\
 \widetilde a_{2(s+1)}   & \cdots& \widetilde a_{2n}\\
  \vdots&\ddots&\vdots \\
 \widetilde a_{p(s+1)}  & \cdots& \widetilde a_{pn}\\} {(n-s) \,\rm columns}.
  \tikz[remember picture,overlay]   \draw[dashed,dash pattern={on 4pt off 2pt}] ([xshift=0.5\tabcolsep,yshift=7pt]a.north) -- ([xshift=0.5\tabcolsep,yshift=-2pt]b.south);\tikz[remember picture,overlay]   \draw[dashed,dash pattern={on 4pt off 2pt}] ([xshift=0.5\tabcolsep,yshift=7pt]e.north) -- ([xshift=0.5\tabcolsep,yshift=-2pt]f.south);
\end{equation}
Let $U_{i_1\cdots i_p}^P$ be the open subset of $\mathbb P(N_{\mathcal V_{(p,0)}/G(p,n)})$ over $U_{i_1\cdots i_p}$  equipped with  local coordinates
\begin{equation}\label{trr12}
  \left(\left(\cdots,\widetilde x_{uv},\cdots\right)_{\substack{1\leq u\leq p,\,1\leq v\leq s\\v\neq i_1,\cdots,i_p}}\,,\,\, [\cdots, \widetilde a_{ij},\cdots]_{1\leq i\leq p,\,s+1\leq j\leq n}\right)=:(\widetilde X,\big[\widetilde A\,\big])_{i_1\cdots i_p}\,,
\end{equation}
where $[\cdots, \widetilde a_{ij},\cdots]$ are the homogeneous coordinates for the fiber. 

It is easy to verify that $(X,[A\,])$ and $\left(\widetilde X,[\widetilde A\,]\right)_{i_1\cdots i_p}$ represent the same point of  $\mathbb P(N_{\mathcal V_{(p,0)}/G(p,n)})$  if and only if the matrices $M_{(X,A)}$ and $M_{\left(\widetilde X,\widetilde A\right)_{i_1\cdots i_p}}$ represent the same point of $G(p,n)$ up to the $\mathbb C^*$-action $\psi_{s,p,n}$, that is,  there exists a matrix $W\in GL(p,\mathbb C)$ and  $\lambda\in\mathbb C^*$ such that 
\begin{equation}
    M_{(X,A)}=W\cdot M_{\left(\widetilde X,\widetilde A\right)_{i_1\cdots i_p}}\cdot \left(\begin{matrix} I_{s\times s}&0\\ 0&\lambda\cdot I_{(n-s)\times (n-s)}\\ \end{matrix} \right)\,.
\end{equation}
\smallskip

Similarly,  we can assign local coordinate charts for $\mathbb P(N_{\mathcal V_{(p-r,r)}/G(p,n)})$ as follows. 

Let $n-s<p<s$. Define an open subset $\widehat U$ of $\mathcal V_{(p-r,r)}\cong G(s+p-n,s)$  by

\begin{equation}\label{nbc2}
\small
\widehat U:=\left\{\underbracedmatrixl{
0 &0& \cdots&0 \\
 0 &0& \cdots&0 \\
  \vdots&\vdots&\ddots&\vdots \\
  0 &0& \cdots&0 \\ \hdashline[4pt/2pt]
  x_{(n-s+1)1} &x_{(n-s+1)2}& \cdots&x_{(n-s+1)(n-p)}  \\
  x_{(n-s+2)1} &x_{(n-s+2)2}& \cdots&x_{(n-s+2)(n-p)}  \\
  \vdots&\vdots&\ddots&\vdots \\ 
  x_{p1} &x_{p2}& \cdots&x_{p(n-p)} \\}{(n-p)\,\rm columns}\hspace{-.22in}\begin{matrix}
  &\hfill\tikzmark{g}\\
  \\
  \\
  \\
  \\
  \\
  \\
  \\
  &\hfill\tikzmark{h}
  \end{matrix}\,\underbracedmatrix{  0  &0& \cdots & 0\\
 0  &0 & \cdots&0\\
  \vdots&\vdots&\ddots&\vdots \\  0  &0 & \cdots& 0\\ \hdashline[4pt/2pt]
 1 & & & \\
   & 1& &\\
  &&\ddots& \\& & &1\\}
   {(s+p-n)\,\rm columns}\hspace{-.2in}
\begin{matrix}
  &\hfill\tikzmark{e}\\
  \\
  \\
  \\
   \\
  \\
  \\
  \\
  &\hfill\tikzmark{f}\end{matrix}\hspace{-.1in}
  \begin{matrix}
  &\hfill\tikzmark{a}\\
  \\ \\
  \\
  \\
  \\
  \\
  \\
  &\hfill\tikzmark{b}\end{matrix}\,
  \underbracedmatrixr{
 1 & & & \\
   & 1& &\\
  &&\ddots& \\& & &1\\ \hdashline[4pt/2pt]
  0  &0& \cdots & 0\\
 0  &0 & \cdots&0\\
  \vdots&\vdots&\ddots&\vdots \\
 0  &0 & \cdots& 0\\} {(n-s) \,\rm columns}\right\}.
  \tikz[remember picture,overlay]   \draw[dashed,dash pattern={on 4pt off 2pt}] ([xshift=0.5\tabcolsep,yshift=7pt]a.north) -- ([xshift=0.5\tabcolsep,yshift=-2pt]b.south);\tikz[remember picture,overlay]   \draw[dashed,dash pattern={on 4pt off 2pt}] ([xshift=0.5\tabcolsep,yshift=7pt]e.north) -- ([xshift=0.5\tabcolsep,yshift=-2pt]f.south);\tikz[remember picture,overlay]   \draw[dashed,dash pattern={on 4pt off 2pt}] ([xshift=0.5\tabcolsep,yshift=7pt]g.north) -- ([xshift=0.5\tabcolsep,yshift=-2pt]h.south);
\end{equation}
Equip  the open subset $\widehat U^N\subset N_{\mathcal V_{(p,0)}/G(p,n)}$ over $U$ with local coordinates 
\begin{equation}
   \left((\cdots,x_{uv},\cdots)_{\substack{n-s+1\leq u\leq p\\1\leq v\leq n-p}}\,,\,\, (\cdots, a_{ij},\cdots)_{\substack{1\leq i\leq n-s\\1\leq j\leq n-p}}\right)=:(X,A)\,.
\end{equation}
Identify $\widehat U^N$ with an open subset of $G(p,n)$ by $ \widehat M:\widehat U^N \rightarrow G(p,n)$
where
\begin{equation}\label{nin2}
\widehat M_{(X,A)}:=\underbracedmatrixl{
a_{11} &a_{12}& \cdots&a_{1(n-p)} \\
a_{21} &a_{22}& \cdots&a_{2(n-p)} \\
  \vdots&\vdots&\ddots&\vdots \\
  a_{(n-s)1} &a_{(n-s)2}& \cdots&a_{(n-s)(n-p)} \\ \hdashline[4pt/2pt]
  x_{(n-s+1)1} &x_{(n-s+1)2}& \cdots&x_{(n-s+1)(n-p)}  \\
  x_{(n-s+2)1} &x_{(n-s+2)2}& \cdots&x_{(n-s+2)(n-p)}  \\
  \vdots&\vdots&\ddots&\vdots \\ 
  x_{p1} &x_{p2}& \cdots&x_{p(n-p)} \\}{(n-p)\,\rm columns}\hspace{-.22in}\begin{matrix}
  &\hfill\tikzmark{g}\\
  \\
  \\
  \\
  \\
  \\
  \\
  \\
  &\hfill\tikzmark{h}
  \end{matrix}\,\underbracedmatrix{  0  &0& \cdots & 0\\
 0  &0 & \cdots&0\\
  \vdots&\vdots&\ddots&\vdots \\  0  &0 & \cdots& 0\\ \hdashline[4pt/2pt]
 1 & & & \\
   & 1& &\\
  &&\ddots& \\& & &1\\}
   {(s+p-n)\,\rm columns}\hspace{-.2in}
\begin{matrix}
  &\hfill\tikzmark{e}\\
  \\
  \\
  \\
   \\
  \\
  \\
  \\
  &\hfill\tikzmark{f}\end{matrix}\hspace{-.1in}
  \begin{matrix}
  &\hfill\tikzmark{a}\\
  \\ \\
  \\
  \\
  \\
  \\
  \\
  &\hfill\tikzmark{b}\end{matrix}\,
  \underbracedmatrixr{
 1 & & & \\
   & 1& &\\
  &&\ddots& \\& & &1\\ \hdashline[4pt/2pt]
  0  &0& \cdots & 0\\
 0  &0 & \cdots&0\\
  \vdots&\vdots&\ddots&\vdots \\
 0  &0 & \cdots& 0\\} {(n-s) \,\rm columns}.
  \tikz[remember picture,overlay]   \draw[dashed,dash pattern={on 4pt off 2pt}] ([xshift=0.5\tabcolsep,yshift=7pt]a.north) -- ([xshift=0.5\tabcolsep,yshift=-2pt]b.south);\tikz[remember picture,overlay]   \draw[dashed,dash pattern={on 4pt off 2pt}] ([xshift=0.5\tabcolsep,yshift=7pt]e.north) -- ([xshift=0.5\tabcolsep,yshift=-2pt]f.south);\tikz[remember picture,overlay]   \draw[dashed,dash pattern={on 4pt off 2pt}] ([xshift=0.5\tabcolsep,yshift=7pt]g.north) -- ([xshift=0.5\tabcolsep,yshift=-2pt]h.south);
\end{equation}
It is clear the the above identification gives a natural birational map between the normal bundle $N_{\mathcal V_{(p-r,r)}/G(p,n)}$ and the grassmannian $G(p,n)$. 
Let $\widehat U^P$ be the open subset of $\mathbb P(N_{\mathcal V_{(p-r,r)}/G(p,n)})$ over $\widehat U$ equipped with local coordinates
\begin{equation}\label{trr2}
   \left((\cdots,x_{uv},\cdots)_{\substack{n-s+1\leq u\leq p\\1\leq v\leq n-p}}\,,\,\, [\cdots, a_{ij},\cdots]_{\substack{1\leq i\leq n-s\\1\leq j\leq n-p}}\right)=:(X,[A\,])\,,
\end{equation}
where $[\cdots, a_{ij},\cdots]_{1\leq i\leq n-s,\,1\leq j\leq n-p}$ are the homogeneous coordinates for the fiber. 

For a sequence of integers $1\leq i_1<i_2 <\cdots<i_{s+p-n}\leq s$, define $\widehat U_{i_1\cdots i_{s+p-n}}\subset\mathcal V_{(p-r,r)}$ by

\begin{equation}
\footnotesize
\left\{ \underbracedmatrixl{
  \cdots &0&0& \cdots&0&0& \cdots&0&0& \cdots  \\
  \ddots&\vdots&\vdots&\ddots&\vdots&\vdots&\ddots&\vdots&\vdots&\ddots\\
  \cdots &0&0& \cdots&0&0& \cdots&0&0& \cdots \\ \hdashline[4pt/2pt]
 \cdots &\widetilde x_{(n-s+1)(i_1-1)}&1& \cdots&\widetilde x_{(n-s+1)(i_2-1)}&0& \cdots&\widetilde x_{(n-s+1)(i_{s+p-n}-1)}&0& \cdots  \\
  \cdots &\widetilde x_{(n-s+2)(i_1-1)}&0& \cdots&\widetilde x_{(n-s+2)(i_2-1)}&1& \cdots&\widetilde x_{(n-s+2)(i_{s+p-n}-1)}&0& \cdots  \\
  \ddots&\vdots&\vdots&\ddots&\vdots&\vdots&\ddots&\vdots&\vdots&\ddots\\
  \cdots &\widetilde x_{p(i_1-1)}&0& \cdots&\widetilde x_{1(i_2-1)}&0& \cdots&\widetilde x_{p(i_{s+p-n}-1)}&1& \cdots \\}{s\,\rm columns}\hspace{-.16in}
\begin{matrix}
  &\hfill\tikzmark{e}\\
  \\
  \\
  \\
  \\
  \\
  \\
  \\
  &\hfill\tikzmark{f}\end{matrix}\hspace{-.1in}
  \begin{matrix}
  &\hfill\tikzmark{a}\\
  \\
  \\
  \\
  \\
  \\
  \\
  \\
  &\hfill\tikzmark{b}\end{matrix}
  \underbracedmatrixr{ 1  & & \\
   &\ddots& \\& &1\\ \hdashline[4pt/2pt]
  0  & \cdots & 0\\
 0   & \cdots&0\\
  \vdots&\ddots&\vdots \\
 0   & \cdots& 0\\} {(n-s) \,\rm columns}\right\}.
  \tikz[remember picture,overlay]   \draw[dashed,dash pattern={on 4pt off 2pt}] ([xshift=0.5\tabcolsep,yshift=7pt]a.north) -- ([xshift=0.5\tabcolsep,yshift=-2pt]b.south);\tikz[remember picture,overlay]   \draw[dashed,dash pattern={on 4pt off 2pt}] ([xshift=0.5\tabcolsep,yshift=7pt]e.north) -- ([xshift=0.5\tabcolsep,yshift=-2pt]f.south);
\end{equation}
Equip  the open subset $\widehat U_{i_1\cdots i_{s+p-n}}^N\subset N_{\mathcal V_{(p-r,r)}/G(p,n)}$ over $\widehat U_{i_1\cdots i_p}$  with local coordinates
\begin{equation}\label{trr21}
   \left(\left(\cdots,\widetilde x_{uv},\cdots\right)_{\substack{n-s+1\leq u\leq p\\1\leq v\leq s\\v\neq i_1,\cdots,i_{s+p-n}}}\,,\,\, (\cdots, \widetilde a_{ij},\cdots)_{\substack{1\leq i\leq n-s\\1\leq j\leq s\\j\neq i_1,\cdots,i_{s+p-n}}}\right)=:(\widetilde X,\widetilde A)_{i_1\cdots i_{s+p-n}}\,.
\end{equation}
Identify $\widehat U_{i_1\cdots i_{s+p-n}}^N$ with an open subset of $G(p,n)$ by $\widehat M: \widehat U_{i_1\cdots i_{s+p-n}}^N \rightarrow G(p,n)$ where 
\begin{equation}\label{tri2}
\small
\widehat M_{(\widetilde X,\widetilde A)_{i_1\cdots i_{s+p-n}}}:=\underbracedmatrixl{
  \cdots &a_{1(i_1-1)}&0&  \cdots&a_{1(i_{s+p-n}-1)}&0& \cdots  \\
   \cdots &a_{2(i_1-1)}&0&  \cdots&a_{2(i_{s+p-n}-1)}&0& \cdots  \\
  \ddots&\vdots&\vdots&\ddots&\vdots&\vdots&\ddots\\
  \cdots &a_{(n-s)(i_1-1)}&0&  \cdots&a_{(n-s)(i_{s+p-n}-1)}&0& \cdots \\ \hdashline[4pt/2pt]
 \cdots &\widetilde x_{(n-s+1)(i_1-1)}&1&  \cdots&\widetilde x_{(n-s+1)(i_{s+p-n}-1)}&0& \cdots  \\
  \cdots &\widetilde x_{(n-s+2)(i_1-1)}&0&  \cdots&\widetilde x_{(n-s+2)(i_{s+p-n}-1)}&0& \cdots  \\
  \ddots&\vdots&\vdots&\ddots&\vdots&\vdots&\ddots\\
  \cdots &\widetilde x_{p(i_1-1)}&0&  \cdots&\widetilde x_{p(i_{s+p-n}-1)}&1& \cdots \\}{s\,\rm columns}\hspace{-.16in}
\begin{matrix}
  &\hfill\tikzmark{e}\\
  \\
  \\
  \\
  \\
  \\
  \\
  \\
  &\hfill\tikzmark{f}\end{matrix}\hspace{-.1in}
  \begin{matrix}
  &\hfill\tikzmark{a}\\
  \\
  \\
  \\
  \\
  \\
  \\
  \\
  &\hfill\tikzmark{b}\end{matrix}
  \underbracedmatrixr{ 1  & &&\\
  &1  &&\\
   &&\ddots& \\&& &1\\ \hdashline[4pt/2pt]
  0  &0& \cdots & 0\\
 0   & 0&\cdots&0\\
  \vdots&\vdots&\ddots&\vdots \\
 0   & 0&\cdots& 0\\} {(n-s) \,\rm columns}.
  \tikz[remember picture,overlay]   \draw[dashed,dash pattern={on 4pt off 2pt}] ([xshift=0.5\tabcolsep,yshift=7pt]a.north) -- ([xshift=0.5\tabcolsep,yshift=-2pt]b.south);\tikz[remember picture,overlay]   \draw[dashed,dash pattern={on 4pt off 2pt}] ([xshift=0.5\tabcolsep,yshift=7pt]e.north) -- ([xshift=0.5\tabcolsep,yshift=-2pt]f.south);
\end{equation}
Let $\widehat U_{i_1\cdots i_{s+p-n}}^P$ be the open subset of $\mathbb P(N_{\mathcal V_{(p-r,r)}/G(p,n)})$ over $\widehat U_{i_1\cdots i_{s+p-n}}$ equipped with local coordinates
\begin{equation}\label{trr3}
   \left(\left(\cdots,\widetilde x_{uv},\cdots\right)_{\substack{n-s+1\leq u\leq p\\1\leq v\leq s\\v\neq i_1,\cdots,i_{s+p-n}}}\,,\,\, [\cdots, \widetilde a_{ij},\cdots]_{\substack{1\leq i\leq n-s\\1\leq j\leq s\\j\neq i_1,\cdots,i_{s+p-n}}}\right)=:(\widetilde X,\big[\widetilde A\,\big])_{i_1\cdots i_{s+p-n}}\,.
\end{equation}

Similarly,  $(X,[A\,])$ and $\left(\widetilde X,[\widetilde A\,]\right)_{i_1\cdots i_{s+p-n}}$ represent the same point of  $\mathbb P(N_{\mathcal V_{(p-r,r)}/G(p,n)})$  if and only if the matrices $\widehat M_{(X,A)}$ and $\widehat M_{\left(\widetilde X,\widetilde A\right)_{i_1\cdots i_{s+p-n}}}$ represent the same point of $G(p,n)$ up to the $\mathbb C^*$-action $\psi_{s,p,n}$, that is,  there exists a matrix $W\in GL(p,\mathbb C)$ and  $\lambda\in\mathbb C^*$ such that 
\begin{equation}
    \widehat M_{(X,A)}=W\cdot \widehat M_{\left(\widetilde X,\widetilde A\right)_{i_1\cdots i_{s+p-n}}}\cdot \left(\begin{matrix} I_{s\times s}&0\\ 0&\lambda\cdot I_{(n-s)\times (n-s)}\\ \end{matrix} \right)\,.
\end{equation}
\smallskip

Next assume that $p<n-s\leq s$.  Define  $\widehat U\subset\mathcal V_{(p-r,r)}\cong G(p,n-s)$ by
\begin{equation}
\widehat U:=\left\{\underbracedmatrixl{ 0  &0& \cdots & 0\\
 0  &0 & \cdots&0\\
  \vdots&\vdots&\ddots&\vdots \\
 0  &0 & \cdots& 0\\
  }{s\,\rm columns}\hspace{-.22in}\begin{matrix}
  &\hfill\tikzmark{g}\\
  \\
  \\
  \\
  &\hfill\tikzmark{h}
  \end{matrix}\hspace{-.1in}
  \begin{matrix}
  &\hfill\tikzmark{a}\\
  \\
  \\
  \\
  &\hfill\tikzmark{b}\end{matrix}\underbracedmatrix{ 1 & & & \\
   & 1& &\\
  &&\ddots& \\& & &1\\}
   {p\,\rm columns}\hspace{-.2in}
\begin{matrix}
  &\hfill\tikzmark{e}\\
  \\
  \\
  \\
  &\hfill\tikzmark{f}\end{matrix}
  \underbracedmatrixr{x_{1(s+p+1)} &x_{1(s+p+2)}& \cdots&x_{1n}  \\
  x_{2(s+p+1)} &x_{2(s+p+2)}& \cdots&x_{2n}  \\
  \vdots&\vdots&\ddots&\vdots \\
  x_{p(s+p+1)} &x_{p(s+p+2)}& \cdots&x_{pn} \\} {(n-s-p) \,\rm columns}\right\}.
  \tikz[remember picture,overlay]   \draw[dashed,dash pattern={on 4pt off 2pt}] ([xshift=0.5\tabcolsep,yshift=7pt]a.north) -- ([xshift=0.5\tabcolsep,yshift=-2pt]b.south);\tikz[remember picture,overlay]   \draw[dashed,dash pattern={on 4pt off 2pt}] ([xshift=0.5\tabcolsep,yshift=7pt]e.north) -- ([xshift=0.5\tabcolsep,yshift=-2pt]f.south);\tikz[remember picture,overlay]   \draw[dashed,dash pattern={on 4pt off 2pt}] ([xshift=0.5\tabcolsep,yshift=7pt]g.north) -- ([xshift=0.5\tabcolsep,yshift=-2pt]h.south);
\end{equation}
Equip  the open subset $\widehat U^N\subset N_{\mathcal V_{(p-r,r)}/G(p,n)}$ over $\widehat U$ with local coordinates
\begin{equation}\label{nbc3}
   \left((\cdots,x_{uv},\cdots)_{\substack{1\leq u\leq p\\s+p+1\leq v\leq n}}\,,\,\, (\cdots, a_{ij},\cdots)_{\substack{1\leq i\leq p\\1\leq j\leq s}}\right)=:(X,A)\,.
\end{equation}
Identify  $\widehat U^N$ with an open subset of $G(p,n)$ by $ \widehat M:\widehat  U^N \rightarrow G(p,n)$ where
\begin{equation}\label{nin3}
\widehat M_{(X,A)}:=\underbracedmatrixl{a_{11}  &a_{12} & \cdots&a_{1(s-1)} & a_{1s}\\
 a_{21}  &a_{22} & \cdots&a_{2(s-1)} & a_{2s}\\
  \vdots&\vdots&\ddots&\vdots&\vdots \\
 a_{p1}  &a_{p2} & \cdots&a_{p(s-1)} & a_{ps}\\
  }{s\,\rm columns}\hspace{-.22in}\begin{matrix}
  &\hfill\tikzmark{g}\\
  \\
  \\
  \\
  &\hfill\tikzmark{h}
  \end{matrix}\hspace{-.1in}
  \begin{matrix}
  &\hfill\tikzmark{a}\\
  \\
  \\
  \\
  &\hfill\tikzmark{b}\end{matrix}\underbracedmatrix{ 1 & & & \\
   & 1& &\\
  &&\ddots& \\& & &1\\}
   {p\,\rm columns}\hspace{-.2in}
\begin{matrix}
  &\hfill\tikzmark{e}\\
  \\
  \\
  \\
  &\hfill\tikzmark{f}\end{matrix}
  \underbracedmatrixr{x_{1(s+p+1)} &x_{1(s+p+2)}& \cdots&x_{1n}  \\
  x_{2(s+p+1)} &x_{2(s+p+2)}& \cdots&x_{2n}  \\
  \vdots&\vdots&\ddots&\vdots \\
  x_{p(s+p+1)} &x_{p(s+p+2)}& \cdots&x_{pn} \\} {(n-s-p) \,\rm columns}.
  \tikz[remember picture,overlay]   \draw[dashed,dash pattern={on 4pt off 2pt}] ([xshift=0.5\tabcolsep,yshift=7pt]a.north) -- ([xshift=0.5\tabcolsep,yshift=-2pt]b.south);\tikz[remember picture,overlay]   \draw[dashed,dash pattern={on 4pt off 2pt}] ([xshift=0.5\tabcolsep,yshift=7pt]e.north) -- ([xshift=0.5\tabcolsep,yshift=-2pt]f.south);\tikz[remember picture,overlay]   \draw[dashed,dash pattern={on 4pt off 2pt}] ([xshift=0.5\tabcolsep,yshift=7pt]g.north) -- ([xshift=0.5\tabcolsep,yshift=-2pt]h.south);
\end{equation}
It is clear the the above identification gives a  birational map between  $N_{\mathcal V_{(p-r,r)}/G(p,n)}$ and $G(p,n)$. 
Let $\widehat U^P$ be the open subset of $\mathbb P(N_{\mathcal V_{(p-r,r)}/G(p,n)})$ over $U$ equipped with local coordinates
\begin{equation}\label{trr41}
   \left((\cdots,x_{uv},\cdots)_{\substack{1\leq u\leq p\\s+p+1\leq v\leq n}}\,,\,\, [\cdots, a_{ij},\cdots]_{\substack{1\leq i\leq p\\1\leq j\leq s}}\right)=:(X,[A\,])\,.
\end{equation}

For a sequence of integers $s+1\leq i_1<\cdots<i_p\leq n$, define $\widehat U_{i_1\cdots i_p}\subset\mathcal V_{(p-r,r)}$ by

\begin{equation}
\widehat U_{i_1\cdots i_p}:=\left\{\underbracedmatrixl{ 0  &0& \cdots & 0\\
 0  &0 & \cdots&0\\
  \vdots&\vdots&\ddots&\vdots \\
 0  &0 & \cdots& 0\\}{s\,\rm columns}\hspace{-.22in}
\begin{matrix}
  &\hfill\tikzmark{e}\\
  \\
  \\
  \\
  &\hfill\tikzmark{f}\end{matrix}\hspace{-.1in}
  \begin{matrix}
  &\hfill\tikzmark{a}\\
  \\
  \\
  \\
  &\hfill\tikzmark{b}\end{matrix}\,
  \underbracedmatrixr{ 
 \cdots &\widetilde x_{1(i_1-1)}&1& \cdots&\widetilde x_{1(i_2-1)}&0& \cdots&\widetilde x_{1(i_p-1)}&0& \cdots  \\
  \cdots &\widetilde x_{2(i_1-1)}&0& \cdots&\widetilde x_{1(i_2-1)}&1& \cdots&\widetilde x_{2(i_p-1)}&0& \cdots  \\
  \ddots&\vdots&\vdots&\ddots&\vdots&\vdots&\ddots&\vdots&\vdots&\ddots\\
  \cdots &\widetilde x_{p(i_1-1)}&0& \cdots&\widetilde x_{1(i_2-1)}&0& \cdots&\widetilde x_{p(i_p-1)}&1& \cdots \\} {(n-s) \,\rm columns}\right\}.
  \tikz[remember picture,overlay]   \draw[dashed,dash pattern={on 4pt off 2pt}] ([xshift=0.5\tabcolsep,yshift=7pt]a.north) -- ([xshift=0.5\tabcolsep,yshift=-2pt]b.south);\tikz[remember picture,overlay]   \draw[dashed,dash pattern={on 4pt off 2pt}] ([xshift=0.5\tabcolsep,yshift=7pt]e.north) -- ([xshift=0.5\tabcolsep,yshift=-2pt]f.south);
\end{equation}
Equip the open subset $\widehat U_{i_1\cdots i_p}^N\subset N_{\mathcal V_{(p-r,r)}/G(p,n)}$ over $\widehat U_{i_1\cdots i_p}$ with  local coordinates
\begin{equation}
   \left(\left(\cdots,\widetilde x_{uv},\cdots\right)_{\substack{1\leq u\leq p\\s+1\leq v\leq n\\v\neq i_1,\cdots,i_p}}\,,\,\, (\cdots, \widetilde a_{ij},\cdots)_{\substack{1\leq i\leq p\\1\leq j\leq s}}\right)=:(\widetilde X,\widetilde A)_{i_1\cdots i_p}\,.
\end{equation}
Identify $\widehat U_{i_1\cdots i_p}^N$ with an open subset of $G(p,n)$ by $\widehat M:\widehat U_{i_1\cdots i_p}^N \rightarrow G(p,n)$ where 
\begin{equation}\label{tri3}
\widehat  M_{\left(\widetilde X,\widetilde A\right)_{i_1\cdots i_p}}:=\underbracedmatrixl{
  \widetilde a_{11}  & \cdots& \widetilde a_{1s}\\
 \widetilde a_{21}   & \cdots& \widetilde a_{2s}\\
  \vdots&\ddots&\vdots \\
 \widetilde a_{p1}  & \cdots& \widetilde a_{ps}\\} {s\,\rm columns}\hspace{-.22in}
\begin{matrix}
  &\hfill\tikzmark{e}\\
  \\
  \\
  \\
  &\hfill\tikzmark{f}\end{matrix}\hspace{-.1in}
  \begin{matrix}
  &\hfill\tikzmark{a}\\
  \\
  \\
  \\
  &\hfill\tikzmark{b}\end{matrix}\,
  \underbracedmatrixr{\cdots &\widetilde x_{1(i_1-1)}&1& \cdots&\widetilde x_{1(i_p-1)}&0& \cdots  \\
  \cdots &\widetilde x_{2(i_1-1)}&0& \cdots&\widetilde x_{2(i_p-1)}&0& \cdots  \\
  \ddots&\vdots&\vdots&\ddots&\vdots&\vdots&\ddots\\
  \cdots &\widetilde x_{p(i_1-1)}&0& \cdots&\widetilde x_{p(i_p-1)}&1& \cdots \\} {(n-s) \,\rm columns}.
  \tikz[remember picture,overlay]   \draw[dashed,dash pattern={on 4pt off 2pt}] ([xshift=0.5\tabcolsep,yshift=7pt]a.north) -- ([xshift=0.5\tabcolsep,yshift=-2pt]b.south);\tikz[remember picture,overlay]   \draw[dashed,dash pattern={on 4pt off 2pt}] ([xshift=0.5\tabcolsep,yshift=7pt]e.north) -- ([xshift=0.5\tabcolsep,yshift=-2pt]f.south);
\end{equation}
Let $\widehat U_{i_1\cdots i_p}^P$ be the corresponding open subset of $\mathbb P(N_{\mathcal V_{(p-r,r)}/G(p,n)})$ with the following local coordinates
\begin{equation}\label{trr5}
   \left(\left(\cdots,\widetilde x_{uv},\cdots\right)_{\substack{1\leq u\leq p\\s+1\leq v\leq n\\v\neq i_1,\cdots,i_p}}\,,\,\, [\cdots, \widetilde a_{ij},\cdots]_{\substack{1\leq i\leq p\\1\leq j\leq s}}\right)=:(\widetilde X,\big[\widetilde A\,\big])_{i_1\cdots i_p}\,.
\end{equation}

$(X,[A\,])$ and $\left(\widetilde X,[\widetilde A\,]\right)_{i_1\cdots i_p}$ represent the same point of  $\mathbb P(N_{\mathcal V_{(p-r,r)}/G(p,n)})$  if and only if the matrices $\widehat M_{(X,A)}$ and $\widehat M_{\left(\widetilde X,\widetilde A\right)_{i_1\cdots i_p}}$ represent the same point of $G(p,n)$ up to the $\mathbb C^*$-action $\psi_{s,p,n}$.

\section{\texorpdfstring{Calculation of the intersection numbers }{rr}}\label{section:inters}

{\bf\noindent Proof of Lemma \ref{i2}.} The proof is the same as  Lemma \ref{i1}. 
It is easy to verify that 
\begin{equation}
    \begin{split}
        &(1)\,\,\,\,{\rm the\,\, projection\,\, of\,\,} \zeta_j^0\,\,{\rm  to}\,\,\mathbb {CP}^{N_{p,n}}\,\,{\rm is\,\, a\,\, point}\,;\\
        &(2)\,\,\,\,{\rm the\,\, projection\,\, of\,\,} \zeta_j^0\,\,{\rm  to}\,\,\mathbb {CP}^{N^{j-1}_{s,p,n}}\,\,{\rm is\,\, a\,\, line}\,;\\
        &(3)\,\,\,\,{\rm the\,\, projection\,\, of\,\,} \zeta_j^0\,\,{\rm  to}\,\,\mathbb {CP}^{N^k_{s,p,n}}\,\,{\rm is\,\, a\,\, point\,\,for\,\,}0\leq k\leq p\,\,{\rm and\,\,}k\neq j-1.\\
    \end{split}
\end{equation}
Therefore, 
\begin{equation}
\begin{split}
&(1)\,\,\,\,\zeta_j^0\cdot R^*_{s,p,n}((\mathcal O_{G(p,n)}(1))=0\,;\\
&(2)\,\,\,\,\zeta_j^0\cdot H_{j-1}\big|_{\mathcal T_{s,p,n}}=1\,;\\
&(3)\,\,\,\,\zeta_j^0\cdot H_k\big|_{\mathcal T_{s,p,n}}=0\,\,{\rm for}\,\,0\leq k\leq p\,\,{\rm and}\,\,k\neq j-1.\\
\end{split}
\end{equation}

We can complete the proof of Lemma \ref{i2} by Lemma \ref{partialline}.\,\,\,\,$\endpf$
\medskip

{\bf\noindent Proof of Lemma \ref{i3}.} The proof is the same as in Lemma \ref{i1}. 
For simplicity, we will only prove when $u=k$ for  the proof of the case $v=s+k$ is the same.

We can show that 
\begin{equation}
    \begin{split}
        &(1)\,\,\,\,{\rm the\,\, projection\,\, of\,\,} \zeta^{0,k}_{k,v}\,\,{\rm  to}\,\,\mathbb {CP}^{N_{p,n}}\,\,{\rm is\,\, a\,\, point}\,;\\
        &(2)\,\,\,\,{\rm the\,\, projection\,\, of\,\,} \zeta^{0,k}_{k,v}\,\,{\rm  to}\,\,\mathbb {CP}^{N^{j}_{s,p,n}}\,\,{\rm is\,\, a\,\, line}\,\,{\rm for}\,\,k\leq j\leq\, {\rm min}\,\{v-s-1,r\}\,;\\
        &(3)\,\,\,\,{\rm the\,\, projection\,\, of\,\,} \zeta^{0,k}_{k,v}\,\,{\rm  to}\,\,\mathbb {CP}^{N^j_{s,p,n}}\,\,{\rm is\,\, a\,\, point\,\,for\,\,}0\leq j\leq k-1\,\,{\rm and\,\,}\\
        &\,\,\,\,\,\,\,\,\,\,\,\,\,\,v-s\leq j\leq r.\\
    \end{split}
\end{equation}
Therefore, 
\begin{equation}
\begin{split}
&(1)\,\,\,\,\zeta^{0,k}_{k,v}\cdot R^*_{s,p,n}((\mathcal O_{G(p,n)}(1))=0;\\
&(2)\,\,\,\,\zeta^{0,k}_{k,v}\cdot H_{j}\big|_{\mathcal T_{s,p,n}}=1\,\,{\rm for}\,\,k\leq j\leq\, {\rm min}\,\{v-s-1,r\};\\
&(3)\,\,\,\,\zeta^{0,k}_{k,v}\cdot H_j\big|_{\mathcal T_{s,p,n}}=0\,\,{\rm for}\,\,0\leq j\leq k-1\,\,{\rm and\,\,}v-s\leq j\leq r.\\
\end{split}
\end{equation}

We can complete the proof of Lemma \ref{i3} by Lemma \ref{partialline}.\,\,\,\,$\endpf$
\medskip

{\bf\noindent Proof of Lemma \ref{i4}.} The proof is the same as in Lemma \ref{i1}. We can show that 
\begin{equation}
    \begin{split}
        &(1)\,\,\,\,{\rm the\,\, projection\,\, of\,\,} \delta_{m_1,m_2}^0\,\,{\rm  to}\,\,\mathbb {CP}^{N_{p,n}}\,\,{\rm is\,\, a\,\, line}\,;\\
        &(2)\,\,\,\,{\rm the\,\, projection\,\, of\,\,} \delta_{m_1,m_2}^0\,\,{\rm  to}\,\,\mathbb {CP}^{N^{j}_{s,p,n}}\,\,{\rm is\,\, a\,\, line}\,\,{\rm for}\,\,0\leq j\leq \min\{m_1-1,r\}\,;\\
        &(3)\,\,\,\,{\rm the\,\, projection\,\, of\,\,} \delta_{m_1,m_2}^0\,\,{\rm  to}\,\,\mathbb {CP}^{N^j_{s,p,n}}\,\,{\rm is\,\, a\,\, point\,\,for\,\,}m_1\leq j\leq r.\\
    \end{split}
\end{equation}
Therefore, 
\begin{equation}
\begin{split}
&(1)\,\,\,\,\delta_{m_1,m_2}^0\cdot R^*_{s,p,n}((\mathcal O_{G(p,n)}(1))=1\,;\\
&(2)\,\,\,\,\delta_{m_1,m_2}^0\cdot H_{j}\big|_{\mathcal T_{s,p,n}}=1\,\,{\rm for}\,\,\min\{m_1-1,r\}\,;\\
&(3)\,\,\,\,\delta_{m_1,m_2}^0\cdot H_j\big|_{\mathcal T_{s,p,n}}=0\,\,{\rm for}\,\,m_1\leq j\leq r.\\
\end{split}
\end{equation}

We can complete the proof of Lemma \ref{i4} by Lemma \ref{partialline}.\,\,\,\,$\endpf$
\medskip

{\bf\noindent Proof of Lemma \ref{i5}.} The proof is the same as in Lemma \ref{i2}. 

For $2\leq j\leq r-l$, it is easy to verify 
\begin{equation}
    \begin{split}
        &(1)\,\,\,\,{\rm the\,\, projection\,\, of\,\,} \zeta_j^l\,\,{\rm  to}\,\,\mathbb {CP}^{N_{p,n}}\,\,{\rm is\,\, a\,\, point}\,;\\
        &(2)\,\,\,\,{\rm the\,\, projection\,\, of\,\,} \zeta_j^l\,\,{\rm  to}\,\,\mathbb {CP}^{N^{l+j-1}_{s,p,n}}\,\,{\rm is\,\, a\,\, line}\,;\\
        &(3)\,\,\,\,{\rm the\,\, projection\,\, of\,\,} \zeta_j^l\,\,{\rm  to}\,\,\mathbb {CP}^{N^k_{s,p,n}}\,\,{\rm is\,\, a\,\, point\,\,for\,\,}0\leq k\leq p\,\,{\rm and\,\,}k\neq l+j-1.\\
    \end{split}
\end{equation}
Therefore, 
\begin{equation}
\begin{split}
&(1)\,\,\,\,\zeta_j^l\cdot R^*_{s,p,n}((\mathcal O_{G(p,n)}(1))=0\,;\\
&(2)\,\,\,\,\zeta_j^l\cdot H_{l+j-1}\big|_{\mathcal T_{s,p,n}}=1\,;\\
&(3)\,\,\,\,\zeta_j^l\cdot H_k\big|_{\mathcal T_{s,p,n}}=0\,\,{\rm for}\,\,0\leq k\leq p\,\,{\rm and}\,\,k\neq l+j-1.\\
\end{split}
\end{equation}

Similarly, when $r-l+2\leq j\leq r$ we have 
\begin{equation}
    \begin{split}
        &(1)\,\,\,\,{\rm the\,\, projection\,\, of\,\,} \zeta_j^l\,\,{\rm  to}\,\,\mathbb {CP}^{N_{p,n}}\,\,{\rm is\,\, a\,\, point}\,;\\
        &(2)\,\,\,\,{\rm The\,\, projection\,\, of\,\,} \zeta_j^l\,\,{\rm  to}\,\,\mathbb {CP}^{N^{r-j+1}_{s,p,n}}\,\,{\rm is\,\, a\,\, line}\,;\\
        &(3).\,\,\,\,{\rm the\,\, projection\,\, of\,\,} \zeta_j^l\,\,{\rm  to}\,\,\mathbb {CP}^{N^k_{s,p,n}}\,\,{\rm is\,\, a\,\, point\,\,for\,\,}0\leq k\leq p\,\,{\rm and\,\,}k\neq r-j+1.\\
    \end{split}
\end{equation}
Therefore, 
\begin{equation}
\begin{split}
&(1)\,\,\,\,\zeta_j^l\cdot R^*_{s,p,n}((\mathcal O_{G(p,n)}(1))=0\,;\\
&(2)\,\,\,\,\zeta_j^l\cdot H_{r-j+1}\big|_{\mathcal T_{s,p,n}}=1\,;\\
&(3)\,\,\,\,\zeta_j^l\cdot H_k\big|_{\mathcal T_{s,p,n}}=0\,\,{\rm for}\,\,0\leq k\leq p\,\,{\rm and}\,\,k\neq r-j+1.\\
\end{split}
\end{equation}

Noticing that $\zeta_j^l$ is disjoint with $D^-_{1},D^-_{2},\cdots,D^-_{l},D^+_1,D^+_2,\cdots,D^+_{r-l}$, we  can conclude (\ref{zei5}) and (\ref{zer5}) by  Lemma \ref{partialline}. We complete the proof of Lemma \ref{i5}.
\,\,\,\,$\endpf$
\medskip

{\bf\noindent Proof of Lemmas \ref{li1} and \ref{li2}.} The proof is the same as in Lemma \ref{i3}. 

First assume that $1\leq k\leq r-l$.
For simplicity, we will only give the proof when $u=l+k$ for  the proof of the case $v=s+l+k$ is the same. We can show 
\begin{equation}
    \begin{split}
        &(1)\,\,\,\,{\rm the\,\, projection\,\, of\,\,} \zeta^{l,k}_{l+k,v}\,\,{\rm  to}\,\,\mathbb {CP}^{N_{p,n}}\,\,{\rm is\,\, a\,\, point}\,;\\
        &(2)\,\,\,\,{\rm the\,\, projection\,\, of\,\,} \zeta^{l,k}_{l+k,v}\,\,{\rm  to}\,\,\mathbb {CP}^{N^{j}_{s,p,n}}\,\,{\rm is\,\, a\,\, line}\,\,{\rm for}\,\,l+k\leq j\leq\, {\rm min}\,\{v-s-1,r\}\,;\\
        &(3)\,\,\,\,{\rm the\,\, projection\,\, of\,\,} \zeta^{l,k}_{l+k,v}\,\,{\rm  to}\,\,\mathbb {CP}^{N^j_{s,p,n}}\,\,{\rm is\,\, a\,\, point\,\,for\,\,}0\leq j\leq l+k-1\,\,{\rm and\,\,}\\
        &\,\,\,\,\,\,\,\,\,\,\,\,\,\,v-s\leq j\leq r.\\
    \end{split}
\end{equation}
Therefore, 
\begin{equation}
\begin{split}
&(1)\,\,\,\,\zeta^{l,k}_{l+k,v}\cdot R^*_{s,p,n}((\mathcal O_{G(p,n)}(1))=0\,;\\
&(2)\,\,\,\,\zeta^{l,k}_{l+k,v}\cdot H_{j}\big|_{\mathcal T_{s,p,n}}=1\,\,{\rm for}\,\,l+k\leq j\leq\, {\rm min}\,\{v-s-1,r\}\,;\\
&(3)\,\,\,\,\zeta^{l,k}_{l+k,v}\cdot H_j\big|_{\mathcal T_{s,p,n}}=0\,\,{\rm for}\,\,0\leq j\leq l+k-1\,\,{\rm and\,\,}v-s\leq j\leq r.\\
\end{split}
\end{equation}

When $r-l+1\leq k\leq r$, without loss of generality, we may assume that $u=r-k+1$. Similarly, we have 
\begin{equation}
    \begin{split}
        &(1)\,\,\,\,{\rm the\,\, projection\,\, of\,\,} \zeta^{l,k}_{r+1-k,v}\,\,{\rm  to}\,\,\mathbb {CP}^{N_{p,n}}\,\,{\rm is\,\, a\,\, point}\,;\\
        &(2)\,\,\,\,{\rm the\,\, projection\,\, of\,\,} \zeta^{l,k}_{r+1-k,v}\,\,{\rm  to}\,\,\mathbb {CP}^{N^{j}_{s,p,n}}\,\,{\rm is\,\, a\,\, line}\,\,{\rm for}\,\,{\rm max}\,\{v-s+p,0\}\leq j\leq r-k\,;\\
        &(3)\,\,\,\,{\rm the\,\, projection\,\, of\,\,} \zeta^{l,k}_{r+1-k,v}\,\,{\rm  to}\,\,\mathbb {CP}^{N^j_{s,p,n}}\,\,{\rm is\,\, a\,\, point\,\,for\,\,}0\leq j\leq v-s+p-1\,\,{\rm and\,\,}\\
        &\,\,\,\,\,\,\,\,\,\,\,r-k+1\leq j\leq r.\\
    \end{split}
\end{equation}
Therefore, 
\begin{equation}
\begin{split}
&(1)\,\,\,\,\zeta^{l,k}_{r+1-k,v}\cdot R^*_{s,p,n}((\mathcal O_{G(p,n)}(1))=0\,;\\
&(2)\,\,\,\,\zeta^{l,k}_{r+1-k,v}\cdot H_{j}\big|_{\mathcal T_{s,p,n}}=1\,\,{\rm for}\,\,{\rm max}\,\{v-s+p,0\}\leq j\leq r-k\,;\\
&(3)\,\,\,\,\zeta^{l,k}_{r+1-k,v}\cdot H_j\big|_{\mathcal T_{s,p,n}}=0\,\,{\rm for}\,\,0\leq j\leq v-s+p-1\,\,{\rm and\,\,}r-k+1\leq j\leq r.\\
\end{split}
\end{equation}
We can complete the proof of Lemmas \ref{li1} and \ref{li2} by Lemma \ref{partialline}.
\,\,\,\,$\endpf$
\medskip

{\bf\noindent Proof of Lemma \ref{li4}.} The proof is the same as in Lemma \ref{i4}. We can show that 
\begin{equation}
    \begin{split}
        &(1)\,\,\,\,{\rm the\,\, projection\,\, of\,\,} \delta_{m_1,m_2}^l\,\,{\rm  to}\,\,\mathbb {CP}^{N_{p,n}}\,\,{\rm is\,\, a\,\, line}\,;\\
        &(2)\,\,\,\,{\rm the\,\, projection\,\, of\,\,} \delta_{m_1,m_2}^l\,\,{\rm  to}\,\,\mathbb {CP}^{N^{j}_{s,p,n}}\,\,{\rm is\,\, a\,\, line}\,\,{\rm for}\,\,l-m_2+1\leq j\leq l+m_1-1\,;\\
        &(3)\,\,\,\,{\rm the\,\, projection\,\, of\,\,} \delta_{m_1,m_2}^l\,\,{\rm  to}\,\,\mathbb {CP}^{N^j_{s,p,n}}\,\,{\rm is\,\, a\,\, point\,\,for\,\,}0\leq j\leq l-m_2\,\,{\rm\,\,or}\\
        &\,\,\,\,\,\,\,\,\,\,\,\,\,\,l+m_1\leq j\leq r.\\
    \end{split}
\end{equation}
Therefore, 
\begin{equation}
\begin{split}
&(1)\,\,\,\,\delta_{m_1,m_2}^l\cdot R^*_{s,p,n}((\mathcal O_{G(p,n)}(1))=1\,;\\
&(2)\,\,\,\,\delta_{m_1,m_2}^l\cdot H_{j}\big|_{\mathcal T_{s,p,n}}=1\,\,{\rm for}\,\,l-m_2+1\leq j\leq l+m_1-1\,;\\
&(3)\,\,\,\,\delta_{m_1,m_2}^l\cdot H_j\big|_{\mathcal T_{s,p,n}}=0\,\,{\rm for\,\,}0\leq j\leq l-m_2\,\,{\rm or}\,\, l+m_1\leq j\leq r.
\end{split}
\end{equation}
We can complete the proof of Lemma \ref{li4} by Lemma \ref{partialline}.\,\,\,\,$\endpf$
\medskip

{\bf\noindent Proof of Lemma \ref{li5}.} The proof is the same as in Lemma \ref{li4}. We can show that 
\begin{equation}
    \begin{split}
        &(1)\,\,\,\,{\rm the\,\, projection\,\, of\,\,} \Delta_{m_1,m_2}^l\,\,{\rm  to}\,\,\mathbb {CP}^{N_{p,n}}\,\,{\rm is\,\, a\,\, line}\,;\\
        &(2)\,\,\,\,{\rm the\,\, projection\,\, of\,\,} \Delta_{m_1,m_2}^l\,\,{\rm  to}\,\,\mathbb {CP}^{N^{j}_{s,p,n}}\,\,{\rm is\,\, a\,\, line}\,\,{\rm for}\,\,l-m_2+1\leq j\leq l+m_1-1\,;\\
        &(3)\,\,\,\,{\rm the\,\, projection\,\, of\,\,} \Delta_{m_1,m_2}^l\,\,{\rm  to}\,\,\mathbb {CP}^{N^j_{s,p,n}}\,\,{\rm is\,\, a\,\, point\,\,for\,\,}0\leq j\leq l-m_2\,\,{\rm\,\,or}\\
        &\,\,\,\,\,\,\,\,\,\,\,\,\,\,l+m_1\leq j\leq r.\\
    \end{split}
\end{equation}
Therefore, 
\begin{equation}
\begin{split}
&(1)\,\,\,\,\Delta_{m_1,m_2}^l\cdot R^*_{s,p,n}((\mathcal O_{G(p,n)}(1))=1\,;\\
&(2)\,\,\,\,\Delta_{m_1,m_2}^l\cdot H_{j}\big|_{\mathcal T_{s,p,n}}=1\,\,{\rm for}\,\,l-m_2+1\leq j\leq l+m_1-1.\\
&(3).\,\,\,\,\Delta_{m_1,m_2}^l\cdot H_j\big|_{\mathcal T_{s,p,n}}=0\,\,{\rm for\,\,}0\leq j\leq l-m_2\,\,{\rm or}\,\, l+m_1\leq j\leq r.
\end{split}
\end{equation}
We can complete the proof of Lemma \ref{li4}.\,\,\,\,$\endpf$
\medskip

{\bf\noindent Proof of Lemma \ref{ri1}.} The proof is the same as in Lemma \ref{i2}. 
It is easy to verify that 
\begin{equation}
    \begin{split}
        &(1)\,\,\,\,{\rm the\,\, projection\,\, of\,\,} \zeta_j^r\,\,{\rm  to}\,\,\mathbb {CP}^{N_{p,n}}\,\,{\rm is\,\, a\,\, point}\,;\\
        &(2)\,\,\,\,{\rm the\,\, projection\,\, of\,\,} \zeta_j^r\,\,{\rm  to}\,\,\mathbb {CP}^{N^{r-j+1}_{s,p,n}}\,\,{\rm is\,\, a\,\, line}\,;\\
        &(3)\,\,\,\,{\rm the\,\, projection\,\, of\,\,} \zeta_j^r\,\,{\rm  to}\,\,\mathbb {CP}^{N^k_{s,p,n}}\,\,{\rm is\,\, a\,\, point\,\,for\,\,}0\leq k\leq p\,\,{\rm and\,\,}k\neq r-j+1.\\
    \end{split}
\end{equation}
Therefore, 
\begin{equation}
\begin{split}
&(1)\,\,\,\,\zeta_j^r\cdot R^*_{s,p,n}((\mathcal O_{G(p,n)}(1))=0\,;\\
&(2)\,\,\,\,\zeta_j^r\cdot H_{r-j+1}\big|_{\mathcal T_{s,p,n}}=1\,;\\
&(3)\,\,\,\,\zeta_j^r\cdot H_k\big|_{\mathcal T_{s,p,n}}=0\,\,{\rm for}\,\,0\leq k\leq p\,\,{\rm and}\,\,k\neq r-j+1.\\
\end{split}
\end{equation}
Since $\zeta_j^r$ is disjoint with $D^-_1,D^-_2,\cdots,D^-_{r}$, we can complete the proof of Lemma \ref{i2} by Lemma \ref{partialline}.\,\,\,\,$\endpf$

\medskip

{\bf\noindent Proof of Lemma \ref{ri2}.} The proof is the same as in Lemmas \ref{li1} and \ref{li2}.

Without loss of generality, we may assume that $u=r-k+1$. Then,
\begin{equation}
    \begin{split}
        &(1)\,\,\,\,{\rm the\,\, projection\,\, of\,\,} \zeta^{r,k}_{k,v}\,\,{\rm  to}\,\,\mathbb {CP}^{N_{p,n}}\,\,{\rm is\,\, a\,\, point}\,;\\
        &(2)\,\,\,\,{\rm the\,\, projection\,\, of\,\,} \zeta^{r,k}_{k,v}\,\,{\rm  to}\,\,\mathbb {CP}^{N^{j}_{s,p,n}}\,\,{\rm is\,\, a\,\, line}\,\,{\rm for}\,\,{\rm max}\,\{v-s+p,0\}\leq j\leq r-k\,;\\
        &(3)\,\,\,\,{\rm the\,\, projection\,\, of\,\,} \zeta^{r,k}_{k,v}\,\,{\rm  to}\,\,\mathbb {CP}^{N^j_{s,p,n}}\,\,{\rm is\,\, a\,\, point\,\,for\,\,}0\leq j\leq v-s+p-1\\
        &\,\,\,\,\,\,\,\,\,\,\,\,\,{\rm or\,\,}r-k+1\leq j\leq r.\\
    \end{split}
\end{equation}
Therefore, 
\begin{equation}
\begin{split}
&(1)\,\,\,\,\zeta^{r,k}_{k,v}\cdot R^*_{s,p,n}((\mathcal O_{G(p,n)}(1))=0\,;\\
&(2)\,\,\,\,\zeta^{r,k}_{k,v}\cdot H_{j}\big|_{\mathcal T_{s,p,n}}=1\,\,{\rm for}\,\,{\rm max}\,\{v-s+p,0\}\leq j\leq r-k\,;\\
&(3)\,\,\,\,\zeta^{r,k}_{k,v}\cdot H_j\big|_{\mathcal T_{s,p,n}}=0\,\,{\rm for}\,\,0\leq j\leq v-s+p-1\,\,{\rm or\,\,}r-k+1\leq j\leq r.\\
\end{split}
\end{equation}
We can complete the proof of Lemma \ref{ri2} by Lemma \ref{partialline}.
\,\,\,\,$\endpf$
\medskip

{\bf\noindent Proof of Lemma \ref{ri3}.} The proof is the same as in Lemma \ref{i4}. We can show that 
\begin{equation}
    \begin{split}
        &(1).\,\,\,\,{\rm the\,\, projection\,\, of\,\,} \delta^r_{m_1,m_2}\,\,{\rm  to}\,\,\mathbb {CP}^{N_{p,n}}\,\,{\rm is\,\, a\,\, line}\,;\\
        &(2).\,\,\,\,{\rm the\,\, projection\,\, of\,\,} \delta^r_{m_1,m_2}\,\,{\rm  to}\,\,\mathbb {CP}^{N^{j}_{s,p,n}}\,\,{\rm is\,\, a\,\, line}\,\,{\rm for}\,\,r-m_2+1\leq j\leq r\,;\\
        &(3).\,\,\,\,{\rm the\,\, projection\,\, of\,\,} \delta^r_{m_1,m_2}\,\,{\rm  to}\,\,\mathbb {CP}^{N^j_{s,p,n}}\,\,{\rm is\,\, a\,\, point\,\,for\,\,}0\leq j\leq r-m_2.\\
    \end{split}
\end{equation}
Therefore, 
\begin{equation}
\begin{split}
&(1)\,\,\,\,\delta^r_{m_1,m_2}\cdot R^*_{s,p,n}((\mathcal O_{G(p,n)}(1))=1\,;\\
&(2)\,\,\,\,\delta^r_{m_1,m_2}\cdot H_{j}\big|_{\mathcal T_{s,p,n}}=1\,\,{\rm for}\,\,r-m_2+1\leq j\leq r\,;\\
&(3)\,\,\,\,\delta^r_{m_1,m_2}\cdot H_j\big|_{\mathcal T_{s,p,n}}=0\,\,{\rm for}\,\,0\leq j\leq r-m_2.\\
\end{split}
\end{equation}
We can complete the proof of Lemma \ref{ri3} by Lemma \ref{partialline}.\,\,\,\,$\endpf$
\medskip

{\bf\noindent Proof of Lemma \ref{ri4}.} The proof is the same as in Lemma \ref{ri3}. We can show that 
\begin{equation}
    \begin{split}
        &(1)\,\,\,\,{\rm the\,\, projection\,\, of\,\,} \Delta^r_{m_1,m_2}\,\,{\rm  to}\,\,\mathbb {CP}^{N_{p,n}}\,\,{\rm is\,\, a\,\, line}\,;\\
        &(2)\,\,\,\,{\rm the\,\, projection\,\, of\,\,} \Delta^r_{m_1,m_2}\,\,{\rm  to}\,\,\mathbb {CP}^{N^{j}_{s,p,n}}\,\,{\rm is\,\, a\,\, line}\,\,{\rm for}\,\,r-m_2+1\leq j\leq r\,;\\
        &(3)\,\,\,\,{\rm the\,\, projection\,\, of\,\,} \Delta^r_{m_1,m_2}\,\,{\rm  to}\,\,\mathbb {CP}^{N^j_{s,p,n}}\,\,{\rm is\,\, a\,\, point\,\,for\,\,}0\leq j\leq r-m_2.\\
    \end{split}
\end{equation}
Therefore, 
\begin{equation}
\begin{split}
&(1)\,\,\,\,\Delta^r_{m_1,m_2}\cdot R^*_{s,p,n}((\mathcal O_{G(p,n)}(1))=1\,;\\
&(2)\,\,\,\,\Delta^r_{m_1,m_2}\cdot H_{j}\big|_{\mathcal T_{s,p,n}}=1\,\,{\rm for}\,\,r-m_2+1\leq j\leq r\,;\\
&(3)\,\,\,\,\Delta^r_{m_1,m_2}\cdot H_j\big|_{\mathcal T_{s,p,n}}=0\,\,{\rm for}\,\,0\leq j\leq r-m_2.\\
\end{split}
\end{equation}
We can complete the proof of Lemma \ref{ri4}.\,\,\,\,$\endpf$
\medskip

\section{\texorpdfstring{Rigidity lemmas for the symmetry groups }{rr}}\label{section:rigid}

\subsection{\texorpdfstring{Proof of Lemma \ref{sigma2}}{rr}}\label{section:rigidmc} Denote by  Con$(\mathcal T_{s,p,n})\subset A_{n-1}(\mathcal T_{s,p,n})\otimes_{\mathbb Z}\mathbb Q$ the cone  of effective divisors of $\mathcal T_{s,p,n}$. By Lemmas \ref{gb} and \ref{cone}, Con$(\mathcal T_{s,p,n})$ is generated by 
\begin{equation}
    \left\{B_0,\cdots,B_r,D^+_1,\cdots,D^+_{r},D^-_1,\cdots,D^-_{r} \right\}.
\end{equation}
Recalling (\ref{bst}), (\ref{bstt0}), (\ref{bsttr}),   we can determine $\mathfrak G$ the set of extremal rays  of Con$(\mathcal T_{s,p,n})$ as follows.
When $p\neq n-s$ and $p<s$, 
\begin{equation}
    \mathfrak G:=\left\{B_0,B_r,D^+_1, \cdots,D^+_{r},D^-_1,\cdots,D^-_{r} \right\};
\end{equation}
when $n-s=p<s$,  
\begin{equation}
    \mathfrak G:=\left\{B_0,B_r,D^+_1,\cdots, D^+_{r},D^-_1,\cdots,D^-_{r-1} \right\};
\end{equation}
when $n-s=p=s$, 
\begin{equation}
    \mathfrak G:=\left\{B_0,B_r,D^+_1,\cdots, D^+_{r-1},D^-_1,\cdots,D^-_{r-1} \right\}.
\end{equation}
Since $\sigma^*$ preserves Con$(\mathcal T_{s,p,n})$, it induces a permutation of $\mathfrak G$. 


Next we prove Lemma \ref{sigma2} base on a case by case argument. For simplicity of notation, we denote $(R_{s,p,n})^*\left(\mathcal O_{G(p,n)}(1)\right)$ by $H$ in the following.

\smallskip

{\bf\noindent Case 1 ($r\geq 2$, $p\neq n-s$ and $p<s$).}
By (\ref{bst}) we have that
\begin{equation}
    (R_{s,p,n})^*\left(\mathcal O_{ G(p,n)}(1)\right)=B_0+\sum_{i=1}^{r}(r+1-i)\cdot D^+_i=B_r+\sum_{i=1}^r(r+1-i)\cdot D^-_i\,;
\end{equation}
hence, 
\begin{equation}\label{b=b}
    \sigma^*(B_0)=\sigma^*(B_r)+\sum_{i=1}^r(r+1-i)\cdot \sigma^*(D^-_i)-\sum_{i=1}^{r}(r+1-i)\cdot \sigma^*(D^+_i)\,.
\end{equation}
By Lemma \ref{picb}, it is clear that (\ref{b=b}) has integer coefficients with respect to the basis of the Picard group $\left\{H,D^-_1, \cdots, D^-_r,D^+_1,\cdots,D^+_r \right\}$.

If $\sigma^*(B_0)\in\left\{D^-_1, \cdots, D^-_r,D^+_1,\cdots,D^+_r\right\}$, we have the following possibilities, by checking the coefficient of $H$ in (\ref{b=b}).
\begin{enumerate}[label=(\Alph*).]
    \item For a certain $1\leq j\leq r$, $\sigma^*(D_j^-)=B_0$ and $\sigma^*(D_j^+)=B_r$.
    \item For a certain $1\leq j\leq r$, $\sigma^*(D_j^+)=B_0$ and $\sigma^*(D_j^-)=B_r$.
    \item $\sigma^*(B_r)=B_r$ and $\sigma^*(D_r^+)=B_0$.
    \item $\sigma^*(B_r)=B_0$ and $\sigma^*(D_r^+)=B_r$.
\end{enumerate}

For Case (A), the coefficient of $D^+_1$ is at most $-r(r-j)$ in the following quantity,
\begin{equation}
  \sum_{i=1}^r(r+1-i)\cdot \sigma^*(D^-_i)-\sum_{i=1}^{r}(r+1-i)\cdot \sigma^*(D^+_i)\,.
\end{equation}
Since $r\geq 2$, we have  $j=r$ by (\ref{b=b}). Then,  by (\ref{b=b}) again we have that
\begin{equation}
    \begin{split}
        &\sigma^*(D_r^-)=B_0\,,\,\,\,\sigma^*(D_r^+)=B_r\,,\,\,\,\sigma^*(B_0)=D_r^-\,,\,\,\,\sigma^*(B_r)=D_r^+\,,\\
        &\sigma^*(D_j^+)=D_j^-\,,\,\,\,\sigma^*(D_j^-)=D_j^+\,\,{\rm for\,\,}1\leq j\leq r-1\,.
    \end{split}
\end{equation}
Since $\sigma^*(K_{\mathcal T_{s,p,n}})=K_{\mathcal T_{s,p,n}}$, we conclude that $n=2s$, and hence $\sigma^*=({\rm USD})^*$.

For Case (B), in the same way as above, the only possibility is that \begin{equation}
    \begin{split}
        &\sigma^*(D_r^+)=B_0\,,\,\,\,\sigma^*(D_r^-)=B_r\,,\,\,\,\sigma^*(B_0)=D_r^+\,,\,\,\,\sigma^*(B_r)=D_r^-\,,\\
        &\sigma^*(D_j^-)=D_j^-\,,\,\,\,\sigma^*(D_j^+)=D_j^+\,\,{\rm for\,\,}1\leq j\leq r-1\,.
    \end{split}
\end{equation}
Notice that the coefficient of $H$ in $K_{\mathcal T_{s,p,n}}$ is $-n$. When $n-s<p$,  the coefficient of $H$ in $\sigma^*(K_{\mathcal T_{s,p,n}})$ is $-2r=-2(n-s)>-2p>-n$; when $p<n-s$, the coefficient of $H$ in $\sigma^*(K_{\mathcal T_{s,p,n}})$ is $-2p>-s-(n-s)=-n$. Both  lead to a contradiction. 

For Case (C),   by (\ref{b=b}) we can show that $\sigma^*$ permutes $D_r^+$ and $B_0$, and fixes all the other divisors. Then the coefficient of $H$ in $\sigma^*(K_{\mathcal T_{s,p,n}})$ is $s-p-n>-n$ which is a contradiction.

For Case (D),  by (\ref{b=b}) we have that
\begin{equation}
    \begin{split}
        &\sigma^*(B_r)=B_0\,,\,\,\,\sigma^*(D_r^+)=B_r\,,\,\,\,\sigma^*(D_r^-)=D_r^+\,,\,\,\,\sigma^*(B_0)=D_r^-\,,\\
        &\sigma^*(D_j^+)=D_j^-\,,\,\,\,\sigma^*(D_j^-)=D_j^+\,\,{\rm for\,\,}1\leq j\leq r-1\,.
    \end{split}
\end{equation}
It is easy to verify that $\sigma^*(H)=H$, and hence $\sigma^*(K_{\mathcal T_{s,p,n}})\neq K_{\mathcal T_{s,p,n}}$. This is a contradiction.

Next, we assume that  $\sigma^*(B_0)\in\left\{B_0, B_r\right\}$. By checking the coefficient of $H$ in (\ref{b=b}), we have the following possibilities. \begin{enumerate}[label=(\alph*).]
    \item  $\sigma^*(B_0)=B_r$ and $\sigma^*(D_r^-)=B_0$.
    \item $\sigma^*(B_0)=B_r$ and $\sigma^*(B_r)=B_0$.
    \item $\sigma^*(B_0)=B_0$ and $\sigma^*(D^-_r)=B_r$.
    \item $\sigma^*(B_0)=B_0$ and $\sigma^*(B_r)=B_r$.
\end{enumerate}

For Case (a), by (\ref{b=b}) we can conclude that \begin{equation}
    \begin{split}
        &\sigma^*(B_0)=B_r\,,\,\,\,\sigma^*(D_r^-)=B_0\,,\,\,\,\sigma^*(D_r^+)=D_r^-\,,\,\,\,\sigma^*(B_r)=D_r^+\,,\\
        &\sigma^*(D_j^+)=D_j^-\,,\,\,\,\sigma^*(D_j^-)=D_j^+\,\,{\rm for\,\,}1\leq j\leq r-1\,.
    \end{split}
\end{equation}
Then $\sigma^*(H)=H$, and hence $\sigma^*(K_{\mathcal T_{s,p,n}})\neq K_{\mathcal T_{s,p,n}}$. This is a contradiction.

For Case (b), by (\ref{b=b}) we can conclude that \begin{equation}
    \begin{split}
        &\sigma^*(B_r)=B_0\,,\,\,\,\sigma^*(B_0)=B_r\,,\,\,\sigma^*(D_j^+)=D_j^-\,,\,\,\,\sigma^*(D_j^-)=D_j^+\,\,{\rm for\,\,}1\leq j\leq r\,.
    \end{split}
\end{equation}
Since $\sigma^*(K_{\mathcal T_{s,p,n}})=K_{\mathcal T_{s,p,n}}$, we conclude that either $n=2p$ and $\sigma^*=({\rm DUAL})^*$, or  $n=2s$ and $\sigma^*=({\rm USD})^*$.

For Case (c),  by (\ref{b=b}) we can show that $\sigma^*$ permutes $D_r^-$ and $B_r$, and fixes all the other divisors. Since $\sigma^*(K_{\mathcal T_{s,p,n}})=K_{\mathcal T_{s,p,n}}$, we have that $p=n-s$ which contradicts the assumption.

For Case (d),  by (\ref{b=b}) we can show that $\sigma^*$ is the identity map.

\smallskip

{\bf\noindent Case 2 ($r\geq 2$ and $n-s=p<s$).} By (\ref{bst}) and (\ref{bsttr}) we have that
\begin{equation}\label{b=b2}
    \sigma^*(B_0)=\sigma^*(B_r)+\sum_{i=1}^{r-1}(r+1-i)\cdot \sigma^*(D^-_i)-\sum_{i=1}^{r}(r+1-i)\cdot \sigma^*(D^+_i)\,.
\end{equation}

If $\sigma^*(B_0)\in\left\{D^-_1, \cdots, D^-_{r-1},D^+_1,\cdots,D^+_r\right\}$, we have the following possibilities by checking the coefficients of $H$ in (\ref{b=b2}).
\begin{enumerate}[label=(\Alph*).]
    \item For a certain $1\leq j\leq r-1$, $\sigma^*(D_j^-)=B_0$ and $\sigma^*(D_j^+)=B_r$.
    \item For a certain $1\leq j\leq r-1$, $\sigma^*(D_j^+)=B_0$ and $\sigma^*(D_j^-)=B_r$.
    \item $\sigma^*(B_r)=B_r$ and $\sigma^*(D_r^+)=B_0$.
    \item $\sigma^*(B_r)=B_0$ and $\sigma^*(D_r^+)=B_r$.
\end{enumerate}

For Case (A),  the coefficient of $D^+_1$ is at most $-r(r-j)\leq -r$ in the following quantity,
\begin{equation}\label{dtru}
  \sum_{i=1}^{r-1}(r+1-i)\cdot \sigma^*(D^-_i)-\sum_{i=1}^{r}(r+1-i)\cdot \sigma^*(D^+_i)\,.
\end{equation}
This leads to a contradiction.

For Case (B), the coefficient of $D^-_1$ is at most $-r(r-j)\leq -r$ in (\ref{dtru}) which leads to a contradiction.

For Case (C),  by (\ref{b=b2}) we can show that $\sigma^*$ permutes $D_r^+$ and $B_0$, and fixes all the other divisors. However, then  $\sigma^*(K_{\mathcal T_{s,p,n}})-K_{\mathcal T_{s,p,n}}=(s-p)(B_0-D^+_r)$ which is a contradiction.

For Case (D),  we can derive a contradiction by computing the sum of coefficients in (\ref{b=b2}). More precisely,   write (\ref{b=b2}) over the basis $\left\{H,D^-_1, \cdots, D^-_{r-1},D^+_1,\cdots,D^+_r \right\}$; the sum of the coefficients on the left hand side is $1$, but that of the right hand side is $-1$. This is impossible.

Next, we assume that  $\sigma^*(B_0)\in\left\{B_0, B_r\right\}$. Checking the coefficient of $H$ in (\ref{b=b2}), we have the following possibilities. \begin{enumerate}[label=(\alph*).]
    \item  $\sigma^*(B_0)=B_0$ and $\sigma^*(B_r)=B_r$.
    \item $\sigma^*(B_0)=B_r$ and $\sigma^*(B_r)=B_0$.
\end{enumerate}

For Case (a),  $\sigma^*$ is the identity map by (\ref{b=b2}).

For Case (b), we can derive a contradiction by computing the sum of coefficients in (\ref{b=b2}).

\smallskip

{\bf\noindent Case 3 ($r\geq 2$ and $n-s=p=s$).} By (\ref{bstt0}) and (\ref{bsttr}) we have that
\begin{equation}\label{b=b4}
    \sigma^*( B_0)=\sigma^*(B_r)+\sum_{i=1}^{r-1}(r+1-i)\cdot \sigma^*(D^-_i)-\sum_{i=1}^{r-1}(r+1-i)\cdot \sigma^*(D^+_i)\,.
\end{equation}

If $\sigma^*(B_0)\in\left\{D^-_1, \cdots, D^-_{r-1},D^+_1,\cdots,D^+_{r-1}\right\}$, we have the following possibilities.
\begin{enumerate}[label=(\Alph*).]
    \item For a certain $1\leq j\leq r-1$, $\sigma^*(D_j^+)=B_r$ and $\sigma^*(D_j^-)=B_0$.
    \item For a certain $1\leq j\leq r-1$, $\sigma^*(D_j^-)=B_r$ and $\sigma^*(D_j^+)=B_0$.
\end{enumerate}

For Case (A),  the coefficient of $D^+_1$ is at most $-r(r-j)\leq -r$ in the following quantity,
\begin{equation}\label{dtru4}
  \sum_{i=1}^{r-1}(r+1-i)\cdot \sigma^*(D^+_i)-\sum_{i=1}^{r}(r+1-i)\cdot \sigma^*(D^-_i)\,.
\end{equation}
This leads to a contradiction.

For Case (B), the coefficient of $D^-_1$ is at most $-r(r-j)\leq -r$ in (\ref{dtru4}) which leads to a contradiction.

If  $\sigma^*(B_0)\in\left\{B_0, B_r\right\}$, we have the following possibilities. \begin{enumerate}[label=(\alph*).]
    \item $\sigma^*(B_r)=B_r$ and $\sigma^*(B_0)=B_0$.
    \item  $\sigma^*(B_r)=B_0$ and $\sigma^*(B_0)=B_r$.
\end{enumerate}

For Case (a),  $\sigma^*$ is the identity map by (\ref{b=b4}).

For Case (b), we can conclude that $\sigma^*=({\rm USD})^*$.
\smallskip

{\bf\noindent Case 4 ($r=p=n-s=1$).}
If $n=2$, Lemma \ref{sigma2} holds trivially.

When $n\geq 3$, we have  that $B_0=H-D^+_1$, $B_1=H$, and that
\begin{equation}
    K_{\mathcal T_{s,p,n}}=-n\cdot H+(n-2)\cdot D^+_1=-n\cdot B_0-2\cdot D^+_1.
\end{equation}
Hence, $\sigma^*$ is the identity map.

\smallskip

{\bf\noindent Case 5 ($r=p=1$ and $2\leq n-s\leq s$).} 
We have that $B_0=H-D^+_1$, $B_1=H-D^-_1$; $\sigma^*$ is a permutation of $\left\{B_0,B_1,D^+_1, D^-_1\right\}$; $K_{\mathcal T_{s,p,n}}=-n\cdot H+(s-1)\cdot D^+_1+(n-s-1)\cdot D^-_1$.

Notice that
\begin{equation}\label{r=11}
    \begin{split}
        K_{\mathcal T_{s,p,n}}&=-n\cdot (B_1+D^-_1)+(s-1)\cdot D^+_1+(n-s-1)\cdot D^-_1\\
        &=-n\cdot B_1+(s-1)\cdot D^+_1+(-s-1)\cdot D^-_1\,.
    \end{split}
\end{equation}
Since $\sigma^*(K_{\mathcal T_{s,p,n}})=K_{\mathcal T_{s,p,n}}$, we can conclude that $\sigma^*$ permutes $B_0$ and $B_1$. If $\sigma^*(B_0)=B_0$, then  $\sigma^*$ is the identity map.
If $\sigma^*(B_0)=B_1$ , then   $n=2s$ and $\sigma^*=({\rm USD})^*$.

\smallskip

{\bf\noindent Case 6 ($r=n-s=1$ and $2\leq p<s$).}
We have that $B_0=H-D^+_1$, $B_1=H-D^-_1$; $\sigma^*$ is a permutation of $\left\{B_0,B_1,D^+_1, D^-_1\right\}$; $K_{\mathcal T_{s,p,n}}=-n\cdot H+(n-p-1)\cdot D^+_1+(p-1)\cdot D^-_1$.

Notice that
\begin{equation}\label{r=13}
    \begin{split}
        K_{\mathcal T_{s,p,n}}&=-n\cdot (B_1+D^-_1)+(n-p-1)\cdot D^+_1+(p-1)\cdot D^-_1\\
        &=-n\cdot B_0+(n-p-1)\cdot D^+_1+(-n+p-1)\cdot D^-_1\,;
    \end{split}
\end{equation}
Then $\sigma^*$ permutes $B_0$ and $B_1$. If $\sigma^*(B_0)=B_0$, then  $\sigma^*$ is the identity map.
If $\sigma^*(B_0)=B_1$ , then  $n=2p$ and $\sigma^*=({\rm DUAL})^*$.

\medskip

We complete the proof of Lemma \ref{sigma2}.
\,\,\,$\endpf$

\subsection{\texorpdfstring{Proof of Lemma \ref{fibpre}}{rr}}\label{section:rigidtc} If $p=n-s$, Lemma \ref{fibpre} follows from Lemma \ref{pre}. Without loss of generality,  we may assume that $2p\leq n$, $p\neq n-s$ and $2\leq p<s$ in the following.

Recall the two fibration structures on $\mathcal M_{s,p,n}$ defined by (\ref{checkr}) and (\ref{checkr1}). Since $\sigma^*$ is the identity map on the Picard group of $\mathcal M_{s,p,n}$, similarly  to Lemma \ref{des} we can derive  the following.
\begin{enumerate}[label=(\alph*)]
    \item $\sigma$ induces an automorphism $\Sigma^1$ (resp. $\Sigma^2$) of $\mathbb P(N_{\mathcal V_{(p,0)}/G(p,n)})$ (resp. $\mathbb P(N_{\mathcal V_{(p-r,r)}/G(p,n)})$).
    \item $\Sigma^1$ (resp. $\Sigma^2$)  maps a fiber of $\kappa^1_{s,p,n}$ (resp. $\kappa^2_{s,p,n}$) to another fiber; $\sigma$ induces an automorphism $\widehat\Sigma^1$ (resp. $\widehat\Sigma^2$) of the base $\mathcal  V_{(p,0)}$ (resp. $\mathcal V_{(p-r,r)}$).
\end{enumerate}
Equivalently, we have the following commutative diagrams.
\begin{equation}
\begin{array}{ccccc}
\vspace{.02in}
\mathcal M_{s,p,n}&\xrightarrow{\,\,\,\,\,\sigma\,\,\,\,\,} &\mathcal M_{s,p,n}&&\\
\vspace{-.03in}
\Big\downarrow\llap{$\scriptstyle \tau^1_{s,p,n}\,\,\,\,\,$}&&\Big\downarrow\rlap{$\scriptstyle \tau^1_{s,p,n}\,\,\,\,\,$}&\\ \vspace{.02in}
\mathbb P(N_{\mathcal V_{(p,0)}/G(p,n)})&\xrightarrow{\,\,\,\,\,\Sigma^1\,\,\,\,\,}&\mathbb P(N_{\mathcal V_{(p,0)}/G(p,n)})\\
\vspace{-.03in}
\Big\downarrow\llap{$\scriptstyle \kappa^1_{s,p,n}\,\,\,\,\,$}&&\Big\downarrow\rlap{$\scriptstyle \kappa^1_{s,p,n}\,\,\,\,\,$}&\\
\mathcal V_{(p,0)}&\xrightarrow{\,\,\,\,\,\widehat\Sigma^1\,\,\,\,\,}&\mathcal V_{(p,0)}\\
\end{array}
\end{equation}
and
\begin{equation}
    \begin{array}{ccccc}
\vspace{.02in}
\mathcal M_{s,p,n}&\xrightarrow{\,\,\,\,\,\sigma\,\,\,\,\,} &\mathcal M_{s,p,n}&&\\
\vspace{-.03in}
\Big\downarrow\llap{$\scriptstyle \tau^2_{s,p,n}\,\,\,\,\,$}&&\Big\downarrow\rlap{$\scriptstyle \tau^2_{s,p,n}\,\,\,\,\,$}&\\ \vspace{.02in}
\mathbb P(N_{\mathcal V_{(p-r,r)}/G(p,n)})&\xrightarrow{\,\,\,\,\,\Sigma^2\,\,\,\,\,}&\mathbb P(N_{\mathcal V_{(p-r,r)}/G(p,n)})\\
\vspace{-.03in}
\Big\downarrow\llap{$\scriptstyle \kappa^2_{s,p,n}\,\,\,\,\,$}&&\Big\downarrow\rlap{$\scriptstyle \kappa^2_{s,p,n}\,\,\,\,\,$}&\\
\mathcal V_{(p-r,r)}&\xrightarrow{\,\,\,\,\,\widehat\Sigma^2\,\,\,\,\,}&\mathcal V_{(p-r,r)}\\
\end{array}\,.
\end{equation}

By Proposition \ref{fiber} and Example \ref{ccolli}, we can conclude that each fiber of $\kappa^1_{s,p,n}$ (resp. $\kappa^2_{s,p,n}$) is isomorphic to $P(M_{p\times (n-s)})$ (resp. $P(M_{p\times s})$ when $p<n-s$ or  $P(M_{(n-s)\times (n-p)})$ when $p>n-s$); each fiber of $\kappa^1_{s,p,n} \circ\tau^1_{s,p,n}$ (resp. $\kappa^2_{s,p,n}\circ\tau^2_{s,p,n}$) is isomorphic to the variety of complete collineations $\widetilde P(M_{p\times (n-s)})\cong\mathcal M_{p,p,n-s+p}$ (resp. $\widetilde P(M_{p\times s})\cong\mathcal M_{p,p,s+p}$ when $p<n-s$ or  $\widetilde P(M_{(n-s)\times (n-p)})\cong\mathcal M_{n-s,n-s,2n-s-p}$ when $p>n-s$). 

Since $\sigma^*$ is the identity map and hence the exceptional divisors are $\sigma$-invariant, $\Sigma^1$ (resp. $\Sigma^2$) preserves the ranks of the matrices in $P(M_{p\times (n-s)})$ (resp. $P(M_{p\times s})$ when $p<n-s$ or  $P(M_{(n-s)\times (n-p)})$ when $p>n-s$). Similarly to Lemma \ref{pre}, we can derive the following.
\begin{enumerate}[label=(\alph*).,start=3]
    \item The restriction of $\Sigma^1$ (resp. $\Sigma^2$) to each fiber of $\kappa^1_{s,p,n}$ (resp. $\kappa^2_{s,p,n}$) can be given by an element of $PGL(p,\mathbb C)\times PGL(n-s,\mathbb C)$ (resp. $PGL(p,\mathbb C)\times PGL(s,\mathbb C)$ when $p<n-s$, or  $PGL(n-s,\mathbb C)\times PGL(n-p,\mathbb C)$  when $p>n-s$). 
\end{enumerate}

In the following, we will prove that $\sigma\in PGL(s,\mathbb C)\times PGL(n-s,\mathbb C)$ by checking Properties (a), (b), (c) in the two fibration structures. 
\smallskip

{\bf\noindent Case 1 ($2\leq p<n-s\leq s$ and $s\neq 2p$).}
Recall that the automorphism group of $\mathcal V_{(p,0)}\cong G(p,s)$ is $PGL(s,\mathbb C)$. Hence,  there is an element $\eta\in PGL(s,\mathbb C)\times PGL(n-s,\mathbb C)$ such that the composition $\eta\circ\sigma$ induces the identity map on the base $\mathcal V_{(p,0)}$. Without loss of generality, we may assume that $\widehat \Sigma^1$ is the identity map.

Recall the local coordinate charts defined by (\ref{nbc1}), (\ref{nbc3}), (\ref{trr1}), (\ref{trr41}) in Appendix \ref{section:projbc} for the normal bundles $N_{\mathcal V_{(p,0)}/G(p,n)}$,  $N_{\mathcal V_{(p-r,r)}/G(p,n)}$, and the projective bundles $\mathbb P(N_{\mathcal V_{(p,0)}/G(p,n)})$, $\mathbb P(N_{\mathcal V_{(p-r,r)}/G(p,n)})$  respectively.
Similarly to (\ref{nbc1}), we can equip the  open set $U^N$ of  $N_{\mathcal V_{(p,0)}/G(p,n)}$ with local coordinates
\begin{equation}\label{coorw}
\small
    W:=\left(\begin{matrix}
    w_{11}&w_{12}&\cdots&w_{1(s-p)}\\
    w_{21}&w_{22}&\cdots&w_{2(s-p)}\\
    \cdots&\cdots&\cdots&\cdots\\
    w_{p1}&w_{p2}&\cdots&w_{p(s-p)}\\
    \end{matrix}\right),\,\,\,\,(X,Z):=\left(\begin{matrix}
    x_{1(s+1)}&\cdots&x_{1(s+p)}&z_{1(s+p+1)}&\cdots&z_{1n}\\
    x_{2(s+1)}&\cdots&x_{2(s+p)}&z_{2(s+p+1)}&\cdots&z_{2n}\\
    \cdots&\cdots&\cdots&\cdots&\cdots&\cdots\\
    x_{p(s+1)}&\cdots&x_{p(s+p)}&z_{p(s+p+1)}&\cdots&z_{pn}\\
    \end{matrix}\right);
\end{equation}
identify $U^N$ with an open set of $G(p,n)$ by
\begin{equation}\label{coorem}
  \big(W,(X,Z)\big)\mapsto   \left(W\,\hspace{-0.15in}\begin{matrix}
  &\hfill\tikzmark{c1}\\
  &\hfill\tikzmark{d1}
  \end{matrix}\,\,\,\, I_{p\times p}\hspace{-0.13in}\begin{matrix}
  &\hfill\tikzmark{c2}\\
  &\hfill\tikzmark{d2}
  \end{matrix}\hspace{-0.1in}\begin{matrix}
  &\hfill\tikzmark{c3}\\
  &\hfill\tikzmark{d3}
  \end{matrix}\,\,\,\, 
    X\hspace{-0.12in}\begin{matrix}
  &\hfill\tikzmark{c4}\\
  &\hfill\tikzmark{d4}
  \end{matrix}\,\,\,\,Z\,\right)\,.
  \tikz[remember picture,overlay]   \draw[dashed,dash pattern={on 4pt off 2pt}] ([xshift=0.5\tabcolsep,yshift=7pt]c1.north) -- ([xshift=0.5\tabcolsep,yshift=-2pt]d1.south);\tikz[remember picture,overlay]   \draw[dashed,dash pattern={on 4pt off 2pt}] ([xshift=0.5\tabcolsep,yshift=7pt]c2.north) -- ([xshift=0.5\tabcolsep,yshift=-2pt]d2.south);\tikz[remember picture,overlay]   \draw[dashed,dash pattern={on 4pt off 2pt}] ([xshift=0.5\tabcolsep,yshift=7pt]c3.north) -- ([xshift=0.5\tabcolsep,yshift=-2pt]d3.south);\tikz[remember picture,overlay]   \draw[dashed,dash pattern={on 4pt off 2pt}] ([xshift=0.5\tabcolsep,yshift=7pt]c4.north) -- ([xshift=0.5\tabcolsep,yshift=-2pt]d4.south);
\end{equation}
Equip  $U^P$ with local coordinates $(W,[X,Z])$ where $U^P$ is the corresponding open subset of  $\mathbb P(N_{\mathcal V_{(p,0)}/G(p,n)})$ and $[X,Z]$ are the homogeneous coordinates such that for $\lambda\in\mathbb C^*$
\begin{equation}
\footnotesize
    \left[\begin{matrix}
    x_{1(s+1)}&\cdots&x_{1(s+p)}&z_{1(s+p+1)}&\cdots&z_{1n}\\
    x_{2(s+1)}&\cdots&x_{2(s+p)}&z_{2(s+p+1)}&\cdots&z_{2n}\\
    \cdots&\cdots&\cdots&\cdots\\
    x_{p(s+1)}&\cdots&x_{p(s+p)}&z_{p(s+p+1)}&\cdots&z_{pn}\end{matrix}\right]=\lambda\cdot\left[\begin{matrix}
    x_{1(s+1)}&\cdots&x_{1(s+p)}&z_{1(s+p+1)}&\cdots&z_{1n}\\
    x_{2(s+1)}&\cdots&x_{2(s+p)}&z_{2(s+p+1)}&\cdots&z_{2n}\\
    \cdots&\cdots&\cdots&\cdots&\cdots&\cdots\\
    x_{p(s+1)}&\cdots&x_{p(s+p)}&z_{p(s+p+1)}&\cdots&z_{pn}\end{matrix}\right]\,.
\end{equation}

Similarly to (\ref{nbc3}), we can equip the open set $\widehat  U^N$ of $N_{\mathcal V_{(p-r,r)}/G(p,n)}$ with the 
local coordinates
\begin{equation}\label{coorw1}
\small
(\widetilde W,\widetilde X):=\left(\begin{matrix}
\widetilde w_{11}&\cdots&\widetilde w_{1(s-p)}&\widetilde x_{1(s-p+1)}&\cdots&\widetilde x_{1s}\\
\widetilde w_{21}&\cdots&\widetilde w_{2(s-p)}&\widetilde x_{2(s-p+1)}&\cdots&\widetilde x_{2s}\\
\cdots&\cdots&\cdots&\cdots&\cdots&\cdots\\
\widetilde w_{p1}&\cdots&\widetilde w_{p(s-p)}&\widetilde x_{p(s-p+1)}&\cdots&\widetilde x_{ps}\\
\end{matrix}\right),\,\,\,\,\,\widetilde Z:=\left(\begin{matrix}
\widetilde z_{1(s+p+1)}&\cdots&\widetilde z_{1n}\\
\widetilde z_{2(s+p+1)}&\cdots&\widetilde z_{2n}\\
\cdots&\cdots&\cdots\\
\widetilde z_{p(s+p+1)}&\cdots&\widetilde z_{pn}\\
\end{matrix}\right);
\end{equation}
identify $\widehat U^N$ with an open set of $G(p,n)$ by
\begin{equation}\label{coorew1}
  \big((\widetilde W,\widetilde X),\widetilde Z)\big)\mapsto   \left(\widetilde W\,\hspace{-0.15in}\begin{matrix}
  &\hfill\tikzmark{c1}\\
  &\hfill\tikzmark{d1}
  \end{matrix}\,\,\,\, 
  \widetilde X
  \hspace{-0.13in}\begin{matrix}
  &\hfill\tikzmark{c2}\\
  &\hfill\tikzmark{d2}
  \end{matrix}\hspace{-0.1in}\begin{matrix}
  &\hfill\tikzmark{c3}\\
  &\hfill\tikzmark{d3}
  \end{matrix}\,\,\,\, 
    I_{p\times p}\hspace{-0.12in}\begin{matrix}
  &\hfill\tikzmark{c4}\\
  &\hfill\tikzmark{d4}
  \end{matrix}\,\,\,\,\widetilde Z\,\right)\,.
  \tikz[remember picture,overlay]   \draw[dashed,dash pattern={on 4pt off 2pt}] ([xshift=0.5\tabcolsep,yshift=7pt]c1.north) -- ([xshift=0.5\tabcolsep,yshift=-2pt]d1.south);\tikz[remember picture,overlay]   \draw[dashed,dash pattern={on 4pt off 2pt}] ([xshift=0.5\tabcolsep,yshift=7pt]c2.north) -- ([xshift=0.5\tabcolsep,yshift=-2pt]d2.south);\tikz[remember picture,overlay]   \draw[dashed,dash pattern={on 4pt off 2pt}] ([xshift=0.5\tabcolsep,yshift=7pt]c3.north) -- ([xshift=0.5\tabcolsep,yshift=-2pt]d3.south);\tikz[remember picture,overlay]   \draw[dashed,dash pattern={on 4pt off 2pt}] ([xshift=0.5\tabcolsep,yshift=7pt]c4.north) -- ([xshift=0.5\tabcolsep,yshift=-2pt]d4.south);
\end{equation}
Equip $\widehat U^P$ with the local coordinates $([\widetilde W,\widetilde X],\widetilde Z)$ where $[\widetilde W,\widetilde X]$ are the corresponding homogeneous coordinates.

In terms of the above local coordinates for $U^P$, we can conclude by Property (c) the following. For fixed $W$   there are matrices $h(W)\in SL(p,\mathbb C)$ and $g(W)\in SL(n-s,\mathbb C)$ such that 
\begin{equation}
    \Sigma^1\big((W,[X,Z])\big)=\big(W,[h(W)\cdot (X,Z)\cdot g(W)]\big)\,,
\end{equation}
where the product is the matrix multiplication.
Since the choice of $h(W)$ and $g(W)$  is unique up to a finite cover (for the natural maps $SL(m,\mathbb C)\rightarrow PGL(m,\mathbb C)$ are finite),  locally we can take $h(W)$ and $g(W)$ as matrix-valued holomorphic functions in the variables $W$. Moreover, we can assume that $h(W)$ and $g(W)$ are  well-defined on $\mathbb C^{p(s-p)}$  for which is simply connected.

Notice that the isomorphism $\mathcal L_{s,p,n}$ defined in Definition \ref{fi} induces a birational map from $\mathbb P(N_{\mathcal V_{(p-r,r)}/G(p,n)})$ to $\mathbb P(N_{\mathcal V_{(p,0)}/G(p,n)})$; by identifying the normal bundles with $G(p,n)$ birationally,   there is a natural birational map $\mathfrak L_{s,p,n}:N_{\mathcal V_{(p-r,r)}/G(p,n)}\dashrightarrow N_{\mathcal V_{(p,0)}/G(p,n)}$ which induces $\mathcal L_{s,p,n}$ on the projective bundles. We can write explicitly $\mathfrak L_{s,p,n}$ in terms of the above local coordinates as follows. $\mathfrak L_{s,p,n} \left(\big((\widetilde W,\widetilde X),\widetilde Z\big)\right)=\big(W,(X,Z)\big)$ where $\big((\widetilde W,\widetilde X),\widetilde Z\big)$ and $\big(W,(X,Z)\big)$  represent the same point in $G(p,n)$  under the identifications (\ref{coorem}) and (\ref{coorew1}); equivalently,
\begin{equation}\label{coorew3}
    \left\{\begin{array}{l}
       \widetilde X=X^{-1}   \\
        \widetilde W=X^{-1}\cdot W\\
       \widetilde Z=X^{-1}\cdot Z\\ 
    \end{array}\right.\,\,\,\,{\rm and}\,\,\,\,\,\, \left\{\begin{array}{l}
        X=\widetilde X^{-1}   \\
         W=\widetilde X^{-1}\cdot\widetilde W\\
        Z=\widetilde X^{-1}\cdot \widetilde Z\\ 
    \end{array}\right.\,.
\end{equation}

Through $\mathfrak L_{s,p,n}$ we can derive an explicit form of the automorphism $\Sigma^2$ of $\mathbb P(N_{\mathcal V_{(p-r,r)}/G(p,n)})$ in terms of $h$, $g$. In particular, $\widehat \Sigma^2$ is given in matrix representatives by a right multiplication by $g$, that is,
\begin{equation}
  \widehat \Sigma^2 \left(\left( I_{p\times p} \,\hspace{-0.15in}\begin{matrix}
  &\hfill\tikzmark{c1}\\
  &\hfill\tikzmark{d1}
  \end{matrix}\,\,\,\, 
  \widetilde Z\right)\right)=\left(
    I_{p\times p}\hspace{-0.12in}\begin{matrix}
  &\hfill\tikzmark{c4}\\
  &\hfill\tikzmark{d4}
  \end{matrix}\,\,\,\,\widetilde Z\,\right)\cdot g\left(\widetilde X^{-1}\cdot\widetilde W\right)\,.
  \tikz[remember picture,overlay]   \draw[dashed,dash pattern={on 4pt off 2pt}] ([xshift=0.5\tabcolsep,yshift=7pt]c1.north) -- ([xshift=0.5\tabcolsep,yshift=-2pt]d1.south);\tikz[remember picture,overlay]   \draw[dashed,dash pattern={on 4pt off 2pt}] ([xshift=0.5\tabcolsep,yshift=7pt]c4.north) -- ([xshift=0.5\tabcolsep,yshift=-2pt]d4.south);
\end{equation}
By Property (b) the automorphism $ \widehat \Sigma^2 $ is independent of the variables $\widetilde X, \widetilde W$; hence, we can take an element $({\rm Id},\tau)\in PGL(s,\mathbb C)\times PGL(n-s,\mathbb C)$ such that $({\rm Id},\tau)\circ\sigma$ induces the identity map on the base $\mathcal V_{(p-r,r)}$. Noticing that $({\rm Id},\tau)$ induces the identity on the base $\mathcal V_{(p,0)}$, without loss of generality, we may assume that $\widehat\Sigma^1$, $\widehat\Sigma^2$ are both identity, and that $g$ is the identity matrix. Then $\Sigma^2$ takes the following form in terms of the local coordinates $\big([\widetilde W,\widetilde X],\widetilde Z\big)$.
\begin{equation}
    \Sigma^2\left([\widetilde W,\widetilde X],\widetilde Z\right)=\left(\left[\widetilde X\cdot h^{-1}\cdot\widetilde X^{-1}\cdot \widetilde W\,\,,\,\,\,\widetilde X\cdot h^{-1}\right],\widetilde Z\right)\,.
\end{equation}

It is clear that to prove Lemma \ref{fibpre}, it suffices to show that $h$ can be taken as the identity matrix. Similarly by Property (c)  we can find matrix-valued holomorphic functions $E(\widetilde Z)\in SL(p,\mathbb C)$ and $F(\widetilde Z)\in SL(s,\mathbb C)$ such that
\begin{equation}\label{easy1}
E(\widetilde Z)\cdot [\widetilde W,\widetilde X]\cdot F(\widetilde Z)=\Sigma^2\left([\widetilde W,\widetilde X],\widetilde Z\right)=\left[\widetilde X\cdot h^{-1}\cdot\widetilde X^{-1}\cdot \widetilde W\,\,,\,\,\,\widetilde X\cdot h^{-1}\right]\,.
\end{equation}

{\bf\noindent Claim.} We can choose $E(\widetilde Z)$ and $F(\widetilde Z)$ as (matrix-valued) constant  functions.
\smallskip

{\bf \noindent Proof of Claim.} Write $E$ and $F$ in components as 
\begin{equation}
    E:=\left(\begin{matrix}
    e_{11}&e_{12}&\cdots&e_{1p}\\
    e_{21}&e_{22}&\cdots&e_{2p}\\
    \cdots&\cdots&\cdots&\cdots\\
    e_{p1}&e_{p2}&\cdots&e_{pp}\\
    \end{matrix}\right)\,\,\,\,{\rm and}\,\,\,\,\,\,F:=\left(\begin{matrix}
    f_{11}&f_{12}&\cdots&f_{1s}\\
    f_{21}&f_{22}&\cdots&f_{2s}\\
    \cdots&\cdots&\cdots&\cdots\\
    f_{s1}&f_{s2}&\cdots&f_{ss}\\
    \end{matrix}\right)\,.
\end{equation}
For $1\leq i\leq p$ and $1\leq j\leq s$, denote by $B_{ij}$ the $(i,j)^{th}$ entry of the matrix $E(\widetilde Z)\cdot [\widetilde W,\widetilde X]\cdot F(\widetilde Z)$. There are integers $1\leq c\leq p$ and $1\leq d\leq s$ such that $e_{1c}\not\equiv 0$ and $f_{d1}\not\equiv 0$. Since the right hand side of (\ref{easy1}) is independent of the variables $\widetilde Z$,   the following quotient is a rational function independent of $\widetilde Z$ for each $1\leq a\leq p$ and $1\leq b\leq s$. 
\begin{equation}\label{bak}
    \frac{B_{ab}}{B_{11}}=\frac{\sum\limits_{i=1}^{p}\sum\limits_{j=1}^{s-p}\left(e_{ai}(\widetilde Z) f_{jb}(\widetilde Z)\right)\widetilde W_{ij}+\sum\limits_{i=1}^{p}\sum\limits_{j=s-p+1}^{s}\left(e_{ai}(\widetilde Z) f_{jb}(\widetilde Z)\right)\widetilde X_{ij}}{\sum\limits_{i=1}^{p}\sum\limits_{j=1}^{s-p}\left(e_{1i}(\widetilde Z) f_{j1}(\widetilde Z)\right)\widetilde W_{ij}+\sum\limits_{i=1}^{p}\sum\limits_{j=s-p+1}^{s}\left(e_{1i}(\widetilde Z) f_{j1}(\widetilde Z)\right)\widetilde X_{ij}}\,\,\,\,\,.
\end{equation}
We can show that(\ref{bak}) is not a constant function  if $a\neq 1$ or $b\neq 1$ for the matrices $E$ and $F$ are non-degenerate. Then we can conclude that the following quotient is a constant function.
\begin{equation}
    \frac{e_{ai}(\widetilde Z) f_{jb}(\widetilde Z)}{e_{1c}(\widetilde Z) f_{d1}(\widetilde Z)}\,,\,\,1\leq a,i\leq p,\,1\leq b,j\leq s.
\end{equation}
Hence, there are constants $\alpha_{ij}$ (resp. $\beta_{ij}$) for $1\leq i,j\leq p$ (resp. $1\leq i,j\leq s$)  such that
\begin{equation}
e_{ij}(\widetilde Z)=\beta_{ij}\cdot e_{1c}(\widetilde Z)\,\,\,\,\,\left({\rm resp.}\,\,f_{ij}(\widetilde Z)=\alpha_{ij}\cdot f_{d1}(\widetilde Z)\right)\,.
\end{equation}

We complete the proof of Claim.\,\,\,\,\,$\endpf$
\smallskip

We assume that $E,F$ are constant matrices in the following. Write $F$ in blocks as 
$F=\left(\begin{matrix}
      F_{11}&F_{12}\\
      F_{21}&F_{22}
    \end{matrix}\right)\,$ where $F_{22}$ is  $p\times p$. Then (\ref{easy1}) yields
\begin{equation}\label{asy1}
\small
\begin{split}
\left[E\cdot \left(\widetilde W\cdot F_{11}\right.\right.&\left.\left. +\widetilde X\cdot F_{21}\right)\,\,,\,\,\,E\cdot \left(\widetilde W\cdot F_{12}+\widetilde X\cdot F_{22}\right)\right]=\left[\widetilde X\cdot h^{-1}\cdot\widetilde X^{-1}\cdot \widetilde W\,\,,\,\,\,\widetilde X\cdot h^{-1}\right]\,,
\end{split}
\end{equation}
and hence
\begin{equation}
\small
\hspace{-.15in}\left(E\cdot \left(\widetilde W\cdot F_{12}+\widetilde X\cdot F_{22}\right)\right)^{-1}\left(E\cdot \left(\widetilde W\cdot F_{11}+\widetilde X\cdot F_{21}\right)\right)=\left(\widetilde X\cdot h^{-1}\right)^{-1}\left(\widetilde X\cdot h^{-1}\cdot\widetilde X^{-1}\cdot \widetilde W\right),    
\end{equation}
\begin{equation}\label{const}
\small
\widetilde W\cdot F_{11}+\widetilde X\cdot F_{21}=\widetilde W\cdot F_{12}\cdot\widetilde X^{-1}\cdot \widetilde W+\widetilde X\cdot F_{22}\cdot\widetilde X^{-1}\cdot \widetilde W\,.
\end{equation}
Since $F$ is a constant matrix, (\ref{const}) holds if only if
$F$ is a scalar, that is,
$F=\lambda\cdot I_{s\times s}$ for a certain $\lambda\in\mathbb C^*$. Substituting it into (\ref{asy1}) we conclude that $E\cdot\widetilde X=\widetilde X\cdot h^{-1}$ and hence $E\cdot X^{-1}=X^{-1}\cdot h^{-1}$. Since $E$ is a constant matrix, we can conclude that $h(W)\equiv\mu\cdot I_{p\times p}$ for a certain $\mu\in\mathbb C^*$.

We conclude that $\sigma\in PGL(s,\mathbb C)\times PGL(n-s,\mathbb C)$.

\smallskip

{\bf\noindent Case 2 ($2\leq p<n-s$ and $s=2p$).}
Recall that the automorphism group of $\mathcal V_{(p,0)}\cong G(p,2p)$ is generated by $PGL(2p,\mathbb C)$ and the dual map (\ref{gdual}). If $\widehat\Sigma^1$ or $\widehat\Sigma^2$ is the identity map, we can prove similarly to Case 1. 

Without loss of generality, we may assume that $\widehat\Sigma^1$ is given by the dual map (\ref{gdual}).

Take the local coordinate charts as in (\ref{coorw}) and (\ref{coorw1}). By Property (c) we can choose matrix-valued holomoprhic functions $h(W)\in SL(p,\mathbb C)$ and $g(W)\in SL(n-s,\mathbb C)$ such that that
\begin{equation}
    \Sigma^1\big((W,[X,Z])\big)=\big(-(W^{-1})^T,\,\,[h(W)\cdot (X,Z)\cdot g(W)]\big)\,.
\end{equation}
Similarly to Case 1, we can verify that the induced automorphism $\widehat \Sigma^2$ of $\mathcal V_{(p-r,r)}$ is given  by
\begin{equation}
  \widehat \Sigma^2 \left(\left( I_{p\times p} \,\hspace{-0.15in}\begin{matrix}
  &\hfill\tikzmark{c1}\\
  &\hfill\tikzmark{d1}
  \end{matrix}\,\,\,\, 
  \widetilde Z\right)\right)=\left(
    I_{p\times p}\hspace{-0.12in}\begin{matrix}
  &\hfill\tikzmark{c4}\\
  &\hfill\tikzmark{d4}
  \end{matrix}\,\,\,\,\widetilde Z\,\right)\cdot g\left(\widetilde X^{-1}\cdot\widetilde W\right)\,.
  \tikz[remember picture,overlay]   \draw[dashed,dash pattern={on 4pt off 2pt}] ([xshift=0.5\tabcolsep,yshift=7pt]c1.north) -- ([xshift=0.5\tabcolsep,yshift=-2pt]d1.south);\tikz[remember picture,overlay]   \draw[dashed,dash pattern={on 4pt off 2pt}] ([xshift=0.5\tabcolsep,yshift=7pt]c4.north) -- ([xshift=0.5\tabcolsep,yshift=-2pt]d4.south);
\end{equation}
Then there is an an element $({\rm Id},\tau)\in PGL(s,\mathbb C)\times PGL(n-s,\mathbb C)$ such that $({\rm Id},\tau)\circ\sigma$ induces the identity map on $\mathcal V_{(p-r,r)}$.

Now we can repeat the argument in Case 1 by switching the places of $\Sigma^1$ and   $\Sigma^2$ therein. 

We thus conclude that $\sigma\in PGL(s,\mathbb C)\times PGL(n-s,\mathbb C)$.

\smallskip

{\bf\noindent Case 3 ($2\leq p$ and $n-s<p<s$).}  The proof is similar to Case 1. First notice that  $s\neq 2(s+p-n)$ for $2p\leq n$.  We may thus assume that $\sigma$ induces the identity map on the base $\mathcal V_{(p-r,r)}\cong G(s+p-n,s)$.

Recall the local coordinate charts defined by (\ref{nbc1}), (\ref{nbc2}), (\ref{trr1}), (\ref{trr2}) in Appendix \ref{section:projbc} for the normal bundles $N_{\mathcal V_{(p,0)}/G(p,n)}$,  $N_{\mathcal V_{(p-r,r)}/G(p,n)}$, and the projective bundles $\mathbb P(N_{\mathcal V_{(p,0)}/G(p,n)})$, $\mathbb P(N_{\mathcal V_{(p-r,r)}/G(p,n)})$  respectively. Similarly to (\ref{nbc1}), we equip the open set $U^N$ of $N_{\mathcal V_{(p,0)}/G(p,n)}$ with  local coordinates
\begin{equation}\label{c34}
\small
 \hspace{-.15in}     W_1:=\left(\begin{matrix}
     w_{11}&  w_{12} &\cdots& w_{1(s-p)}\\
     w_{21}& w_{22}&\cdots& w_{2(s-p)}\\
    \cdots&\cdots&\cdots&\cdots\\
     w_{(n-s)1} & w_{(n-s)2} &\cdots& w_{(n-s)(s-p)}\\
    \end{matrix}\right),\,\,\,\,\,\,\,\,\, W_2:= \left(\begin{matrix}
     w_{1(s+1)} &\cdots&  w_{1n} \\
     w_{2(s+1)} & \cdots &  w_{2n}\\
    \cdots&\cdots&\cdots\\
     w_{(n-s)(s+1)} &\cdots&  w_{(n-s)n}\\
    \end{matrix}\right),
\end{equation}
\begin{equation*}
\small
     Z_1:= \left(\begin{matrix}
      z_{(n-s+1)1}&   z_{(n-s+1)2}&\cdots&  z_{(n-s+1)(s-p)}\\
    z_{(n-s+2)1}&  z_{(n-s+2)2}&\cdots&  z_{(n-s+2)(s-p)}\\
    \cdots&\cdots&\cdots&\cdots\\
      z_{p1}&  z_{p2}&\cdots&  z_{p(s-p)}\\
    \end{matrix}\right),\,\,\,
     Z_2:=\left(
    \begin{matrix}
       z_{(n-s+1)(s+1)}&\cdots& z_{(n-s+1)n}\\
       z_{(n-s+2)(s+1)}&\cdots& z_{(n-s+2)n}\\
      \cdots&\cdots&\cdots\\
       z_{p(s+1)}&\cdots& z_{pn}
    \end{matrix}\right);
\end{equation*}
identify $U^N$ with an open set of $G(p,n)$ by
\begin{equation}\label{c36}
\big(( W_1, W_2) ,( Z_1, Z_2)\big)\mapsto\,\,\,\,\left(\begin{matrix}
   W_1&I_{(n-s)\times(n-s)}&0& W_2\\  Z_1&0&I_{(s+p-n)\times(s+p-n)}& Z_2
\end{matrix}
\right)\,\,.
\end{equation}
Let ${U^P}$ be the corresponding open set of  $\mathbb P(N_{\mathcal V_{(p,0)}/G(p,n)})$.
Equip ${U^P}$ with the local coordinates $\left( W_1,  X_1, \left[\begin{matrix}
      W_2\\
      Z_2
\end{matrix}\right]\right)\,$ where $\left[\begin{matrix}
        W_2\\
        Z_2
     \end{matrix}\right]\,$
are the corresponding homogeneous coordinates.

Similarly to (\ref{nbc2}), we equip the open set $\widehat U^N$ of  $N_{\mathcal V_{(p-r,r)}/G(p,n)}$ with local coordinates
\begin{equation}\label{c31}
\small
   \hspace{-.25in} \widetilde W_1:=\left(\begin{matrix}
   \widetilde  w_{11}&\widetilde w_{12}&\cdots&\widetilde w_{1(s-p)}\\
    \widetilde w_{21}&\widetilde w_{22}&\cdots&\widetilde w_{2(s-p)}\\
    \cdots&\cdots&\cdots&\cdots\\
    \widetilde w_{(n-s)1}&\widetilde w_{(n-s)2}&\cdots&\widetilde w_{(n-s)(s-p)}\\
    \end{matrix}\right),\,\,\,\widetilde W_2:=\left(\begin{matrix}
  \widetilde  w_{1(s-p+1)}&\cdots&\widetilde w_{1(n-p)}\\
   \widetilde w_{2(s-p+1)} &\cdots&\widetilde w_{2(n-p)}\\
    \cdots&\cdots&\cdots\\
  \widetilde  w_{(n-s)(s-p+1)}&\cdots&\widetilde w_{(n-s)(n-p)}\\
    \end{matrix}\right),\\
\end{equation}
\begin{equation*}
\small
 \widetilde    Z_1:=\left(\begin{matrix}
    z_{(n-s+1)1}&\widetilde z_{(n-s+1)2}&\cdots&\widetilde z_{(n-s+1)(s-p)}\\
    z_{(n-s+2)1}&\widetilde z_{(n-s+2)2}&\cdots&\widetilde z_{(n-s+2)(s-p)}\\
    \cdots&\cdots&\cdots&\cdots\\
  \widetilde   z_{p1}&\widetilde z_{p2}&\cdots&\widetilde z_{p(s-p)}\\
    \end{matrix}\right),\,\,\,
   \widetilde  Z_2:=\left(
    \begin{matrix}
    \widetilde   z_{(n-s+1)(s-p+1)}&\cdots&\widetilde z_{(n-s+1)(n-p)}\\
     \widetilde  z_{(n-s+2)(s-p+1)}&\cdots&\widetilde z_{(n-s+2)(n-p)}\\
      \cdots&\cdots&\cdots\\
     \widetilde  z_{p(s-p+1)}&\cdots&\widetilde z_{p(n-p)}
    \end{matrix}\right);
\end{equation*}
identify $\widehat U^N$ with an open set of $G(p,n)$ by
\begin{equation}\label{c33}
\big((\widetilde W_1,\widetilde W_2),(\widetilde Z_1,\widetilde Z_2)\big)\mapsto\,\,\,\,\left(\begin{matrix}
 \widetilde  W_1&\widetilde W_2&0&I_{(n-s)\times(n-s)}\\ \widetilde Z_1&\widetilde Z_2&I_{(s+p-n)\times(s+p-n)}&0
\end{matrix}
\right)\,\,.
\end{equation}
Denote by $\big([\widetilde W_1,\widetilde W_2],(\widetilde Z_1,\widetilde Z_2)\big)$ the local coordinates of $U^P$ where  $U^P$ is the corresponding open set of $\mathbb P(N_{\mathcal V_{(p-r,r)}/G(p,n)})$ and $[\widetilde W_1,\widetilde W_2]$ are the homogeneous coordinates.

By identifications (\ref{nin1}), (\ref{nin2}), there is  a natural birational map $\mathfrak L_{s,p,n}$ from $N_{\mathcal V_{(p-r,r)}/G(p,n)}$ to $N_{\mathcal V_{(p,0)}/G(p,n)}$ such that it induces the birational map $\mathcal L_{s,p,n}$ on the projective bundles.  In the above local coordinate charts, $\mathfrak L_{s,p,n}\left( \big((\widetilde W_1,\widetilde W_2) ,(\widetilde Z_1,\widetilde Z_2)\big)\right)=\big((W_1,W_2) ,(Z_1,Z_2)\big)$ where 
\begin{equation}\label{lawe}
    \left\{\begin{array}{l}
       \widetilde W_2=W_2^{-1}   \\
        \widetilde W_1=W_2^{-1}\cdot W_1\\
       \widetilde Z_1=Z_1-Z_2\cdot W_2^{-1}\cdot W_1\\
       \widetilde Z_2=-Z_2\cdot W_2^{-1}\\ 
    \end{array}\right.\,\,\,\,{\rm and}\,\,\,\,\,\, \left\{\begin{array}{l}
        W_2=\widetilde W_2^{-1}   \\
     W_1=\widetilde W_2^{-1}\cdot \widetilde W_1\\ Z_1=\widetilde Z_1-\widetilde Z_2\cdot \widetilde W_2^{-1}\cdot \widetilde W_1\\
      Z_2=- \widetilde  Z_2\cdot  \widetilde  W_2^{-1}\\ 
    \end{array}\right.\,.
\end{equation}

Similarly to Case 1, we can derive  holomorphic functions  $E(\widetilde Z_1,\widetilde Z_2)\in SL(n-s,\mathbb C)$ and $F(\widetilde Z_1,\widetilde Z_2):=\left(\begin{matrix}
  F_{11}&F_{12}\\
  F_{21}&F_{22}
\end{matrix}\right)\in SL(n-p,\mathbb C)$
such  that
\begin{equation}
\begin{split}
&\Sigma^2\big(\big([\widetilde W_1,\widetilde W_2],(\widetilde Z_1,\widetilde Z_2)\big)\big)=\big(\big([E\cdot(\widetilde W_1,\widetilde W_2)\cdot F],(\widetilde Z_1,\widetilde Z_2)\big)\big)\\
&=\big(\left[E\cdot (\widetilde W_1\cdot F_{11}+\widetilde W_2\cdot F_{21}),E\cdot (\widetilde W_1\cdot F_{12}+\widetilde W_2\cdot F_{22})\right],(\widetilde Z_1,\widetilde Z_2)\big)\,.
\end{split}
\end{equation}
Through $\mathfrak L_{s,p,n}$ we can derive $\widehat \Sigma^1$  explicitly as follows.
\begin{equation}
\small
\widehat \Sigma^1 \left(\left(\begin{matrix}
   W_1&I_{(n-s)\times(n-s)}&0\\  Z_1&0&I_{(s+p-n)\times(s+p-n)}\end{matrix}\right)\right)=\left(\begin{matrix}
  G&I_{(n-s)\times(n-s)}&0\\ K&0&I_{(s+p-n)\times(s+p-n)}\end{matrix}\right)\,
\end{equation}
where 
\begin{equation}\label{lawe2}
    G:=\left( W_1\cdot F_{12}+ F_{22}\right)^{-1} \cdot \left( W_1\cdot F_{11}+F_{21} \right)\,
\end{equation}
and
\begin{equation}
\begin{split}
   K:&=Z_1- Z_2\cdot  W_2^{-1}\cdot  W_1+ Z_2\cdot  W_2^{-1}\cdot \left( W_1\cdot F_{12}+ F_{22}\right)^{-1}\cdot \left( W_1\cdot F_{11}+ F_{21}\right)\\
   &= Z_1- Z_2\cdot W_2^{-1}\cdot  (W_1-G)\,.
\end{split}
\end{equation}
By Property (b) $\widehat \Sigma^1$ is independent of the variables $ W_2$ and $ Z_2$; hence matrices $G$ and $K$ are holomorphic functions in the variables $ W_1$ and $ Z_1$. We can thus conclude that $G= W_1$ and $K= Z_1$. Substituting $G= W_1$ into (\ref{lawe2}), we have that
\begin{equation}
     \left( W_1\cdot F_{12}+ F_{22}\right)\cdot W_1=\left( W_1\cdot F_{11}+F_{21} \right)\,,
\end{equation}
or equivalently,
\begin{equation}
\small
\left(\widetilde W_2^{-1}\cdot \widetilde W_1\cdot F_{12}(\widetilde Z_1,\widetilde Z_2)+ F_{22}(\widetilde Z_1,\widetilde Z_2)\right)\cdot \widetilde W_2^{-1}\cdot\widetilde W_1= \left(\widetilde W_2^{-1}\cdot \widetilde W_1\cdot F_{11}(\widetilde Z_1,\widetilde Z_2)+F_{21}(\widetilde Z_1,\widetilde Z_2) \right)\,.
\end{equation}
Then it is easy to verify that $F(\widetilde Z_1,\widetilde Z_2)\equiv\lambda(\widetilde Z_1,\widetilde Z_2)\cdot I_{(n-p)\times (n-p)}$ where $\lambda(\widetilde Z_1,\widetilde Z_2)$ is a holomorphic scalar function.

Without loss of generality, we can assume that $F(\widetilde Z_1,\widetilde Z_2)\equiv I_{(n-p)\times (n-p)}$. In what follows, we will show that $E(\widetilde Z_1,\widetilde Z_2)$ is a constant matrix. 

Computation yields that $\Sigma^1$ acts on the fiber by
\begin{equation}
    \left[\begin{matrix}
      W_2\\
      Z_2
\end{matrix}\right]\,\,\,\mapsto\,\,\,\,\,\left[\begin{matrix}
      W_2\cdot E^{-1}\\
      Z_2\cdot E^{-1}
\end{matrix}\right]\,.
\end{equation}
Then $E(\widetilde Z_1,\widetilde Z_2)=E\left( Z_1- Z_2\cdot  W_2^{-1}\cdot  W_1,-  Z_2\cdot  W_2^{-1}\right)$ is independent of $ W_2$ and $ Z_2$. Define new matrix-valued variable $L=(l_{ij}):= Z_2\cdot  W_2^{-1}$;  then $E(\widetilde Z_1,\widetilde Z_2)=E\left( Z_1-L\cdot  W_1,-L\right)$ is independent of $L$. 

Denote by $e_{\alpha\beta}(\widetilde Z_1,\widetilde Z_2)$ the $(\alpha,\beta)^{th}$ entry of $E(\widetilde Z_1,\widetilde Z_2)$; denote the components of the variables $\widetilde Z_1$, $\widetilde Z_2$, $W_1$ by $\left((\widetilde Z_1)_{ij}\right)$, $\widetilde Z_2$ by $\left((\widetilde Z_2)_{ij}\right)$, $\big(( W_1)_{ij}\big)$ respectively.  Taking the partial derivatives with respect to $l_{ij}$, we have that
\begin{equation}
    -\sum_{l=1}^{n-s}\frac{\partial e_{\alpha\beta}(\widetilde Z_1,\widetilde Z_2)}{\partial(\widetilde Z_1)_{il}}\cdot \left(W_1\right)_{jl}-\frac{\partial e_{\alpha\beta}(\widetilde Z_1,\widetilde Z_2)}{\partial(\widetilde Z_2)_{ij}}\equiv 0\,,\,\,1\leq \alpha,\beta\leq p. 
\end{equation}
Then all partial derivatives of $E$ vanish and hence it is a constant matrix-valued function.

Therefore, we conclude that $\sigma\in PGL(s,\mathbb C)\times PGL(n-s,\mathbb C)$. \smallskip

As a conclusion, we complete the proof of Lemma \ref{fibpre}.\,\,\,$\endpf$

\subsection{An alternative proof of Lemma \ref{fibpre}.}\label{section:rigidtc2} We provide an alternative proof of  Lemma \ref{fibpre} which is based on the semi-positivity of $K_{\mathcal T_{s,p,n}}$ and the Pl\"ucker relations. 

Let $\sigma$ be an automorphism of $\mathcal M_{s,p,n}$ which induces the identity map on the Picard group. By Remark \ref{mcls}, we can show that   $\sigma$ is induced by an automorphism of the ambient space $\mathbb {CP}^{N_{p,n}}\times\mathbb {CP}^{N^0_{s,p,n}}\times\cdots\times\mathbb {CP}^{N^p_{s,p,n}}$.

Similarly to Lemma \ref{des}, $\sigma$ induces an automorphism of $\mathbb P(N_{\mathcal V_{(p,0)}/G(p,n)})$ as follows.
\vspace{-.1in}
\begin{equation}
\begin{array}{ccccc}
\mathcal M_{s,p,n}&\xrightarrow{\,\,\,\sigma\,\,\,} &\mathcal M_{s,p,n}&&\\
\Big\downarrow\llap{$\scriptstyle  \tau^1_{s,p,n}\,\,\,\,\,$}&&\Big\downarrow\rlap{$\scriptstyle \tau^1_{s,p,n}\,\,\,\,\,$}&&\\
\mathbb P(N_{\mathcal V_{(p,0)}/G(p,n)})&\xrightarrow{\,\,\,\Sigma\,\,\,}&\mathbb P(N_{\mathcal V_{(p,0)}/G(p,n)})&\xhookrightarrow{\,\,\,i\,\,\,}&\mathbb{CP}^{N_{p,n}}\times\mathbb{CP}^{N^1_{s,p,n}}\\
\end{array}\,\,\,\,\,\,\,\,\,\,\,\,\,\,\,\,.
\end{equation}
Moreover,  $\Sigma$ maps each fiber of $\kappa^1_{s,p,n}:\mathbb P(N_{\mathcal V_{(p,0)}/G(p,n)})\rightarrow \mathcal V_{(p,0)}$ to another fiber;   each fiber of $\kappa^1_{s,p,n}$ is isomorphic to  $P(M_{p\times (n-s)})$ and $\Sigma$ preserves the ranks of the matrices in $P(M_{p\times (n-s)})$.

\smallskip

{\bf\noindent Case 1 ($2\leq p<s$ and $s\neq 2p$).} We can assume that $\Sigma$ preserves the fibers of $\kappa^1_{s,p,n}$  by composing $\sigma$ with an element of $PGL(s,\mathbb C)\times PGL(n-s,\mathbb C)$.  Then the following diagram commutes.
\begin{equation}\label{st}
\begin{tikzcd}
&\mathbb P(N_{\mathcal V_{(p,0)}/G(p,n)})\arrow{r}{\,\,\Sigma}\arrow[hookrightarrow]{d}{\,\,i}  &\mathbb P(N_{\mathcal V_{(p,0)}/G(p,n)})\arrow[hookrightarrow]{d}{\,\,i}\\
&\mathbb{CP}^{N_{p,n}}\times\mathbb{CP}^{N^1_{s,p,n}}\arrow{r}{\,\,({\rm Id},\,\tau)\,\,\,\,\,\,} &\mathbb{CP}^{N_{p,n}}\times\mathbb{CP}^{N^1_{s,p,n}} \\
\end{tikzcd}\vspace{-20pt}\,
\end{equation} 
where $\tau$ is an automorphism of $\mathbb{CP}^{N^1_{s,p,n}}$.

In the following, we will show that
$\tau\in PGL(s,\mathbb C)\times PGL(n-s,\mathbb C)$. 

Recall the local coordinates $(X,[A\,])$  defined by (\ref{trr1}) in Appendix \ref{section:projbc} for the open subset $U^P$ of $\mathbb P(N_{\mathcal V_{(p,0)}/G(p,n)})$. Then locally the embedding $i$ can be given by 
\begin{equation}\label{coordinatesemb}
 i\left((X,[A\,])\right)=   \big[\cdots,P_I(M_{(X,A)}),\cdots\big]_{I\in\mathbb I_{p,n}}\times\big[\cdots,P_I(M_{(X,A)}),\cdots\big]_{I\in\mathbb I^1_{s,p,n}}\,
\end{equation}
where $P_I$ is the Pl\"ucker coordinate function and $M_{(X,A)}$ is the matrix defined by (\ref{mnin1}).

To distinguish the source and target $\mathbb P(N_{\mathcal V_{(p,0)}/G(p,n)})$ in (\ref{st}), we  use the above notation for the source and the following for the target. Denote by $\widetilde U^P$ the open subset $U^P$ in the target;  denote the local coordinates of $\widetilde U^P$ by 
\begin{equation}\label{yuv}
   (\cdots,y_{uv},\cdots)_{1\leq u\leq p,\,1\leq v\leq s}\times [\cdots, b_{ij},\cdots]_{1\leq i\leq p,\,s+1\leq j\leq n}=:(Y,[B])\,.
\end{equation}

Since $\tau$ is a projective linear transformation of $\mathbb{CP}^{N^1_{s,p,n}}$, we can write $\Sigma$ in the above local coordinates charts as follows. 
\begin{equation}\label{cij}
\begin{split}
&\,\,\,y_{uv}\left((Y,B)\right)=x_{uv}\left((X,[A])\right)\,,\,\,\,1\leq u\leq p\,,\,\,1\leq v\leq s\,,\\
&\left[\cdots,P_{I^{\prime}}\left(M_{(Y,[B])}\right),\cdots\right]_{I^{\prime}\in\mathbb I^1_{s,p,n}}=\left[\cdots,\sum_{I\in\mathbb I^1_{s,p,n}}C_{I^{\prime}}^{I}\cdot P_I(M_{(X,A)}),\cdots\right]_{I^{\prime}\in\mathbb I^1_{s,p,n}}\,,
\end{split}
\end{equation}
where $(Y,[B])=\Sigma\big((X,[A])\big)$ and $C^I_{I^{\prime}}\in\mathbb C$.

For simplicity, we introduce the following notation. $I_0:=(s,s-1,\cdots,s-p+1)$. 
For $1\leq r\leq p$ and $1\leq t\leq s-p$, 
\begin{equation}
I_{rt}:=(s,s-1,\cdots,s-p+r+1,\widehat{s-p+r},s-p+r-1\cdots,s-p+1,t)\,.
\end{equation}
For $1\leq r\leq p$ and $s+1\leq t\leq n$, 
\begin{equation}
I^*_{rt}:=(t,s,s-1,\cdots,s-p+r+1,\widehat{s-p+r},s-p+r-1,\cdots,s-p+1)\,.
\end{equation}
For $s-p+1\leq \alpha_1,\alpha_2\leq s$ and $\alpha_1\neq\alpha_2$, $1\leq \beta\leq s-p$ and $s+1\leq \gamma\leq n$,  
\begin{equation}
I^{\alpha_1\alpha_2}_{\beta\gamma}:=(\gamma,s,s-1,\cdots,\alpha_1+1,\widehat{\alpha_1},\alpha_1-1,\cdots,\alpha_2+1,\widehat{\alpha_2},\alpha_2-1,\cdots,s-p+1,\beta)\,.
\end{equation}

Computation yields that 
\begin{equation}
P_{I_0}(M_{(X,A)})=1\,, \,\,P_{I_{rt}}(M_{(X,A)})=(-1)^{r-1}\cdot  x_{rt}\,, \,\,
 P_{I^*_{rt}}(M_{(X,A)})=(-1)^{p-r}\cdot a_{rt}\,, \,\,
\end{equation}
and
\begin{equation}
\small
P_{I^{\alpha_1\alpha_2}_{\beta\gamma}}(M_{(X,[A])})=\left\{\begin{array}{cc}
   (-1)^{p+\alpha_1-\alpha_2}\left( x_{(\alpha_1-s+p)\beta}\cdot a_{(\alpha_2-s+p)\gamma}-x_{(\alpha_2-s+p)\beta}\cdot a_{(\alpha_1-s+p)\gamma}\right)  &{\rm if\,\,} \alpha_1>\alpha_2 \\
(-1)^{p-1+\alpha_1-\alpha_2}\left( x_{(\alpha_1-s+p)\beta}\cdot a_{(\alpha_2-s+p)\gamma}-x_{(\alpha_2-s+p)\beta}\cdot a_{(\alpha_1-s+p)\gamma}\right)     & {\rm if\,\,}\alpha_1<\alpha_2\\
\end{array}\right..
\end{equation}

We make the following claim.
\smallskip

{\bf\noindent Claim.}  Let $C_{I^{\prime}}^I$ be the constants in (\ref{cij}). Assume that $C_{I^*_{ij}}^I\neq0$ for a certain index $I\in\mathbb I^1_{s,p,n}$ and integers $1\leq i\leq p$, $s+1\leq j\leq n$. Then  $I=I^*_{rt}$ where $1\leq r\leq p$ and $s+1\leq t\leq n$.  
\smallskip

{\bf\noindent Proof of Claim.} Suppose  that $C_{I^*_{ij}}^{I^*}\neq0$ for integers $1\leq i\leq p$, $s+1\leq j\leq n$, and an index  $I^*:=(i_1,\cdots,i_p)\in\mathbb I^1_{s,p,n}\backslash\left\{I^*_{rt}\big|1\leq r\leq p,s+1\leq t\leq n\right\}$.  Then $i_p\leq s-p$. Take an integer $1\leq i^{\prime}\leq p$ such that $i^{\prime}\not\in\{i,i_2-s+p,i_3-s+p,\cdots, i_{p-1}-s+p\}$.

Recall that the grassmannian $G(p,n)$ is a subvariety of $\mathbb {CP}^{N_{p,n}}$ defined by the Pl\"ucker relations.
In particular, we have 
\begin{equation}\label{plcuker}
  P_{I^{(s-p+i^{\prime})(s-p+i)}_{i_pj}}\cdot P_{I_0}=\left\{\begin{array}{cc}
     P_{I_{i^{\prime}i_p}}\cdot P_{I^*_{ij}}- P_{I_{ii_p}}\cdot P_{I^*_{i^{\prime}j}}  & {\rm if\,\,} i<i^{\prime}\\
     -\left(P_{I_{i^{\prime}i_p}}\cdot P_{I^*_{ij}}- P_{I_{ii_p}}\cdot P_{I^*_{i^{\prime}j}} \right) & {\rm if\,\,} i>i^{\prime}
  \end{array}\right.\,.
\end{equation}
Since the image $\Sigma(U^P)$ satisfies the same equation, substituting (\ref{cij}) into (\ref{plucker}) we have that
\begin{equation}
\begin{split}
 \sum_{I\in\mathbb I^1_{s,p,n}}C_{I^{(s-p+i^{\prime})(s-p+i)}_{i_pj}}^{I}&\cdot P_I(M_{(X,A)})=\pm x_{i^{\prime}i_p}\cdot \left(\sum_{I\in\mathbb I^1_{s,p,n}}C_{I^*_{ij}}^{I} P_I(M_{(X,A)})\right)\\
 &\mp x_{ii_p}\cdot \left(\sum_{I\in\mathbb I^1_{s,p,n}}C_{I^*_{i^{\prime}j}}^{I} P_I(M_{(X,A)})\right)\,.
\end{split}
\end{equation}
Then the following equality holds as polynomials.
\begin{equation}\label{contr}
\begin{split}
 C_{I^*_{ij}}^{I^*}&\cdot x_{i^{\prime}i_p} P_{I^*}(M_{(X,A)})+x_{i^{\prime}i_p}\cdot \left(\sum_{I\in\mathbb I^1_{s,p,n}\,,\,\,I\neq I^*}C_{I^*_{ij}}^{I} P_I(M_{(X,A)})\right)\\
 &=\pm\sum_{I\in\mathbb I^1_{s,p,n}}C_{I^{(s-p+i^{\prime})(s-p+i)}_{i_pj}}^{I} P_I(M_{(X,A)})+ x_{ii_p}\cdot \left(\sum_{I\in\mathbb I^1_{s,p,n}}C_{I^*_{i^{\prime}j}}^{I} P_I(M_{(X,A)})\right)\,.
\end{split}
\end{equation}
Notice that the first term on the left hand side of (\ref{contr}) contains a nonzero monomial with a factor $x_{i^{\prime}i_p}^2$; however, there is no nonzero monomial with a factor $x_{i^{\prime}i_p}^2$ on the right hand side. This is a contradiction.

We complete the proof of Claim.\,\,\,$\endpf$
\smallskip

By Claim, we have that
\begin{equation}
b_{i^{\prime}j^{\prime}}\left(\Sigma(X,[A])\right)=\sum_{i=1}^p\sum_{j=s+1}^n(-1)^{i-i^{\prime}}\cdot C_{I^*_{i^{\prime}j^{\prime}}}^{I^*_{ij}}\cdot a_{ij}\,\,,\,\,\,{\,\,\,\,\,\,C_{I^*_{i^{\prime}j^{\prime}}}^{I^*_{ij}}}\in\mathbb C\,.\\
\end{equation}
Since $\Sigma$ preserves the ranks of the matrices in $P(M_{p(n-s)})$, similarly to Lemma \ref{pre}, we can show that there is an element $\sigma_l\in PGL(p,\mathbb C)$ and an element $\sigma_r\in PGL(n-s,\mathbb C)$ such that
\begin{equation}
\left(b_{i^{\prime}j^{\prime}}\left(\Sigma(X,[A])\right)\right)=\sigma_l\cdot (a_{ij})\cdot \sigma_r\,.\\
\end{equation}
Composed with a $PGL(n-s,\mathbb C)$-action, we can assume that the automorphism  $\Sigma$ is given by the left action $\sigma_l$.

We claim that $b_{ij}\left(\Sigma(X,[A])\right)= C_{I^*_{ij}}^{I^*_{ij}}\cdot a_{ij}$ for each  $1\leq i\leq p$, $s+1\leq j\leq n$. Otherwise, there are integers $1\leq i\leq p$, $1\leq i^{\prime}\leq p$, $i\neq i^{\prime}$, $s+1\leq j\leq n$ such that $C_{I^*_{ij}}^{I^*_{i^{\prime}j}}\neq 0$. However, then 
$P_{I^{(s-p+i)(s-p+i^{\prime})}_{(s-p)j}}\left(\Sigma(X,[A])\right)$ contains a monomial with a factor $a_{i^{\prime}j}\cdot x_{i^{\prime}(s-p)}$ which can not be a linear combination of $P_I(M_{(X,[A])})$. This is a contradiction.  
 
Similarly, we can conclude that $C_{I^*_{ij}}^{I^*_{ij}}=C_{I^*_{i^{\prime}j}}^{I^*_{i^{\prime}j}}$ for $1\leq i,i^{\prime}\leq p$. Then $\sigma_l$ is trivial.

We complete the proof of Lemma \ref{fibpre} when $2p\neq s$.
\smallskip

{\bf\noindent Case 2 ($2\leq p$ and $s=2p$).} Recall that the automorphism group of $G(p,2p)$ is generated by $PGL(2p,\mathbb C)$ and the dual automorphism. Similarly  to Case 2 in the first proof of Lemma \ref{fibpre} in Appendix \ref{section:rigidmc}, we can show that $\Sigma$ induces an automorphism of $G(p,2p)$ in $PGL(2p,\mathbb C)$. Then the same argument in  Case 1 applies here.
\medskip

We complete the proof of Lemma \ref{fibpre}.\,\,\,$\endpf$

\section{\texorpdfstring{Calculation for Delcroix's criterion\\ by Kexing Chen}{rr}}\label{section:nc}
This appendix {\footnote {Contact information for the author of this appendix: Kexing Chen: Shaoxing Vocational \& Technical College, 
526 Shanyin Rd, Shaoxing, Zhejiang, China, 312010. Email: chenkx$ @$sxvtc.com.}} provides certain  detailed calculation for Delcroix's criterion.
\begin{lemma}\label{numcal}Assume that $2\leq n-s\leq s<n$. Then 
(\ref{crt0}) holds if and only if $n=2s$.
\end{lemma}
{\bf\noindent Proof of Lemma \ref{numcal}.}
We make the following change of variables,
\begin{equation}
       y_1=x_1  \,,\,\,\,\,
      y_2=2x_2-x_1+2 \,.
\end{equation}
Then,
\begin{equation}
\begin{split}
&\rho\left(x_1,x_2\right)=y_2^2\left(\frac{-y_1+ y_2}{2} +n-s\right)^{n-s-2}\left(\frac{-y_1-y_2} {2}+n-s\right)^{n-s-2}\\
&\hspace{.5in}\cdot\left(\frac{y_1-y_2} {2}+s\right)^{s-2}\left(\frac{y_1+y_2} {2}+s\right)^{s-2}=:\eta\left(y_1,y_2\right)\,;\\
\end{split}
\end{equation}
the corresponding domain is transformed to
\begin{equation}
    \Omega:=\left\{(y_1,y_2)\in\mathbb R^2\,\big|\,-1\leq y_1\leq 1\,,\,\,y_2\geq 0\,,\,\,y_2\leq 4+y_1\,,\,\,y_2\leq 4-y_1\,\right\}\,.
\end{equation}
Hence, (\ref{crt0}) is equivalent to that
\begin{equation}\label{crt1}
0=\displaystyle\int_{\Omega}y_1\cdot\eta\left(y_1,y_2\right)dy_1dy_2\,,
\end{equation}
and that
\begin{equation}\label{crt2}
 0<\displaystyle\int_{\Omega}(y_2-2)\cdot\eta\left(y_1,y_2\right)dy_1dy_2\,.
\end{equation}

Denote by $\Omega^+$ the following subdomain of $\Omega$
\begin{equation}
    \Omega^+:=\left\{(y_1,y_2)\in\mathbb R^2\,\big|\,0\leq y_1\leq 1\,,\,\,y_2\geq 0\,,\,\,y_2\leq 4+y_1\,,\,\,y_2\leq 4-y_1\,\right\}\,.
\end{equation}
Then (\ref{crt1}) is equivalent to 
\begin{equation}
0=\displaystyle\int_{\Omega^+}y_1y_2^2\cdot \left(F^1_{s,n}(y_1,y_2)-F^2_{s,n}(y_1,y_2)\right)\cdot dy_1dy_2\,,
\end{equation}
where
\begin{equation}
\footnotesize
\begin{split}
&F^1_{s,n}(y_1,y_2)=\left(\frac{-y_1+ y_2}{2} +n-s\right)^{n-s-2}\left(\frac{-y_1-y_2} {2}+n-s\right)^{n-s-2}\left(\frac{y_1-y_2} {2}+s\right)^{s-2}\left(\frac{y_1+y_2} {2}+s\right)^{s-2}\,,\\
&F^2_{s,n}(y_1,y_2)=\left(\frac{y_1+ y_2}{2} +n-s\right)^{n-s-2}\left(\frac{y_1-y_2} {2}+n-s\right)^{n-s-2}\left(\frac{-y_1-y_2} {2}+s\right)^{s-2}\left(\frac{-y_1+y_2} {2}+s\right)^{s-2}\,.\\
\end{split}
\end{equation}
Notice that $F^1_{s,n}(y_1,y_2)$ and $F^2_{s,n}(y_1,y_2)$ are non-negative for $s\geq 2$.

Computation yields that
\begin{equation}
\begin{split}
&\frac{\partial}{\partial y_1}\left\{\log\left\{\frac{F^1_{s,n}(y_1,y_2)}{F^2_{s,n}(y_1,y_2)}\right\}\right\}=\left(\frac{-\frac{1}{2}\left(n-s-2\right)}{n-s-\frac{1}{2}y_1+\frac{1}{2}y_2}+\frac{-\frac{1}{2}\left(n-s-2\right)}{n-s-\frac{1}{2}y_1-\frac{1}{2}y_2}+\frac{\frac{1}{2}(s-2)}{s+\frac{1}{2}y_1-\frac{1}{2}y_2}\right.\\
&\,\,\,\,\,\,\,\,\,\,\,\,\,\,\,\,\,\,\,\,\,\left.+\frac{\frac{1}{2}(s-2)}{s+\frac{1}{2}y_1+\frac{1}{2}y_2}\right)-\left(\frac{\frac{1}{2}\left(n-s-2\right)}{n-s+\frac{1}{2}y_1+\frac{1}{2}y_2}+\frac{\frac{1}{2}\left(n-s-2\right)}{n-s+\frac{1}{2}y_1-\frac{1}{2}y_2}\right.\\
&\,\,\,\,\,\,\,\,\,\,\,\,\,\,\,\,\,\,\,\,\,\,\,\,\,\,\,\,\left.+\frac{-\frac{1}{2}(s-2)}{s-\frac{1}{2}y_1-\frac{1}{2}y_2}+\frac{-\frac{1}{2}(s-2)}{s-\frac{1}{2}y_1+\frac{1}{2}y_2}\right)\geq 0\,,\,\,\,(y_1,y_2)\in\Omega^+.\\
\end{split}
\end{equation}
Then the minimum of $\log\left\{\frac{F^1_{s,n}(y_1,y_2)}{F^2_{s,n}(y_1,y_2)}\right\}$ appears at $y_1=0$. Plugging $y_1=0$ into $\log\left\{\frac{F^1_{s,n}(y_1,y_2)}{F^2_{s,n}(y_1,y_2)}\right\}$, we can conclude that
\begin{equation}
  F^1_{s,n}(y_1,y_2)-F^2_{s,n}(y_1,y_2)\geq 0 \,,\,\,\,(y_1,y_2)\in\Omega^+.
\end{equation}

It is easy to verify that $ F^1_{s,n}(y_1,y_2)-F^2_{s,n}(y_1,y_2)\equiv 0$ if and only if $n=2s$. 
\smallskip

Next we will verify that (\ref{crt2}) holds when $n=2s$. First notice that (\ref{crt2}) is equivalent to
\begin{equation}
\footnotesize
\begin{split}
 0&<\displaystyle\int_{\Omega}(y_2-2)y_2^2\left(\left(\frac{-y_1+ y_2}{2} +s\right)\left(\frac{-y_1-y_2} {2}+s\right)\left(\frac{y_1-y_2} {2}+s\right)\left(\frac{y_1+y_2} {2}+s\right)\right)^{s-2}dy_1dy_2 \\
&=2\displaystyle\int_{\Omega^+}(y_2-2)y_2^2\left(\left(\frac{-y_1+ y_2}{2} +s\right)\left(\frac{-y_1-y_2} {2}+s\right)\left(\frac{y_1-y_2} {2}+s\right)\left(\frac{y_1+y_2} {2}+s\right)\right)^{s-2}dy_1dy_2\\ 
&=2s^{4s-8}\displaystyle\int_{\Omega^+}(y_2-2)y_2^2\left(1-\frac{2y_1^2+2y_2^2}{4s^2}+\frac{(y_1^2-y_2^2)^2}{(4s^2)^2}\right)^{s-2}dy_1dy_2\,.\\
\end{split}
\end{equation}

In what follows, we will show that
\begin{equation}\label{crt3}
\begin{split}
 0<\displaystyle\int_{\Omega^+}(y_2-2)y_2^2\left(1-\frac{2y_1^2+2y_2^2}{4s^2}+\frac{(y_1^2-y_2^2)^2}{(4s^2)^2}\right)^{s-2}dy_1dy_2\,.\\
\end{split}
\end{equation}

First notice that
\begin{equation}
\small
\begin{split}
&\displaystyle\int_{\Omega^+}(y_2-2)y_2^2\left(1-\frac{2y_1^2+2y_2^2}{4s^2}+\frac{(y_1^2-y_2^2)^2}{(4s^2)^2}\right)^{s-2}dy_1dy_2=\displaystyle\int_{\Omega^+}(y_2-2)y_2^2\,dy_1dy_2\\
&\,\,\,\,\,+\displaystyle\int_{\Omega^+}(y_2-2)y_2^2\left\{\left(1-\frac{2y_1^2+2y_2^2}{4s^2}+\frac{(y_1^2-y_2^2)^2}{(4s^2)^2}\right)^{s-2}-1\right\}dy_1dy_2\\
&\geq\displaystyle\int_{\Omega^+}(y_2-2)y_2^2\,dy_1dy_2-\sup_{(y_1,y_2)\in\Omega^+}\rule[-.16in]{0.01in}{.45in}\left(1-\frac{2y_1^2+2y_2^2}{4s^2}+\frac{(y_1^2-y_2^2)^2}{(4s^2)^2}\right)^{s-2}-1\,\rule[-.16in]{0.01in}{.45in}\cdot\displaystyle\int_{\Omega^+}|y_2-2|y_2^2\,dy_1dy_2\\
\end{split}
\end{equation}
Calculation yields that
\begin{equation}
\begin{split}
\displaystyle\int_{\Omega^+}(y_2-2)y_2^2\,dy_1dy_2=\frac{593}{60}\,\,\,\,\,\,{\rm and}\,\,\,\,\,\,\displaystyle\int_{\Omega^+}|y_2-2|y_2^2\,dy_1dy_2=\frac{251}{20}\,.
\end{split}
\end{equation}
Recall the following inequality
\begin{equation}
(1-x)^{\frac{1}{x}-1}\geq e^{-1}\,\,,\,\,\,0<x<1;
\end{equation}
then for $(y_1,y_2)\in\Omega^+$,
\begin{equation}
1\geq \left(1-\frac{2y_1^2+2y_2^2}{4s^2}+\frac{(y_1^2-y_2^2)^2}{(4s^2)^2}\right)^{s-2}\geq \left(1-\frac{8}{s^2}\right)^{s-2}\geq e^{-\frac{8(s-2)}{s^2-8}}\,.
\end{equation}
Therefore to prove (\ref{crt3}) it suffices to show that
\begin{equation}
\frac{593}{60}-\rule[-.1in]{0.01in}{.3in}\, e^{-\frac{8(s-2)}{s^2-8}}-1\,\rule[-.1in]{0.01in}{.3in}\, \frac{251}{20}\geq 0\,,
\end{equation}
or equivalently
\begin{equation}
\frac{8(s-2)}{s^2-8}\leq\log\left(\frac{753}{160}\right)\,.
\end{equation}
It is clear that when $s\geq 5$ the above inequality holds.

For small $s$ we compute the integral in (\ref{crt3}) by Mathematica as follows. When $s=2$
\begin{equation}
\begin{split}
 \displaystyle\int_{\Omega^+}(y_2-2)y_2^2\left(1-\frac{2y_1^2+2y_2^2}{4s^2}+\frac{(y_1^2-y_2^2)^2}{(4s^2)^2}\right)^{s-2}dy_1dy_2=\frac{593}{60}\,;\\
\end{split}
\end{equation}
when $s=3$
\begin{equation}
\begin{split}
 \displaystyle\int_{\Omega^+}(y_2-2)y_2^2\left(1-\frac{2y_1^2+2y_2^2}{4s^2}+\frac{(y_1^2-y_2^2)^2}{(4s^2)^2}\right)^{s-2}dy_1dy_2=\frac{43301123}{9797760}\,;\\
\end{split}
\end{equation}
when $s=4$
\begin{equation}
\begin{split}
 \displaystyle\int_{\Omega^+}(y_2-2)y_2^2\left(1-\frac{2y_1^2+2y_2^2}{4s^2}+\frac{(y_1^2-y_2^2)^2}{(4s^2)^2}\right)^{s-2}dy_1dy_2=\frac{196456943526409}{45343781683200}\,.\\
\end{split}
\end{equation}

We complete the proof of Lemma \ref{numcal}.
\,\,\,\,$\endpf$

\begin{lemma}\label{numcal2}Assume that $2\leq n-s\leq s<n$. Then, \begin{equation}
\small
\begin{split}
\displaystyle\int_{-1}^1x(2x+2)^2(n-s+1+x)^{n-s-2}(n-s-1+x)^{n-s-2}(s+1+x)^{s-2}(s-1-x)^{s-2}\,dx>0\,.\\
\end{split}
\end{equation}
\end{lemma}
{\bf\noindent Proof of Lemma \ref{numcal2}.} The proof is the same as in Lemma \ref{numcal}. Write
\begin{equation}
\small
\begin{split}
&\displaystyle\int_{-1}^1x(2x+2)^2(n-s+1+x)^{n-s-2}(n-s-1+x)^{n-s-2}(s+1+x)^{s-2}(s-1-x)^{s-2}\,dx\\
&\,\,=(n-s+1)^{n-s-2}(n-s-1)^{n-s-2}(s+1)^{s-2}(s-1)^{s-2}\\
&\,\,\,\,\,\,\,\,\,\,\,\,\cdot\displaystyle\int_{-1}^1x(2x+2)^2\left(1-\frac{2x+x^2}{(n-s)^2-1}\right)^{n-s-2}\left(1-\frac{2x+x^2}{s^2-1}\right)^{s-2}\,dx\,.\\
\end{split}
\end{equation}
To prove Lemma \ref{numcal2} it suffices to show that
\begin{equation}
\displaystyle\int_{-1}^1x(2x+2)^2\left(1-\frac{2x+x^2}{(n-s)^2-1}\right)^{n-s-2}\left(1-\frac{2x+x^2}{s^2-1}\right)^{s-2}\,dx>0\,.
\end{equation}

Notice that for $-1\leq x\leq 1$,
\begin{equation}
e^{\frac{-3}{(n-s)^2-4}}\leq 1-\frac{2x+x^2}{(n-s)^2-1}\leq e^{\frac{1}{(n-s)^2-1}}\,,
\end{equation}
and
\begin{equation}
e^{\frac{-3}{s^2-4}}\leq 1-\frac{2x+x^2}{s^2-1}\leq e^{\frac{1}{s^2-1}}\,.
\end{equation}
Hence for $-1\leq x\leq 1$,
\begin{equation}
\begin{split}
&e^{-\frac{3}{s+2}-\frac{3}{n-s+2}}\leq\left(1-\frac{2x+x^2}{(n-s)^2-1}\right)^{n-s-2}\left(1-\frac{2x+x^2}{s^2-1}\right)^{s-2}\leq e^{\frac{s-2}{s^2-1}+\frac{n-s-2}{(n-s)^2-1}}\\
\end{split}
\end{equation}
Then 
\begin{equation}\label{crit3}
\small
\begin{split}
&\displaystyle\int_{-1}^1x(2x+2)^2\left(1-\frac{2x+x^2}{(n-s)^2-1}\right)^{n-s-2}\left(1-\frac{2x+x^2}{s^2-1}\right)^{s-2}\,dx\\
&\geq\displaystyle\int_{-1}^1x(2x+2)^2\,dx-\left(\displaystyle\int_{-1}^1|x|(2x+2)^2\,dx\right)\cdot\sup\left\{1-e^{-\frac{3}{s+2}-\frac{3}{n-s+2}},\,e^{\frac{s-2}{s^2-1}+\frac{n-s-2}{(n-s)^2-1}}-1\right\}\,.\\
\end{split}
\end{equation}

Computation yields that
\begin{equation}\label{int2}
    \displaystyle\int_{-1}^1x(2x+2)^2=\frac{16}{3}\,\,\,\,\,{\rm and}\,\,\,\,\,\,\displaystyle\int_{-1}^1|x|(2x+2)^2\,dx=6\,.
\end{equation}
Hence, to guarantee the right hand side of (\ref{crit3}) is strictly positive, it suffice to have
\begin{equation}\label{int3}
\frac{3}{s+2}+\frac{3}{n-s+2}<\log 9\,,
\end{equation}
and
\begin{equation}\label{int4}
\frac{s-2}{s^2-1}+\frac{n-s-2}{(n-s)^2-1}<\log\left(\frac{17}{9}\right)\,.
\end{equation}

By (\ref{int2}) Lemma \ref{numcal2} holds when $n-s=s=2$; when $n-s\geq 3$ or $n-s=2$ and  $s\geq 3$, (\ref{int3}) and (\ref{int4}) hold.

We complete the proof of Lemma \ref{numcal2}.\,\,\,\,$\endpf$
\medskip

By a change of variables $y_i=x_i+\left(s-\frac{n}{2}+i-1\right)(p+1-i)$ in (\ref{criterion}), we can derive the following.
\begin{lemma}\label{criterion4}
Assume that $2\leq  n-s\leq p\leq s$ $(r=n-s)$. Let $x_1=\left(\frac{n}{2}-p\right)(n-s)$ and $x_{r+1}=0$.  Define a polytope 
$\Delta$ by
\begin{equation}
\Delta=\left\{(x_2,\cdots,x_r)\in\mathbb{R}^{r-1}\,\rule[-.16in]{0.01in}{.45in}\small{\begin{matrix}
\,\,0\leq x_i\leq 1+\left(\frac{n}{2}-p+i-1\right)(n-s+1-i)\,\,{\rm for}\,\,2\leq i\leq r\,;\\
x_{j-1}-2x_{j}+x_{j+1}\leq 0\,\,{\rm for}\,\,2\leq j\leq r\,
\end{matrix}}\right\}\,.
\end{equation} 
Define a density function $\rho(x_2,\cdots,x_r)$ by
\begin{equation}
\small
\rho:=\left(\prod\limits_{1\leq i<j\leq r}((x_{i+1}-x_{i})-(x_{j+1}-x_{j}))^2\right)\left(\prod\limits_{i=1}^r\left(\frac{n}{2}+(x_{i+1}-x_{i})\right)^{s+p-n}\left(\frac{n}{2}-(x_{i+1}-x_{i})\right)^{s-p}\right)\,.
\end{equation}
Then there exists a K\"ahler-Einstein metric on $\mathcal M_{s,p,n}$ if and only if for $2\leq k\leq r$,
\begin{equation}\label{scri}
    \frac{\displaystyle\int_{\Delta}\left(x_k\cdot\rho(x_2,\cdots,x_r)\right)\, dx_2dx_3\cdots dx_r}{\displaystyle\int_{\Delta}\rho(x_2,\cdots,x_r)\, dx_2dx_3\cdots dx_r}>\left(\frac{n}{2}-p+i-1\right)(n-s+1-i)\,.
\end{equation}
\end{lemma}
Using software Mathematica,  we have the following.
\begin{lemma}\label{critm1}
There exists a K\"ahler-Einstein metric on
$\mathcal M_{p,p,2p}$ if $p=4,5$.
\end{lemma}
{\bf\noindent Proof of Lemma \ref{critm1}.}
When $p=4$,
\begin{equation}
\Delta=\left\{(x_2,x_3,x_4)\in\mathbb{R}^{3}\,\rule[-.16in]{0.01in}{.45in}\small{\begin{matrix}
\,\,0\leq x_2\leq 4,\,\,0\leq x_3\leq 5,\,\,{\rm and}\,\,0\leq x_4\leq 4\,;\\
-2x_{2}+x_{3}\leq 0,\,\,
x_{2}-2x_{3}+x_{4}\leq 0,\,\,{\rm and}\,\,x_{3}-2x_{4}\leq 0\\
\end{matrix}}\right\}\,;
\end{equation} 
\begin{equation}
\begin{split}
\rho&=(2x_2-x_3)^2(x_2-x_4+x_3)^2(x_2+x_4)^2(x_2-2x_3+x_4)^2(x_3-x_2+x_4)^2(2x_4-x_3)^2\,.
\end{split}
\end{equation}
Then,
\begin{equation}
\begin{split}
&\displaystyle\int_{\Delta}\rho\, dx_2dx_3dx_4=\frac{2243664235225939}{567567000}\approx 3.95313\times 10^6\,, \\
&\displaystyle\int_{\Delta}\left(x_2\cdot\rho\right)\, dx_2dx_3dx_4=\frac{55382785289338434218971}{4067390354227200}\approx 1.36163\times10^7\,,  \\
&\displaystyle\int_{\Delta}\left(x_3\cdot\rho\right)\, dx_2dx_3dx_4=\frac{5416920544038914803}{305543145600}\approx 1.77288\times10^7\,,   \\
&\displaystyle\int_{\Delta}\left(x_4\cdot\rho\right)\, dx_2dx_3dx_4=\frac{55382785289338434218971}{4067390354227200} \approx 1.36163\times10^7  \\
\end{split}
\end{equation}
Then
\begin{equation}
\begin{split}   &\frac{\displaystyle\int_{\Delta}\left(x_2\cdot\rho(x_2,x_3,x_4)\right)\, dx_2dx_3dx_4}{\displaystyle\int_{\Delta}\rho(x_2,x_3,x_4)\, dx_2dx_3dx_4}=\frac{\displaystyle\int_{\Delta}\left(x_4\cdot\rho(x_2,x_3,x_4)\right)\, dx_2dx_3dx_4}{\displaystyle\int_{\Delta}\rho(x_2,x_3,x_4)\, dx_2dx_3dx_4}\approx 3.4>3\,,\\
&\hspace{1.2in}\frac{\displaystyle\int_{\Delta}\left(x_3\cdot\rho(x_2,x_3,x_4)\right)\, dx_2dx_3dx_4}{\displaystyle\int_{\Delta}\rho(x_2,x_3,x_4)\, dx_2dx_3dx_4}\approx 4.4>4\,.\\
\end{split}
\end{equation}
Therefore, there exists a K\"ahler-Einstein metric on $\mathcal M_{p,p,2p}$ when $p=4$ by (\ref{scri}).
\medskip

When $p=5$,
\begin{equation}
\Delta=\left\{(x_2,x_3,x_4,x_5)\in\mathbb{R}^{4}\,\rule[-.16in]{0.01in}{.45in}\small{\begin{matrix}
\,\,0\leq x_2\leq 5,\,\,0\leq x_3\leq 7,\,\,0\leq x_4\leq 7,\,\,{\rm and}\,\,0\leq x_5\leq 5\,;\\
-2x_{2}+x_{3}\leq 0,\,\,
x_{2}-2x_{3}+x_{4}\leq 0,\\
\,x_{3}-2x_{4}+x_{5}\leq 0,\,\,{\rm and}\,\,x_{4}-2x_{5}\leq 0\\
\end{matrix}}\right\}\,;
\end{equation} 
\begin{equation}
\begin{split}
\rho&=(2x_2-x_3)^2(x_2+x_3-x_4)^2(x_2+x_4-x_5)^2(x_2-2x_3+x_4)^2(x_2-2x_3+x_4)\\
&\,\,\,\,\,\,\,\,\,\,\,\cdot(x_3-x_2+x_4-x_5)^2(x_3-x_2+x_5)^2(x_3-2x_4+x_5)^2(x_4-x_3+x_5)^2(2x_5-x_4)^2\,.
\end{split}
\end{equation}
Then,
\begin{equation}
\small
\begin{split}
&\displaystyle\int_{\Delta}\rho\, dx_2dx_3dx_4dx_5=\frac{57336210099961579033706793911}{74803289175014400}\approx 7.66493\times 10^{11}\,,\\
&\displaystyle\int_{\Delta}\left(x_2\cdot\rho\right)\, dx_2dx_3dx_4dx_5=\frac{1760441835266075955851497040448393302143609}{521057531362788347090042880000}\approx 3.37859\times10^{12}\,,\\
&\displaystyle\int_{\Delta}\left(x_3\cdot\rho\right)\, dx_2dx_3dx_4dx_5=\frac{2567995351960762288549954674509341582094471}{521057531362788347090042880000}\approx 4.92843\times10^{12}\,, \\
&\displaystyle\int_{\Delta}\left(x_4\cdot\rho\right)\, dx_2dx_3dx_4dx_5=\frac{2567995351960762288549954674509341582094471}{521057531362788347090042880000}\approx 4.92843\times10^{12}\,,  \\
&\displaystyle\int_{\Delta}\left(x_5\cdot\rho\right)\, dx_2dx_3dx_4dx_5=\frac{1760441835266075955851497040448393302143609}{521057531362788347090042880000}\approx 3.37859\times10^{12}\,, \\
\end{split}
\end{equation}
Then
\begin{equation}
\begin{split}   &\frac{\displaystyle\int_{\Delta}\left(x_2\cdot\rho\right)\, dx_2dx_3dx_4dx_5}{\displaystyle\int_{\Delta}\rho(x_2,x_3,x_4,x_5)\, dx_2dx_3dx_4}=\frac{\displaystyle\int_{\Delta}\left(x_5\cdot\rho\right)\, dx_2dx_3dx_4dx_5}{\displaystyle\int_{\Delta}\rho\, dx_2dx_3dx_4dx_5}\approx 4.4>4\,,\\
&\frac{\displaystyle\int_{\Delta}\left(x_3\cdot\rho\right)\, dx_2dx_3dx_4dx_5}{\displaystyle\int_{\Delta}\rho\, dx_2dx_3dx_4}=  \frac{\displaystyle\int_{\Delta}\left(x_4\cdot\rho\right)\, dx_2dx_3dx_4dx_5}{\displaystyle\int_{\Delta}\rho\, dx_2dx_3dx_4dx_5}\approx 6.4>6\,.\\
\end{split}
\end{equation}
Therefore, there are K\"ahler-Einstein metric on $\mathcal M_{p,p,2p}$ when $p=5$ by (\ref{scri}).

We complete the proof of Lemma \ref{critm1}.\,\,\,\,$\endpf$

\end{document}